\newcommand{\pointsize}{11pt}

\documentclass[oneside, \pointsize, reqno, fleqn]{amsbook}

\usepackage[
   includehead,
   includefoot,
     left = 1in, 
      top = 1in, 
    right = 1in,
   bottom = 1in
]{geometry}
\usepackage{fancyhdr}
\usepackage{setspace}
\usepackage{calc}
\usepackage[nocompress]{cite} 
\usepackage[pdfborder={0 0 0}, pdfpagemode=UseNone, pdfstartview=FitH]{hyperref} 

\fancyheadoffset[R]{0.5in} 

\fancypagestyle{prelim}{%
   \renewcommand{\headrulewidth}{0pt} 
   \fancyhf{}           
   \pagenumbering{roman}    
   \cfoot{\thepage}       
}

\fancypagestyle{maintext}{%
   \renewcommand{\headrulewidth}{0pt} 
   \pagenumbering{arabic}
   \fancyhf{}
   \cfoot{\thepage}
}

\fancypagestyle{abstract1}{%
   \renewcommand{\headrulewidth}{0pt} 
	 \fancyheadoffset[R]{0in}
   \pagenumbering{arabic}
   \fancyhf{}
}

\fancypagestyle{abstract2}{%
   \renewcommand{\headrulewidth}{0pt} 
   \fancyheadoffset[R]{0in}
	 \pagenumbering{arabic}
   \fancyhf{}
   \rhead{\thepage}
}

\numberwithin{figure}{chapter} 
\numberwithin{table}{chapter}
\numberwithin{equation}{chapter}
\numberwithin{section}{chapter}

\usepackage{graphicx}
\usepackage[font=footnotesize]{caption}
\usepackage[list=true]{subcaption} 

\usepackage{lineno,hyperref}
\usepackage{amsfonts}
\usepackage{amssymb}
\usepackage{multirow}
\usepackage{multicol}
\usepackage{wrapfig}
\usepackage{bm}
\usepackage{amsmath} 
\usepackage{subcaption}
\usepackage{dsfont}
\usepackage{enumitem}
\usepackage{scalerel}

\makeatletter
\def\@tocline#1#2#3#4#5#6#7{\relax
\ifnum #1>\c@tocdepth 
\else
\par \addpenalty\@secpenalty\addvspace{#2}%
\begingroup \hyphenpenalty\@M
\@ifempty{#4}{%
  \@tempdima\csname r@tocindent\number#1\endcsname\relax
}{%
  \@tempdima#4\relax
}%
\parindent\z@ \leftskip#3\relax \advance\leftskip\@tempdima\relax
\rightskip\@pnumwidth plus4em \parfillskip-\@pnumwidth
#5\leavevmode\hskip-\@tempdima #6\nobreak\relax
\ifnum#1<0\hfill\else\dotfill\fi\hbox to\@pnumwidth{\@tocpagenum{#7}}\par
\nobreak
\endgroup
\fi}
\makeatother

\usepackage{titlesec}
\newcommand*{\justifyheading}{\raggedright}

\titleformat{\chapter}[display]
  {\normalfont\bfseries\LARGE\justifyheading}
  {\chaptertitlename~\thechapter}{1pc}
  {{\color{black}\titlerule[2pt]}\vspace{1pc}}
\titleformat{name=\chapter,numberless}[display]
  {\normalfont\bfseries\LARGE}{}{1pc}
  {{\color{black}\titlerule[2pt]}\vspace{1pc}}
  
\titleformat{\section}
  {\normalfont\large\bfseries\justifyheading}{\thesection}{1em}{}
\titleformat{\subsection}
  {\normalfont\large\justifyheading}{\thesubsection}{1em}{}
  \titlespacing*{\chapter}
{0pt}{0pt}{1cm}
  \titlespacing*{\section}
{0pt}{1cm}{0.5cm}
\titlespacing*{\subsection}
{0pt}{1cm}{0.5cm}

\newcommand{\x}{\bm{x}}
\newcommand{\X}{\bm{X}}
\newcommand{\n}{\bm{n}}
\renewcommand{\u}{\bm{u}}
\newcommand{\ub}{\bm{\bar u}}
\newcommand{\F}{\bm{F}}

\newcommand{\sig}{\bm{\sigma}}
\newcommand{\gt}{\bm{\tilde g}}
\newcommand{\g}{\bm{g}}
\newcommand{\G}{\bm{G}}
\newcommand{\T}{\bm{T}}
\newcommand{\p}{\bm{p}}
\newcommand{\A}{\bm{A}}
\newcommand{\Q}{\bm{Q}}
\newcommand{\U}{\bm{U}}
\newcommand{\grad}{\nabla}
\newcommand{\Div}{\nabla \boldsymbol{\cdot}}
\newcommand{\dotp}{\boldsymbol{\cdot}}
\newcommand{\C}{\mathcal{C}}
\renewcommand{\L}{\mathcal{L}}
\newcommand{\jump}[1]{\ensuremath{[\![#1]\!]} }
\newcommand{\squeezeup}{\vspace{-2mm}}

\hypersetup{
  colorlinks   = true, 
  urlcolor     = blue, 
  linkcolor    = blue, 
  citecolor   = blue 
}

\usepackage{titletoc}

\makeatletter
\renewcommand{\tocchapter}[3]{%
  \indentlabel{{ \@ifnotempty{#2}{\bfseries{{#2.\quad#3}}}}}}
\makeatother

\dottedcontents{section}[8em]{}{2em}{0.4pc}
\dottedcontents{subsection}[13em]{}{3em}{0.4pc}  
\dottedcontents{figure}[3.5em]{}{3em}{0.4pc}
\newcounter{lofdepth}
\begin{document}

   \frontmatter

   \pagestyle{prelim}
   
   %
   \fancypagestyle{plain}{%
      \fancyhf{}
      \cfoot{\thepage}
   }%
   \begin{center}
   \null\vfill
   \textbf{%
      The Immersed Boundary Double Layer (IBDL) Method
   }%
   \\
   \bigskip
   By \\
   \bigskip
   BRITTANY JAE LEATHERS \\
   DISSERTATION \\
   \bigskip
   Submitted in partial satisfaction of the requirements for the
   degree of \\
   \bigskip
   DOCTOR OF PHILOSOPHY \\
   \bigskip
   in \\
   \bigskip
   APPLIED MATHEMATICS \\
   \bigskip
   in the \\
   \bigskip
   OFFICE OF GRADUATE STUDIES \\
   \bigskip        
   of the \\
   \bigskip
   UNIVERSITY OF CALIFORNIA \\
   \bigskip
   DAVIS \\
   \bigskip
   Approved: \\
   \bigskip
   \bigskip
   \makebox[3in]{\hrulefill} \\
  Robert Guy, Chair\\
   \bigskip
   \bigskip
   \makebox[3in]{\hrulefill} \\
   Becca Thomases \\
   \bigskip
   \bigskip
   \makebox[3in]{\hrulefill} \\
   Gregory Miller \\
   \bigskip
   Committee in Charge \\
   \bigskip
   2022 \\
   \vfill
\end{center}

   \newpage
	
	 \thispagestyle{empty}
	 \begin{titlepage}
	 \vspace*{50em}
	 \begin{center}
		 \copyright \ Brittany J.\ Leathers, 2022.  All rights reserved.  
	 \end{center}
	 \end{titlepage}
	 \newpage
	 \stepcounter{page}
	
	 \thispagestyle{plain}
	 \vspace*{20em}
	 \begin{center}
	   To Colin
	 \end{center}
	 \newpage
   
   %
   \doublespacing
\let\oldaddcontentsline\addcontentsline
\renewcommand{\addcontentsline}[3]{}
   
   \chapter*{Abstract}
   The Immersed Boundary (IB) method of Peskin (J. Comput. Phys., 1977) is useful for problems that involve fluid-structure interactions or complex geometries. By making use of a regular Cartesian grid that is independent of the geometry, the IB framework yields a robust numerical scheme that can efficiently handle immersed deformable structures. Additionally, the IB method has been adapted to problems with prescribed motion and other PDEs with given boundary data. IB methods for these problems traditionally involve penalty forces which only approximately satisfy boundary conditions, or they are formulated as constraint problems. In the latter approach, one must find the unknown forces by solving an equation that corresponds to a poorly conditioned first-kind integral equation. This operation can therefore require a large number of iterations of a Krylov method, and since a time-dependent problem requires this solve at each step in time, this method can be prohibitively inefficient without preconditioning. In this dissertation, we introduce a new, well-conditioned IB formulation for boundary value problems, which we call the Immersed Boundary Double Layer (IBDL) method. We formulate the method for Poisson, Helmholtz, Brinkman, Stokes, and Navier-Stokes equations and demonstrate its efficiency over the original constraint method. In this double layer formulation, the equation for the unknown boundary distribution corresponds to a well-conditioned second-kind integral equation that can be solved efficiently with a small number of iterations of a Krylov method without preconditioning. Furthermore, the iteration count is independent of both the mesh size and spacing of the immersed boundary points. The method converges away from the boundary, and when combined with a local interpolation, it converges in the entire PDE domain. Additionally, while the original constraint method applies only to Dirichlet problems, the IBDL formulation can also be used for Neumann boundary conditions. 
   
   \newpage
   
   \chapter*{Acknowledgments}
  
  There are many people to whom I owe thanks, but I would like to thank two people in particular. 

Thank you to my advisor, Bob Guy. You have given me a tremendous amount of your time, support, and guidance. You have pushed me to become a better mathematician and researcher. You have shown me patience and kindness when I needed it, and you have taught me what a difference a great mentor can make. 

And thank you to my partner, Joshua Parker. Thank you for letting me run ideas by you, for cooking meals when research took up all my time, and for being there when I needed you. You have supported me in more ways than I can count. You have also introduced me to new ways of thinking and helped me grow as a person. Finally, you have believed in me even when my confidence has faltered, and for this I am exceedingly grateful.

This work was supported in part by NSF grant DMS-1664679 to R.D.G.

   \let\addcontentsline\oldaddcontentsline

       \tableofcontents
  
   \newpage
   
   \let\oldaddcontentsline\addcontentsline
\renewcommand{\addcontentsline}[3]{}
   \listoffigures
 \newpage
\listoftables

   \let\addcontentsline\oldaddcontentsline
\newpage
   
   \mainmatter
   
   \pagestyle{maintext}
   
   %
   \fancypagestyle{plain}{%
      \renewcommand{\headrulewidth}{0pt}
      \fancyhf{}
      \cfoot{\thepage}
   }%
   
   \chapter{Introduction}
   \label{chapter 1}
 
 The Immersed Boundary (IB) method \cite{Peskin77, Peskin02} is a valuable numerical tool for fluid-structure interactions. It was initially developed by Peskin for problems involving elastic, deformable structures, such as those involved in cardiac dynamics \cite{Peskin-72, Peskin-81-heart, mitral-valve, heart-model, PeskinGriffithheart}. However, the robustness and simplicity of the IB method has led to its use in many different applications (see \cite{red-blood-cells, Peskin-Fogelson, pulp}, for just a few). Additionally, since its creation, there have been many developments and variations in the IB method. For instance, the IB framework has been altered to incorporate porous boundaries \cite{porous1,porous2} and model elastic rods with a curvature or twist \cite{gIB1, gIB2}. It has also been adapted to the flow of non-Newtonian fluids \cite{nonNewtonian1,nonNewtonian2} and coupled with internal force mechanisms to model swimming organisms \cite{C-F-C-D,cilia}. In recent years, work has been done to apply the Immersed Boundary framework to problems involving prescribed boundary values \cite{Taira, GriffithDonev, GriffithBhalla, Uhlmann, SuLaiLin, Guylewis}, which is the category of problems that this dissertation addresses. 

The robustness of the IB method comes from its use of two coordinate systems: a Lagrangian system that moves with the structure and a fixed Eulerian system on which the fluid equations are solved. It uses convolutions with discrete delta functions to link these systems together and map forces from the structure to the grid. One can then solve the PDE on a regular Cartesian mesh. This elimination of the physical boundary makes it possible to use a PDE solver that is efficient and independent of the geometry of the structure.

As stated, this dissertation will focus on problems involving prescribed boundary values. In the field of fluid dynamics, this includes the motion of rigid bodies and fluid flow through domains with stationary boundaries or boundaries with prescribed motion. Since the IB method bypasses the need for a conforming mesh, it also has obvious advantages in the broader case of solving PDEs on complex domains. 

When the IB method is applied to deformable structures, the boundary force density is found using a constitutive law from the structure characteristics \cite {Peskin02}, but in the case of rigid bodies, one needs a different way to find or interpret the boundary force. There have been several routes taken to use the IB framework in these situations. One method is to consider the Lagrangian points to be tethered to specified locations by springs. The boundary force is then interpreted as the spring restoring force, penalizing deviations from the prescribed boundary position \cite{Beyer, Goldstein, formal, TeranPeskin}. By using the IB framework, such penalty methods can be efficiently implemented, but they use parameters, such as spring constants, to approximate a rigid limit, and the required magnitude of the appropriate parameter leads to numerically stiff equations that necessitate very small time steps. 

Another approach for applying the IB method to a rigid body problem consists of viewing the boundary force density as a Lagrange multiplier, used to enforce the no-slip boundary condition \cite{Taira}. Chapter \ref{chapter 2} gives a description of this method, and we refer to it as the Immersed Boundary Single Layer (IBSL) method. In this constraint method, the velocity and force are both unknowns in an algebraic system. One way to solve this system is to invert the Schur complement to solve first for the Lagrange multiplier force and then for the velocity. The main disadvantage of such a method is that this operator suffers from poor conditioning \cite{GriffithDonev}. Furthermore, the conditioning of the discrete problem worsens when the Cartesian grid is refined or when the Lagrangian point spacing is refined relative to the grid. To avoid solving the fully discretized constraint problem, most numerical methods using this constraint approach rely on some form of time step splitting \cite{Taira, GriffithBhalla, Uhlmann, SuLaiLin, Guylewis}. Taira and Colonius \cite{Taira}, for example, solve a simpler unconstrained system for an intermediate velocity, use this velocity to find pressure and the unknown boundary force, and then complete a projection step to remove non-divergence-free and slip components of the velocity. However, like the penalty methods, these fractional step methods only satisfy the constraint equations \textit{approximately}, which can result in fluid penetration into a rigid body. Additionally, such a method cannot be used for steady Stokes or other time-independent PDE. 

However, if instead of a time-splitting scheme, one inverts the Schur complement with a Krylov method, the number of iterations required could be very large due to its poor conditioning. The computational cost then becomes prohibitive for a time-dependent problem, when this solve would be required at each time step. One way around this is to design a preconditioner \cite{Ceniceros, guyphilipgriffith, Stein}, such as the physics-based approximation of the Schur complement constructed by Kallemov et al. \cite{GriffithDonev}. A preconditioner can allow for a more efficient solution to the constraint problem, but developing and implementing one is a nontrivial undertaking, and such a preconditioner can involve computing the inverse of a dense matrix, which can be computationally expensive when the number of boundary points is large. Another limitation is that in order to control the conditioning, there is often a requirement that the boundary points be spaced to about twice the grid spacing, which, depending on the application and the Eulerian-Lagrangian coupling scheme, may result in decreased accuracy of the solution \cite{GriffithDonev, Griffithpointspacing, Griffithpressure}. 

In this dissertation, we present a reformulated Immersed Boundary method for prescribed boundary values. Like the Immersed Boundary Single Layer method, our method enforces the boundary conditions \textit{exactly}, but the resulting linear system is very well-conditioned. We can therefore solve it with an unpreconditioned Krylov method with very few iterations. This conditioning does not worsen as we refine the mesh, nor as we tighten the spacing of our boundary points relative to the grid. We are therefore able to avoid the need for a preconditioner altogether.

We formulate this new IB method by utilizing a connection between the IBSL method and a single layer boundary integral equation. Boundary integral methods rely on reformulating a boundary value problem as an integral equation with an unknown density on the boundary. These methods are particularly useful for linear, elliptic, and homogeneous PDE \cite{Pozblue}. One advantage of these methods, which we exploit, is that there are well-conditioned integral representations available. For example, for the Helmholtz or Poisson equation, starting with a broad integral representation, one can derive single and double layer integral representations, the second of which produces an operator with much better conditioning. There are, however, some disadvantages to boundary integral methods. Firstly, once the boundary density is found, it is expensive to evaluate the solution on an entire grid, whereas by using efficient solvers, the IB method can do this quickly. Additionally, in order to directly implement integral methods, one needs an analytical Green's function for each specific PDE and problem domain, making nontrivial exterior boundary conditions complicated to implement. The IB method, on the other hand, does not require a Green's function and can be used on more general domains. The method we present in this paper maintains the flexibility of the IB method while capturing the better conditioning of a double layer integral equation.

The connection between the IBSL method and a regularized single layer integral equation has been identified in a few recent works. Usabiaga et al. \cite{GriffithDonev2} illustrate this connection and discuss that their rigid multiblob method can be seen as a technique for solving a regularized first-kind integral equation for Stokes flow. Eldredge \cite{eldredge} made a more general connection between IB methods and boundary integral equations by extending the form of a PDE to govern a variable that defines a different function for each side of an immersed boundary. The resulting PDE contains jumps in field quantities across the boundary that correspond to the strengths of single and double layer potentials. This general connection allows for solutions on either side of the boundary. However, in this work, we look specifically at the IB method for prescribed boundary value problems for which the PDE domain exists on only one side of the boundary, and we create an IB formulation that corresponds to the use of a double layer potential with an unknown strength in order to take advantage of the better conditioning. 

In this dissertation, we present the new \textit{Immersed Boundary Double Layer (IBDL) method}, which is able to achieve the same order of accuracy as the original IBSL method away from the boundary, while only requiring a small number of iterations of a Krylov method. We derive, implement, and analyze this method for the 2-D Helmholtz, Poisson, Brinkman, Stokes, and Navier-Stokes equations with Dirichlet boundary conditions. While we focus on two dimensions in this dissertation, the method is not restricted to this case. Additionally, while the IBSL method applies only to Dirichlet problems, our new method can also be used for Neumann boundary conditions. 

This dissertation is organized as follows. The first three chapters provide introductory and background material. In Chapter \ref{chapter 2}, we introduce the Immersed Boundary Single Layer method for solving PDEs with Dirichlet boundary conditions. In Chapter \ref{chapter 3}, we provide an introduction to Green's functions and boundary integral equations for the relevant PDEs. 

In Chapter \ref{chapter 4}, we explicitly connect the IBSL method for Poisson and Stokes equations to corresponding single layer integral equations. In Chapter \ref{chapter 5}, we present the form of the Immersed Boundary Double Layer method as it applies to Poisson and Helmholtz equations. We then explicitly demonstrate that the method corresponds to a regularized double layer integral equation. We then discuss a linear interpolation that we use in order to obtain convergence in the solution near the immersed boundary. Additionally, we form and analyze the IBDL method for Neumann boundary conditions. Lastly, the method is tested and analyzed on several scalar PDEs on interior and exterior domains. We demonstrate that the method achieves the same first-order accuracy as the IBSL method, while requiring far fewer iterations of a Krylov method.

The IBSL method has a very direct generalization from scalar PDEs to Stokes equation, in which fluid velocity and pressure are coupled. However, matching the form of the double layer integral representation for Stokes equation requires a non-trivial adaptation of the IBDL method. Chapter \ref{chapter 6} presents the form of the IBDL method in this case and demonstrates its connection to a regularized double layer integral equation. We again test and analyze several PDEs and domains, including Stokes flow past a periodic array of cylinders, for which we compare drag forces against several numerical methods and asymptotic approximations. We additionally implement the IBDL method in the time-dependent Navier-Stokes equation and again compare drag forces with other numerical methods. Finally, in Chapter \ref{chapter 7}, we discuss the benefits and drawbacks of the IBDL method. We also discuss several directions that our future work will take in order to improve and generalize the IBDL method.

   \chapter{The Immersed Boundary method for boundary value problems}
   \label{chapter 2}

In this chapter, we give an introduction to the Immersed Boundary (IB) method for solving elliptic Dirichlet boundary value problems. For a scalar function $u$ and an elliptic operator $\L$, such a problem has the form
\begin{subequations} \label{pde}
\begin{alignat}{2}
& \L u = g \qquad && \text{in } \Omega  \label{pde1}\\
&u=U_b \qquad && \text{on } \Gamma,  \label{pde2}
\end{alignat}
\end{subequations}
where the domain of the PDE is $\Omega$, with boundary $\Gamma=\partial \Omega$. We will let $\C$ be a larger computational domain containing $\Omega$. The method presented is not specific to two dimensions, but for this dissertation, we will take $\C$ to be a two-dimensional torus or periodic box.  We will also assume $\Gamma$ is a smooth one-dimensional curve that does not self-intersect.

The constraint method presented in this chapter has been developed in previous works \cite{Peskin02, Taira, GriffithDonev}, and in this dissertation, we will refer to this method as the Immersed Boundary Single Layer (IBSL) method. The rationale for this name will become clear in Chapter \ref{chapter 4}. It is presented here in order to introduce the notation and concepts that will be used in the presentation of the Immersed Boundary Double Layer (IBDL) method in subsequent chapters. Additionally, since one aim of this dissertation is to improve the efficiency of the IBSL method while preserving its robustness, we will make frequent comparisons to this method.

The first two sections of this chapter will focus on the IBSL method applied to scalar elliptic PDEs and the numerical implementation of the method. In Section \ref{ch 2 stokes}, we then discuss the application of the method to Stokes and Brinkman equations, and in Section \ref{ch 2 results}, we give some numerical results. Lastly, in Section \ref{ch 2 discussion}, we discuss the benefits and limitations of the method in order to further motivate the IBDL method. 

\section{Mathematical description of IBSL method: scalar elliptic equations}\label{ch2 description}

\begin{figure}
\centering
\begin{subfigure}{0.495\textwidth}
\centering
\includegraphics[width=\textwidth]{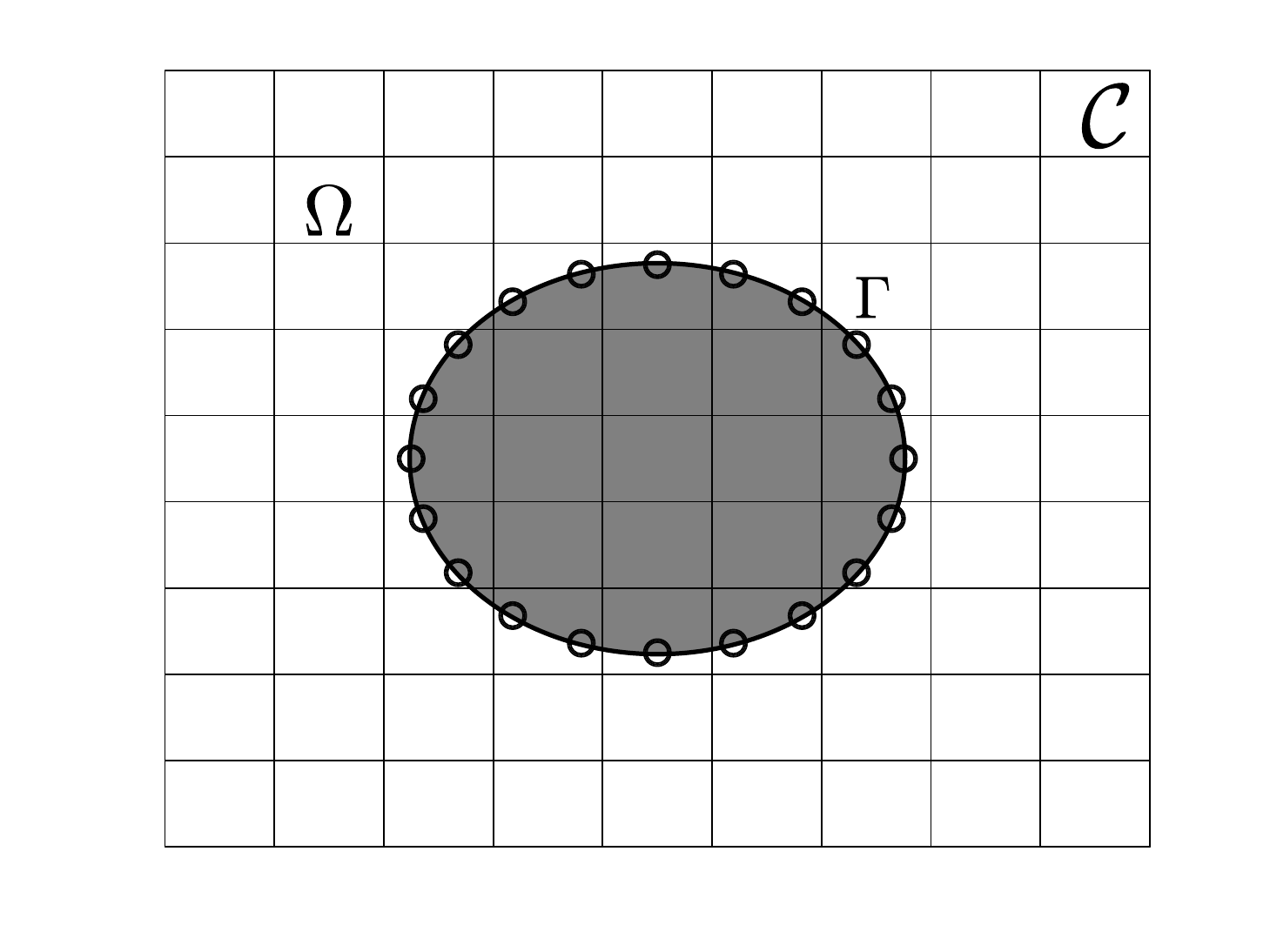}
\subcaption{\normalsize Exterior domain}
\label{exterior domain}
\end{subfigure}
\begin{subfigure}{0.495\textwidth}
\centering
\includegraphics[width=\textwidth]{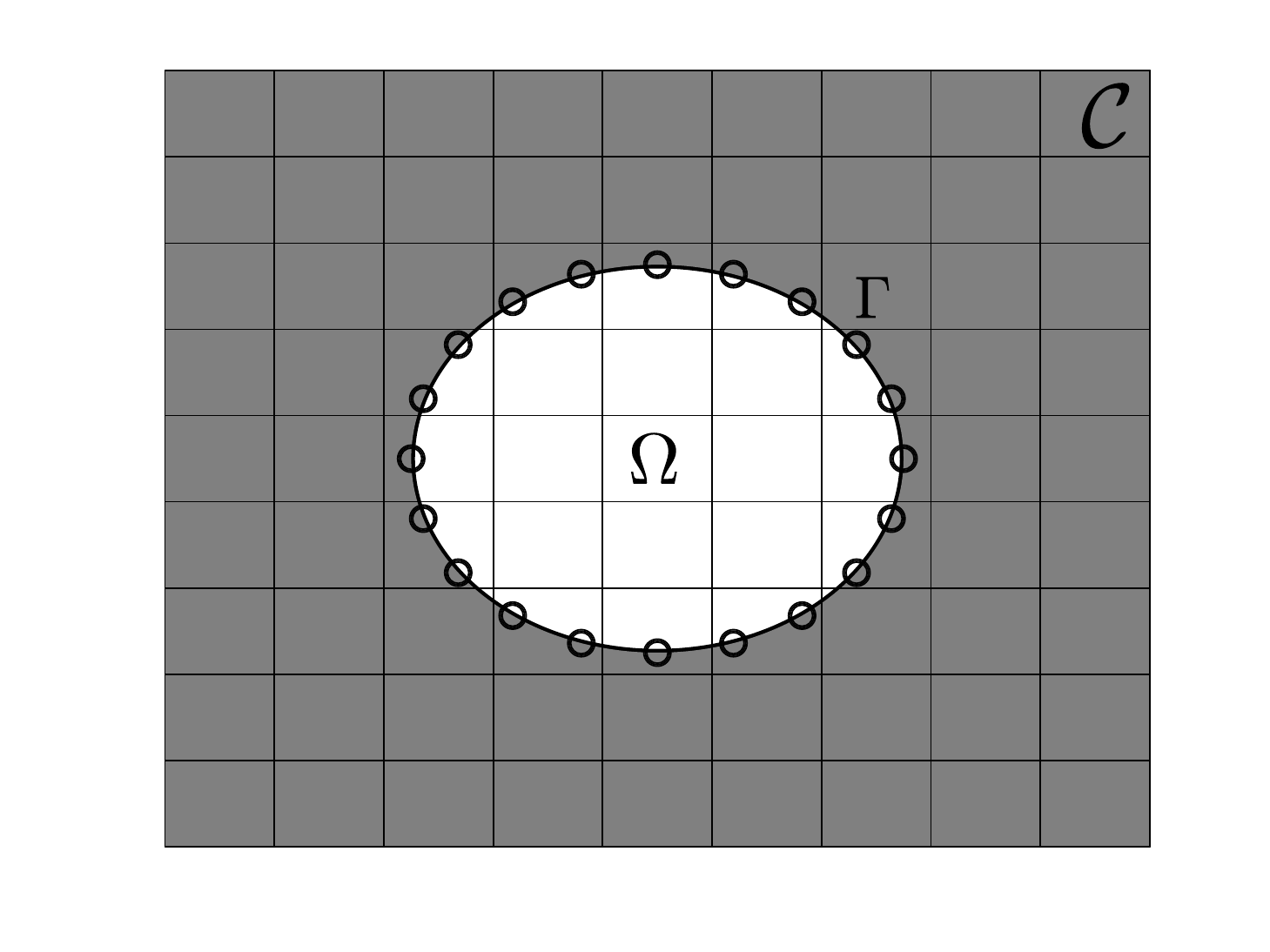}
\subcaption{\normalsize Interior domain}
\label{interior domain}
\end{subfigure}
\caption[Illustration of Eulerian and Lagrangian coordinate systems on two possible domains]{Diagram showing the two IB coordinate systems for an immersed structure $\Gamma$, embedded in a 2-D computational domain, $\C$. The PDE or physical domain, given by $\Omega$, is shown in white, and the non-physical domain is shown in grey. Figure \ref{exterior domain} shows an exterior PDE domain, and Figure \ref{interior domain} shows an interior PDE domain.}\label{IB domains}
\end{figure}

In the Immersed Boundary method, one embeds the PDE domain $\Omega$ into a geometrically simple computational domain $\C$. The method then utilizes a Lagrangian coordinate system located on the immersed boundary, as well as an Eulerian coordinate system for the fluid. Figure \ref{IB domains} illustrates simple discretizations of these coordinate systems for two possible PDE domains. 

In the original IB framework for interactions between fluid and elastic structures, the density of the force exerted by the structure on the fluid is represented by $F$, and it is supported on the boundary $\Gamma$. $F$ is mapped to the Eulerian coordinate system through a convolution with a delta function. The resulting Eulerian force density is incorporated into the PDE, which is then solved in $\C$.  With the elimination of the physical boundary, one can solve the discretized PDE on a regular grid using efficient solvers that are independent of the potentially complex geometry of $\Omega$. After the PDE is solved, the resulting solution $u$ can then be mapped back to the boundary through another convolution with a delta function, and this gives the boundary values.

This dissertation will focus on the IB method for PDEs with \textit{prescribed} boundary values, such as the PDE in Equation \eqref{pde}. In this so-called constraint formulation of the Immersed Boundary method, the force density $F$ is an unknown Lagrange multiplier, used to enforce the boundary condition. For scalar elliptic problems, we will focus on the Helmholtz operator, $\L\equiv \Delta -k^2$, where $k=0$ gives the Laplacian operator. We will use $s$ and $\x$ as our Lagrangian and Eulerian coordinates, respectively, and $\X(s,t)$ gives the location of the boundary $\Gamma$. The inclusion of time, $t$, allows for a moving boundary, but for brevity, we will generally assume a fixed boundary and omit $t$ from our equations. The equations for the continuous IB formulation of Equation \eqref{pde} are then given by
\begin{subequations} \label{constraint IB}
\begin{alignat}{2}
& \L \tilde u(\x) +\int_{\Gamma} F(s)\delta(\x-\X(s, t))ds = \tilde g(\x) \qquad && \text{in } \mathcal{C}  \label{constraint IB1}\\
&\int_{\C}\tilde  u(\x)\delta(\x-\X(s, t))d\x=U_b(s) \qquad && \text{on } \Gamma  . \label{constraint IB2}
\end{alignat}
\end{subequations}
In the above formulation, our unknowns are $F$ and $\tilde u$. Note that $\tilde u(\x)$ and $\tilde g (\x)$ represent extensions of the solution $u$ and the function $g$ from $\Omega$ to the larger computational domain $\C$. Once the solution $\tilde u(x)$ is found in $\C$, the solution to the original PDE is then given by $\tilde u|_{\Omega}$. As such, we will drop the tilde notation on $u(x)$. On the other hand, we can define the extension of $g$ by $\tilde g = g\chi_{\scaleto{\Omega}{4.5pt}} + g_e \chi_{\scaleto{\C \setminus \Omega}{6pt}}$, where $\chi$ is an indicator function. Since $\C\setminus\Omega$ is not in the PDE domain, there is flexibility in our choice for $g_e$. For example, one can simply use a smooth extension of $g$ or use $g_e=0$.

In the IB method, mappings between the coordinate systems are carried out with a regularized delta function $\delta_h$, where the regularization lengthscale, $h$, is generally chosen to be on the order of the discretized grid point spacing. We can then view the first convolution as \textit{spreading} the force density to the nearby grid points, and the second convolution as \textit{interpolating} the solution values onto the boundary points. Therefore, let us define the spread and interpolation operators, respectively, as
\begin{subequations} \label{operators}
\begin{alignat}{1}
&(SF)(\x) =\int_{\Gamma} F(s)\delta_h(\x-\X(s, t))ds  \label{spread operator}\\
&(S^*u)(s)=\int_{\C} u(\x)\delta_h(\x-\X(s, t))d\x.   \label{interpolation operator}
\end{alignat}
\end{subequations}
These operators are adjoint in the sense that 
\begin{equation}\label{adjoint}
\langle SF, u\rangle_{\C} = \langle F, S^*u\rangle_{\Gamma},
\end{equation}
where the inner products are the usual $L^2$ inner products on $\Omega$ and $\Gamma$, respectively. Discretization and numerical implementation will be discussed more in Section \ref{ch 2 numerical implementation}, but note that we will also use $S$ and $S^*$ to refer to the discretized versions of these operators. Similarly, let $\L$ represent our continuous or discretized elliptic differential operator. We then get the following system of equations. 
\begin{subequations} \label{constraint IB lL}
\begin{alignat}{2}
& \L u +SF = \tilde g \qquad && \text{in } \mathcal{C}  \label{constraint IB L1}\\
&S^* u = U_b \qquad && \text{on } \Gamma,  \label{constraint IB L2}
\end{alignat}
\end{subequations}
or
\begin{equation} \label{saddle point equation}
\begin{pmatrix}
\L & S\\
S^* & 0 
\end{pmatrix}
\begin{pmatrix}
u\\
F 
\end{pmatrix}= 
\begin{pmatrix}
\tilde g\\U_b
\end{pmatrix}.
\end{equation}
As indicated in the beginning of the chapter, we will refer to the method described by Equation \eqref{constraint IB lL} as the Immersed Boundary Single Layer (IBSL) method. 


\section{Numerical implementation}\label{ch 2 numerical implementation}

\subsection{Discretization of space and differential operators}\label{2.2 space and L}

Unless otherwise specified, we take the computational domain to be the two-dimensional periodic box, $\C=[-L/2, L/2]^2$, in order to make use of efficient finite difference and Fourier spectral methods to solve Equation \eqref{constraint IB L1}. The computational domain is discretized with a regular Cartesian mesh with $N$ points in each direction, giving us $\Delta x = \Delta y = L/N$. We then have $N^2$ nodes located at $\x_{i,j}=(-L/2+i\Delta x,-L/2+ j\Delta y)$ for integers $0\leq i,j \leq N-1$. We can easily generalize this to a periodic box with different horizontal and vertical lengths and node spacings, but we make these choices for simplicity. 

Since we will be looking at domains of various sizes, we rescale the discrete $L^1(\Omega)$ and $L^2(\Omega)$ norms by the area of $\Omega$, denoted $|\Omega|$. Our refinement studies therefore use the discrete function norms given by 
\begin{subequations}
\begin{alignat}{2}
&||u||_1 &&= \frac{\Delta x \Delta y}{|\Omega|} \sum_{i,j} u(\x_{i,j})\chi_{\scaleto{\Omega}{4.5pt}}\\
&||u||_2 &&= \Bigg(\frac{\Delta x \Delta y}{|\Omega|} \sum_{i,j} (u(\x_{i,j}))^2\chi_{\scaleto{\Omega}{4.5pt}}\Bigg)^{1/2}.
\end{alignat}
\end{subequations}

The IBSL method has flexibility in the choice for discretization of the differential operator. In this dissertation, we utilize both finite difference and Fourier spectral methods. In the case of finite differences, we will use the standard five-point, second-order accurate approximation for the Laplacian, given by
\begin{equation}\label{standard laplacian}
(\Delta u )_{i,j}\approx \frac{1}{\Delta x \Delta y}\Big(u_{i-1,j}+u_{i,j-1}+u_{i+1,j}+u_{i,j+1}-4u_{i,j}\Big),
\end{equation}
where $u_{i,j}=u(\x_{i,j})$. For Fourier spectral methods, the Laplacian is given by 
\begin{equation}
\Delta u(x)=-\mathcal{F}^{-1}(\mathbf{k}\dotp \mathbf{k} \mathcal{F}u(x)),
\end{equation}
where $\mathcal{F}$ is the Discrete Fast Fourier Transform, $\mathcal{F}^{-1}$ is its inverse, and $\mathbf{k}=(k_1,k_2)$ gives the vector of wave-numbers.

\subsection{Discretization of spread and interpolation operators}\label{2.2 spread}

To discretize the boundary $\Gamma$, we use a set of $N_{IB}$ boundary points, given by $\{X(s_i)\}_{i=1}^{N_{IB}}$. Since the locations of $\X(s_i)$ do not coincide with the grid points, we use the spread and interpolation operators, $S$ and $S^*$, to map between the Lagrangian and Eulerian grids. 

\begin{figure}
\centering
\includegraphics[width=0.5\textwidth]{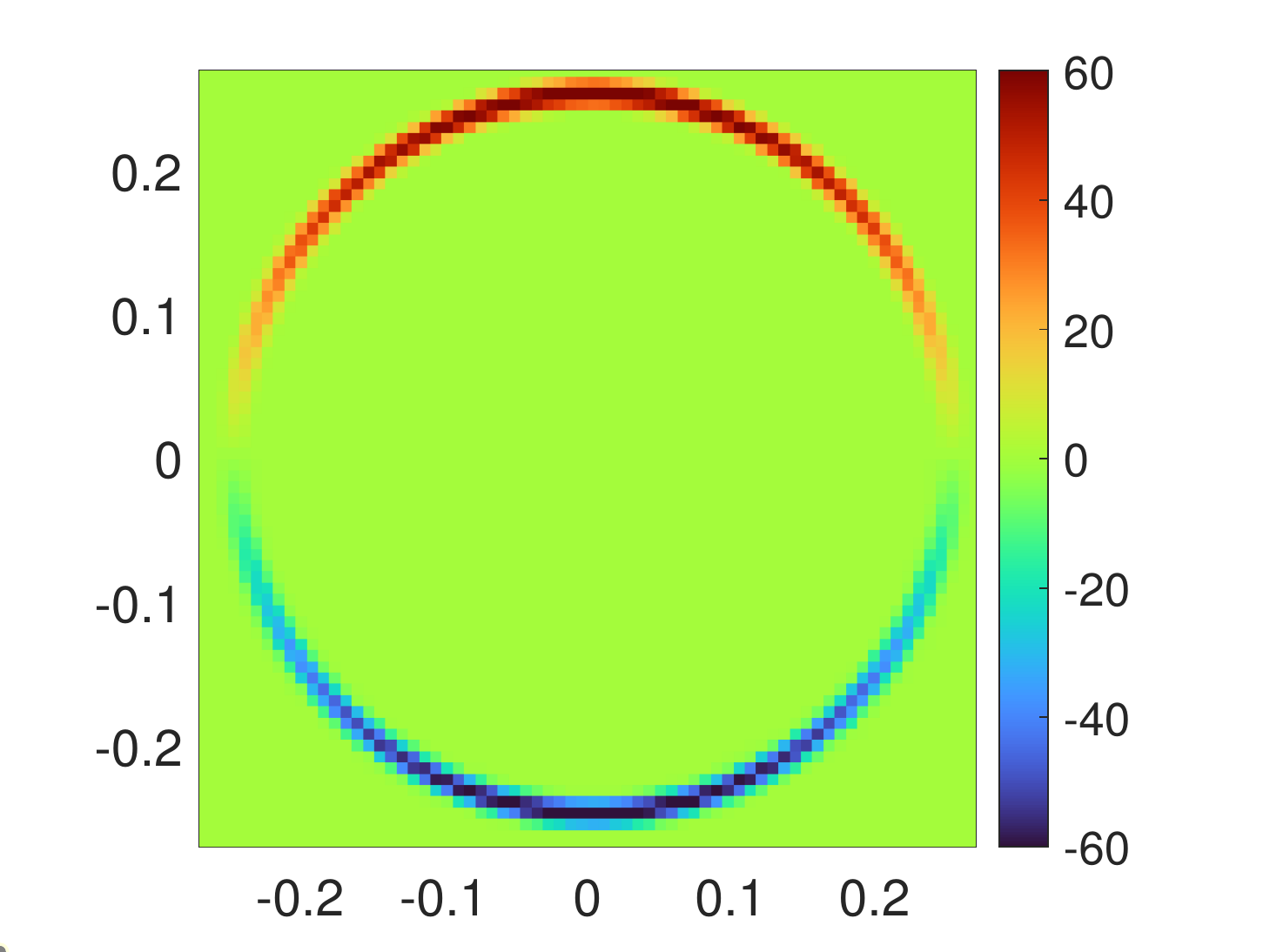}
\caption[Graphical illustration of the action of the discrete spread operator $S$]{Illustration of the spread operator defined in Equation \eqref{spread operator again2}. The boundary $\Gamma$ is a circle of radius 0.25 centered at the origin, and the spread operator is applied to the distribution $F=\sin{\theta}$. The boundary point spacing is given by $\Delta s \approx 2 \Delta x$.}\label{spread operator figure}
\end{figure}

The spread operator maps a singular force distribution supported on $\Gamma$ to an appropriate density on the Eulerian grid. It is formed as an approximate convolution with a regularized two-dimensional delta function. There is flexibility in the choice of the regularized delta function, but a suitable choice will satisfy certain properties \cite{Peskin02, liumori, ontheorder}. In this dissertation, unless otherwise specified, we use the traditional Peskin four-point delta function \cite{Peskin02}, which in two dimensions is given by 
\begin{equation}
\delta_h=\frac{1}{h^2} \phi\bigg(\frac{x}{h}\bigg)\phi\bigg(\frac{y}{h}\bigg),
\end{equation}
for $h=\Delta x = \Delta y$ and 
\begin{equation}
\phi (r)=\begin{cases}
 \frac{1}{8} \big( 3-2|r| + \sqrt{1+4|r|-4|r|^2}\big) & 0\leq |r|\leq 1 \\ 
      \frac{1}{8} \big( 5-2|r| - \sqrt{-7+12|r|-4|r|^2}\big) & 1\leq |r|\leq 2 \\
	0 & |r|\geq 2  .  \end{cases}
\end{equation}
Using this delta function, we define the continuous spread operator as
\begin{equation}
(SF)(\x) =\int_{\Gamma} F(s)\delta_h(\x-\X(s))ds .    \label{spread operator again}
\end{equation}
Then, discretizing the integral gives us the discrete operator,
\begin{equation}
(SF)(\x) =\sum_{i=1}^{N_{IB}} F(s_i)\delta_h(\x-\X(s_i))\Delta s_i .    \label{spread operator again2}
\end{equation}
A graphical illustration of the spread operator is given in Figure \ref{spread operator figure}. Note that this integral is respect to the parameter $s$. It is not necessary that $s$ be an arclength parameter, nor is it necessary that the points be equally spaced, either physically or in terms of the  parameter $s$. However, if we do space the points equally in one of these ways, we will get a spectrally accurate quadrature rule for a smooth periodic function $F$, defined on $\Gamma$. Unless otherwise specified, we will use arclength parametrization and equally spaced points with $\Delta s_i = L_{IB}/N_{IB}$, where $L_{IB}$ is the length of the immersed boundary. Additionally, we will select $N_{IB}$ such that $\Delta s \approx \alpha \Delta x$ for various values of $\alpha$. As will be discussed further in Section \ref{ch 2 results}, a commonly used point spacing for the IBSL method is $\Delta s \approx 2 \Delta x$ \cite{GriffithDonev}. If the exact values of $\Delta s_i$ are unavailable, approximations can be made. In the case of an arclength parametrization, using $\Delta s_i \approx ||\X(s_{i+1})-\X(s_i)||_2$ is sufficient. 

\squeezeup

We lastly look at the the interpolation operator, which maps a discrete function on the Eulerian grid to a discrete function supported on the Lagrangian grid. It is discretized as follows. 
\begin{equation}
(S^*u)(s)=\int_{\C} u(\x) \delta_h(\x-\X(s)) d\x \approx \sum_{x_{i,j}} u(\x_i)\delta_h(\x_i-\X(s)) (\Delta x)(\Delta y),  \label{interpolation again}
\end{equation}
where the sum is over the Cartesian gridpoints, and for a fixed $s$ and a four-point delta function, the sum has at most 16 non-zero terms. 
\subsection{Solution to discrete system for invertible $\mathcal{L}$}\label{2.2 invertible L}

We now look at the discretized system 
\begin{subequations} \label{constraint IB lL again}
\begin{alignat}{2}
& \L u +SF = \tilde g \qquad && \text{in } \mathcal{C}  \label{constraint IB L1 again}\\
&S^* u = U_b \qquad && \text{on } \Gamma .  \label{constraint IB L2 again}
\end{alignat}
\end{subequations}

In the case that the differential operator is invertible, such as for $\L=\Delta -k^2$, for $k\neq0$, we can invert the operator to obtain
\begin{equation}
u=-\L^{-1}SF+\L^{-1}\tilde g  .
\end{equation}
Then, by applying the interpolation operator and using Equation \eqref{constraint IB L2 again}, we obtain
\begin{equation}
-(S^* \L^{-1}S)F=U_b -S^* \L^{-1} \tilde g .  \label{schur complement again}
\end{equation}
We can therefore solve the saddle point problem in Equation \eqref{constraint IB lL again} by first solving Equation \eqref{schur complement again} for $F$ and then obtaining $u$ from Equation \eqref{constraint IB L1 again}. Here, the operator that must be inverted, $-(S^* \L^{-1}S)$, is called the Schur complement of the system. We do not form the Schur complement matrix explicitly but instead create a routine that performs the action of the operator, and we solve for $F$ with a Krylov method. Since the Schur complement is a symmetric operator, we use \texttt{minres}, to a tolerance of $10^{-8}$, unless otherwise specified.


\subsection{Solution to discrete system for Poisson equation}\label{2.2 Laplacian L}

In the case that the differential operator is the periodic Laplacian, $\L=\Delta$, the solution method outlined in the previous section must be adjusted to account for the nullspace of $\L$. The boundary value problem given by 
\begin{subequations} \label{laplaces pde}
\begin{alignat}{2}
& \Delta u = g \qquad && \text{in } \Omega  \label{l pde1}\\
&u=U_b \qquad && \text{on } \Gamma  \label{l pde2}
\end{alignat}
\end{subequations}
has a unique solution in $\Omega$, but, by using the IB framework, we embed the PDE into a periodic computational domain on which $\Delta$ is not invertible. Stein et al. \cite{Stein} present a process for solving Equation \eqref{constraint IB lL again} in this case. Since we will require a similar method when applying the IBDL method to Stokes equations, we present a version of this method here. 

We first decompose the solution $u$ as 
\begin{equation}
u=u_0+\bar u \label{u0+ubar}, 
\end{equation}
where $u_0$ has mean $0$ on $\C$ and $\bar u = \frac{1}{|\C|}\int_{\C} u d\x $ is a constant giving the mean value of $u$ on the computational domain. Note that a constant function spans the nullspace of the periodic Laplacian.

The PDE on the computational domain $\C$ that results from the IBSL formulation is then given by
\begin{equation}
\Delta u_0 + SF  =\tilde g, \label{laplaces with ib}
\end{equation}
where we have used that $\Delta \bar u  =0$. To derive the solvability condition, we integrate Equation \eqref{laplaces with ib} over a general computational domain $\C$. Let us use $\partial \C$ as the boundary of $\C$ and $\n_{\C}$ as the unit normal on $\partial \C$, to distinguish it from $\n$, which we will later use as the unit normal on the immersed boundary $\Gamma$.  We note that since $S$ is defined using a regularized delta function, our functions are smooth, and we can use the divergence theorem. This gives us 
\begin{equation}
\int_{\partial \C} \grad u_0 \dotp \n_{\C} dl(\x) + \int_{\C} S F d\x = \int_{\C} \tilde g d\x,  \label{first step after divergence theorem}
\end{equation}
where we use $dl(\x)$ to denote an integral with respect to arclength. If we take $\C$ to be the periodic box used in this work, the first term disappears, and we get the following solvability condition for the periodic computational domain.
\begin{equation}
\int_{\C}S Fd\x=\int_{\C} \tilde g d\x  \label{solvability constraint}
\end{equation}

Equation \eqref{solvability constraint} then provides the constraint needed to solve for the additional unknown, $\bar u$. We can then discretize the previous equations and summarize the discrete IBSL system as  
\begin{subequations} \label{constraint IB poisson}
\begin{alignat}{1}
&\Delta u_0 +SF = \tilde g  \label{constraint IB L1 poisson}\\
&S^* u_0+S^*\bar u \mathds{1}_{N^2} = U_b   \label{constraint IB L2 poisson}\\
&(\Delta x\Delta y) \mathds{1}_{N^2}^\intercal S F = (\Delta x\Delta y) \mathds{1}_{N^2}^\intercal \tilde g  , \label{constraint IB L3 poisson}
\end{alignat}
\end{subequations}
where $\mathds{1}_{N^2}$ denotes a vector of length $N^2$ consisting of all ones. Since we have ensured that the solvability constraint is satisfied, the solution to Equation \eqref{constraint IB L1 poisson} can now be found, and the solution is unique up to an additive constant. We then let $\Delta_0^{-1}$ denote the operation that inverts the Laplacian by returning a solution with mean $0$ on $\C$. Then we have $\Delta_0^{-1}\Delta u_0=u_0$. Using this to invert the Laplacian in  Equation \eqref{constraint IB L1 poisson} and then applying $S^*$, we get
\begin{equation}
S^* u_0+S^*\Delta_0^{-1}SF=S^*\Delta_0^{-1} \tilde g.
\end{equation}
Then, we can replace $S^* u_0$ using Equation \eqref{constraint IB L2 poisson}, and we get
\begin{equation}
U_b-S^*\bar u \mathds{1}_{N^2}+S^*\Delta_0^{-1}SF=S^*\Delta_0^{-1} \tilde g.
\end{equation}
We then have the following system with unknowns, $F$ and $\bar u$:
\begin{equation} \label{poisson schur system}
\begin{pmatrix}
-S^*\Delta_0^{-1}S & S^*\mathds{1}_{N^2}\\
(\Delta x\Delta y) \mathds{1}_{N^2}^\intercal S & 0 
\end{pmatrix}
\begin{pmatrix}
F\\
\bar u 
\end{pmatrix}= 
\begin{pmatrix}
U_b-S^*\Delta_0^{-1}\tilde g\\(\Delta x\Delta y) \mathds{1}^\intercal_{N^2} \tilde g
\end{pmatrix}.
\end{equation}

We can then use properties of the spread and interpolation operators to simplify this system. First, we note that the interpolation of a constant function will return a constant function on $\Gamma$, giving us
\begin{equation}
S^* \bar u  \mathds{1}_{N^2} = \bar u \mathds{1}_{N_{IB}} , \label{interp of constant}
\end{equation}
where $\mathds{1}_{N_{IB}}$ is a vector of length $N_{IB}$ consisting of all ones. We can use this to change the operator acting on $\bar u$ in Equation \eqref{poisson schur system}. Next, letting $u=1$ in the adjoint property given in Equation \eqref{adjoint}, we see that
\begin{equation}
\int_{\C} SF d\x = \int_{\Gamma} F ds,
\end{equation}
or, in discretized terms, 
\begin{equation}
(\Delta x \Delta y) \mathds{1}_{N^2}^\intercal SF = \Delta s \mathds{1}_{N_{IB}}^\intercal F.
\end{equation}
The system can now be written as
\begin{equation} \label{poisson schur system again}
\begin{pmatrix}
-S^*\Delta_0^{-1}S & \mathds{1}_{N_{IB}}\\
\Delta s\mathds{1}_{N_{IB}}^\intercal & 0 
\end{pmatrix}
\begin{pmatrix}
F\\
\bar u 
\end{pmatrix}= 
\begin{pmatrix}
U_b-S^*\Delta_0^{-1}\tilde g\\(\Delta x\Delta y)\mathds{1}^\intercal_{N^2} \tilde g
\end{pmatrix}.
\end{equation} 
Therefore, we can make this system symmetric by dividing the second equation by $\Delta s$, and our system becomes
\begin{equation} \label{poisson schur system 2}
\begin{pmatrix}
-S^*\Delta_0^{-1}S & \mathds{1}_{N_{IB}}\\
 \mathds{1}_{N_{IB}}^\intercal  & 0 
\end{pmatrix}
\begin{pmatrix}
F\\
\bar u 
\end{pmatrix}= 
\begin{pmatrix}
U_b-S^*\Delta_0^{-1}\tilde g\\\frac{\Delta x \Delta y}{\Delta s} \mathds{1}^\intercal \tilde g
\end{pmatrix}.
\end{equation}

Since the operator is symmetric, we again solve for $F$ and $\bar u$ using \texttt{minres}. We can then obtain $u_0$ from Equation \eqref{constraint IB L1 poisson} and add $\bar u$ to obtain our final solution, 
\begin{equation}
u=-\Delta_0^{-1} SF+\Delta_0^{-1}\tilde g +\bar u.
\end{equation}


\section{Mathematical description of IBSL method: Brinkman and Stokes equations}\label{ch 2 stokes}

In this section, we apply the Immersed Boundary Single Layer method to boundary value problems of the form 
\begin{subequations} \label{stokes and brinkman pde}
\begin{alignat}{2}
& \mu \Delta \u -k^2 \u -\grad p= \g \qquad && \text{in } \Omega  \label{s and b pde 1}\\
& \grad \dotp \u = 0 \qquad && \text{in } \Omega  \label{s and b pde 2}\\
&\u=\U_b \qquad && \text{on } \Gamma,  \label{s and b pde 3}
\end{alignat}
\end{subequations}
where in this dissertation, we will take $\u$ and $\g$ to be two-dimensional vector-valued functions. For $k=0$, this is Stokes equation for viscous flow with zero Reynolds number, and $k\neq 0 $ gives the Brinkman equation, which can be used for problems involving porous media \cite{Brinkman}. The Brinkman equation will also arise in our implicit-explicit discretization of the Navier-Stokes equation, discussed in Chapter \ref{chapter 6}. The system of equations that results from applying the IBSL method to Equation \eqref{stokes and brinkman pde} is 
\begin{subequations} \label{ibsl brinkman 1}
\begin{alignat}{2}
&\mu  \Delta \u -k^2 \u - \grad p +S\F =  \gt \qquad && \text{in } \mathcal{C}  \label{ibsl brinkman 1 1}\\
& \grad \dotp \u =0 \qquad && \text{in } \mathcal{C}  \label{ibsl brinkman 1 2}\\
&S^* \u = \U_b \qquad && \text{on } \Gamma.  \label{ibsl brinkman 1 3}
\end{alignat}
\end{subequations} 
We can see that the IBSL method generalizes in a straightforward way to the Stokes equation. Note that $\F$ is now a two-dimensional force density on the boundary $\Gamma$. Additionally, the spread and interpolation operators act element-wise on the discrete vector-valued functions. For example, $S\F=(S\F_1, S\F_2) $. In the case that $k=0$, the nullspace of the periodic Laplacian will again require a more careful consideration similar to that seen in Section \ref{2.2 Laplacian L}. We will therefore handle the two cases separately following a brief description of the discretization of the differential operators.

\subsection{Discretization of differential operators}\label{2.3 space and op}

The discretization of space and the spread and interpolation operators is identical to that described in Section \ref{ch 2 numerical implementation}. We again use both finite difference and Fourier spectral methods for discretizing the differential operators. For Fourier spectral methods and two dimensions, we compute the derivative with multi-index $\alpha = (\alpha_1, \alpha_2)$ as 
\begin{equation}
u^{\alpha}(\x)= \mathcal{F}^{-1}\Bigg( i^{|\alpha |} k_1^{\alpha_1} k_2^{\alpha_2}\mathcal{F}u(x)\Bigg), 
\end{equation}
where $\mathcal{F}$ and $\mathcal{F}^{-1}$ give the Discrete Fast Fourier Transform and its inverse, and $\mathbf{k}=(k_1,k_2)$ gives the vector of wave-numbers. 

For finite difference methods, we use centered, second-order accurate approximations for first derivatives, and we again use the standard five-point, second-order accurate approximation for the Laplacian. However, as will be described in the following section, in order to solve for pressure, we will need to invert the operator given by $\grad \dotp \grad$, which will not correspond to the standard five-point Laplacian given in Equation \eqref{standard laplacian}. It instead corresponds to a Laplacian with a wider stencil, given by 
\begin{equation}
(\grad \dotp \grad u)_{i,j} \approx \frac{1}{4\Delta x \Delta y}\Big(u_{i-2,j}+u_{i,j-2}+u_{i+2,j}+u_{i,j+2}-4u_{i,j}\Big),
\end{equation}
This will result in an operator with a four-dimensional nullspace of constants on each of 4 subgrids. Therefore, in the following section, when we say that the mean of the pressure will be fixed at $0$, in the case of this finite difference method, the mean will be fixed at $0$ on each of the 4 subgrids. 

A staggered-grid, or MAC, discretization is a common way to avoid the de-coupling of these 4 subgrids, as well as to improve volume conservation in Immersed Boundary fluid applications \cite{Hornung, GriffithBhalla, Griffithaortic, firstMAC, Griffithvolume}. In this dissertation, we do not make use of a staggered grid in order to maintain consistency between our implementation of the IBSL and IBDL methods, and as will be discussed further in Section \ref{6.4 space}, generalizing this discretization to the IBDL method is non-trivial and merits future exploration.

\subsection{Solution to discrete system for Brinkman equation}\label{2.3 Brinkman}

The IBSL formulation of the Brinkman equation with Dirichlet boundary conditions is again given by
\begin{subequations} \label{ibsl brinkman}
\begin{alignat}{2}
& \mu \Delta \u -k^2 \u - \grad p +S\F =  \gt \qquad && \text{in } \mathcal{C}  \label{ibsl brinkman 1}\\
& \grad \dotp \u =0 \qquad && \text{in } \mathcal{C}  \label{ibsl brinkman 2}\\
&S^* \u = \U_b \qquad && \text{on } \Gamma,  \label{ibsl brinkman 3}
\end{alignat}
\end{subequations}
or
\begin{equation} \label{saddle point equation}
\begin{pmatrix}
\L & -\grad & S\\
\grad \dotp & 0 & 0\\
S^* & 0 & 0\\
\end{pmatrix}
\begin{pmatrix}
\u \\
p \\
\F
\end{pmatrix}= 
\begin{pmatrix}
 \gt\\0\\ \U_b
\end{pmatrix}, 
\end{equation}
where $\L= \mu \Delta -k^2 \mathds{I} $, and $\mathds{I} $ is the identity operator. 

Taking the divergence of Equation \eqref{ibsl brinkman 1} and using the incompressibility of $\u$, we get
\begin{equation}
-\Delta p + \Div S\F =  \Div \gt. \label{div of eqn 1}
\end{equation}
Equation \eqref{div of eqn 1} is solvable on the periodic domain, and we invert the periodic Laplacian with the operator $\Delta_0^{-1}$, which was described in Section \ref{2.2 Laplacian L}. We therefore fix its mean value on $\C$ to be $0$. Completing this inversion, we get
\begin{equation}
p=-\Delta_0^{-1} \Div (\gt-S\F). 
\end{equation}
Using this expression for $p$ in Equation \eqref{ibsl brinkman 1}, we get
\begin{equation}
\L \u =(\gt - S\F)-\grad \Delta_0^{-1}\Div (\gt - S\F).\label{after p is in}
\end{equation}
Let us denote the operator that projects onto divergence-free fields as $\mathbb{P}=\mathds{I} - \grad \Delta_0^{-1} \Div$. Equation \eqref{after p is in} then becomes
\begin{equation}
\L \u =\mathbb{P} (\gt - S\F).\label{after p is in 2}
\end{equation}
Inverting $\L$ and applying the interpolation operator $S^*$, we get the following equation for $\F$:
\begin{equation}\label{brink schur}
-\Big(S^*\L^{-1}\mathbb{P}S\Big)\F=\U_b-S^*\L^{-1}\mathbb{P}\gt.
\end{equation}
As we did in Section \ref{2.2 invertible L}, we solve this equation for $\F$ using \texttt{minres}, and then we solve Equations \eqref{ibsl brinkman 1} and \eqref{ibsl brinkman 2} for the velocity $\u$ and pressure $p$. 

\subsection{Solution to discrete system for Stokes equation}\label{2.3 stokes}

In the case of Stokes equation, for which $k=0$, we must again adjust the method of solution since the periodic Laplacian is not invertible. We can follow a process similar to that in Section \ref{2.2 Laplacian L} to find the solvability constraint for Equation \eqref{ibsl brinkman 1}. We again begin by decomposing the solution $\u$ as 
\begin{equation}
\u=\u_0+\ub,
\end{equation}
where $\u_0$ has mean $0$ on $\C$, and $\ub=(\bar u, \bar v)$ gives the mean values of the horizontal and vertical components of velocity, $u$ and $v$, respectively. By taking $\L=\mu \Delta$ and integrating Equation \eqref{after p is in} over a general computational domain $\C$, we get
\begin{equation}
\mu \int_{\C} \Delta \u_0 d\x = \int_{\C} \Big(\gt - S\F - \grad \Delta_0^{-1} \Div (\gt - S\F) \Big)d \x.
\end{equation}
Using the divergence theorem where applicable, we get
\begin{equation}
\mu \int_{\partial \C} \grad \u_0 \dotp \n_{\C} dl(\x) = \int_{\C} \gt d\x  -\int_{\C} S \F d\x  -  \grad\Delta_0^{-1}\int_{\partial \C} \gt\dotp \n_{\C} dl(\x)+ \grad\Delta_0^{-1} \int_{\partial \C}S\F\dotp \n_{\C} dl(\x)  .  \label{first step after divergence theorem}
\end{equation}
If we assume $\Gamma$ is away from the boundary of $\C$, the last integral on the right-hand-side vanishes  due to the compact support of the integrand. If we take $\C$ to be a periodic box, the left-hand-side also vanishes. If the PDE domain $\Omega$ is an exterior domain, then $\gt$ on $\partial \C$ is $\g$, the periodic forcing function in Equation \eqref{s and b pde 1}. If, on the other hand, $\Omega$ is an interior domain, since we have flexibility in choosing our extended function $\gt$, we can let $\gt\chi_{\scaleto{\C \setminus \Omega}{6pt}}=\g_e=0$. In either of these cases, the third integral on the right-hand-side vanishes. Using that $\int_{\C} S\F d\x=\int_{\Gamma} \F ds$, we again get the solvability constraint
\begin{equation}
\int_{\Gamma} \F ds = \int_{\C} \gt d\x.
\end{equation}
Therefore, using a process similar to that of Section \ref{2.2 Laplacian L}, we obtain the following system for unknowns $\F$ and $\ub$.
\begin{equation}\label{stokes schur}
\begin{pmatrix}
-\frac{1}{\mu}S^*\Delta_0^{-1}\mathbb{P}S & \mathds{1}_{N_{IB}}\\
 \mathds{1}_{N_{IB}}^\intercal  & 0 
\end{pmatrix}
\begin{pmatrix}
\F\\
\ub 
\end{pmatrix}= 
\begin{pmatrix}
\U_b-\frac{1}{\mu}S^*\Delta_0^{-1}\mathbb{P}\gt \\ \frac{(\Delta x)^2}{\Delta s} \mathds{1}^\intercal \gt
\end{pmatrix}.
\end{equation}
We again solve this system using \texttt{minres}, and after obtaining $\F$ and $\ub$, we get $\u$ by solving Equation \eqref{ibsl brinkman 1} for $\u_0$ and then adding $\ub$ to obtain our final solution, 
\begin{equation}
\u=\frac{1}{\mu}\Delta_0^{-1}\mathbb{P}\Big(\gt-S\F\Big)+\ub.
\end{equation}


\section{IBSL method: results}\label{ch 2 results}

In this section, we provide some numerical results for the Immersed Boundary Single Layer method presented in this chapter. We draw particular attention to the high number of iterations of \texttt{minres} required to solve the systems found in Equations \eqref{schur complement again}, \eqref{poisson schur system}, \eqref{brink schur}, and \eqref{stokes schur}. These saddle point problems resulting from Dirichlet boundary conditions are poorly conditioned \cite{Benzisaddlepoint}. Therefore, as mentioned in the introduction, much work has been devoted to forming preconditioners \cite{GriffithDonev, Ceniceros, guyphilipgriffith, Stein}. In this section, we will see that without proper preconditioning, the iteration counts can become extremely large. With an eye toward time-dependent PDE, such as the Navier-Stokes equation, in which the saddle point system would be solved in each time step, this can be prohibitively time consuming. Additionally, preconditioning itself can be computationally expensive, as it often involves inverting a dense matrix. This poor conditioning is one of the main motivations for the development of the Immersed Boundary Double Layer (IBDL) method, which is able to achieve the same order of accuracy as the IBSL method while avoiding these high iteration counts and therefore eliminating the need for preconditioning. This method will be introduced for scalar elliptic PDEs in Chapter \ref{chapter 5} and for Stokes equation in Chapter \ref{chapter 6}. Comparable tests will be performed with the IBDL method to compare to the results presented here.


\subsection{Helmholtz equation}\label{2.4 Poisson}

Our first test problem is the PDE 
\begin{subequations} \label{hh pde}
\begin{alignat}{2}
& \Delta u -  u = 0 \qquad && \text{in } \Omega  \label{hh pde 1}\\
&u=\sin{2\theta} \qquad && \text{on } \Gamma,  \label{hh pde 2}
\end{alignat}
\end{subequations}
where $\Omega$ is the interior of a circle of radius 0.25, centered at the origin. The analytical solution is given by 
\begin{equation}
u=\frac{I_2(r)\sin{2\theta}}{I_2(0.25)}, 
\end{equation}
where $I_2$ is the first-kind modified Bessel function of order 2. Our computational domain here is the periodic box $[-0.5, 0.5]^2$, and we use equally spaced boundary points with $\Delta s \approx \alpha \Delta x$ for various values of $\alpha$. We use a finite difference method, and the solutions are computed for grid sizes ranging from $N=2^6$ to $2^{12}$. 

\begin{figure}
\centering
\includegraphics[width=0.5\textwidth]{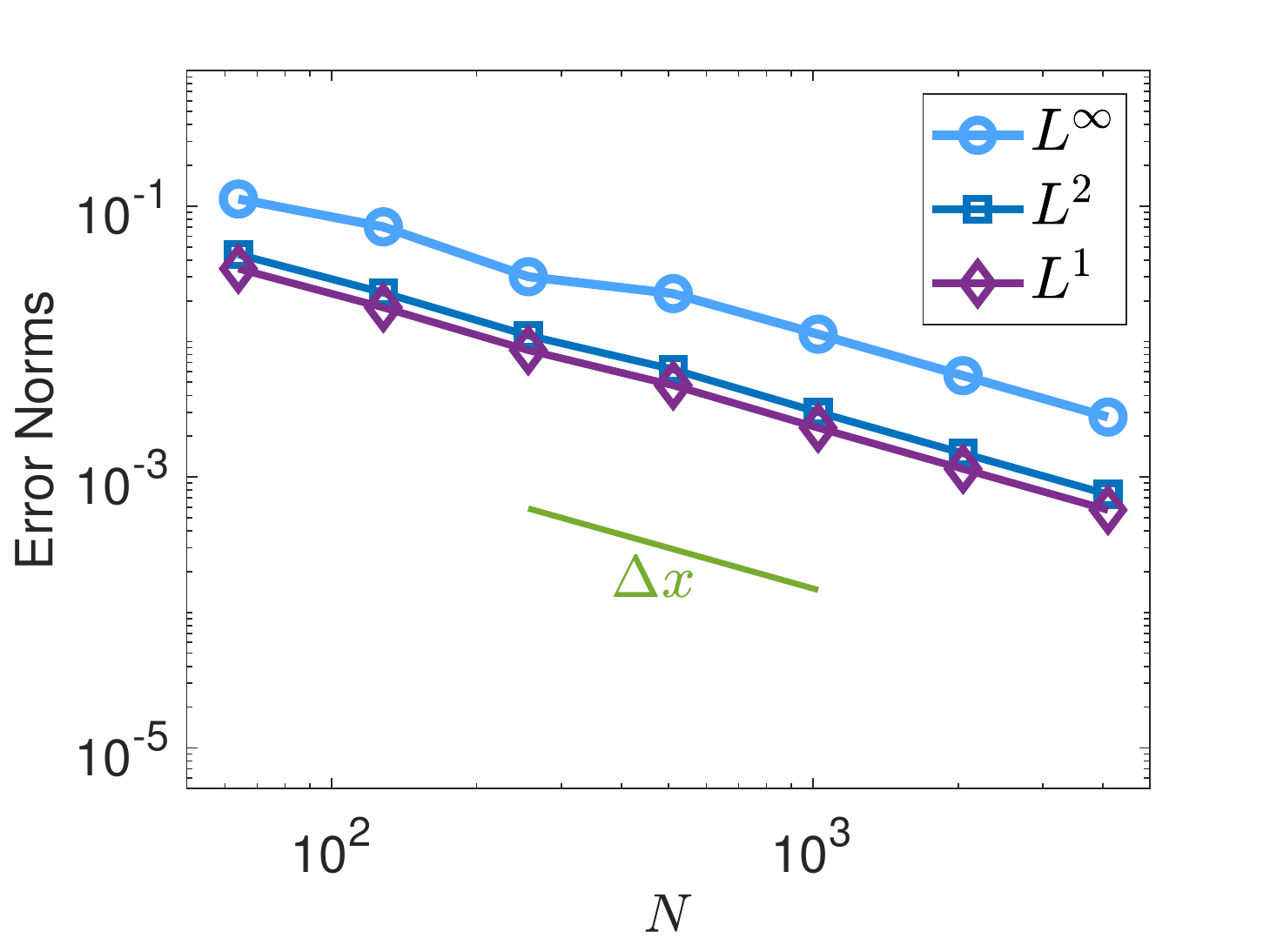}
\caption[Refinement study for solutions to the Helmholtz equation \eqref{hh pde} using the Immersed Boundary Single Layer (IBSL) method]{Refinement studies for solutions to Equation \eqref{hh pde} found using the IBSL method with finite differences. The computational domain is the periodic box $[-0.5, 0.5]^2$, $\Omega$ is the interior of a circle of radius 0.25, and the prescribed boundary values are given by $U_b=\sin{2\theta}$. The boundary point spacing is $ \Delta s  \approx 0.75\Delta x$, and the plot shows the $L^1$, $L^2$, and $L^{\infty}$ norms of the absolute errors. In green, we denote the slope corresponding to first-order convergence.}\label{IBSL interior refinement 1}
\end{figure}

Table \ref{iteration table 1} gives the iteration counts for using \texttt{minres}, to a tolerance of $10^{-8}$, without preconditioning, to solve Equation \eqref{hh pde} with the IBSL method. We can see that the the iteration counts increase as we refine our grid spacing. These counts would be prohibitively large in the case of a time-dependent problem, in which this solve would be require at each time step. The source of these large iteration counts will be discussed in Chapters \ref{chapter 3} and \ref{chapter 4}. We also see the iteration counts increase as we refine the spacing of our immersed boundary points relative to the grid. For this reason, in order to make the IBSL method practical, it is common to use $\Delta s \approx 2 \Delta x$ and to make use of preconditioning \cite{Ceniceros, GriffithDonev, guyphilipgriffith, Stein}. Figure \ref{IBSL interior refinement 1} demonstrates the first-order convergence of the solution, which is characteristic of the IBSL method. In Chapter \ref{chapter 5}, we present a method that maintains this order of convergence while avoiding these high iteration counts.

\begin{figure}
\begin{center}
 \begin{tabular}{||c | c | c | c | c ||} 
 \hline
 \multicolumn{5}{||c||}{Iteration Counts - Circular Boundary} \\
 \hline
  $\Delta x$ &$\Delta s \approx 2\Delta x $& $\Delta s \approx  1.5 \Delta x$ &$ \Delta s \approx  1 \Delta x$& $\Delta s \approx  0.75 \Delta x$  \\
 \hline
$2^{-6}$ &   17  &    45   &  249  &   491  \\
$2^{-7} $&   36   &   57    &  691  &  752  \\
$2^{-8 }$&   49  &    71   & 1233   &  3057      \\
$2^{-9}$&  61    &    110   & 1922   &   5829      \\
$ 2^{-10}$& 68   &    144   & 1936   &    8084    \\
$ 2^{-11}$ & 95  &  234   &  4364 &   9335  \\
$ 2^{-12}$  & 142 &303     &  4535  &   10589    \\
  \hline
  \end{tabular}
\captionof{table}[Number of iterations of \texttt{minres} to solve the Helmholtz equation \eqref{hh pde} using the Immersed Boundary Single Layer (IBSL) method]{Number of iterations of \texttt{minres}, with tolerance $10^{-8}$, needed to solve the Helmholtz equation \eqref{hh pde}  with the IBSL method, using finite differences, on the periodic computational domain $[-0.5, 0.5]^2$, where the circular boundary of radius 0.25 has prescribed boundary values given by $U_b=\sin{2\theta}$.  }
\label{iteration table 1}
\end{center}
\end{figure}

\subsection{Brinkman equation} \label{2.4 Stokes}

Our second test problem is the Brinkman PDE 
\begin{subequations} \label{brink pde}
\begin{alignat}{2}
& \Delta \u -  \u -\grad p = \g \qquad && \text{in } \Omega  \label{brink pde 1}\\
& \Div \u=0 \qquad && \text{in } \Omega  \label{brink pde 2}\\
&u=\U_b \qquad && \text{on } \Gamma.  \label{brink pde 3}
\end{alignat}
\end{subequations}
We use an analytical solution given by
\begin{subequations} \label{brink soln}
\begin{alignat}{2}
& u=e^{\sin{x}}\cos{y} \\
& v=-\cos{x}\hspace{0.1cm} e^{\sin{x}}\sin{y}\\
& p= e^{\cos{y}}\label{p eqn}
\end{alignat}
\end{subequations}
to determine the boundary values, $\U_b=\u|_{\Gamma}$, and obtain the forcing function given by 
\begin{equation} \label{brink g}
\g =\begin{pmatrix} e^{\sin{x}}\cos{y}\hspace{0.1cm} (\cos^2{x}-\sin{x}-2)\\
-\cos{x} \hspace{0.1cm} e^{\sin{x}}\sin{y}\hspace{0.1cm}(\cos^2{x}-3\sin{x}-3)+\sin{y}\hspace{0.1cm}e^{\cos{y}}\end{pmatrix}.
\end{equation}
Here $\Omega$ is the interior of a circle of radius 0.75, centered at the origin, and the computational domain is the periodic box, $\C=[-1, 1]^2$. We use equally spaced boundary points with $\Delta s \approx \alpha \Delta x$ for various values of $\alpha$. We use a Fourier spectral method, and the solutions are computed for grid sizes ranging from $N=2^6$ to $2^{12}$. 

\begin{figure}
\centering
\begin{subfigure}{0.5\textwidth}
\centering
\includegraphics[width=\textwidth]{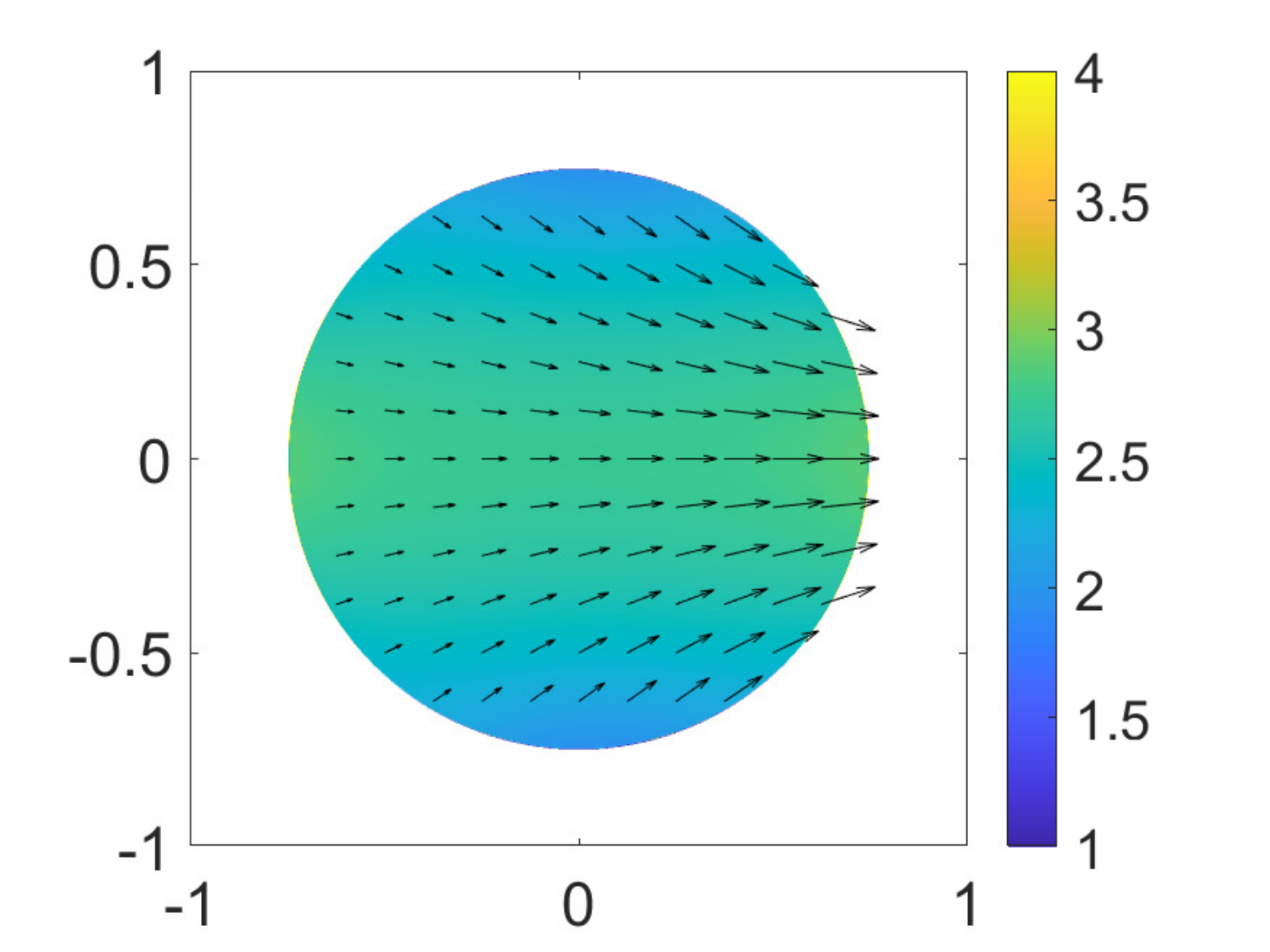}
\caption{\normalsize  Solution}
\label{IBSL brink plot}S
\end{subfigure}

\begin{subfigure}{0.395\textwidth}
\centering
\includegraphics[width=\textwidth]{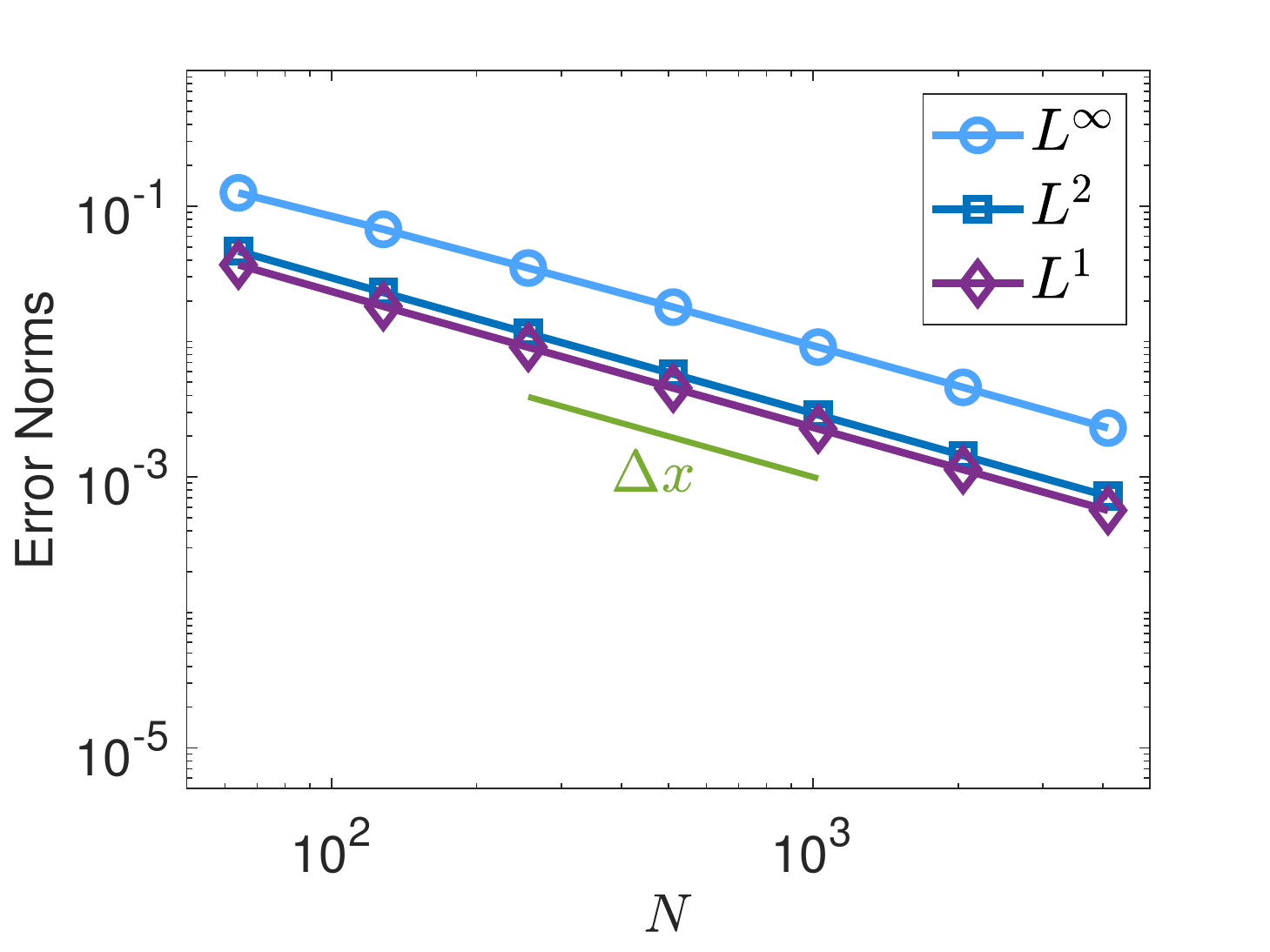}
\caption{\normalsize  $u$}
\label{IBSL brink refinement u}
\end{subfigure}
\begin{subfigure}{0.395\textwidth}
\centering
\includegraphics[width=\textwidth]{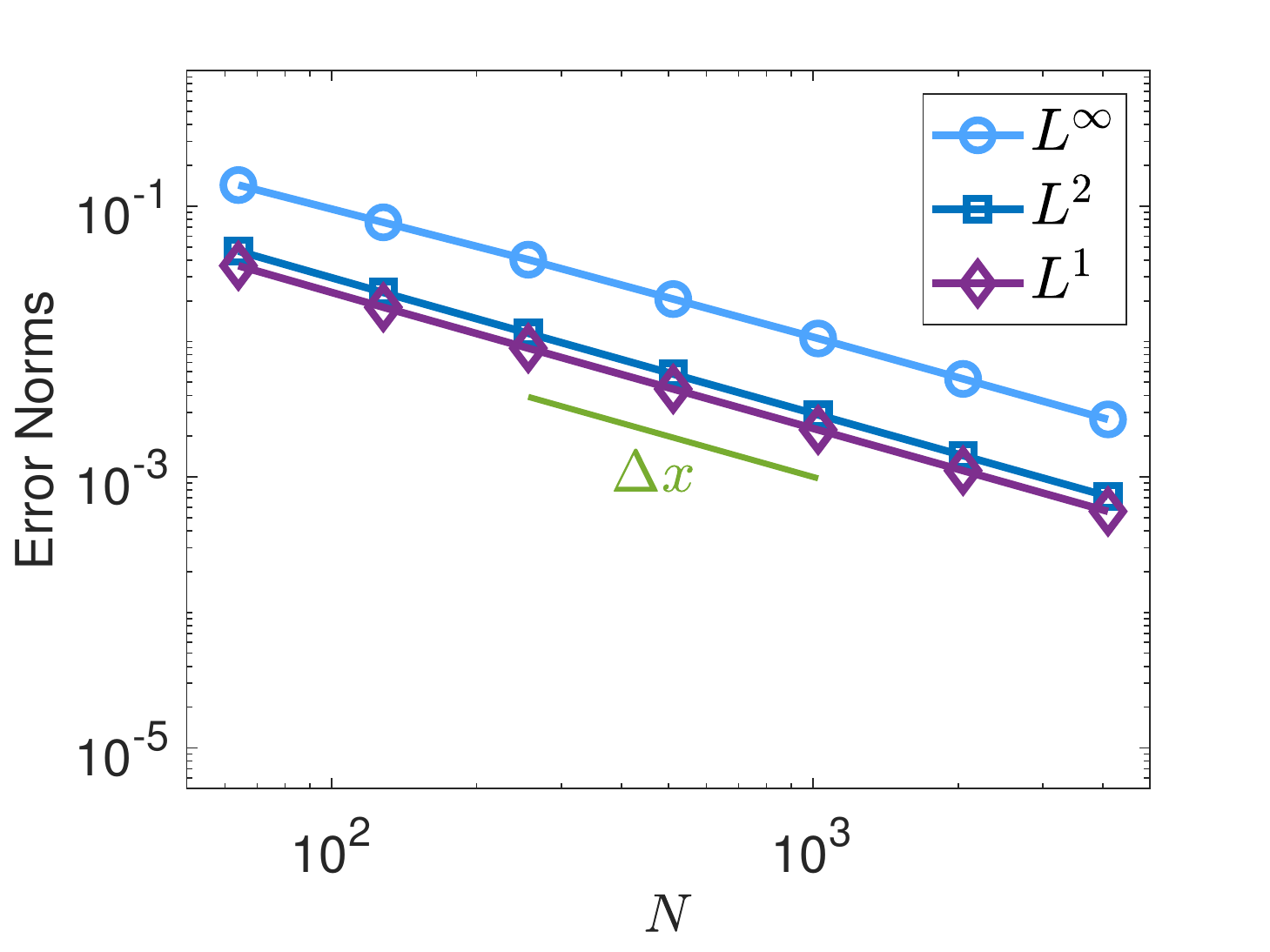}
\caption{\normalsize  $v$}
\label{IBSL brink refinement v}
\end{subfigure}
\caption[Solution plot and refinement studies for solution to the Brinkman equation \eqref{brink pde} using the IBSL method]{Solution plot and refinement studies, showing absolute errors,  for solution to Equation \eqref{brink pde} found using the IBSL method. The computational domain is the periodic box $[-1, 1]^2$, and the boundary point spacing is $ \Delta s  \approx 1.5\Delta x$. Figure \ref{IBSL brink plot} shows the velocity solution with vectors and pressure solution with color. It is found using the grid size $N=2^{10}$.}\label{IBSL interior refinement}
\end{figure}

Table \ref{iteration table 2} gives the iteration counts for using \texttt{minres}, without preconditioning, to solve Equation \eqref{brink pde} with the IBSL method. Note that for $\Delta s\approx 0.75 \Delta x$, \texttt{minres} stagnated without reaching the tolerance of $10^{-8}$, and the iteration counts in the last column give the points at which the method stagnated. We can again see that the iteration counts increase as we refine our grid spacing or the spacing of our immersed boundary points relative to the grid. Figure \ref{IBSL brink plot} provides a plot of the solution on $\Omega$ for $N=2^{10}$ and $\Delta s=1.5 \Delta x$. The velocity vectors are plotted every $2^6$ meshwidths, and the color shows the pressure. As discussed in Section \ref{2.3 Brinkman}, the pressure is unique only up to an additive constant. We therefore solve for the pressure that has mean $0$ on $\C$. We then find the difference between the computed pressure and analytical solution in Equation \eqref{p eqn} at the origin, and we add this constant to the computed pressure for the plot. Figures \ref{IBSL brink refinement u}-\ref{IBSL brink refinement v} demonstrate the first-order convergence of the velocity solutions. We return to this problem in Chapter \ref{chapter 6} to compare these results to those of the IBDL method.

\begin{figure}
\begin{center}
 \begin{tabular}{||c | c | c | c | c ||} 
 \hline
 \multicolumn{5}{||c||}{Iteration Counts - Circular Boundary} \\
 \hline
  $\Delta x$ &$\Delta s \approx 2\Delta x $& $\Delta s \approx  1.5 \Delta x$ &$ \Delta s \approx  1 \Delta x$& $\Delta s \approx  0.75 \Delta x$  \\ [0.5ex] 
 \hline
$2^{-5}$ & 123    &   296    & 2558   &   16190 \hspace{0.1cm} stag. \\
$2^{-6} $&    185  & 472      &  5578  &  47863 \hspace{0.1cm} stag. \\
$2^{-7 }$& 257    &   678    &  14438  &  95009  \hspace{0.1cm} stag.    \\
$2^{-8}$&  351    &  825     &  17475  &    68861  \hspace{0.1cm} stag.   \\
$ 2^{-9}$& 478   &   1137    &  17789  &    62697  \hspace{0.1cm} stag.  \\
$ 2^{-10}$ & 622  &   1620  &  20333 &  75910\hspace{0.1cm} stag.  \\
$ 2^{-11}$  & 821&   2180  &  27467 &    79699 \hspace{0.1cm} stag. \\
  \hline
  \end{tabular}
\captionof{table}[Number of iterations of \texttt{minres} to solve the Brinkman equation \eqref{brink pde} using the IBSL method]{Number of iterations of \texttt{minres}, with tolerance $10^{-8}$, needed to solve Equation \eqref{brink pde} with the IBSL method, using a Fourier spectral method, on the periodic computational domain $[-1, 1]^2$. The word \textit{stag.} is used to indicate that the method stagnated at this number of iterations without converging to the desired tolerance.  }
\label{iteration table 2}
\end{center}
\end{figure}

\section{Discussion}\label{ch 2 discussion}

The Immersed Boundary Single Layer method is extremely robust and flexible. One can use the discretization and PDE solver of their choice without needing it to conform to the geometry of the boundary. It also has some key advantages over the boundary integral methods that will be discussed in the next chapter. For example, it can be easily adapted to various exterior domains, such as the periodic domain we utilize in this dissertation. It can also be easily used with nonhomogeneous PDE without greatly increasing the computing work required, and the solution can be calculated on an entire domain mesh very quickly. Additionally, this method can be used for Navier-Stokes and other nonlinear problems by treating the nonlinearity explicitly in time. For these reasons, we strive to maintain the overall Immersed Boundary framework in our new method. However, we seek to improve on two areas. Firstly, the IBSL method cannot be used for PDEs with Neumann boundary conditions, and secondly, the conditioning of the Schur complement can make this method inefficient, especially in time-dependent problems. Therefore, we seek to remedy these two areas by reformulating the method using concepts from boundary integral equations.

   \chapter{Boundary integral equations}
   \label{chapter 3}


Boundary integral methods can also be used to solve the PDEs in Chapter \ref{chapter 2} by reformulating boundary value problems as integral equations using the appropriate Green's functions. In this chapter, we give an introduction to these boundary integral equations and the relevant properties that we will use when analyzing the IBSL method and developing the IBDL method. We begin by looking at boundary integral equations for Laplace's equation, and then we do the same for Stokes equation. Much of the material in this chapter can be found in boundary integral textbooks, such as \cite{Pozred} and \cite {Pozblue}, but we present it here in order to establish notation and motivation for the form of the IBDL method.

\section{Green's functions}\label{ch 3 Greens}

We begin by looking at Laplace's equation, 
\begin{equation}\label{laplaces equation}
\Delta u = 0  \hspace{0.6cm}\text{in } \Omega . 
\end{equation}
A Green's function for this PDE is a function $G(\x,\x_0$) that 
satisfies\begin{equation}
\Delta G(\x,\x_0) = - \delta (\x-\x_0) \hspace{0.6cm} \text{in } \Omega ,   \label{greens function}
\end{equation}
where $\x_0$ is the location of a singular point force. One must also specify boundary conditions to determine a particular Green's function. For example, the free-space Green's function for an infinite solution domain in two-dimensions is
\begin{equation}
G^{FS}(\x, \x_0)=-\frac{1}{2\pi}\ln{|\x-\x_0|}. \label{fs greens}
\end{equation} 
This Green's function would be sufficient for solving boundary value problems for interior PDE domains. However, the exterior domains explored in this dissertation are not in free-space, but in a periodic box. However, we can use this function to analyze behavior near the singular point because all Green's functions for Laplace's equation exhibit the same singular behavior as $\x\longrightarrow \x_0$. In fact, any Green's function can be decomposed as
\begin{equation}
G(\x,\x_0)= G^{FS}(\x,\x_0) + H^C(\x,\x_0), \label{decompose G}
\end{equation}
where $H^C(\x,\x_0)$ is a complementary harmonic function that is non-singular on $\Omega$ \cite{Pozblue}.  

We next derive an integral property of the Green's function derivative that will be important for obtaining our integral equations.  For convenience, let our problem domain $\Omega$ be an interior domain, unless otherwise specified, with boundary $\Gamma$. In this chapter, when we omit explicit discussion of an exterior domain, the steps and conclusions follow similarly as for the interior domain. Taking $\x_0 \notin \Omega$, we can integrate Equation \eqref{greens function} about $\Omega$, and use the divergence theorem to get
\begin{equation}
0=\int_{\Omega} \Delta G(\x,\x_0) d\x = \int_{\Gamma} \grad G(\x,\x_0) \dotp \n(\x) dl(\x) , \label{aaa1}
\end{equation}
where $dl(\x)$ denotes an integral with respect to arclength, and $\n$ is the unit normal pointing out of $\Omega$, regardless of whether $\Omega$ is an exterior or interior domain. 

\begin{figure}
\centering
\begin{subfigure}{0.495\textwidth}
\centering
\includegraphics[width=\textwidth]{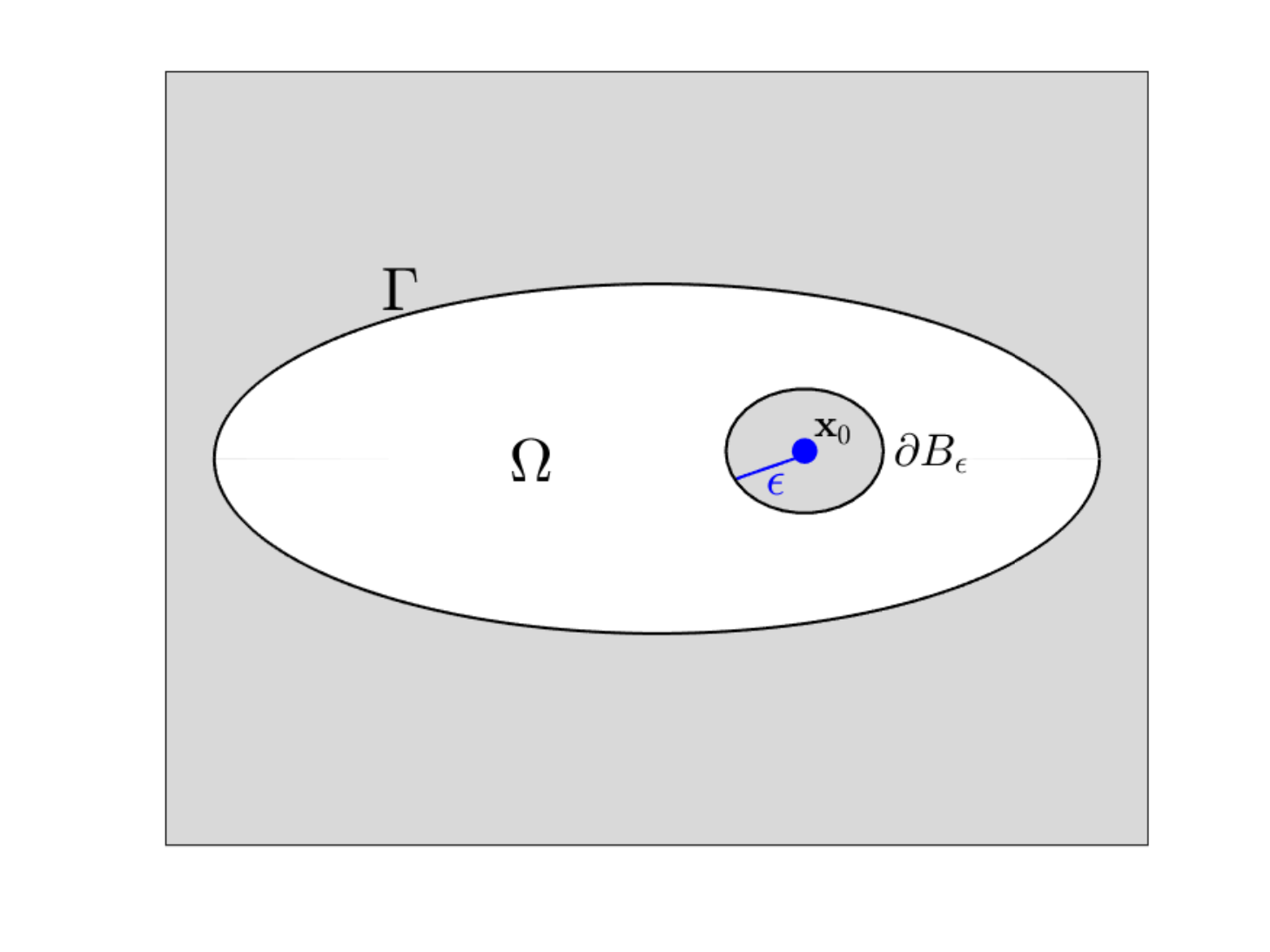}
\subcaption{\normalsize $\x_0\in \Omega$}
\label{inside}
\end{subfigure}
\begin{subfigure}{0.495\textwidth}
\centering
\includegraphics[width=\textwidth]{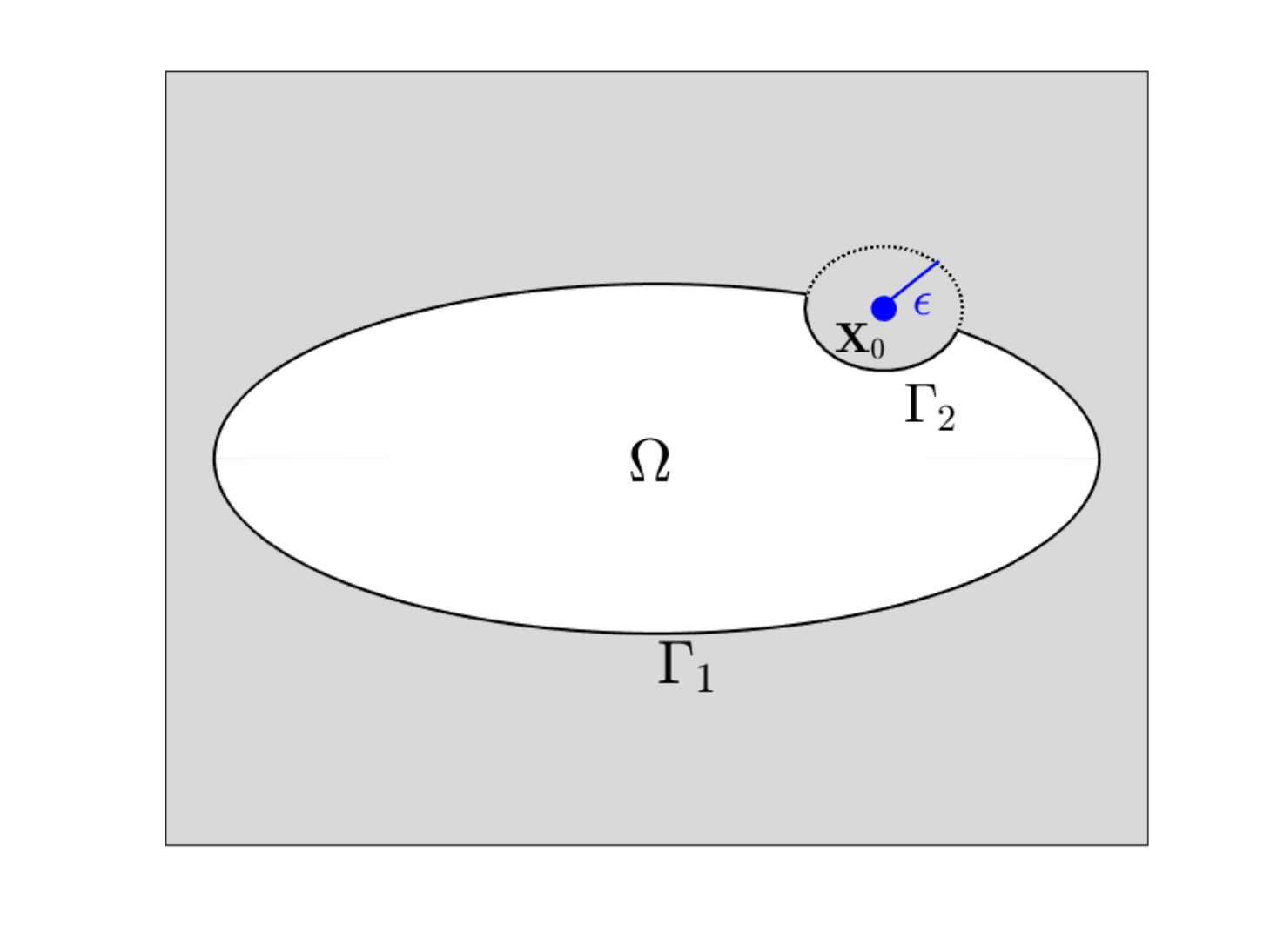}
\subcaption{\normalsize $\x_0 \in \Gamma$}
\label{on bound}
\end{subfigure}
\caption[Altered PDE domains to exclude singular point in two cases]{Figure \ref{inside} illustrates the removal of ball $B_{\epsilon}(\x_0)$ from the PDE domain $\Omega$ for $\x_0\in\Omega$ in order to use the divergence theorem in Equation \eqref{first inside one}. The resulting domain, $\Omega'$, excludes $\x_0$ and has the boundary $\Gamma\cup \partial B_{\epsilon}$. Figure \ref{on bound} illustrates the removal of $\Omega \cap B_{\epsilon}(\X_0)$ from $\Omega$ for $\X_0 \in \Gamma$ The resulting domain, $\Omega'$, excludes $\X_0$ and has the boundary $\Gamma_1\cup \Gamma_2$.}\label{greens pic}
\end{figure}

If, on the other hand, $\x_0\in \Omega$, since the Green's function is not continuously differentiable at $\x=\x_0$, in order to use the divergence theorem, we can remove a circle centered at $\x_0$ with sufficiently small radius $\epsilon$ and use the divergence theorem on the region $\Omega'=\Omega\setminus B_{\epsilon}(\x_0)$. Figure \ref{inside} gives an illustration. Using the decomposition of $G$ in Equation \eqref{decompose G}, we then have 
\begin{equation}
0=\int_{\Omega'} \Delta G(\x,\x_0) d\x = \int_{\Omega'} \Delta G^{FS}(\x,\x_0) d\x + \int_{\Omega'} \Delta H^C(\x,\x_0) d\x = \int_{\Omega'} \Delta G^{FS}(\x,\x_0) d\x ,
\end{equation}
where we have used that $H^C$ is a harmonic function. Using the divergence theorem then gives us 
\begin{equation}
0=\int_{\Gamma} \grad G^{FS}(\x,\x_0) \dotp \n(\x) dl(\x) +\int_{\partial B_{\epsilon}} \grad G^{FS}(\x,\x_0) \dotp \n(\x) dl(\x), \label{first inside one}
\end{equation}
where both unit normal vectors are those that point out of $\Omega'$. The free-space Green's function derivative is
\begin{equation}
\grad G^{FS}(\x,\x_0)=-\frac{1}{2\pi |\x-\x_0|^2}\begin{bmatrix} x-x_0 \\ y-y_0\end{bmatrix}.
\end{equation}
Additionally, the unit normal on $\partial B_{\epsilon}$ is 
\begin{equation}
\n=-\frac{\x-\x_0}{|\x-\x_0|}.
\end{equation}
Then the second integral is 
\begin{equation}
\int_{\partial B_{\epsilon}} \grad G^{FS}(\x,\x_0) \dotp \n(\x) dl(\x) =\frac{1}{2\pi}\int_{\partial B_{\epsilon}} \frac{|\x-\x_0|^2}{|\x-\x_0|^3} dl(\x) =  \frac{1}{2\pi\epsilon}\int_{\partial B_{\epsilon}} dl(\x) = 1. 
\end{equation}
Then Equation \eqref{first inside one} becomes
\begin{equation}
-1 = \int_{\Gamma} \grad G^{FS}(\x,\x_0) \dotp \n(\x) dl(\x) .
\end{equation}
In order to find this property for the general Green's function, note that, using that $H^C$ is harmonic and non-singular on $\Omega$, we have
\begin{equation}
0=\int_{\Omega} \Delta H^C(\x,\x_0)d\x= \int_{\Gamma} \grad H^C(\x,\x_0)\dotp \n(\x) dl(\x). 
\end{equation}
Therefore, we have 
\begin{equation}
\int_{\Gamma}\grad G(\x,\x_0) \dotp \n(\x)dl(\x) = \int_{\Gamma}\grad G^{FS}(\x,\x_0) \dotp \n(\x)dl(\x) + \int_{\Gamma}\grad H^C(\x,\x_0) \dotp \n(\x)dl(\x)= -1, \label{inside one}
\end{equation}
for $\x_0 \in \Omega$. 

We next take the singularity to be on the boundary $\Gamma$. Since in the IBSL method, we denoted points on the boundary with a capital $\X(s)$, in order to maintain consistency and also distinguish from the case when $\x_0$ is not on the boundary, we will call this singularity $\X_0$. In order to use the divergence theorem, we then remove from $\Omega$ a circle centered at $\X_0$ with radius $\epsilon$ and again call this new region $\Omega'$. An illustration of this is given in Figure \ref{on bound}. Let $\Gamma_1$ denote the region of $\Gamma$ that borders $\Omega'$ (i.e. the portion of $\Gamma$ that excludes the region near $\X_0$), and let $\Gamma_2=\partial \Omega'\setminus \Gamma_1$. Then, since $\X_0 \notin \Omega'$, we can use the divergence theorem to get 
\begin{equation}
0=\int_{\Gamma_1 \cup \Gamma_2} \grad G(\x, \X_0)\dotp \n(\x) dl(\x) = \int_{\Gamma_1 } \grad G(\x, \X_0)\dotp \n(\x) dl(\x) +\int_{ \Gamma_2} \grad G(\x, \X_0)\dotp \n(\x) dl(\x) .
\end{equation}
We can again use the decomposition of $G$ into the free-space Green's function and the complementary harmonic function in order to use the free-space Green's function derivative. For brevity, we omit the details as they follow the same process as used for $\x_0\in\Omega$. We can therefore use the form of the free-space Green's function and note that since $\Gamma$ is smooth, as $\epsilon \longrightarrow 0$, the second integral approaches the integral over a semicircle of radius $\epsilon$. Then, this integral simplifies to $1/2$ , and we get
\begin{equation}
-\frac12=\int_{\Gamma_1} \grad G(\x, \X_0)\dotp \n(\x) dl(\x) \label{gam1int}.
\end{equation}
Let us define the principal-value integral, denoted $PV$, as the integral for which $\X_0$ is placed exactly on $\Gamma$. It is computed by integrating over $\Gamma \setminus B_{\epsilon}(\X_0)$, in which a disk around $\X_0$ has been removed from the boundary. The $PV$ integral then comes from taking the limit as $\epsilon$ goes to 0. If we take $\epsilon\longrightarrow 0$ for the definition of $\Gamma_1$ in Equation \eqref{gam1int}, we get exactly this principal-value integral, so we have
\begin{equation}
-\frac{1}{2}=\int_{\Gamma}^{PV} \grad G(\x, \X_0)\dotp \n(\x) dl(\x).\label{aaa3}
\end{equation}
Therefore, to summarize Equations \eqref{aaa1}, \eqref{inside one}, and \eqref{aaa3}, we can give the following integral identity for the Green's function derivative: 
\begin{equation} \label{greens deriv prop 1}
\int_{\Gamma}\grad G(\x, \x_0)\dotp \n(\x) dl(\x) = 
    \begin{cases}
       -1, & \text{if } \x_0\in \Omega\\
       -\frac12, & \text{if } \x_0\in \Gamma\\
       0, & \text{if } \x_0 \notin \bar \Omega
        \end{cases}, 
\end{equation}
where $\bar \Omega=\Omega \cup \Gamma$. We can also rewrite Equation \eqref{greens deriv prop 1} as 
\begin{equation} \label{greens deriv prop 2}
\int_{\Gamma} \grad G(\x, \x_0)\dotp \n(\x) dl(\x) = 
    \begin{cases}
       \int_{\Gamma}^{PV} \grad G(\x, \X_0)\dotp \n(\x) dl(\x)-\frac12 & \text{if } \x_0\in \Omega\\
      \int_{\Gamma}^{PV} \grad G(\x, \X_0)\dotp \n(\x) dl(\x)+\frac12 & \text{if } \x_0\notin \bar \Omega ,
        \end{cases}
\end{equation}
and this version of the property will be used in Section \ref{ch 3 integral eqns}.

\section{Integral representations for solutions to Laplace's equation}\label{ch 3 reps}

In this section, we derive integral representations for solutions to Equation \eqref{laplaces equation}. We begin with Green's second identity. When applied to the Green's function and a solution to the PDE, $u$, the identity simplifies to
\begin{equation}
u(\x) \delta(\x-\x_0)= \grad \dotp \Big(G(\x, \x_0) \grad u(\x) - u(\x) \grad G(\x,\x_0)\Big). \label{greens second identity}
\end{equation} 
Let $\x_0 \in \Omega $. Then, in order to use the divergence theorem, we can again remove from $\Omega$ a ball of radius $\epsilon$, centered at the singularity. By following similar steps to those in the previous section, we can then get an integral representation of $u$ given by
\begin{equation}
u(\x_0) = \int_{\Gamma} G(\x, \x_0) \grad u(\x) \dotp \n (\x) dl(\x) - \int_{\Gamma} u(\x) \grad G(\x,\x_0) \dotp \n(\x) dl(\x), \label{integral rep}
\end{equation} 
where $\n$ is the unit normal of $\Gamma$, pointing out of $\Omega$. 

The first term of Equation \eqref{integral rep} is called a \emph{single layer potential}, which can be viewed as a distribution of point forces on the boundary, and the second term is called a \emph{double layer potential}, which can be viewed as a distribution of point force dipoles on the boundary. Notice that the strength of the single layer potential is given by the boundary distribution of the normal derivative of $u$, and the strength of the double layer potential is given by the boundary values of $u$. 

There are, however, other integral representations of the solution $u$ in $\Omega$. We can, for instance, represent $u$ solely with a single layer potential or a double layer potential. Let us consider two solutions: $u^{Int}$ that satisfies the PDE on $\Omega$ and $u^{Ext}$ that satisfies the PDE on $\Omega^{Ext}$, the region exterior to $\Gamma$. Let these two solutions share the same values on $\Gamma$. In order to find a single layer representation for $u^{Int}(\x_0)$ for $\x_0\in\Omega$ we need two equations. Firstly, we use Equation \eqref{integral rep} for $u^{Int}$. Secondly, by the same process that gave us Equation \eqref{integral rep}, we can get a similar equation for $u^{Ext}$, where the corresponding normal vector points outward from $\Omega^{Ext}$, and the left-hand side vanishes because $\x_0\notin \Omega^{Ext}$. By negating the expression in order to use the same unit normal as that in Equation \eqref{integral rep}, we get the following equation for $u^{Ext}$:
\begin{equation}
0= -\int_{\Gamma} G(\x, \x_0) \grad u^{Ext}(\x) \dotp \n(\x) dl(\x) + \int_{\Gamma} u^{Ext}(\x) \grad G(\x,\x_0) \dotp \n(\x) dl(\x). \label{integ rep 2}
\end{equation} 
Then, by adding Equation \eqref{integral rep}, for $u^{Int}$, and Equation \eqref{integ rep 2},  we get the \emph{generalized single layer integral representation} for $u^{Int}$ given by 
\begin{equation}
u^{Int} (\x_0) = \int_{\Gamma}  \sigma (\x)  G(\x, \x_0)dl(\x)  ,        \label{single layer rep}
\end{equation}
where the strength of the single layer potential is given by the jump in normal derivatives across the boundary, or
\begin{equation}
\sigma = \grad(u^{Int}-u^{Ext}) \dotp \n \Big|_{\Gamma}.
\end{equation} 

If we consider two solutions, $u^{Int}$ and $u^{Ext}$ that, instead of matching values on the boundary, match normal derivatives on the boundary, we get the \emph{generalized double layer integral representation} for $u^{Int}$ given by 
\begin{equation}
u^{Int} (\x_0) = \int_{\Gamma} \gamma(\x) \grad G(\x, \x_0)\dotp \n(\x) dl(\x)  ,        \label{double layer rep}
\end{equation}
where the strength of the double layer potential is given by the jump in solution values across the boundary, or
\begin{equation}
\gamma = (u^{Ext}-u^{Int})\Big|_{\Gamma}.
\end{equation}

\section{Integral equations for Laplace's equation with Dirichlet boundary conditions}\label{ch 3 integral eqns}

In this section, we introduce the integral equations that arise from using the integral representations from the previous section to solve the Dirichlet Laplace's boundary value problem given by
\begin{subequations} \label{laplaces bvp}
\begin{alignat}{2}
& \Delta u = 0\qquad && \text{in } \Omega  \label{laplaces bvp 1}\\
&u=U_b \qquad && \text{on } \Gamma.  \label{laplaces bvp 2}
\end{alignat}
\end{subequations}

In order to use our boundary data, we must take $\x_0 \longrightarrow \Gamma$, and we will again denote a point on the boundary as $\X_0$. The single layer potential seen in Equation \eqref{single layer rep} is continuous as $\x_0$ approaches and crosses the boundary \cite{Pozblue}, but the double layer potential is not. We therefore seek the limit of Equation \eqref{double layer rep} for $\x_0 \longrightarrow \Gamma$. We can first split up the integral using the potential strength at $\x_0$, or
\begin{equation}
\begin{split}
&\lim_{\x_0\rightarrow \Gamma} \int_{\Gamma}\gamma(\x)\grad G(\x, \x_0)\dotp \n(\x) dl(\x) \\
&=\lim_{\x_0\rightarrow \Gamma}\Bigg( \int_{\Gamma}\Big(\gamma(\x)-\gamma(\x_0)\Big)\grad G(\x, \x_0)\dotp \n(\x) dl(\x) +\gamma(\x_0) \int_{\Gamma} \grad G(\x, \x_0)\dotp \n(\x) dl(\x) \Bigg).
\end{split}
\end{equation}
Then we can use the identity given in Equation \eqref{greens deriv prop 2} to rewrite this as
\begin{multline}
\lim_{\x_0\rightarrow \Gamma}\Bigg( \int_{\Gamma}\Big(\gamma(\x)-\gamma(\x_0)\Big)\grad G(\x, \x_0)\dotp \n(\x) dl(\x)\\
+\gamma(\x_0)\Bigg( \int_{\Gamma}^{PV} \grad G(\x, \X_0)\dotp \n(\x) dl(\x)-\frac12\Bigg) \Bigg).
\end{multline}
Finally, we distribute, take the limit, and recombine the integrals to obtain
\begin{equation}
\lim_{\x_0\rightarrow \Gamma} \int_{\Gamma}\gamma(\x)\grad G(\x, \x_0)\dotp \n(\x) dl(\x) = \int_{\Gamma}^{PV}\gamma(\x)\grad G(\x, \X_0)\dotp \n(\x) dl(\x) -\frac{1}{2}\gamma(\X_0). \label{dl limit}
\end{equation}

Using the limit given in Equation \eqref{dl limit} for the double layer representation and the boundary condition in Equation \eqref{laplaces bvp 2}, the single and double layer integral equations are given by
\begin{equation}
U_b(\X_0) = \int_{\Gamma}  \sigma (\x)  G(\x, \X_0)dl(\x)  ,      \label{single layer eqn}
\end{equation}
\vspace{-0.3cm}
\begin{equation}
U_b (\X_0) = \int_{\Gamma}^{PV} \gamma(\x) \grad G(\x, \X_0)\dotp \n(\x) dl(\x) - \frac12 \gamma(\X_0), \label{double layer eqn}
\end{equation}
where the unknown quantities in Equations \eqref{single layer eqn} and \eqref{double layer eqn} are the potential strengths, $\sigma$ and $\gamma$, respectively. Note that the integral equations for a Dirichlet Helmholtz boundary value problem are the same as those for Laplace's equation, with the exception of the Green's function.

\section{Green's functions for Stokes equation}\label{ch 3 Stokes Greens}
In this section, we look at Stokes equation for incompressible viscous flow, 
\begin{subequations}\label{homog stokes pde}
\begin{alignat}{2}
& \mu \Delta \u -\grad p=\grad \dotp \sig = 0\qquad && \text{in } \Omega  \label{homog stokes pde 1}\\
&\Div \u =0  \qquad &&  \text{in } \Omega , \label{homog stokes pde 2}
\end{alignat}
\end{subequations}
where $\u$ is velocity, $p$ is pressure, $\sig$ is the stress tensor, and $\mu$ is the viscosity of the fluid. We will again focus on the case where $\u$ is a two-dimensional vector-valued function. In this case, we define the Green's function to be the solution to the singularly forced Stokes equation, 
\begin{subequations}\label{greens homog stokes pde}
\begin{alignat}{2}
& \mu \Delta \u^{G} -\grad p^{G}=\Div \sig^G= -\delta(\x-\x_0)\mathbf{c} \qquad && \text{in } \Omega  \label{greens homog stokes pde 1}\\
&\Div \u^{G} =0  \qquad &&  \text{in } \Omega , \label{greens homog stokes pde 2}
\end{alignat}
\end{subequations}
where $\mathbf{c}$ is a 2-D constant vector, and we use the superscript $G$ to distinguish the singularly forced solution from a general solution to Stokes equation. Then, letting $\G$ be the Green's function and $\p$ and $\T$ be the pressure vector and the stress tensor corresponding to the Green's function, we can write the solution to Equation \eqref{greens homog stokes pde} as 
\begin{subequations}\label{greens funct}
\begin{alignat}{2}
& u_i(\x) =\frac{1}{\mu}G_{ij}(\x,\x_0)c_j \label{greens funct 1}\\
&p(\x)=p_i(\x,\x_0)c_j \label{greens funct 2}\\
&\sigma_{ik}(\x)=T_{ijk}(\x,\x_0)c_j \label{greens funct 3}.
\end{alignat}
\end{subequations}
Note that we use Einstein notation to indicate summation over repeated indices. We will use this notation often for sections on Stokes equation. As discussed in Section \ref{ch 3 Greens}, the Green's function depends on the boundary conditions enforced on the solution domain, but the free-space Green's function provides the form of any such Green's function near the singular point, so we can use it to establish the properties of the Green's function. In free-space, the Green's function and corresponding stress tensor are called the Stokeslet and Stresslet, respectively. If we use $\hat{\x}=(\x-\x_0)$, then in two dimensions, the Stokeslet, Stresslet, and correponsing pressure vector are given by 
\begin{subequations}\label{fs greens funct}
\begin{alignat}{2}
& G_{ij}^{FS}(\x,\x_0)=-\frac{\delta_{ij}\ln{|\hat \x|}}{4\pi}+\frac{\hat{x_i}\hat{x_j}}{4\pi |\hat \x|^2} \label{fs greens funct 1}\\
&p_i^{FS}(\x,\x_0)=\frac{\hat{x_i}}{2\pi |\hat \x|^2} \label{fs greens funct 2}\\
&T_{ijk}^{FS}(\x,\x_0)=-\frac{\hat{x_i}\hat{x_j}\hat{x_k}}{\pi |\hat \x|^4} \label{fs greens funct 3},
\end{alignat}
\end{subequations}
where $\delta_{ij}$ is the Kronecker delta, which is $1$ when $i=j$ and $0$ otherwise. A general Green's function, stress tensor, and pressure vector can be decomposed into the free-space functions and complementary functions that are non-singular in $\Omega$. Near the singularity, all Green's functions therefore exhibit the same leading order behavior as the free-space Green's function.

We will now derive an identity for the stress tensor $\T$ that is similar to the identity for the Green's function derivative given in Equation \eqref{greens deriv prop 2}. We start with the singularly forced Stokes equation and use the expression for $\sig$ to get 
\begin{equation}
\frac{\partial }{\partial x_{k}}\Big(T_{ijk}c_j\Big)= -\delta(\x-\x_0)c_i .
\end{equation}
In two dimensions, this is two equations, but by reorganizing each one to eliminate the arbitrary constants, $c_i$, we get four equations summarized by
\begin{equation}
\frac{\partial T_{ijk}}{\partial x_{k}}= -\delta_{ij}\delta(\x-\x_0).
\end{equation}
For $\x_0\notin \Omega$, we can use the divergence theorem to get
\begin{equation}
\int_{\Gamma} T_{ijk}(\x,\x_0)n_k(\x) dl(\x)=0.
\end{equation}
Then for the case of $\x_0 \in \Omega$, as we did in Section \ref{ch 3 Greens}, we can use the altered domain in Figure \ref{inside} and a decomposition of the stress tensor to find that
\begin{equation}
\int_{\Gamma} T_{ijk}(\x,\x_0)n_k(\x) dl(\x)=-\delta_{ij}.
\end{equation}

Lastly, we look at the case in which the singularity, $\X_0$ is on the boundary $\Gamma$. We use the domain in Figure \ref{on bound}, created by removing from $\Omega$ a circle centered at $\X_0$ with radius $\epsilon$. Then, since $\X_0 \notin \Omega'$, we can use the divergence theorem to get 
\begin{equation}
0=\int_{\Gamma_1 \cup \Gamma_2} T_{ijk}(\x,\x_0)n_k(\x) dl(\x) = \int_{\Gamma_1 }T_{ijk}(\x,\x_0)n_k(\x) dl(\x) +\int_{ \Gamma_2}T_{ijk}(\x,\x_0)n_k(\x) dl(\x) .
\end{equation}
Then, by using the decomposition of the stress tensor, we can again use the form of the Stresslet to calculate the integral. Additionally, as $\epsilon\longrightarrow 0$, the second integral approaches the integral over a semicircle. Then this integral simplifies to
\begin{equation}
-\frac{1}{\pi\epsilon^3}\int_0^{\pi} x_ix_jx_kn_k d\theta =  \frac{1}{\pi\epsilon^2}\int_0^{\pi} x_ix_j d\theta =\frac12 \delta_{ij}.
\end{equation}
We therefore get the principal-value integral, 
\begin{equation}
-\frac{1}{2}\delta_{ij}=\int_{\Gamma}^{PV} T_{ijk}(\x,\x_0)n_k(\x) dl(\x),
\end{equation}
and we have the following integral identity for the stress tensor: 
\begin{equation} \label{stresslet prop 1}
\int_{\Gamma}T_{ijk}(\x,\x_0)n_k(\x)dl(\x) = 
    \begin{cases}
       -\delta_{ij}, & \text{if } \x_0\in \Omega\\
       -\frac12\delta_{ij}, & \text{if } \x_0\in \Gamma\\
       0, & \text{if } \x_0 \notin \bar \Omega,
        \end{cases}
\end{equation}
where $\bar \Omega=\Omega \cup \Gamma$. This property is again used to derive the integral equations presented in Section \ref{ch 3 stokes integ eqns}.

\section{Integral representations for solutions to Stokes equation}\label{ch 3 stokes reps}

\subsection{Full integral representation}\label{3.5 full}
In this section, we derive the full integral representation for a solution to Equation \eqref{homog stokes pde}. We begin by deriving the reciprocal relation for Stokes equation, which plays the same role that Green's second identity did in Section \ref{ch 3 reps}. Let $\u$, $p$, $\sig$ and $\u'$, $p'$, $\sig'$ define two solutions to Equation \eqref{homog stokes pde}. Then the product rule gives us 
\begin{equation}
u_i'\frac{\partial \sigma_{ij}}{\partial x_j}=\frac{\partial  (u_i'\sigma_{ij} ) }{\partial  x_j }-\sigma_{ij}\frac{\partial   u_i' }{\partial  x_j }\label{sig thing}
\end{equation}
Next, we use 
\begin{equation}
\sigma_{ij}=-\delta_{ij}p+\mu\Bigg(\frac{\partial u_i }{\partial x_j}+\frac{\partial u_j}{\partial x_i}\Bigg),
\end{equation}
to replace $\sigma_{ij}$ in the second term of Equation \eqref{sig thing}. Then we simplify and use that $\Div \u'=0$ to get 
\begin{equation}
u_i' \frac{\partial \sigma_{ij}}{\partial x_j}= \frac{\partial  (u_i'\sigma_{ij} ) }{\partial  x_j }-\mu\Bigg(\frac{\partial u_i }{\partial x_j}+\frac{\partial u_j}{\partial x_i}\Bigg) \frac{\partial   u_i' }{\partial  x_j }.
\end{equation}
Switching the flows to form a similar expression for $u_i(\partial \sigma_{ij}'/\partial x_j)$ and subtracting them, we get 
\begin{equation}
u_i' \frac{\partial \sigma_{ij}}{\partial x_j}- u_i \frac{\partial \sigma_{ij}'}{\partial x_j}= \frac{\partial}{\partial x_j}\Big(u_i'\sigma_{ij}-u_i\sigma_{ij}'\Big).\label{before recip}
\end{equation}
Lastly, since they are both solutions to Stokes equations, $\sigma_{ij}/\partial x_j=\sigma_{ij}'/\partial x_j=0$, and we get the reciprocal relation:
\begin{equation}
0= \frac{\partial}{\partial x_j}\Big(u_i'\sigma_{ij}-u_i\sigma_{ij}'\Big).\label{recip}
\end{equation}

Now, instead of using two solutions to Equation \eqref{homog stokes pde}, we replace the second flow with the flow from the singularly forced Stokes equation. Then Equation \eqref{before recip} becomes 
\begin{equation}
u_i^G \frac{\partial \sigma_{ij}}{\partial x_j}- u_i \frac{\partial \sigma_{ij}^G}{\partial x_j} = \frac{\partial}{\partial x_j}\Big(u_i^G\sigma_{ij}-u_i\sigma_{ij}^G\Big),
\end{equation}
or 
\begin{equation}
u_i\delta(\x-\x_0) c_i = \frac{\partial}{\partial x_j}\Bigg(\frac{1}{\mu}G_{im}c_m\sigma_{ij}-u_iT_{imj}c_m\Bigg).
\end{equation}
This is a scalar equation, but by reorganizing to eliminate the arbitrary constants, $c_i$, we get two equations summarized by
\begin{equation}
u_i(\x)\delta(\x-\x_0)= \frac{\partial}{\partial x_j}\Bigg(\frac{1}{\mu}G_{ki}(\x,\x_0)\sigma_{kj}(\x)-u_k(\x)T_{kij}(\x,\x_0)\Bigg).
\end{equation}
By again first making the necessary adjustment to the domain to ensure the differentiability of the integrand, we can integrate over the altered domain and use the divergence theorem to obtain the integral representation of $u$. By also making a switch of indices for a cleaner expression, we get 
\begin{equation}
u_j(\x_0)=\frac{1}{\mu}\int_{\Gamma}G_{ij}(\x,\x_0)\sigma_{ik}(\x)n_k(\x)dl(\x)-\int_{\Gamma}u_i(\x)T_{ijk}(\x,\x_0)n_k(\x)dl(\x).\label{stokes rep}
\end{equation}

The first term of Equation \eqref{stokes rep} is the single layer potential, and its strength is given by the boundary traction, $\sig\dotp \n$. The second term is the double layer potential, and its strength is given by the boundary values.

\subsection{Single layer integral representation}\label{3.5 SL}

In a similar fashion as we saw in Section \ref{ch 3 reps}, we can define two solutions, one on the interior domain $\Omega$ and another on the region exterior to $\Gamma$. If we let these two solutions have equal values on boundary, we can take Equation \eqref{stokes rep} for $\u^{Int} $ and a similar equation for $\u^{Ext}$,
\begin{equation}
0=-\frac{1}{\mu}\int_{\Gamma}G_{ij}(\x,\x_0)\sigma_{ik}^{Ext}(\x)n_k(\x)dl(\x)+\int_{\Gamma}u_i^{Ext}(\x)T_{ijk}(\x,\x_0)n_k(\x)dl(\x),
\end{equation}
where $\n$ is pointing out of $\Omega$. Then, by adding these equations, we get the single layer integral equation for $u^{Int}$ given by 
\begin{equation}
u_j^{Int} (\x_0) = \frac{1}{\mu}\int_{\Gamma}    G_{ij}(\x, \x_0) \Phi_i(\x) dl(\x)  ,        \label{stokes single layer rep}
\end{equation}
where the strength of the single layer potential is given by the jump in traction across the boundary, or
\begin{equation}
\bm{\Phi} = (\sig^{Int}-\sig^{Ext})\dotp \n\Big|_{\Gamma}.
\end{equation} 

\textbf{Force and Torque. } It will be useful for subsequent chapters to calculate the net force and net torque exerted on a closed curve. As we have established, the boundary distribution, $\bm{\Phi}$ is supported on the boundary $\Gamma$. Let $\widetilde \Gamma$ be another curve that encloses $\Gamma$. Then the net force on $\widetilde \Gamma$ is given by 
\begin{equation}
\bm{B}=\int_{\widetilde \Gamma} \sig(\x)\dotp\n(\x) dl(\x),
\end{equation}
where we use $\bm{B}$ instead of $\F$ in order to avoid conflict with our use of $\F$ as the IBSL constraint force. Let $R$ be the area enclosed by $\widetilde \Gamma$. Then we can use the divergence theorem to find 
\begin{equation}
0=\int_R \Div\sig(\x) d\x= \int_{\widetilde \Gamma}\sig(\x)\dotp \n(\x) dl(\x) + \int_{\Gamma}\jump{\sig \dotp \n} dl(\x) .\label{force derivation}
\end{equation}
Therefore, the net force is 
\begin{equation}
 \bm{B} =- \int_{\Gamma}\jump{\sig \dotp \n} dl(\x) .
\end{equation}
In the case where $\bm{Phi}$ is defined as we have in this section, the net force is then
\begin{equation}
 \bm{B}=- \int_{\Gamma}\bm{\Phi}(\x) dl(\x) .
\end{equation}

Then, the net torque on $\widetilde \Gamma$, in two dimensions, is given by  
\begin{equation}
\bm{L}=\int_{\widetilde \Gamma} \x \times (\sig \dotp \n)dl(\x)=\epsilon_{ijk} \int_{\widetilde \Gamma} x_j\sigma_{km}n_mdl(\x),
\end{equation}
where $\epsilon_{ijk}$ is the Levi-Civita symbol, which in 2-D gives $1$ for an even permutation of $(1,2)$, $-1$ for an odd permutation, and $0$ otherwise. Before deriving the expression for the net torque, let us note that 
\begin{equation}
 \epsilon_{ijk} \int_R\frac{\partial }{\partial x_m}(x_j\sigma_{km})d\x= \epsilon_{ijk} \int_Rx_j\frac{\partial \sigma_{km} }{\partial x_m}d\x+ \epsilon_{ijk} \int_R\delta_{mj} \sigma_{km}d\x = 0 ,
\end{equation}
where in the penultimate expression, the first integral is $0$ since $\Div \sig=0$ and the second integral is $0$ by the symmetry of $\sig$. Then, we can use the divergence theorem to get 
\begin{equation}
 0=\epsilon_{ijk} \int_R\frac{\partial }{\partial x_m}(x_j\sigma_{km})d\x=\epsilon_{ijk} \int_{\widetilde \Gamma} x_j\sigma_{km}n_mdl(\x)+\epsilon_{ijk} \int_{ \Gamma} x_j\jump{\sigma_{km}n_m}dl(\x). \label{torque derivation}
\end{equation}
Therefore the net torque exerted on $\widetilde \Gamma$ is 
\begin{equation}
 \bm{L}=-\epsilon_{ijk} \int_{\Gamma} x_j\jump{\sigma_{km}n_m}dl(\x). 
\end{equation}
Note that in two dimensions, this quantity can be viewed as a scalar. As with the net force, in terms of the single layer potential strength, this is equivalent to
\begin{equation}
L=- \int_{\Gamma} \x \times \bm{\Phi} (\x) dl(\x).
\end{equation}

\subsection{Double layer integral representation}\label{3.5 DL}

If we follow the the same process as the previous section but require that the interior and exterior solutions have equal boundary tractions instead of equal boundary values, then we get the double layer integral equation for $u^{Int}$ given by 
\begin{equation}
u_j^{Int} (\x_0) = \int_{\Gamma}   \Psi_i(\x) T_{ijk}(\x,\x_0)n_k(\x) dl(\x)  ,        \label{stokes double layer rep}
\end{equation}
where the strength of the double layer potential is given by the jump in velocity across the boundary, or
\begin{equation}
\bm{\Psi} = (\u^{Ext}-\u^{Int})\Big|_{\Gamma}.
\end{equation} 

We now seek to find the net force and net torque on a closed curve $\widetilde \Gamma$ surrounding $\Gamma$. We can still use the derivations in Equations \eqref{force derivation} and \eqref{torque derivation}, but since the double layer representation was formed by assuming no jump in the traction across $\Gamma$, we get that the net force and net torque on $\widetilde \Gamma$ are $0$. This shows us that the double layer potential alone is incapable of representing flows  with a net force or net torque. For a domain interior to $\Gamma$, this does not cause a problem since a solution to Stokes equation on this domain has no net force, as seen by 
\begin{equation}
\bm{B}=\int_{\Gamma}\sig\dotp \n dl(\x)= \int_{\Omega} \Div \sig d\x = 0,
\end{equation}
and similarly for the net torque. However, on an exterior domain, we require the ability to represent a flow with net force and net torque. In order to achieve this, we can supplement the double layer potential with a flow producing a finite force and torque on $\Gamma$. The resulting representation is called a compound or completed double layer representation. For a review of ways to do this, one can see Section 4.7 in \cite{Pozred}. We briefly present two such methods here. 

One solution, presented by Power and Miranda \cite{powermiranda, power}, is to add a pair of singularities to the interior of $\Omega$ using a Stokeslet, or a point force, and a rotlet, or a point torque. A second solution, which we will make explicit use of in Chapter \ref{chapter 6}, was presented by Hsiao and Kress \cite{hsiaokress} and Hebeker \cite{hebeker}. It consists of adding a flow created by a single layer potential whose strength is a constant multiple of the double layer potential strength. If we let this constant be $\eta/\mu$, where $\eta$ is an arbitrary positive constant, then our integral representation is given by 
\begin{equation}
u_j (\x_0) = \int_{\Gamma}   \Psi_i(\x) T_{ijk}(\x,\x_0)n_k(\x) dl(\x)+\frac{\eta}{\mu} \int_{\Gamma} G_{ij}(\x,\x_0)\Psi_i(\x)dl(\x)  .       \label{stokes completed double layer rep}
\end{equation}
Using what we found in Section \ref{3.5 SL}, we then have the net force and net torque as
\begin{subequations}
\begin{alignat}{2}
&\bm{B}=-\int_{ \Gamma}\bm{\Psi}(\x) dl(\x)\\
&L= -\int_{ \Gamma} \x \times \bm{\Psi} (\x) dl(\x).
\end{alignat}
\end{subequations}

\section{Integral equations for Stokes equation with Dirichlet boundary conditions}\label{ch 3 stokes integ eqns}

In this section, we introduce the integral equations that arise from using the representations from the previous section to solve the Stokes equation with Dirichlet boundary conditions given by
\begin{subequations}\label{homog stokes pde again }
\begin{alignat}{2}
& \mu\Delta \u -\grad p=\grad \dotp \sig = 0\qquad && \text{in } \Omega  \label{homog stokes pde again 1}\\
&\Div \u =0  \qquad &&  \text{in } \Omega \label{homog stokes pde again 2}\\
& \u=\U_b \qquad &&  \text{on } \Gamma . \label{homog stokes pde again 3}
\end{alignat}
\end{subequations}

We omit the derivation, as it is similar to that in Section \ref{ch 3 integral eqns}, but using the property of $\T$ given in Equation \eqref{stresslet prop 1}, we get the single and double layer integral equations given by 
\begin{equation}
U_{b_j}(\X_0) = \frac{1}{\mu}\int_{\Gamma}  G_{ij}(\x, \X_0) \Phi_i(\x) dl(\x)  ,      \label{stokes single layer eqn}
\end{equation}
\vspace{-0.3cm}
\begin{equation}
U_{b_j} (\X_0) = \int_{\Gamma}^{PV} \Psi_i(\x) T_{ijk}(\x,\X_0)n_k(\x) dl(\x)   - \frac12 \Psi_j(\X_0), \label{stokes double layer eqn 3}
\end{equation}
where the unknown quantities in Equations \eqref{stokes single layer eqn} and \eqref{stokes double layer eqn 3} are the potential strengths, $\bm{\Phi}$ and $\bm{\Psi}$, respectively. Note that again, the integral equations for the Brinkman equation have the same form. 
For the reasons discussed in Section \ref{3.5 DL}, the double layer integral equation above is only valid for flows with no net force or torque. For other exterior flows, we can use the completed double layer integral equation given by 
\begin{equation}
U_{b_j} (\X_0) = \int_{\Gamma}   \Psi_i(\x) T_{ijk}(\x,\X_0)n_k(\x) dl(\x)+\frac{\eta}{\mu} \int_{\Gamma} G_{ij}(\x,\X_0)\Psi_i(\x)dl(\x)  - \frac12 \Psi_j(\X_0),        \label{stokes completed double layer eqn}
\end{equation}
which comes from taking the limit as $\x_0\longrightarrow \Gamma$ in Equation \eqref{stokes completed double layer rep}.

\section{Discussion}\label{ch 3 discussion}

Boundary integral methods can be used to get very accurate solutions to homogeneous linear boundary value problems by reformulating them as integral equations with one fewer dimension. However, the efficiency is greatly reduced when the problem is nonhomogeneous as this introduces an integral over the entire PDE domain. It can also be expensive to use the potential strength to obtain solution values at all points in a discretized domain. Additionally, by requiring the analytical form of the Green's function and its derivative, it can be difficult to generalize to non-free-space exterior domains, such as the periodic box used in this dissertation. Therefore, as discussed in Section \ref{ch 2 discussion}, we seek to use the Immersed Boundary framework instead of explicitly using boundary integral methods. However, boundary integral methods have the advantage that well-conditioned integral operators exist, and this makes the step of solving for the unknown boundary distribution much more efficient. 

Let us examine the integral operators for Laplace's equation to explore this conditioning issue. The single layer integral equation,
\begin{equation}
U_b(\X_0) = \int_{\Gamma}  \sigma (\x)  G(\x, \X_0)dl(\x)  ,      \label{single layer eqn again}
\end{equation}
has the form of a Fredholm integral of the first kind, and the double layer integral equation, 
\begin{equation}
U_b (\X_0) = \int_{\Gamma}^{PV} \gamma(\x) \grad G(\x, \X_0)\dotp \n(\x) dl(\x) - \frac12 \gamma(\X_0), \label{double layer eqn again}
\end{equation}
has the form of a Fredholm integral of the second kind \cite{atkinson}. The operators given by  
\begin{subequations} \label{integral operators}
\begin{alignat}{1}
&K_1 \sigma = \int_{\Gamma}  \sigma (\x)  G(\x, \x_0)dl(\x)    \label{sl integral operator}\\
& K_2\gamma = \int_{\Gamma} \gamma(\x) \grad G(\x, \x_0)\dotp \n(\x) dl(\x)    \label{dl integral operator}
\end{alignat}
\end{subequations}
have eigenvalues in the interval $(-0.5, 0.5)$, and the only limit point of the eigenvalues is $0$ \cite{coltonkress, spectralproperties}. The condition number, defined as the ratio of the largest and smallest eigenvalues, is therefore infinite in the continuous case and large in the discretized case. On the other hand, Equation \eqref{double layer eqn again} can be rewritten as
\begin{equation}
\Big(K_2 - \frac12 \mathds{I}\Big)\gamma = U_b, 
\end{equation}
where $\mathds{I}$ is the identity operator. Shifting the operator in this way shifts the eigenvalues, and the only limit point for the eigenvalues becomes -1/2. Therefore, the condition number of the discretized operator is finite and does not grow with refinement of the discretization. The better conditioning of this operator is the characteristic of the double layer representation that we exploit to form our new Immersed Boundary Double Layer method.

   \chapter{Connection between IBSL method and boundary integral methods}
   \label{chapter 4}


In this chapter, we relate the Immersed Boundary Single Layer method as described in Chapter \ref{chapter 2} to the single layer integral equations described in Chapter \ref{chapter 3}. As discussed in Chapter \ref{chapter 1}, this connection has been established in several recent works \cite{GriffithDonev2, eldredge}, but here, we present this connection explicitly from the IBSL constraint systems seen in Equation \eqref{constraint IB lL} for Helmholtz and Poisson equations and in Equation \eqref{ibsl brinkman} for Brinkman and Stokes equations. By establishing this connection, we can see that the large iteration counts seen in Chapter \ref{chapter 2} can be explained by the poor conditioning of the integral operator, as discussed at the end of Chapter \ref{chapter 3}. We will then propose that a solution to the conditioning problem is to formulate an Immersed Boundary version of a double layer integral equation. Such formulations will be presented in Chapters \ref{chapter 5} and \ref{chapter 6}. 

For this chapter and the derivation of the IBDL method, we will simplify the presentation by focusing on the homogeneous case, or $g=0$. However, this connection holds for the inhomogeneous case as well, and we use the new method on several inhomogenous problems in Chapters \ref{chapter 5} and \ref{chapter 6}.

\section{Poisson and Helmholtz equations}\label{ch 4 helmholtz}
We begin by making this connection explicit in the case of the Poisson and Helmholtz equations. For simplicity, we assume that the differential operator $\L=\Delta -k^2$ is invertible, either by focusing on Helmholtz or assuming a computational domain $\C$ for which $\Delta$ is invertible. 

\subsection{Regularized Green's function}\label{4.1 reg Greens}
Let us first define a regularized Green's function by 
\begin{equation}
G_h(\x,\x_0) \equiv G * \delta_h = \int_{\C} G(\mathbf{y}, \x_0)\delta_h(\x-\mathbf{y})d\mathbf{y}. \label{Gh}
\end{equation}
The linearity of the differential operator then gives us that 
\begin{equation}
\L G_h(\x,\x_0)=-\delta_h(\x-\x_0).    \label{Gh2}
\end{equation}
In this way, $G_h$ can be seen both as a regularization of a Green's function through a convolution with a smoothing function and as the solution to the PDE when the forcing function is a regularized delta function.

\subsection{IBSL method as a regularized single layer integral equation}\label{4.1 connect}

The homogeneous Immersed Boundary Single Layer system is given by 
\begin{subequations} \label{constraint IB lL again again }
\begin{alignat}{2}
& \L u +SF = 0 \qquad && \text{in } \mathcal{C}  \label{constraint IB L1 again again}\\
&S^* u = U_b \qquad && \text{on } \Gamma . \label{constraint IB L2 again again}
\end{alignat}
\end{subequations}
Starting with Equation \eqref{constraint IB L1 again again} and using the definition of the spread operator, we have 
\begin{equation}
\L u = - SF = - \int_{\Gamma} F(s)\delta_h(\x-\X(s)) ds. \label{connection1}
\end{equation}
Inverting the operator and using Equation \eqref{Gh2}, we get
\begin{equation}
u(\x) =  \int_{\Gamma} F(s)G_h(\x,\X(s)) ds \label{connection2}
\end{equation}
for $\x\in\Omega\setminus \Gamma$. Then, Equation \eqref{constraint IB L2 again again} gives us 
\begin{equation}
U_b(s') =\int_{\C} u(\x)\delta_h(\x-\X(s')) d\x,  \label{connection3}
\end{equation}
where we are using $s'$ to distinguish from our previous variable of integration or, alternatively, viewing $\X(s')$ as selecting a particular boundary point. Combining this with Equation \eqref{connection2} and changing the order of integration, we get
\begin{equation}
U_b(s') =\int_{\Gamma} F(s) \int_{\C} G_h(\x,\X(s)) \delta_h(\x-\X(s')) d\x ds. \label{connection4}
\end{equation}
Notice that the second integral is equivalent to $G_h * \delta_h$, evaluated at $\X(s')$. Therefore, by denoting a twice-regularized Green's function by $G_{hh}$, this becomes
\begin{equation}
U_b(s') =\int_{\Gamma} F(s)  G_{hh}(\X(s'),\X(s)) ds.\label{connection5}
\end{equation}
The symmetry of the Green's function, which is preserved through convolutions with the regularized delta function, gives us that $G_{hh}(\X(s'),\X(s)) = G_{hh}(\X(s),\X(s'))$. Using this to switch the arguments of $G_{hh}$ and appropriately redefining $F$ and $U_b$ as functions of $\x$, we get 
\begin{equation}
U_b(\X(s')) =\int_{\Gamma} F(\X(s)) G_{hh}(\X(s),\X(s')) ds. \label{connection6}
\end{equation}

Recall from Equation \eqref{single layer eqn} that the single layer integral equation for $\X_0$ on the boundary is given by 
\begin{equation}
 U_b(\X_0) =\int_{\Gamma} \sigma(\x)  G(\x, \X_0) dl(\x).  \label{connection7}
\end{equation}
 By parametrizing $\Gamma$ by $\X(s)$, this becomes 
 \begin{equation}
 U_b(\X(s')) =\int_{\Gamma} \sigma(\X(s))  G(\X(s),\X(s')) \bigg|\frac{\partial \X(s)}{\partial s}\bigg| ds.  \label{connection8}
\end{equation}
By comparing Equations \eqref{connection6} and \eqref{connection8}, we can associate the IBSL constraint force $F$ with the potential strength modified by the parametrization term, $\big|\partial \X(s)/\partial s\big| $. In other words, in the limit that the regularization width, $h$, approaches $0$, we have
\begin{equation}
F(\X)=\sigma(\X)\bigg|\frac{\partial \X}{\partial s}\bigg|. \label{F eqn}
\end{equation}
Specifically, in the case of an arclength parametrization, where $ \big|\partial \X(s)/\partial s\big| =1$, we see from Chapter \ref{chapter 3} that $F$ gives the jump in the normal derivative of the solution across the boundary. We have now established that the IBSL method is equivalent to a \emph{regularized} single layer integral equation. It is this connection that motivates the name Immersed Boundary Single Layer method in order to distinguish it from the new method we propose in Chapter \ref{chapter 5}.

\section{Stokes equation}\label{ch 4 stokes}

We next look at establishing this connection between the IBSL method and a single layer integral equation for the Stokes equation. Again, for simplicity we assume that the differential operators involved are invertible on the computational domain $\C$. Section \ref{2.3 stokes} discusses how to utilize the IBSL method when they are not invertible.

\subsection{Regularized Green's functions and related identities}\label{4.2 reg Greens}

Since Stokes equation involves the Laplacian, we will require the use of the regularized Green's function for Laplace's equation developed in Section \ref{4.1 reg Greens}. To distinguish it from the Stokes Green's function in this section, we will call this function $G_h^L$. In this section, we will define other regularized Green's functions and derive some identities that will be used in Section \ref{4.2 connect} and in Chapter \ref{chapter 6}.

\textbf{Regularized biharmonic Green's function. } Let us define the Green's function to the biharmonic equation as the function $H(\x,\x_0)$ that satisfies
\begin{equation}
\Delta^2 H = -\delta (\x-\x_0).\label{Hdef}
\end{equation}
We also define the regularized biharmonic Green's function by 
\begin{equation}
H_h(\x,\x_0) \equiv H * \delta_h = \int_{\C} H(\mathbf{y}, \x_0)\delta_h(\x-\mathbf{y})d\mathbf{y}. \label{Hh}
\end{equation}
The linearity of the differential operator again gives us that 
\begin{equation}
\Delta^2 H_h(\x,\x_0)=-\delta_h(\x-\x_0).    \label{Hh2}
\end{equation}
Then using our regularized Laplace's Green's function, we have
\begin{equation}
\Delta^2 H_h(\x,\x_0)=\Delta G_h^L(\x,\x_0),  
\end{equation}
and therefore that 
\begin{equation}
\Delta H_h(\x,\x_0)= G_h^L(\x,\x_0) . \label{G and H relationship}
\end{equation}
Note that we also have this relationship for the exact Green's functions, so 
\begin{equation}
\Delta H(\x,\x_0)= G^L(\x,\x_0) . \label{G and H relationship no h}
\end{equation}

\textbf{Relationship between Stokes and biharmonic Green's functions. } We next seek to find a relationship between the Stokes Green's function $\G$ and the biharmonic Green's function $H$. We start with the singularly forced Stokes equation, given by
\begin{subequations} \label{sing forced stokes}
\begin{alignat}{2}
& \mu \Delta\u - \grad p =-\delta(\x-\x_0)\bm{c}  \label{sing forced stokes 1} \\
& \grad \dotp \u =0 . \label{sing forced stokes 2}
\end{alignat}
\end{subequations}
By taking the divergence of Equation \eqref{sing forced stokes 1} and using the incompressibility of the velocity, we get 
\begin{equation}
\Delta p = \Div \delta(\x-\x_0)\bm{c}.
\end{equation}
Then we can invert the operator to get
\begin{equation}
 p(\x) = -\Div G^L(\x,\x_0)\bm{c}= -\grad G^L(\x,\x_0) \dotp \bm{c}.\label{p thing}
\end{equation}
Using this in Equation \eqref{sing forced stokes 1}, we get
\begin{equation}
\mu \Delta \u = -\grad( \grad G^L(\x,\x_0)\dotp \bm{c})-\delta(\x-\x_0)\bm{c}.
\end{equation}
Again inverting the operator and using Equation \eqref{G and H relationship no h} for the first term and Equation \eqref{Hdef} for the second term, we get
\begin{equation}
\mu \u(\x) = -\grad( \grad H(\x,\x_0)\dotp \bm{c})+\Delta H (\x,\x_0)\bm{c}.
\end{equation}
Switching to Einstein notation, we have 
\begin{equation}
\mu u_i(\x) = - \frac{\partial H}{\partial x_i\partial x_j }(\x,\x_0)c_j+ \frac{\partial H}{\partial x_k\partial x_k}(\x,\x_0)c_i.
\end{equation}
Therefore, since the the solution using the Stokes Green's function is given by 
\begin{equation}
u_i(\x)=\frac{1}{\mu}G_{ij}(\x,\x_0)c_j,    
\end{equation}
we see that 
\begin{equation}
G_{ij}=- \frac{\partial H}{\partial x_i\partial x_j }+ \delta_{ij}  \frac{\partial H}{\partial x_k\partial x_k}.\label{stokes g to other ones}
\end{equation}

\textbf{Relationship between the Stokes pressure vector and Laplace's Green's function. } 
We will now find an identity for the Stokes pressure vector that will be needed in Section \ref{ch 6 IBDL}. From equation \eqref{p thing}, we have 
\begin{equation}
 p(\x) =  -\grad G^L(\x,\x_0) \dotp \bm{c}.
\end{equation}
Since the pressure vector is defined in Equation \eqref{greens funct 2} as $\p$ such that $p(\x)=p_jc_j$, we can see that 
\begin{equation}
p_j =  -\frac{\partial G^L}{\partial x_j}(\x,\x_0).
\end{equation}

\textbf{Regularized Stokes Green's function. } Let us also define the regularized Stokes Green's function and the corresponding regularized pressure vector and stress tensor as 
\begin{subequations}
\begin{equation}
G_{ij}^h(\x,\x_0) \equiv G_{ij} * \delta_h= \int_{\C} G_{ij}(\mathbf{y}, \x_0)\delta_h(\x-\mathbf{y})d\mathbf{y},  \label{Gh stokes}
\end{equation}
\begin{equation}
p_j^h(\x,\x_0) \equiv p_j * \delta_h= \int_{\C} p_j(\mathbf{y}, \x_0)\delta_h(\x-\mathbf{y})d\mathbf{y},  \label{Gh stokes}
\end{equation}
\begin{equation}
T_{ijk}^h(\x,\x_0) \equiv T_{ijk} * \delta_h= \int_{\C} T_{ijk}(\mathbf{y}, \x_0)\delta_h(\x-\mathbf{y})d\mathbf{y},  \label{Gh stokes}
\end{equation}
\end{subequations}

By linearity, the relationships that we have previously shown for the exact functions also hold for the regularized ones, so that 
\begin{equation}
G_{ij}^h=- \frac{\partial H_h}{\partial x_i\partial x_j }+ \delta_{ij}  \frac{\partial H_h}{\partial x_k\partial x_k},\label{stokes g to other ones h}
\end{equation}
and the solution to a Stokes equation that is forced with the regularized delta function will exactly be the solution corresponding to this regularized Green's function. 
Additionally, we have 
\begin{equation}
\p_j^h =  -\frac{\partial G_h^L}{\partial x_j}(\x,\x_0).\label{p to GL}
\end{equation}

\subsection{IBSL method as a regularized single layer integral equation}\label{4.2 connect}
The homogeneous Immersed Boundary Single Layer system for Stokes equation is given by 
\begin{subequations} \label{ibsl brinkman again}
\begin{alignat}{2}
& \mu \Delta\u - \grad p +S\F =  0 \qquad && \text{in } \mathcal{C}  \label{ibsl brinkman 1 again}\\
& \grad \dotp \u =0 \qquad && \text{in } \mathcal{C}  \label{ibsl brinkman 2 again}\\
&S^* \u = \U_b \qquad && \text{on } \Gamma.  \label{ibsl brinkman 3 again}
\end{alignat}
\end{subequations}
We will now explicitly derive the connection to a regularized single layer integral equation. Many of the steps will mimic the steps taken in Section \ref{4.2 reg Greens}. Taking the divergence of Equation \eqref{ibsl brinkman 1 again} and using Equation \eqref{ibsl brinkman 2 again} to eliminate the velocity term, we get 
\begin{equation}
\Delta p =  \Div S \F = \Div \int_{\Gamma} \F(s)\delta_h(\x-\X(s)) ds. 
\end{equation}
We can then invert the operator and use the property of the regularized Green's function for Laplace's equation, given in Equation \eqref{Gh2}. We then have
\begin{equation}
 p(\x) = - \Div \int_{\Gamma} \F(s)G_h^L(\x,\X(s)) ds. 
\end{equation}
Bringing in the divergence and manipulating the expression, we get 
\begin{equation}
 p(\x) = -  \int_{\Gamma} \F(s)\dotp \grad G_h^L(\x,\X(s)) ds. 
\end{equation}
Using this expression for pressure, Equation \eqref{ibsl brinkman 1 again} then gives us 
\begin{equation}
\mu \Delta \u =\grad p - S\F= -\grad  \int_{\Gamma} \F(s)\dotp \grad G_h^L(\x,\X(s)) ds -\int_{\Gamma} \F(s)\delta_h(\x-\X(s)) ds . 
\end{equation}
At this point, we switch to Einstein summation notation, and by bringing in the gradient, we get
\begin{equation}
\mu \Delta u_j =  -\int_{\Gamma}  \frac{\partial^2 G_h^L}{\partial x_i\partial x_j }(\x,\X(s)) F_i(s) ds -\int_{\Gamma} F_j(s)\delta_h(\x-\X(s)) ds . \label{sl part 1}
\end{equation}
Inverting the operator on $u_i$ and using Equation \eqref{G and H relationship} for the first term and Equation \eqref{Hh2} for the second term, we get
\begin{equation}
\mu u_j(\x) =  -\int_{\Gamma}  \frac{\partial^2 H_h}{\partial x_i\partial x_j }(\x,\X(s)) F_i(s) ds +\int_{\Gamma} F_j(s)\frac{\partial^2 H_h}{\partial x_k \partial x_k}(\x,\X(s)) ds . \label{almost there}
\end{equation}
Recognizing the form of the regularized Stokes Green's function from Equation \eqref{stokes g to other ones h}, we have 
\begin{equation}
u_j(\x)= \frac{1}{\mu} \int_{\Gamma} G_{ij}^h (\x, \X(s))F_i(s) ds\label{sl part 2}
\end{equation}
for $\x\in\Omega\setminus \Gamma$. Next, we use Equation \eqref{ibsl brinkman 3 again} to get 
\begin{equation}
\U_b(s') =\int_{\C} \u(\x)\delta_h(\x-\X(s')) dx. 
\end{equation}
We then follow the same steps seen in Section \ref{4.1 connect}. We combine this with Equation \eqref{almost there}, change the order of integration, recognize the presence of $G_{ij}^h * \delta_h$, and denote the twice-regularized Green's function by $G_{ij}^{hh}$. This leaves us with 
\begin{equation}
U_{b_j}(s') =\frac{1}{\mu}\int_{\Gamma} G_{ij}^{hh}(\X(s'),\X(s)) F_i(s) ds. 
\end{equation}
The symmetry of the Green's function, which is preserved through convolutions with the regularized delta function, gives us that $G_{ij}^{hh}(\X(s'),\X(s)) = G_{ij}^{hh}(\X(s),\X(s'))$ \cite{Pozred}. Using this to switch the arguments of $G_{hh}$ and appropriately redefining $\F$ and $\U_b$ as functions of $\x$, we get 
\begin{equation}
U_{b_j}(\X(s')) =\frac{1}{\mu}\int_{\Gamma} G_{ij}^{hh}(\X(s),\X(s')) F_i(s) ds.\label{IBSL stokes one}
\end{equation}

Recall that the single layer integral equation for $\X_0$ on the boundary from Section \ref{ch 3 stokes integ eqns} is given by 
\begin{equation}
U_{b_j}(\X_0) =\frac{1}{\mu}\int_{\Gamma} G_{ij}(\x,\X_0) \Phi_i(\x) dl(\x).  
\end{equation}
By parametrizing $\Gamma$ by $\X(s)$, we get 
\begin{equation}
U_{b_j}(\X(s')) =\frac{1}{\mu}\int_{\Gamma} G_{ij}(\X(s),\X(s')) \Phi_i(\X(s)) \bigg|\frac{\partial \X(s)}{\partial s}\bigg| ds.  \label{bem sl one}
\end{equation}
By comparing Equations \eqref{IBSL stokes one} and \eqref{bem sl one}, in the case of an arclength parametrization, we can associate the IBSL constraint force $\F$ with the potential strength $\bm{\Phi}$, and in the limit that the regularization width, $h$, approaches $0$, we have
\begin{equation}
\F(\X)=\Phi(\X). \label{ stokesF eqn}
\end{equation}
We then see from Chapter \ref{chapter 3} that $\F$ gives the jump in the traction across the boundary. We have now established that the IBSL method for Stokes is equivalent to a \emph{regularized} single layer integral equation.

   \chapter{The Immersed Boundary Double Layer (IBDL) method: Helmholtz and Poisson equations}
   \label{chapter 5}


In this chapter, we introduce the Immersed Boundary Double Layer (IBDL) method for Helmholtz and Poisson equations. The development of this method was motivated by the desire to maintain the robustness of the Immersed Boundary method for boundary value problems, while avoiding the poor conditioning introduced by the form of the IBSL method, namely the first-kind integral equation to which it corresponds. In this chapter, we see that the IBDL method is in fact able to achieve the same first-order convergence as the IBSL method while improving this conditioning drastically by utilizing the form of a second-kind integral equation. In this chapter, we also explore other advantages of the IBDL method, such as the convergence of the potential strength. This convergence also allows us to use the framework of the IBDL method to solve a PDE with Neumann boundary conditions. This is in contrast to the IBSL method, which only applies to Dirichlet problems. 

This chapter is organized as follows. In Section \ref{ch 5 IBDL}, we introduce the method and demonstrate the connection to a regularized double layer integral equation. We also discuss the discontinuity of the resulting solution and how this affects pointwise convergence near the boundary. In Section \ref{ch 5 neumann}, we present the method for using the IBDL method to solve PDEs with Neumann boundary conditions. In Section \ref{ch 5 numerical}, we discuss the numerical implementation, focusing on the areas in which this is different from the IBSL numerical implementation discussed in Section \ref{ch 2 numerical implementation}. In this section, we also discuss the way in which the IBDL method can be used to easily flag points as interior or exterior to the immersed boundary. In Section \ref{ch 5 results}, we compute solutions to Helmholtz and Poisson equations and compare the results to the IBSL method.  In Section \ref{ch 5 max norm}, we further discuss the factors affecting the pointwise convergence of the solution. In Section \ref{ch 5 potential}, we illustrate the convergence of the strength of the potential. Finally, in Section \ref{ch 5 neumann results}, we solve a PDE with Neumann boundary conditions.

\section{Mathematical description of the method}\label{ch 5 IBDL}

We look again at the boundary value problem given by 
\begin{subequations} \label{pde again}
\begin{alignat}{2}
& \L u = g \qquad && \text{in } \Omega  \label{pde1 again}\\
&u=U_b \qquad && \text{on } \Gamma,  \label{pde2 again}
\end{alignat}
\end{subequations}
where $\L=\Delta-k^2$ is the Helmholtz or Laplacian operator. As in Chapter \ref{chapter 2}, we assume that $\Omega$ is a subset of the larger computational domain $\C$ and that $\Gamma$ is a smooth one-dimensional curve parametrized by $\X(s)$. 

\subsection{Formulation}\label{5.1 formulation}
The Immersed Boundary Double Layer formulation for Equation \eqref{pde again} is given by 
\begin{subequations} \label{ibdl}
\begin{alignat}{2}
& \L u +\widetilde S Q = \tilde g \qquad && \text{in } \mathcal{C}  \label{ibdl 1}\\
&S^* u + \frac12 Q = U_b \qquad && \text{on } \Gamma,   \label{ibdl 2}
\end{alignat}
\end{subequations}
where we assume an arclength parametrization of $\Gamma$, and where we define the operator $\widetilde S$ as 
\begin{equation}\widetilde S Q \equiv \grad\dotp (SQ\n), 
\end{equation}
for $\n$, the unit normal pointing out of $\Omega$. Using the definition of $S$, this gives us 
\begin{equation}
(\widetilde S Q )(\x)= \Div \int_{\Gamma} Q(s)\delta_h(\x-\X(s)) ds =  \int_{\Gamma} Q(s)\Div\delta_h(\x-\X(s)) ds.
\end{equation}
The form of $\widetilde S$ then implies that $Q$ gives the strength of a dipole force distribution on the boundary, and this is our Lagrange multiplier used to enforce the boundary condition. Recall that $\tilde g$ represents an extension of the function $g$ from $\Omega$ to the computational domain $\C$, and we define it as $\tilde g \equiv g\chi_{\scaleto{\Omega}{4.5pt}} + g_e \chi_{\scaleto{\C \setminus \Omega}{6pt}}$. For the Helmholtz equation, there is flexibility in our choice for $g_e$. However, we will see in Section \ref{5.3 poisson} that we can sometimes use our choice of $g_e$ to deal with the nullspace of the periodic Laplacian. 

We will next demonstrate that this formulation corresponds to a regularized double layer integral equation. For simplicity, we will consider the homogeneous case, $g=0$. Starting with Equation \eqref{ibdl 1}, we have 
\begin{equation}
\L u = - \widetilde SQ = - \grad \dotp \int_{\Gamma} Q(s)\n(s) \delta_h(\x-\X(s)) ds. \label{ibdl connect1}
\end{equation}
Inverting the operator and using the property of the regularized Green's function given by Equation \eqref{Gh2}, we get
\begin{equation}
u(\x) = \grad \dotp  \int_{\Gamma} Q(s)\n(s)G_h(\x,\X(s)) ds \label{ibdl connect2}
\end{equation}
for $\x\in\Omega\setminus \Gamma$. Bringing in the divergence and manipulating the expression, we get
\begin{equation}
u(\x) =   \int_{\Gamma} Q(s) \grad G_h(\x,\X(s)) \dotp\n(s) ds .\label{ibdl connect3}
\end{equation}

Then, the second equation of the IBDL method, Equation \eqref{ibdl 2} gives us 
\begin{equation}
U_b(s') =\int_{\C} u(\x)\delta_h(\x-\X(s')) dx +\frac12Q(s'). \label{ibdl connect4}
\end{equation}
Combining Equation \eqref{ibdl connect3} with Equation \eqref{ibdl connect4}, changing the order of integration, and again recognizing the presence of $G_{hh} = G_h*\delta_h$, we get 
\begin{equation}
U_b(s') =\int_{\Gamma} Q(s) \grad G_{hh}(\X(s'),\X(s)) \dotp\n(s) ds +\frac12Q(s'). \label{ibdl connect5}
\end{equation}

The odd symmetry of the gradient of the Green's function \cite{Pozblue}, which is again preserved through the convolutions with the regularized delta function, gives us that $\grad G_{hh}(\X(s'),\X(s)) = -\grad G_{hh}(\X(s),\X(s'))$.  Using this to switch the arguments of $\grad G_{hh}$ and appropriately redefining $Q$, $U_b$, and $\n$ as functions of $\x$, we get
\begin{equation}
U_b(\X(s')) =-\int_{\Gamma} Q(\X(s)) \grad G_{hh}(\X(s),\X(s'))\dotp \n(\X(s)) ds + \frac12Q(\X(s')). \label{ibdl connect5}
\end{equation}

Recall that the double layer integral equation for $\X_0$ on the boundary is given by 
\begin{equation}
 U_b(\X_0) =\int_{\Gamma}^{PV} \gamma(\x) \grad G(\x, \X_0) \dotp \n(\x) dl(\x) - \frac12 \gamma(\X_0).  \label{ibdlconnect6}
\end{equation}
Using our arclength parametrization, we can rewrite this as
\begin{equation}
 U_b(\X(s')) =\int_{\Gamma}^{PV} \gamma(\X(s)) \grad G(\X(s), \X(s')) \dotp \n(\X(s)) ds - \frac12 \gamma(\X(s')).  \label{ibdlconnect7}
\end{equation}
By comparing Equations \eqref{ibdl connect5} and \eqref{ibdlconnect7}, we see that in the limit that the regularization width, $h$, approaches $0$, we have 
\begin{equation}
Q(\X)= - \gamma(\X).
\end{equation}
Therefore, we can see that the IBDL method is equivalent to a \emph{regularized} double layer integral equation. Note from Section \ref{ch 3 integral eqns} that $Q$ gives the jump in the solution across the boundary. 

We reiterate that the formulation, specifically Equation \eqref{ibdl 2}, assumes that $X(s)$ gives an arclength parametrization of $\Gamma$. In the IBSL method, if one uses a different paramterization of the boundary, the rescaling of $F$ seen in Equation \eqref{F eqn} accounts for the transformation, and no alteration in the method is needed to find the solution. However, in the IBDL method, since the $1/2$ in Equation \eqref{ibdl 2} is derived using an arclength parametrization, one would need to alter this term in the case that  $  \big|\partial \X(s)/\partial s\big| \neq 1$. This will be discussed further in Section \ref{5.1 param}, but otherwise, we will assume an arclength parametrization throughout this dissertation.

\subsection{Discontinuity of the solution}\label{5.1 discontinuity}

\begin{figure}
\centering
\begin{subfigure}{0.495\textwidth}
\centering
\includegraphics[width=\textwidth]{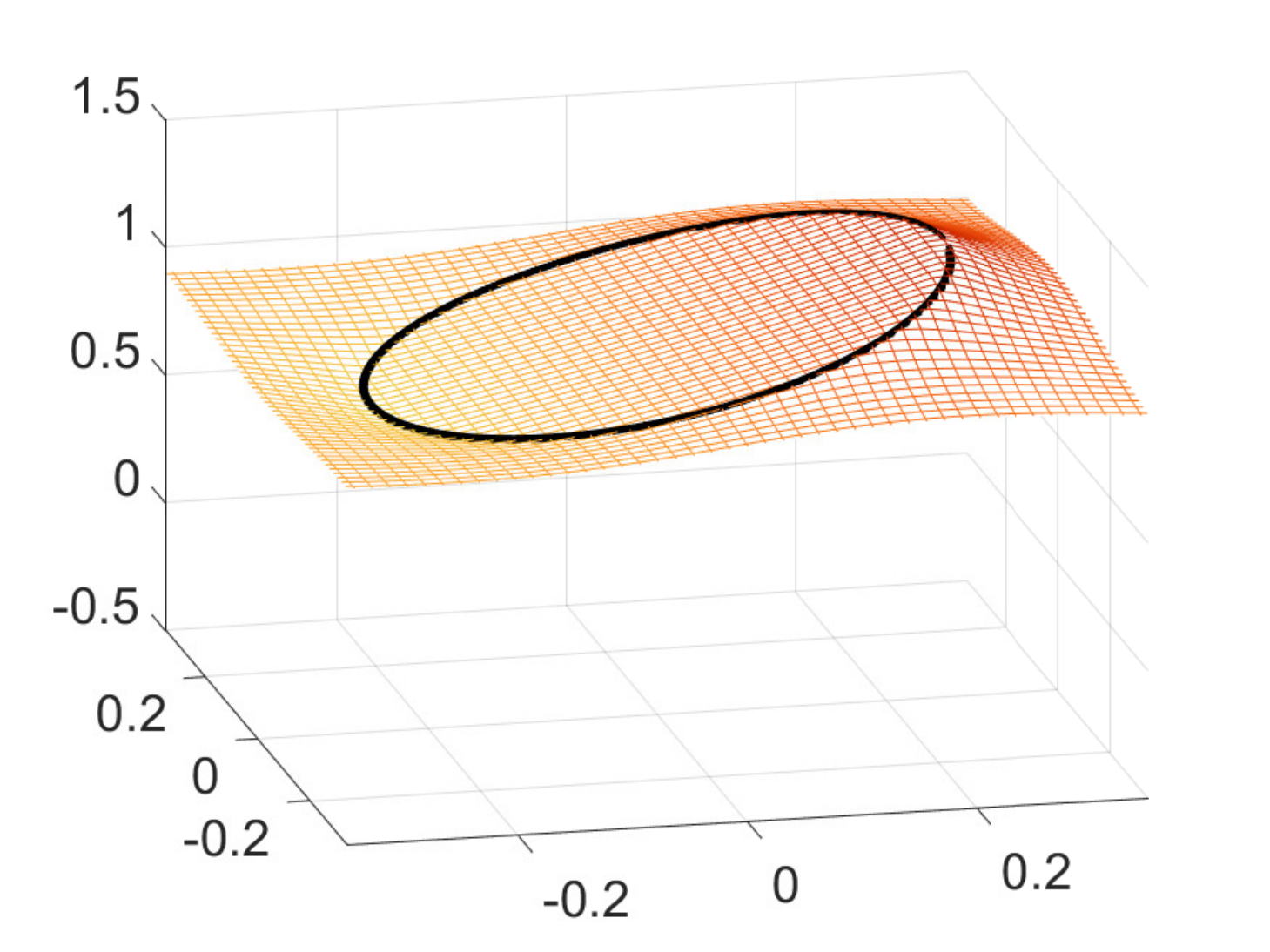}
\caption{\normalsize IBSL Solution}
\label{ibsl solution}
\end{subfigure}
\begin{subfigure}{0.495\textwidth}
\centering
\includegraphics[width=\textwidth]{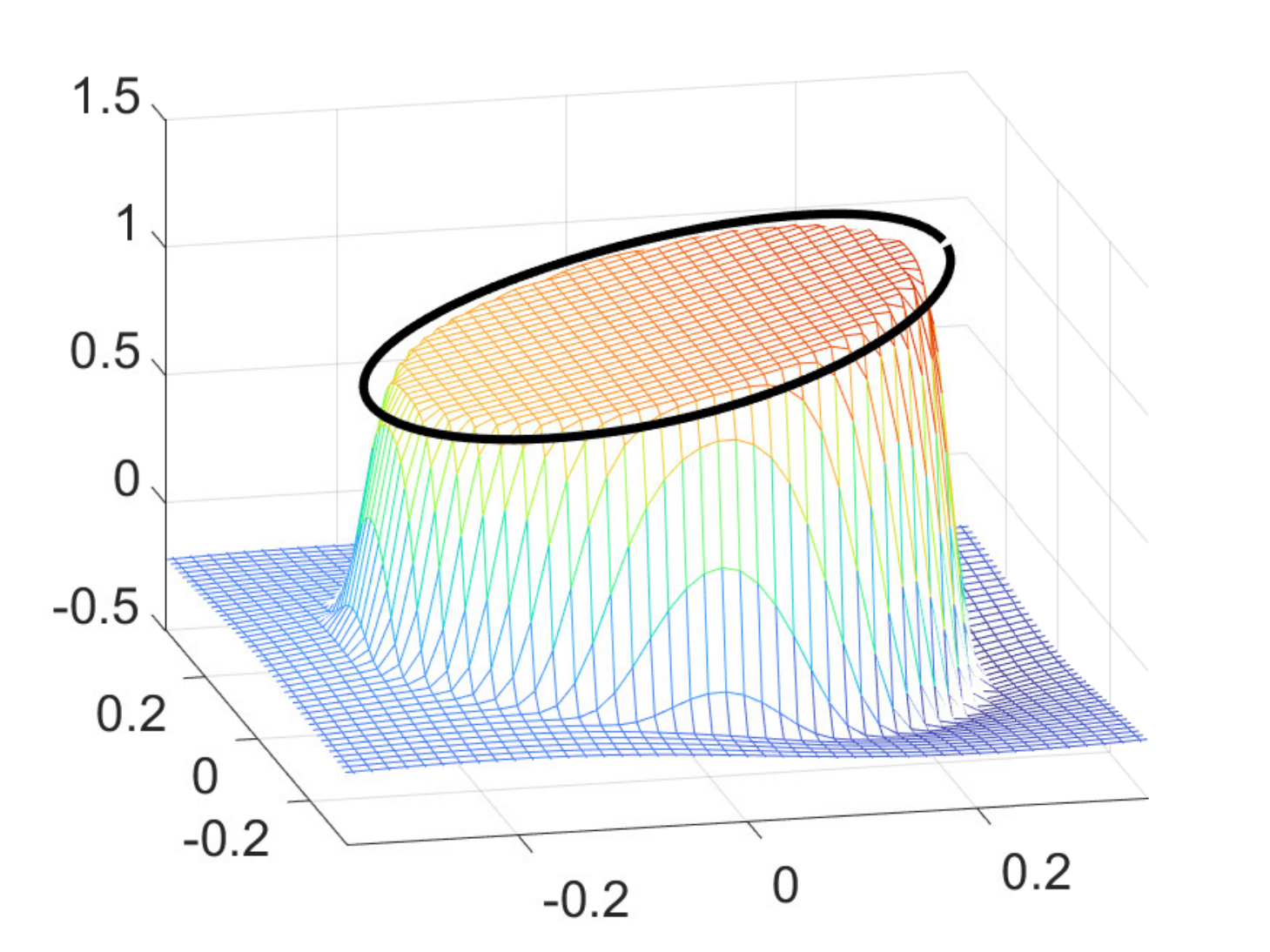}
\caption{\normalsize IBDL Solution}
\label{ibdl solution}
\end{subfigure}
\begin{subfigure}{0.495\textwidth}
\centering
\includegraphics[width=\textwidth]{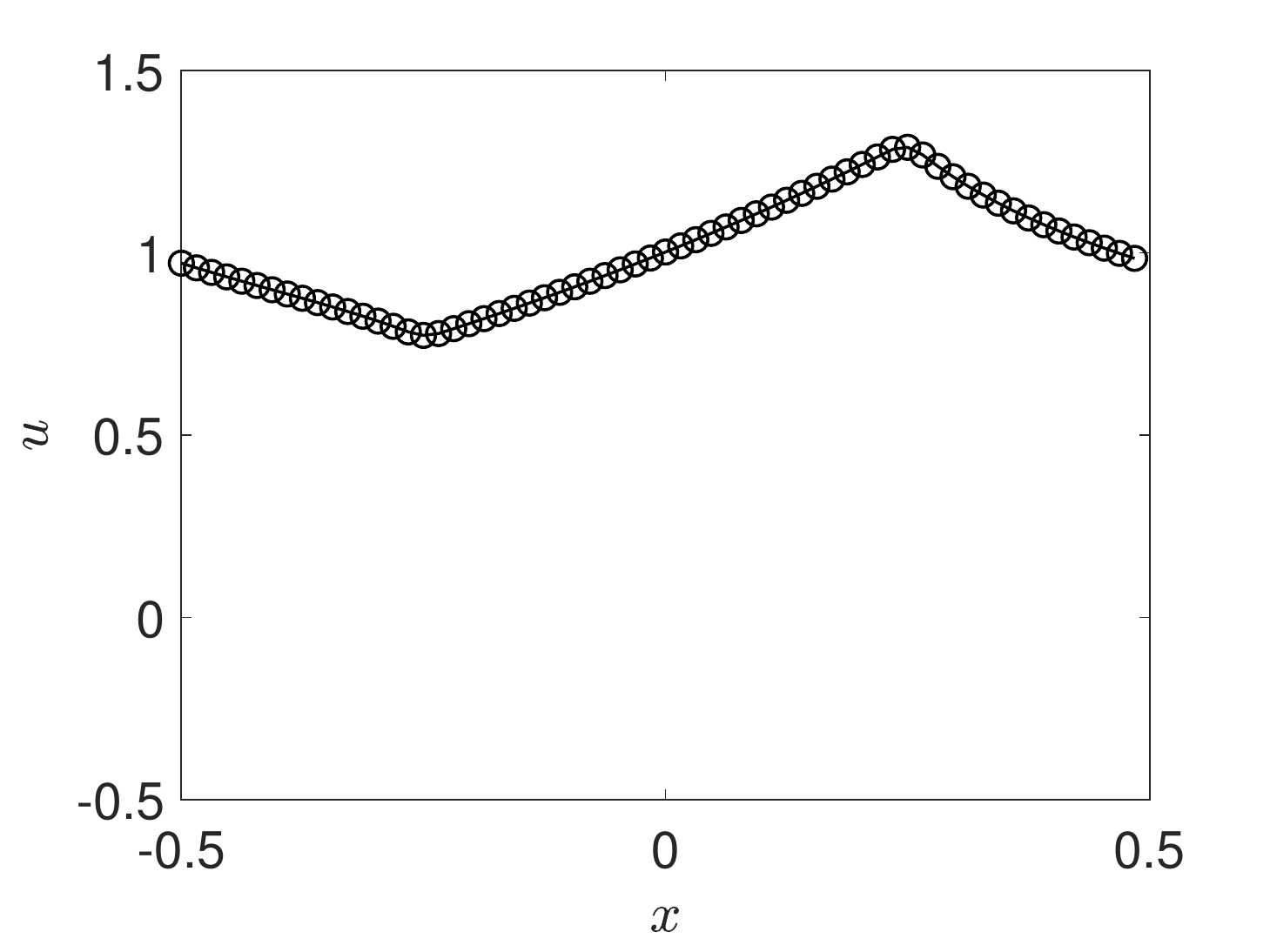}
\caption{\normalsize IBSL Solution}
\label{ibsl slice}
\end{subfigure}
\begin{subfigure}{0.495\textwidth}
\centering
\includegraphics[width=\textwidth]{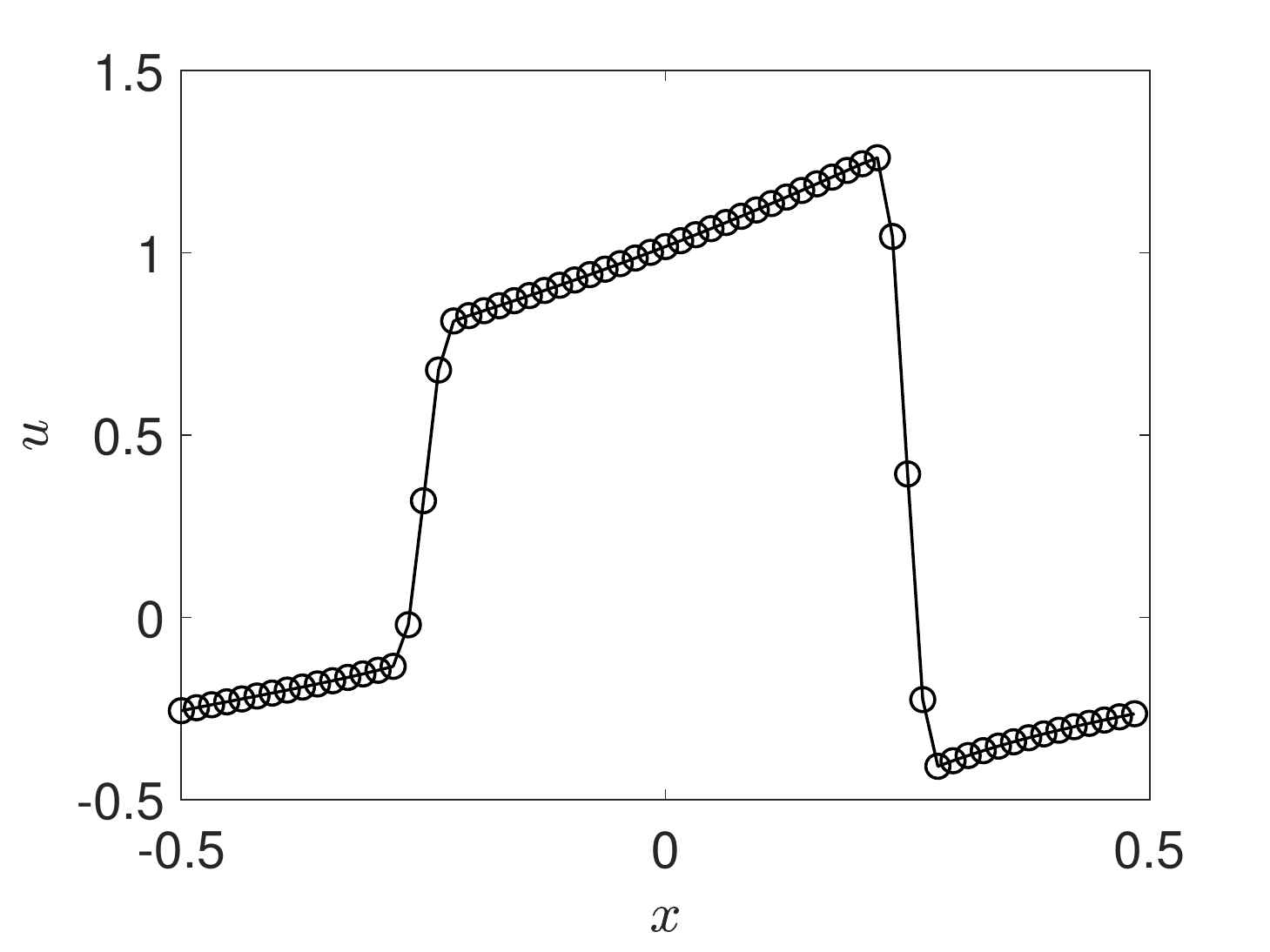}
\caption{\normalsize IBDL Solution}
\label{ibdl slice}
\end{subfigure}
\caption[Solution plots and slices found using the IBSL and IBDL methods to solve Equation \eqref{pde again} demonstrating the discontinuity of the IBDL solution]{Solution plots found using the IBSL and IBDL methods to solve Equation \eqref{pde again} with $k=1$ and $g=0$, and where the circular boundary of radius 0.25 has prescribed boundary values given by $U_b=e^x$. The periodic computational domain is $[-0.5, 0.5]^2$, and the grid and boundary point spacing is $\Delta s \approx \Delta x = 2^{-6}$. Figure \ref{ibsl solution} gives the solution from using the IBSL method, and Figure \ref{ibdl solution}, the IBDL method. $U_b$ is shown with a black curve. Figures \ref{ibsl slice} and \ref{ibdl slice} show solution values for $y=-0.015625$ for the IBSL and IBDL methods, respectively.}\label{solution plots} 
\end{figure}

Since our boundary density $Q$ corresponds to the jump in the solution values across the boundary, clearly our solution $u$ will be discontinuous. In the original IBSL method, the solution was continuous, but the normal derivative was not. Figure \ref{solution plots} illustrates the solutions to Equation \eqref{pde} produced by the IBSL and IBDL methods on the entire computational domain, $\C=[-0.5,0.5]^2$. The PDE domain $\Omega$ is the interior of a circle, and the boundary condition is $U_b=e^x$. We can see that the two solutions match on the portion of $\Omega$ that is about a couple meshwidths away from the boundary. However, the IBDL method gives a solution that is discontinuous across the boundary.

The lack of smoothness in the derivative across the boundary causes the IBSL method to achieve only first-order accuracy. Since the IBDL method will instead yield a discontinuous function, we will not see pointwise convergence near the boundary if we use the solution values given directly from the method in Equation \eqref{ibdl}. We therefore must replace the values of $u$ for grid points \textit{near} the boundary. 

Integral equation methods are generally able to achieve higher accuracy than the Immersed Boundary method. However, this lack of smoothness related to the singularity in the Green's function derivative results in the need to employ analytical techniques in order to achieve the same level of accuracy for points near $\Gamma$. For instance, Beale and Lai use a regularized Green's function and then analytically derive the correction terms \cite{BealeLai}. Kl\"{o}ckner et al. use analytical expansions centered at points several meshwidths from the boundary to evaluate the solution for points nearer to the boundary \cite{QBX}, and Carvalho et al. use asymptotic analysis to approximate the solution at near-boundary points with known boundary data and a nonlocal correction \cite{ShilpaBI}. 

In the Immersed Boundary method, however, one only expects first-order accuracy, and we are able to obtain this level of accuracy with a simple linear interpolation using known boundary values and approximate solution values several meshwidths into the PDE domain. Section \ref{5.3 interp} provides details of the numerical implementation. 

\begin{figure}
\centering
\begin{subfigure}{0.495\textwidth}
\centering
\includegraphics[width=\textwidth]{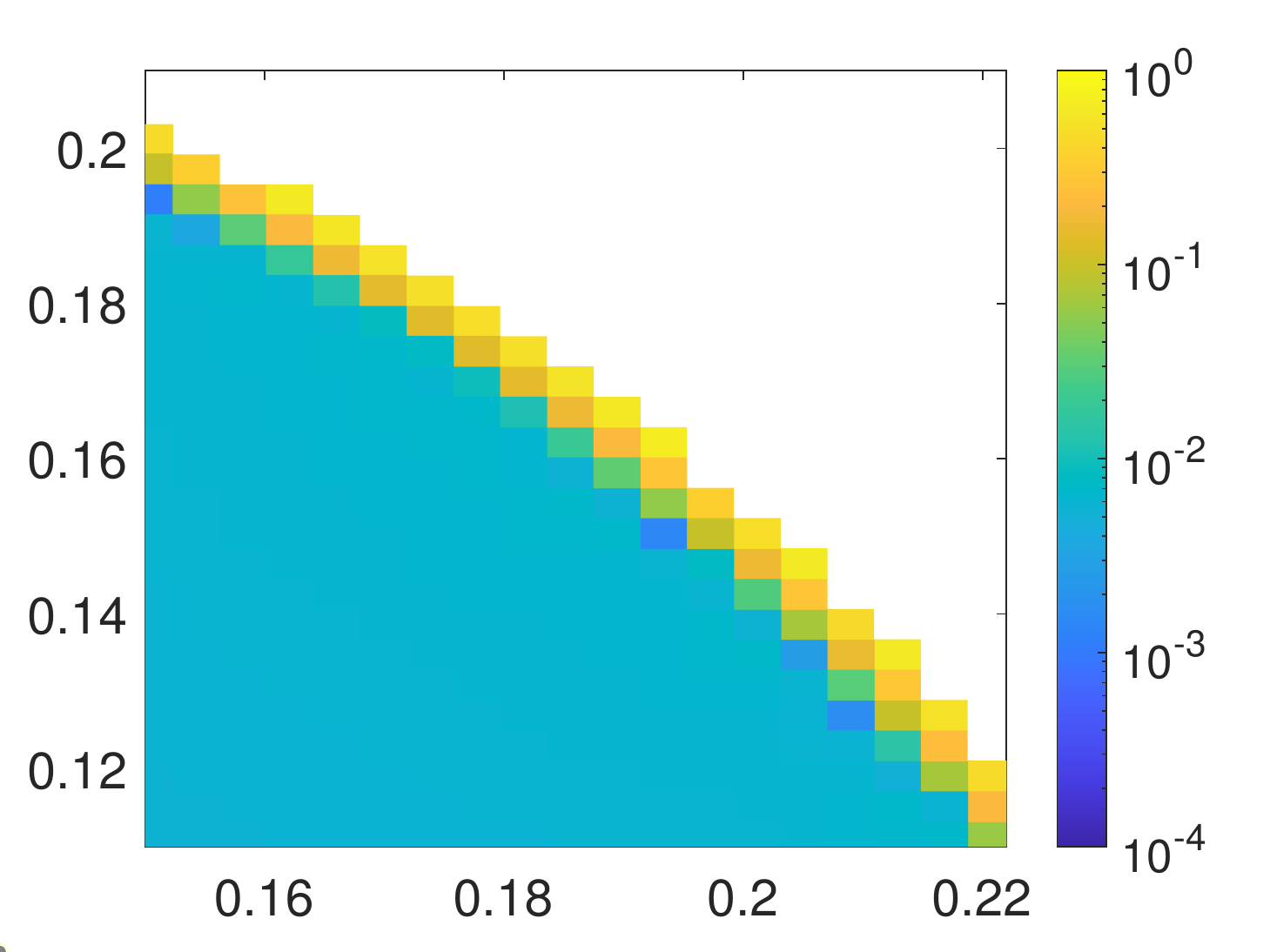}
\caption{\normalsize $\Delta x = 2^{-8}$}
\label{course}
\end{subfigure}
\begin{subfigure}{0.495\textwidth}
\centering
\includegraphics[width=\textwidth]{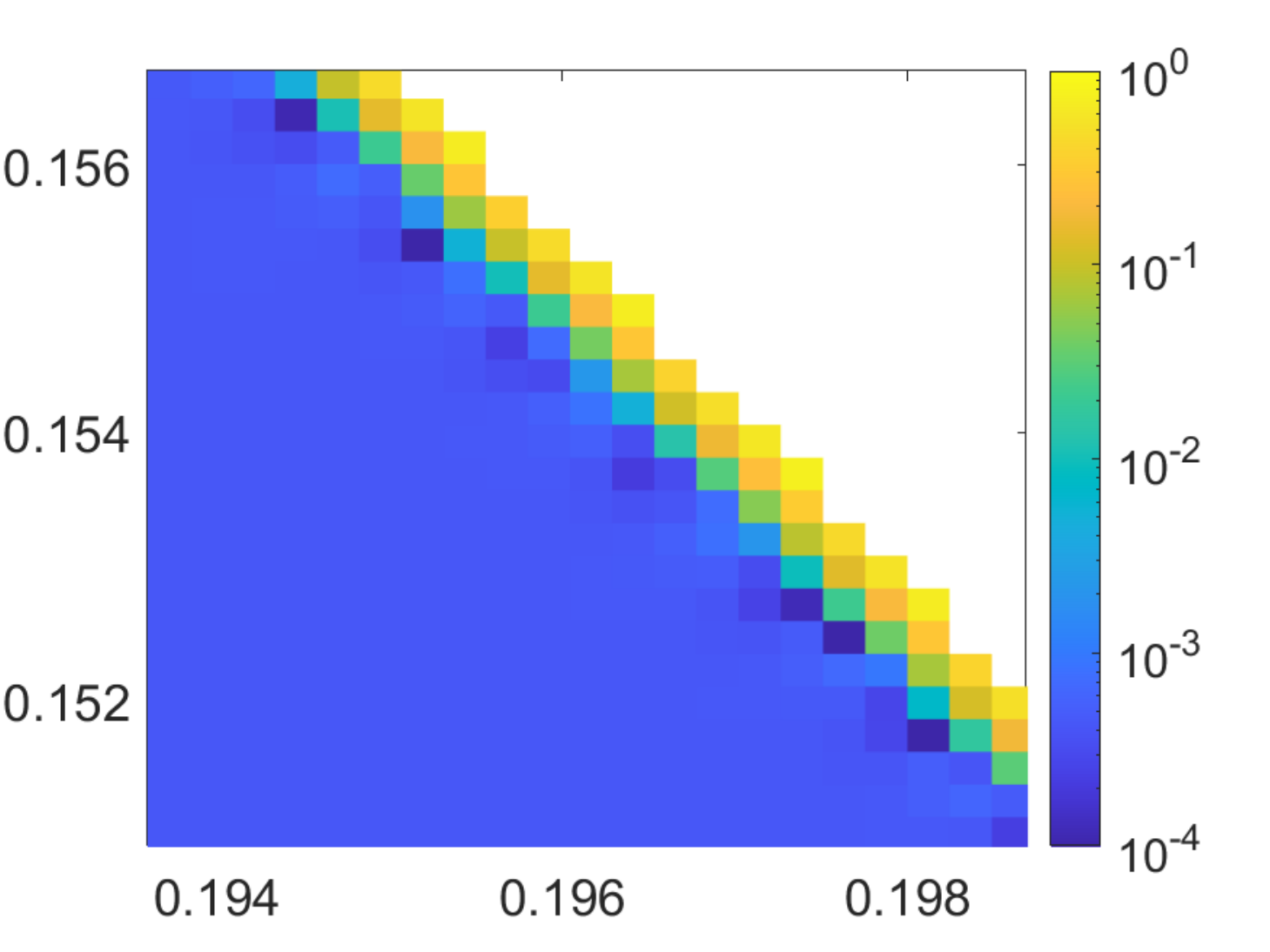}
\caption{\normalsize  $\Delta x = 2^{-12}$}
\label{fine}
\end{subfigure}
\caption[Error plots for the IBDL method without interpolation, showing large errors within 2-3 meshwidths of the boundary]{Plots of errors from using IBDL method to solve Equation \eqref{pde again} with $k=1$ and $g=0$ on the periodic computational domain $[-0.5, 0.5]^2$, where the circular boundary of radius 0.25 has prescribed boundary values given by $U_b=e^x$. The errors are shown for two grid point spacings and with $\Delta s \approx \Delta x$. The figure is zoomed-in to view the large errors that lie in a small region near $\Gamma$.}\label{near boundary zoomed in}
\end{figure}

There are two regions of grid points for which we do not see pointwise convergence. The first region contains the largest errors, and it is illustrated in Figure \ref{near boundary zoomed in} for an interior circular domain with boundary data $U_b=e^x$. We can see that this region remains localized to only a small number of grid points, so that as we refine the grid, the width of this region quickly goes to 0. These large errors are the direct result of the smoothing of the dipole forces, and the number of meshwidths is therefore determined by the support of the regularized delta function used. As discussed in Section \ref{2.2 spread}, we generally use a 4-point delta function, and this results in about $2-3$ meshdwidths of large errors on one side of the boundary, which can be seen in Figure \ref{near boundary zoomed in}. 

Secondly, the numerical method used to discretize the PDE spreads the error from the discontinuity into a larger region near the boundary. Therefore, if we use the interpolation to only correct for the large errors within the first few meshwidths, we will still see a lessening of first-order convergence in the max norm for sufficiently fine grids. The width of the region on which the solution fails to converge pointwise still approaches 0, but the number of meshwidths affected can increase. One can recover pointwise convergence with further interpolation, and in practice, one can generally obtain a maximum error comparable to that of the IBSL method by using interpolation for only a relatively small number of meshwidths. In Section \ref{ch 5 max norm}, we further discuss factors that contribute to the number of meshwidths needed.

The full IBDL method for Poisson and Helmholtz equations can now be described in the following steps:
{\begin{enumerate}[topsep=2ex]
\item Use a Krylov method to solve the following system for $Q$. 
\begin{equation} 
-(S^* \L^{-1}\widetilde S)Q + \frac12 Q=U_b - S^*\L^{-1} \tilde g\label{ibdl schur}
\end{equation}
\item Use $Q$ to solve the following equation for $u$ in $\C$. 
 \begin{equation} 
\L u + \widetilde S Q=\tilde g \label{ibdl 1 again}
\end{equation}
\item Use an interpolation to replace the values of $u$ at $\x \in \Omega$ that are near $\Gamma$. Factors that determine the number of meshwidths to include in this step are discussed in Section \ref{ch 5 max norm}.
\end{enumerate}
}

\subsection{Arbitrary immersed boundary parametrization}\label{5.1 param}

In this dissertation, we utilize an arclength parametrization of the immersed boundary $\Gamma$. Numerically, this means that in the discretization of the spread operator, $\Delta s_i$ is the size of the arclength connecting $\X(s_i)$ and $\X(s_{i+1})$. In practice, this choice is not a limiting factor because even if an analytical arclength parametrization of the immersed boundary curve is not available, one may still approximate this value with the distance $\Delta s_i\approx ||\X(s_{i+1})-\X(s_i)||_2$. However, to be thorough, in this section, we discuss the use of an arbitrary parametrization. 

Let $\X(s)$ be an arbitrary parametrization of $\Gamma$. Using this parametrization and the derivations from Chapter \ref{chapter 3}, the double layer integral equation for Poisson would have the form 
\begin{equation}
U_b(\X(s'))= \int_{\Gamma}^{PV} \gamma(\X(s)) \grad G(\X(s),\X(s'))\dotp\n(\X(s)) \Big|\frac{\partial\X}{ds}\Big|d s-\frac12 \gamma(\X(s')).\label{param1}
\end{equation}

In this case, we can form the IBDL method as
\begin{subequations}\label{ibdl param}
\begin{alignat}{2}
& \L u +\widetilde S Q = \tilde g \qquad && \text{in } \mathcal{C}  \label{ibdl 1 param}\\
&S^* u + \frac{1}{2|\partial\X/\partial s|} Q = U_b \qquad && \text{on } \Gamma.   \label{ibdl 2 param}
\end{alignat}
\end{subequations}
Then, using the same steps as in Section \ref{5.1 formulation}, we would obtain the equation 
\begin{equation}
U_b(\X(s'))= -\int_{\Gamma}^{PV} Q(\X(s)) \grad G_{hh}(\X(s),\X(s'))\dotp\n(\X(s)) ds+ \frac{1}{2|\partial\X/\partial s|}Q(\X(s')).\label{param2}
\end{equation}
By comparing Equations \eqref{param1} and \eqref{param2}, we see that in the limit $h\longrightarrow 0$, we get 
\begin{equation}
Q(\X)=-\gamma(\X) \Big|\frac{\partial\X}{\partial s}\Big|.
\end{equation}

As discussed in Section \ref{ch 3 discussion}, the constant $1/2$ shifts all of the eigenvalues of the integral operator so that the only limit point is away from $0$. This ensures a small, constant condition number as the mesh is refined and therefore small iteration counts for the Krylov method. An arbitrary paramatrization may not shift the eigenvalues in this uniform manner. As such, the iteration counts may increase. In Section \ref{5.4 hh}, we will revisit the Helmholtz equation from Section \ref{2.4 Poisson} given by 
\begin{subequations} 
\begin{alignat}{2}
& \Delta u -  u = 0 \qquad && \text{in } \Omega  \\
&u=\sin{2\theta} \qquad && \text{on } \Gamma, 
\end{alignat}
\end{subequations}
where $\Omega$ is the interior of a circle of radius 0.25, centered at the origin. We will see that the iteration count for solving this equation with the IBDL method using arclength parametrization is about 4-5 iterations. If we instead use the parametrization given by 
\begin{equation}
\begin{pmatrix} x \\y\end{pmatrix} = \begin{pmatrix} 0.25 \cos{(s^2+s)} \\ 0.25\sin{(s^2+s)}\end{pmatrix}, 
\end{equation}
for $0\leq s\leq \frac12(-1+\sqrt{1+8\pi})$, we have 
\begin{equation}
\Big|\frac{\partial \X}{\partial s}\Big|=0.25(2s+1).
\end{equation}
Using immersed boundary points that are equally spaced in the parameter $s$, the iteration count is about 20. The iteration count is larger, but it does remain constant as the mesh is refined. This is likely to be the case with an arbitrary parametrization as the eigenvalues will still be shifted in some manner away from $0$. It is therefore acceptable to implement the IBDL method in this manner, but we then need the values or approximations of $|\partial \X/\partial s|$. To maintain consistency in the formulation of the IBDL method, we will use an arclength parametrization, with equally spaced points, in all other sections of this dissertation.

\section{Neumann boundary conditions}\label{ch 5 neumann}

We now use the IB framework and our connection to boundary integral equations to solve a PDE with Neumann boundary conditions. The original IBSL method is unable to handle such a problem since the solution derivatives are not convergent at the boundary. 

Let us look at the PDE
\begin{subequations} \label{pde neumann}
\begin{alignat}{2}
& \L u = g \qquad && \text{in } \Omega  \label{pde1 neumann}\\
&\frac{\partial u}{\partial n} =V_b \qquad && \text{on } \Gamma,  \label{pde2 neumann}
\end{alignat}
\end{subequations}
where we will first take $g=0$ for simplicity. In Section \ref{ch 3 integral eqns}, we introduced the full integral representation of $u$. For $\x_0 \in \Omega\setminus \Gamma $, we have
\begin{equation}
u(\x_0) = \int_{\Gamma}  \grad u(\x) \dotp \n (\x) G(\x, \x_0)dl(\x) - \int_{\Gamma} u(\x) \grad G(\x,\x_0) \dotp \n(\x) dl(\x). \label{integral rep again}
\end{equation}  
Then, we can take $\x_0\longrightarrow \Gamma$ in Equation \eqref{integral rep again} and denote it as $\X_0$. Using the property of the Green's function derivative from Section \ref{ch 3 Greens}, we get
\begin{equation}
u(\X_0) = \int_{\Gamma}  \grad u(\x) \dotp \n (\x)G(\x, \X_0) dl(\x) - \int_{\Gamma}^{PV} u(\x) \grad G(\x,\X_0) \dotp \n(\x) dl(\x) +\frac{1}{2} u(\X_0). \label{integral eqn again}
\end{equation} 
The first term in Equation \eqref{integral eqn again} is a single layer potential with the strength given by the known boundary derivatives, $V_b$. The second term is a double layer potential whose strength is given by the \textit{unknown} boundary values. If we let $U_b\equiv u\big|_{\Gamma}$ be the unknown boundary values, we can then write the PDE in the IB framework as
\begin{equation}
 \L u +\widetilde S U_b+SV_b = 0 \label{homog neumm 1}
\end{equation}

Then, recall that the second IBDL equation in the Dirichlet case is given by 
\begin{equation}
S^*u+\frac12 Q = U_b. \label{Dir sstar}
\end{equation}
In the Neumann case, since $U_b$ is our unknown boundary distribution, we replace $Q$ in Equation \eqref{Dir sstar} with $U_b$, and we get the corresponding equation in the Neumann case. Allowing for a non-zero $g$, we can therefore write the PDE in the Immersed Boundary framework as 
\begin{subequations} \label{ib neumann}
\begin{alignat}{2}
& \L u +\widetilde S U_b+SV_b = \tilde g \qquad && \text{in } \mathcal{C}  \label{ib neumann 1}\\
&S^* u = \frac12 U_b \qquad && \text{on } \Gamma.   \label{ib neumann 2}
\end{alignat}
\end{subequations}
For Neumann boundary conditions, $U_b$ is the \textit{unknown} potential strength on the boundary, corresponding to the unknown boundary values. 

We can solve this system by first solving the equation 
\begin{equation}
(-S^* \L^{-1}\widetilde S)U_b-\frac12 U_b=S^*\L^{-1}SV_b-S^*\L^{-1}\tilde g \label{Neumann saddle point}
\end{equation} 
for the boundary values, $U_b$, and then obtaining $u$ from Equation \eqref{ib neumann 1}. We can see by comparing to Equation \eqref{ibdl schur} that the operator is similar to that of the IBDL method for Dirichlet boundary conditions, where the only difference is the sign on the $1/2$. Therefore, we again get a well-conditioned problem that can be solved with a small number of iterations of a Krylov method. 

We now provide a more detailed demonstration of the connection between Equation \eqref{integral eqn again} and Equation \eqref{ib neumann}. We use steps similar to those used in Sections \ref{4.1 connect} and \ref{5.1 formulation}.  

Starting with the Equation \eqref{ib neumann 1} and using $\tilde g=0$ for simplicity, we have 
\begin{equation}
\L u =- S V_b - \widetilde S U_b =- \int_{\Gamma} V_b(s) \delta_h(\x-\X(s)) ds - \grad \dotp \int_{\Gamma} U_b(s)\n(s) \delta_h(\x-\X(s)) ds  . \label{neum connect 1}
\end{equation}
Inverting the operator and using that $\L G_h(\x,\x_0)=-\delta_h(\x-\x_0)$, we get
\begin{equation}
u(\x) =  \int_{\Gamma} V_b(s) G_h(\x,\X(s)) ds+\grad \dotp  \int_{\Gamma} U_b(s)\n(s)G_h(\x,\X(s)) ds   \label{neum connect 2}
\end{equation}
for $\x\in\Omega\setminus \Gamma$. Bringing in the divergence and manipulating the expression, we get
\begin{equation}
u(\x) =  \int_{\Gamma} V_b(s) G_h(\x,\X(s)) ds+  \int_{\Gamma} U_b(s) \grad G_h(\x,\X(s)) \dotp\n(s) ds. \label{neum connect 3}
\end{equation}
The second equation of the method, Equation \eqref{ib neumann 2}, gives us 
\begin{equation}
\frac12 U_b(s') =\int_{\C} u(\x)\delta_h(\x-\X(s')) dx. \label{neum connect 4}
\end{equation}
Combining Equation \eqref{neum connect 3} with Equation \eqref{neum connect 4}, changing the order of integration, and recognizing the presence of $G_{hh} = G_h*\delta_h$, we get 
\begin{equation}
\frac12 U_b(s') =\int_{\Gamma} V_b(s) G_{hh}(\X(s'),\X(s)) ds+\int_{\Gamma} U_b(s) \grad G_{hh}(\X(s'),\X(s)) \dotp\n(s) ds . \label{neum connect 5}
\end{equation}
Note that due to the presence of a double layer potential, as discussed in Section \ref{5.1 formulation}, we assume that $\X(s)$ is the arclength parametrization of $\Gamma$ so that $  \big|\partial \X(s)/\partial s\big| =1$.

The symmetry of the Green's function, which is preserved through convolutions with the regularized delta function, gives us that $G_{hh}(\X(s'),\X(s)) = G_{hh}(\X(s),\X(s'))$. Additionally, the odd symmetry of the gradient of the Green's function gives us  $ \grad G_{hh}(\X(s'),\X(s))=- \grad G_{hh}(\X(s),\X(s'))$ \cite{Pozblue}. Using these to switch the arguments of these functions, we get
\begin{equation}
\frac12 U_b(s') =\int_{\Gamma} V_b(s) G_{hh}(\X(s),\X(s')) ds -\int_{\Gamma} U_b(s) \grad G_{hh}(\X(s),\X(s')) \dotp\n(s) ds . \label{neum connect 6}
\end{equation}
Appropriately redefining $U_b$, $V_b$, and $\n$ as functions of $\x$, we get 
\begin{multline}
\frac12 U_b(\X(s')) =\int_{\Gamma} V_b(\X(s)) G_{hh}(\X(s'),\X(s)) ds\\
-\int_{\Gamma} U_b(\X(s)) \grad G_{hh}(\X(s'),\X(s)) \dotp\n(\X(s)) ds . \label{neum connect 7}
\end{multline}
 
Using our arclength parametrization in Equation \eqref{integral eqn again} and combining the $u(\X_0)$ terms, we can write the integral equation for $\X_0$ on the boundary as
\begin{multline}
\frac12 u(\X(s')) =\int_{\Gamma}  \grad u(\X(s)) \dotp \n (\X(s))G(\X(s), \X(s')) ds\\
- \int_{\Gamma}^{PV} u(\X(s)) \grad G(\X(s),\X(s')) \dotp \n(\X(s)) ds  . \label{neum connect 8}
\end{multline}
By comparing Equations \eqref{neum connect 7} and \eqref{neum connect 8}, we can see that the IB formulation presented in Equation \eqref{ib neumann} is equivalent to a regularized integral equation.

\section{Numerical implementation}\label{ch 5 numerical}

We now discuss the numerical implementation of the IBDL method, particularly the areas in which the method differs from the IBSL method. There is no change in the discretization of space, the differential operators, and the basic spread and interpolation operators, and this information can be found in Section \ref{ch 2 numerical implementation}.

\subsection{Discretization of spread operator}\label{5.3 spread}

We will first discuss the implementation of the operator $\widetilde S$ defined by 
\begin{equation}
\widetilde S Q = \grad \dotp (S Q\n). \label{Stilde}
\end{equation}
This operation consists of first multiplying our discretized boundary density $Q$ by the outward unit normal vectors  at the immersed boundary points. We then spread these vectors to the grid and take the divergence. We discuss the implementation of each of these steps below. 

The first component of the operator $\widetilde S$ is the divergence. To compute this, we use either a finite difference or Fourier spectral method to match the method used to discretize the PDE. For finite differences, we use the standard second-order accurate centered difference for derivatives, and for the Fourier spectral method, we use the derivatives discussed in Section \ref{2.3 space and op}. 

Next, we need the unit normal vectors to the immersed boundary. These can be treated as an input to the method, or we can approximate them if they are unavailable. Approximating the curve with linear elements and using the resulting normal vectors is sufficient for maintaining first-order accuracy. For instance, for an interior domain and a counterclockwise curve parametrization, we can use 
\begin{equation}
\n(s_i)= \begin{bmatrix}Y(s_{i+1})-Y(s_{i-1})\\   -\big(X(s_{i+1})-X(s_{i-1})\big)\end{bmatrix} \Bigg/ \big|\big|\X(s_{i+1})-\X(s_{i-1})\big|\big|_2.\label{normal vector approximations}
\end{equation}
For an exterior domain, we would negate this.

We next look at the spread operator $S$, which is defined as
\begin{equation}
(SF)(\x) =\int_{\Gamma} F(s)\delta_h(\x-\X(s))ds .    \label{spread operator again}
\end{equation}
As discussed in Section \ref{ch 2 numerical implementation}, we use the traditional Peskin four-point delta function for $\delta_h$ \cite{Peskin02}, unless otherwise specified.

Then, as mentioned in Section \ref{ch 5 IBDL}, we discretize the integral in Equation \eqref{spread operator again} with respect to arclength, so henceforth, I will treat $s$ as an arclength parameter. With this new assumption on the coordinate $s$, $\Delta s_i$ gives the length of the curve connecting the immersed boundary points, $\X(s_i)$ and $\X(s_{i+1})$. Then, the discrete spread operator is given by
\begin{equation}
(SF)(\x) =\sum_{i=1}^{N_{IB}} F(s_i)\delta_h(\x-\X(s_i))\Delta s_i .   
\end{equation}

Depending on the parametrization, the boundary points may not be equally spaced, but since we will be looking at a periodic distribution on $\Gamma$, choosing equally spaced points will give spectral accuracy in our integral approximation. Unless otherwise specified, we will space the boundary points equally with $\Delta s \approx \alpha \Delta x$ for various values of $\alpha$. Then $\Delta s = L_{IB}/N_{IB}$, where $L_{IB}$ gives the length of the immersed boundary. If the exact values of $\Delta s_i$ are unavailable, approximations can be made. Again, using $\Delta s_i \approx ||\X(s_{i+1})-\X(s_i)||_2$ is sufficient to maintain first-order accuracy.

\subsection{Solution to discrete system for invertible $\L$  }\label{5.3 invertible L}

We now look at the discretized system 
\begin{subequations} \label{ibdl again}
\begin{alignat}{2}
& \L u +\widetilde S Q = \tilde g \qquad && \text{in } \mathcal{C}  \label{ibdl again 1}\\
&S^* u + \frac12 Q = U_b \qquad && \text{on } \Gamma.   \label{ibdl again 2}
\end{alignat}
\end{subequations}

In the case that the differential operator is invertible, such as for $\L=\Delta -k^2$, for $k\neq0$, we can invert the operator to obtain
\begin{equation}
u=-\L^{-1}\widetilde S Q+\L^{-1}\tilde g  .
\end{equation}
Then, by applying the interpolation operator and using Equation \eqref{ibdl again 2}, we obtain
\begin{equation}
-(S^* \L^{-1}\widetilde S)Q + \frac12 Q=U_b - S^*\L^{-1} \tilde g . \label{ibdl schur again}
\end{equation}
We can therefore solve Equation \eqref{ibdl again} by first solving Equation \eqref{ibdl schur again} for $Q$ and then obtaining $u$ from Equation \eqref{ibdl again 1}. Here, the operator that must be inverted, $-(S^* \L^{-1}\widetilde S) + \frac12 \mathds{I}$, is not symmetric, and we therefore solve for $Q$ with \texttt{gmres}, to a tolerance of $10^{-8}$, unless otherwise specified.

\subsection{Solution to discrete system for Poisson equation  }\label{5.3 poisson}

In the case that the differential operator is the periodic Laplacian, $\L=\Delta$, the solution method outlined in the previous section must be adjusted, both for the nullspace of the discrete constraint system and the nullspace of $\L$ . 

\textbf{Nullspace of Constraint System. } In the case of the IBSL method applied to the Poisson equation, the constraint matrix, given by 
\begin{equation}
\begin{pmatrix} \Delta &S\\ S^* &0\end{pmatrix},
\end{equation}
is invertible. By handling the nullspace of $\Delta$ as described in Section \ref{2.2 Laplacian L}, we can therefore recover the unique solution to the constraint system. In the case of the IBDL method, however, the constraint matrix, given by
\begin{equation}
\begin{pmatrix} \Delta &\widetilde S\\ S^* &\frac12\mathbb{I}\end{pmatrix},
\end{equation}
has a one-dimensional nullspace. The analytical nullspace is spanned by 
\begin{subequations}\label{nullspace shit}
\begin{alignat}{2}
&u=\chi_{\scaleto{\C \setminus \Omega}{6pt}}\\
&Q=-1,
\end{alignat}
\end{subequations}
where $\chi_{\scaleto{\C \setminus \Omega}{6pt}}$ gives the indicator function that is $0$ on $\Omega$ and $1$ otherwise. Section \ref{5.3 flagging} demonstrates the use of a similar solution to create an indicator function. According to Equation \eqref{nullspace shit}, since an element of the nullspace would be $0$ on the physical domain, any solution we get would be the same on $\Omega$, and this nullspace would therefore not pose a problem for the IBDL method. We have observed this to be true for interior domains or small exterior domains. However, we have observed that for some boundary value problems on exterior domains  for which the proportion of the periodic box occupied by $\Omega$ is large, the IBDL method outlined does not converge to the correct solution. We conjecture that this may be explained by two facts. Firstly, while Equation \eqref{nullspace shit} gives the \textit{analytical} nullspace, our regularization and discretization results in a different numerical nullspace. Additionally, in the infinite exterior domain case, the decay of the double layer potential means that it is unable to represent solutions to the Poisson equation with arbitrary boundary values \cite{periodic}. As we increase the periodic domain length $L$, we may therefore approach an analytical problem that does not have a solution, and the contribution from an element of our numerical nullspace may be larger. We give a solution to this issue below, but more investigation is needed into the subject of this nullspace.

For large exterior domains, therefore, we may use the completed double layer representation, which was discussed in Section \ref{3.5 DL}, as it relates to Stokes equation. We can formulate this completed representation by the addition of a single layer potential whose strength is a constant multiple of the strength of the double layer potential. The completed system is then given by 
\begin{subequations} \label{ibdl completed}
\begin{alignat}{2}
& \Delta u +\eta S Q+ \widetilde S Q = \tilde g \qquad && \text{in } \mathcal{C}  \label{ibdl 1 completed}\\
&S^* u + \frac12 Q = U_b \qquad && \text{on } \Gamma,  \label{ibdl 2 completed}
\end{alignat}
\end{subequations}
where $\eta$ is a positive constant. Appropriate choices for $\eta$ are discussed in Section \ref{6.5 eta}. We will refer to this form of the IBDL method as the \textit{completed IBDL method}. In this case, the constraint matrix, given by 
\begin{equation}
\begin{pmatrix} \Delta &\widetilde S+\eta S\\ S^* &\frac12\mathbb{I}\end{pmatrix},
\end{equation}
is invertible, and then we need only deal with the nullspace of the periodic Laplacian. In Section \ref{5.4 poisson shit}, we examine a problem for which the completed IBDL method is needed.

\textbf{Nullspace of $\L$. } As discussed in Section \ref{2.2 Laplacian L}, the original boundary value problem on $\Omega$ has a unique solution, but, by using the IB framework, we embed the PDE into a periodic computational domain on which $\Delta$ is not invertible. 

As we did in Section \ref{2.2 Laplacian L}, let us begin by decomposing the solution $u$ as 
\begin{equation}
u=u_0+\bar u,
\end{equation}
where $u_0$ has mean $0$ on $\C$ and $\bar u$ gives the mean value of $u$ on $\C$. The PDE on $\C$ that results from the completed IBDL formulation in Equation \eqref{ibdl completed} is given by 
\begin{equation}
\Delta u + \eta SQ +\grad\dotp (SQ\n)  =\tilde g, \label{laplaces with ibdl}
\end{equation}
where $\tilde g \equiv g\chi_{\scaleto{\Omega}{4.5pt}} + g_e \chi_{\scaleto{\C \setminus \Omega}{6pt}}$ is an extension of $g$. To derive the solvability condition, we integrate Equation \eqref{laplaces with ibdl} over a general computational domain $\C$. Noting again that our spread operator is defined using a regularized delta function, the functions are smooth, and we can use the divergence theorem. Let us again use $\partial \C$ as the boundary of $\C$ and $\n_{\C}$ as the unit normal on $\partial \C$, to distinguish it from $\n$, which we continue to use as the unit normal on the immersed boundary $\Gamma$. We then get
\begin{equation}
\int_{\partial \C} \grad u \dotp \n_{\C} dl(\x)  +\eta \int_{\C} SF d\x+ \int_{\partial \C} S Q \n\dotp \n_{\C} dl(\x) = \int_{\C} \tilde g d\x.  \label{first step after divergence theorem again}
\end{equation}
If $\Gamma$ is away from $\partial \C$, the third term vanishes due to the compact support of the integrand. Then, if we take $\C$ to be the periodic box used in this work, the first term also disappears. Using that $\int_{\C} SF d\x=\int_{\Gamma}Fds$, the solvability condition on the periodic computational domain $\C$ then has the form
\begin{equation}
\eta \int_{\Gamma} Qds=\int_{\C}\tilde g(\x) d\x .\label{poisson completed solv}
\end{equation}
When we are not using the completed formulation, $\eta=0$, and we can easily satisfy this condition with our choice of $ g_e$. For example, we can choose 
\begin{equation}
g_e(\x) =-\frac{1}{|\C\setminus \Omega|} \int_{\Omega} g(\x)d\x, 
\end{equation} 
where we use $|\dotp | $ to denote the area of the domain. With this constraint satisfied, the solution to Equation \eqref{laplaces with ibdl} can be found, and it is unique up to an additive constant. Therefore, since we are not interested in the solution on the nonphysical domain, $\C\setminus \Omega$, we can use the solution with mean $0$ on $\C$. If we select a solution with a different mean, the potential strength $Q$ will adjust to enforce the boundary condition on the physical side of $\Gamma$, maintaining the unique solution on the PDE domain $\Omega$. 

If, on the other hand, we use the completed IBDL method, we can use the constraint in Equation \eqref{poisson completed solv} to find the value of $\bar u$, as we did in the IBSL method in Section \ref{2.2 Laplacian L}. Discretizing the constraint, we get the following system of equations. 
\begin{subequations} \label{ibdl completed again }
\begin{alignat}{2}
& \Delta u_0 +\eta S Q+ \widetilde S Q = \tilde g \\
&S^* u_0 + \bar u \mathds{1}_{N_{IB}}+ \frac12 Q = U_b\\
& \eta(\Delta s) \mathds{1}^\intercal_{N_{IB}}Q = (\Delta x\Delta y)\mathds{1}^\intercal_{N^2} \tilde g.
\end{alignat}
\end{subequations}

Using the operator previously defined to return a solution with mean $0$ on $\C$, we use \texttt{gmres} to solve  the following system of discrete equations for $Q$ and $\bar u$.
\begin{subequations}
\begin{alignat}{2}
&-(S^* \Delta_0^{-1}\widetilde S)Q -\eta(S^*\Delta_0^{-1}S)Q+\bar u \mathds{1}_{N_{IB}} +\frac12 Q=U_b - S^*\Delta_0^{-1} \tilde g\\
&\eta(\Delta s) \mathds{1}^\intercal_{N_{IB}}Q = (\Delta x\Delta y)\mathds{1}^\intercal_{N^2} \tilde g
\end{alignat}
\end{subequations}
We can then obtain $u$ by computing
\begin{equation}
u=-\Delta_0^{-1}(\widetilde S +\eta S)Q +\Delta_0^{-1}\tilde g+\bar u \mathds{1}_{N^2}.
\end{equation}

\subsection{Interpolation for near-boundary points  }\label{5.3 interp}

As discussed in Section \ref{5.1 discontinuity}, we use a linear interpolation to replace solution values within $m_1$ meshwidths of the boundary $\Gamma$. For the Dirichlet problem, the interpolation for an individual gridpoint uses one approximate solution value located about $m_2$ meshwidths away from $\Gamma$ and one value approximated from the known solution values located on $\Gamma$. For the Neumann problem, the solution values located on $\Gamma$ are not known, and therefore they are also approximate. In this section, we discuss the details of this linear interpolation.

\begin{figure}
\centering
\includegraphics[width=0.6\textwidth]{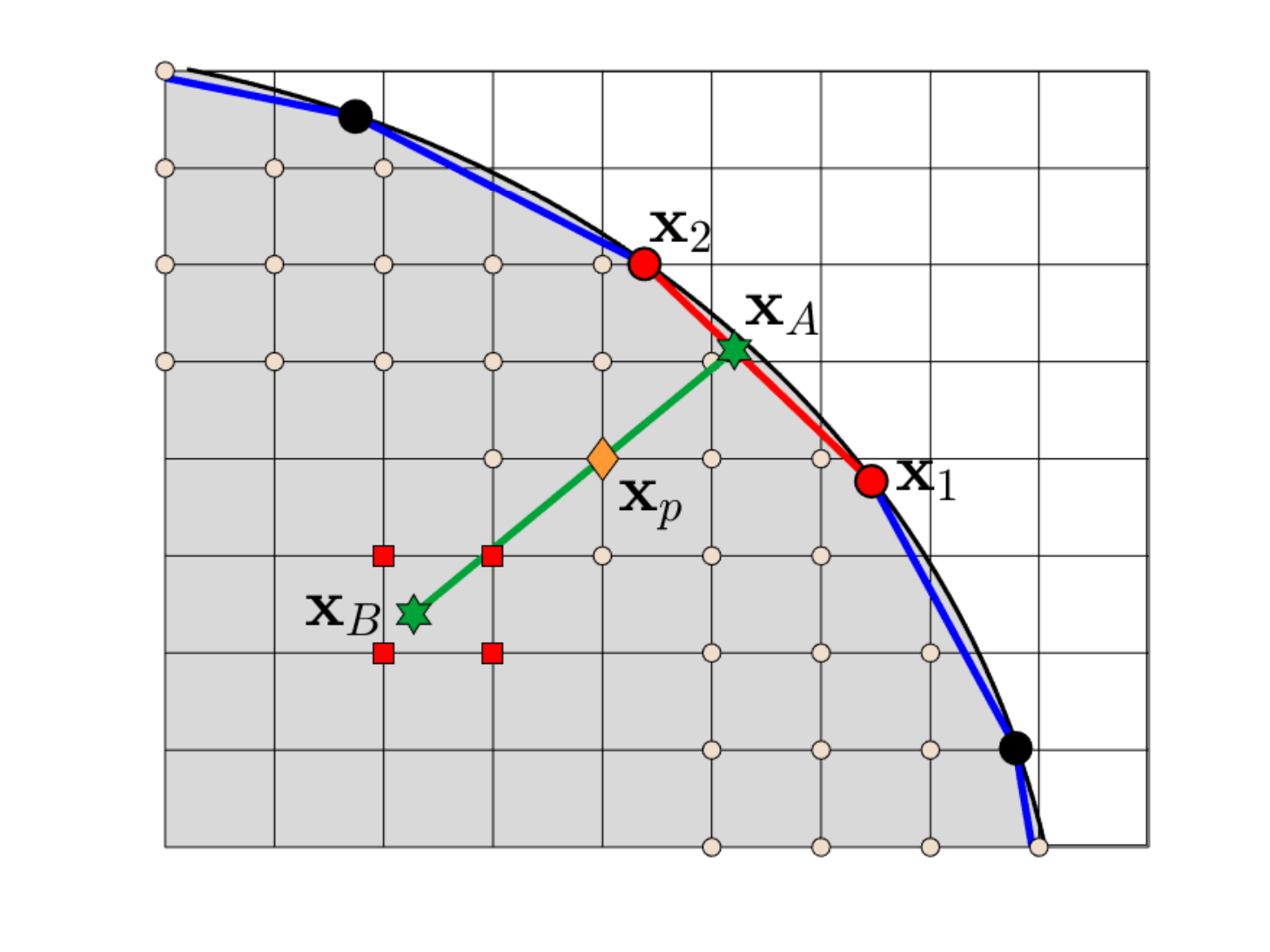}
\caption[Illustration of interpolation step of the IBDL method]{Illustration of the interpolation step of the IBDL method. The PDE domain $\Omega$ is shown in grey. The exact boundary curve is shown in black with black circles indicating immersed boundary points. Lines connecting the IB points are shown in blue. Grid points within $m_1=3$ meshwidths of the boundary are indicated by beige circles. The diagram illustrates the process for interpolating the value for one point, denoted with an orange diamond and labeled $\x_p$. The two closest boundary points are marked with red circles and labeled $\x_1$ and $\x_2$. The closest point on the line segment connecting $\x_1$ and $\x_2$ is called $\x_A$ and is marked with a green pentagram. Also marked with a green pentagram is the point $m_2$ meshwidths away from $\x_A$, and it is labeled $\x_B$. Marked with red squares are the four closest grid points to $\x_B$.   }\label{interp pic}
\end{figure}

Figure \ref{interp pic} illustrates the interpolation for one point, $\x_p$. Given such a point, we first locate the two nearest boundary points, labeled in the picture as $\x_1$ and $\x_2$. We use a projection to find $\x_A$, the point on the segment between $\x_1$ and $\x_2$ that is closest to $\x_p$. It is given by 
\begin{equation}
\x_A=\x_2+\frac{(\x_2-\x_p)\dotp (\x_1-\x_2)}{||\x_1-\x_2||^2}(\x_1-\x_2). 
\end{equation}
If the curvature and point placements are such that this calculation results in a point not located on the line segment between $\x_1$ and $\x_2$, we simply find the third closest boundary point and use the two outermost points to find $\x_A$. 
Next, we use a simple interpolation of the values at $\x_1$ and $\x_2$ to approximate the value at $\x_A$, and this is given by
\begin{equation}
u(\x_A)\approx U_b(\x_1)\frac{||\x_A-\x_2||}{||\x_2-\x_1||} + U_b(\x_2)\frac{||\x_A-\x_1||}{||\x_2-\x_1||} ,
\end{equation}
where the values at $\x_1$ and $\x_2$ are denoted with $U_b$ because these are known boundary values.

We next find the second point of interpolation, located $m_2$ meshwidths away from point $\x_A$. It is labeled in the figure as $\x_B$, and it is given by 
\begin{equation}
\x_B=\x_A+\frac{\x_p-\x_A}{||\x_p-\x_A||}(m_2\Delta x). \label{find B}
\end{equation}
Note that this step assumes that points are labeled as interior or exterior according to the \textit{discretized} shape created by the IB points. One way of doing this is presented in Section \ref{5.3 flagging}. Otherwise, if we are completing this interpolation for a point $\x_p$ that is in the interior of the analytical boundary shape but exterior to the discretized shape, we would need to alter Equation \eqref{find B} to use subtraction instead of addition. 

We then estimate the value for $\x_B$ with a bilinear interpolation of the computed solution values at the four grid points located around $\x_B$, marked in the figure with red squares. Let $\x_{00}$ be the lower left point, $\x_{01}$ be the upper left, and so on, and let $d_x$ and $d_y$ be the $x$ and $y$ distances from $\x_B$ to $\x_{00}$. Then our bilinear interpolation is 
\begin{multline}
u(\x_B)\approx \Big(\frac{d_x}{\Delta x}\Big) \Big(\frac{d_y}{\Delta y}\Big) u(\x_{11}) + \Big(\frac{d_y}{\Delta y}\Big)\Big(1-\frac{d_x}{\Delta x}\Big)u(\x_{01})\\+\Big(\frac{d_x}{\Delta x}\Big)\Big(1-\frac{d_y}{\Delta y}\Big)u(\x_{10}) + \Big(1-\frac{d_x}{\Delta x}-\frac{d_y}{\Delta y}+ \Big(\frac{d_x}{\Delta x}\Big) \Big(\frac{d_y}{\Delta y}\Big) \Big)u(\x_{00}).
\end{multline}
Note that if $\x_B$ is located directly on a grid line, this can be simplified. 

Once we have the approximated values at $\x_B$ and $\x_A$, we again use a simple linear interpolation to approximate $u(\x_p)$, and this is given by
\begin{equation}
u(\x_p)\approx u(\x_A)\frac{||\x_B-\x_p||}{||\x_B-\x_A||} + u(\x_B)\frac{||\x_A-\x_p||}{||\x_B-\x_A||} .
\end{equation}

\subsection{Flagging interior and exterior points using the IBDL method }\label{5.3 flagging}

In order to complete the interpolation step and to isolate the solution on $\Omega$, it is necessary to identify grid points as interior or exterior to $\Gamma$. Depending on the boundary shape, this can be done with the exact parametrization of the curve, but we can also use the IBDL method itself to identify grid points as interior to the \textit{discretized} curve. In this section, we show that by using analytical features of the indicator function for the interior of the immersed boundary $\Gamma$, we arrive at a familiar IBDL expression that enables us to use the framework already set forth in order to flag grid points as either exterior or interior to $\Gamma$. Eldredge also discusses this in Appendix A.2 of \cite{eldredge}, and his work uses such an indicator function as a foundation on which to build PDEs that govern different variables on different sides of the immersed boundary. Here, we present the derivation and application of this mechanism in the context of the IBDL method and also discuss how to apply it to a periodic computational domain. 

Let the indicator function be defined as
\begin{equation}
\chi(\x)=\begin{cases}
      1 & \x \text{ in } \Omega  \\
      0 & \x \text{ in } \C\setminus\Omega, \\
\end{cases}\label{indic}
\end{equation}
where we assume here that $\Omega$ is the interior of $\Gamma$. Let $\psi\in C_c^{\infty}$ be a smooth, compactly supported test function. Then we have
\begin{equation}
\int_{\C} \chi \Delta \psi d\x = \int_{\Omega} \Delta \psi d\x= \int_{\Gamma} \grad \psi \dotp \n dl(\x) , \label{B2}
 \end{equation}
where we first use the definition of $\chi$ and then the divergence theorem. Alternatively, we can use integration by parts and the compact support of $\chi$ to obtain 
\begin{equation}
 \int_{\C} \chi \Delta \psi d\x=-\int_{\C} \grad \chi \dotp \grad \psi d\x  , \label{chi 3}
\end{equation}
where we can view $\grad \chi$ as the weak derivative of $\chi$. By equating the far right-hand side of Equation \eqref{B2} and the right-hand side of Equation \eqref{chi 3}, we see that the weak derivative is the distribution given by
\begin{equation}
\grad \chi = - \delta (\x-\X) \n(\x), 
\end{equation}
where the $\X$ is a point on $\Gamma$, indicating that the the support of $\grad \chi$ is $\Gamma$. This then gives us 
\begin{equation}
- \Delta \chi = -\grad \dotp \grad \chi = \grad \dotp (\delta \n) ,
\end{equation}
and we therefore the PDE for $\chi$ given by
\begin{equation}
 \Delta \chi + \grad \dotp (\delta \n)= 0.
\end{equation}
We can recognize $\grad \dotp (\delta \n)$ as the unregularized version of the term found in the IBDL method, $\widetilde S Q = \grad \dotp( S Q\n)$, where $Q$ is a constant distribution of $1$. This then gives us a way to construct the indicator function $\chi$ by solving the PDE on $\C$ given by
\begin{equation}
\Delta \chi + \widetilde S \n=0,\label{ibdl with 1}
\end{equation} 
using the IBDL method. 

\begin{figure}
\centering
\includegraphics[width=0.6\textwidth]{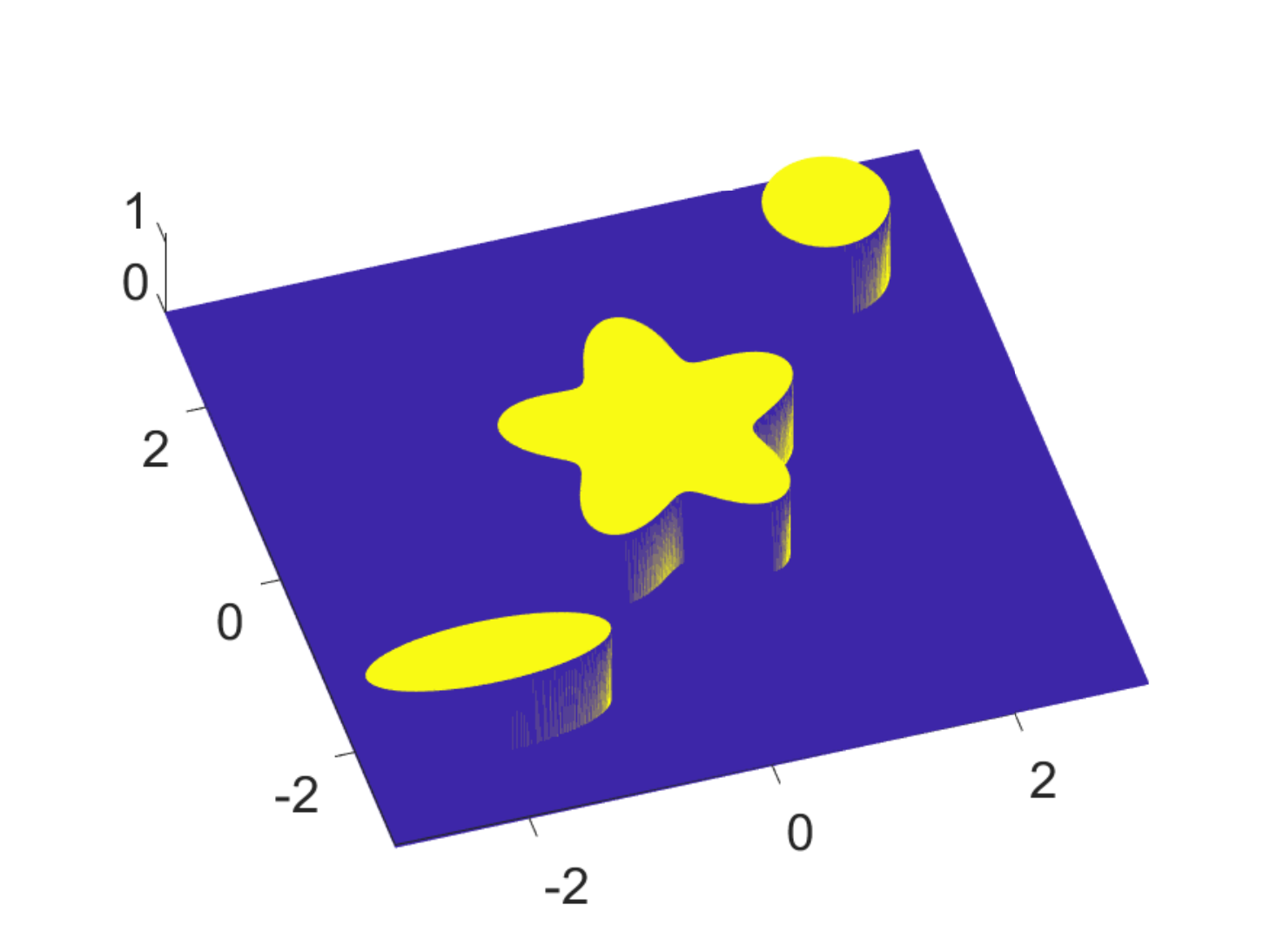}
\caption[Illustration of the IBDL method as a means to create an indicator function for the interior of $\Gamma$]{ Illustration showing the indicator function resulting from the process described in Section \ref{5.3 flagging} when applied to the discretized boundaries of three shapes.}\label{indicator pic}
\end{figure}

Then, if we are using a periodic computational domain, as we do in this dissertation, the solution is only unique up to an additive constant. Therefore, we can first find the solution with mean $0$ on $\C$. Then, to shift the exterior values to $0$, we simply need to subtract the appropriate constant, which can be approximated using a value of $\chi$ far from $\Gamma$. A likely choice is the corner of the periodic domain. Subtracting this value from $\chi$ then results in a regularized indicator function that is approximately equal to $1/2$ on the boundary $\Gamma$. To retrieve the final indicator function, we simply set all values of $\chi$ that are less than $1/2$ equal to $0$ and all those above $1/2$ to 1. Then $\chi$ is the indicator function described in Equation \eqref{indic} that is equal to 1 for all points on the interior of the discretized boundary, $\Gamma$. Figure \ref{indicator pic} gives an illustration of an indicator function resulting from applying this method to a set of three shapes. This process can be applied to multiple non-intersecting shapes at the same time. In other words, the PDE in Equation \eqref{ibdl with 1} needs only to be solved once, with the discrete spread operator applying to a distribution on the boundary points for all of the shapes. 

\section{Results: Dirichlet Helmholtz and Poisson equations}\label{ch 5 results}

In this section, we apply the Immersed Boundary Double Layer method and discuss numerical results. In Section \ref{5.4 hh}, we revisit the Dirichlet Helmholtz problem from Section \ref{2.4 Poisson} and demonstrate the improved efficiency of the IBDL method. In Section \ref{5.4 poisson}, we apply the method to a nonhomogeneous Poisson equation on a complex exterior domain. Lastly, in Section \ref{5.4 poisson shit}, we explore the need for the completed IBDL method for certain exterior domains

\subsection{Helmholtz equation revisited}\label{5.4 hh}

First, we revisit the Helmholtz PDE from Section \ref{2.4 Poisson},
\begin{subequations} \label{pde again again}
\begin{alignat}{2}
& \Delta u -  u = 0 \qquad && \text{in } \Omega  \label{pde1 again again}\\
&u=\sin{2\theta} \qquad && \text{on } \Gamma,  \label{pde2 again again}
\end{alignat}
\end{subequations}
where $\Omega$ is the interior of a circle of radius 0.25, centered at the origin. The analytical solution is given by 
\begin{equation}
u=\frac{I_2(r)\sin{2\theta}}{I_2(0.25)}, 
\end{equation}
where $I_2$ is the first-kind modified Bessel function of order 2. Our computational domain here is the periodic box $[-0.5, 0.5]^2$. For the IBDL method, we replace solution values within $m_1=6$ meshwidths of the boundary and use interpolation points located $m_2=8$ meshwidths from the boundary. We use equally spaced boundary points with $\Delta s \approx \alpha \Delta x$ for various values of $\alpha$. Both methods use a finite difference method, and the solutions are computed for grid sizes ranging from $N=2^6$ to $2^{12}$. 

\begin{figure}
\begin{center}
 \begin{tabular}{||c | c | c | c | c | c | c | c | c||} 
 \hline
 \multicolumn{9}{||c||}{Iteration Counts - Circular Boundary} \\
 \hline
  &\multicolumn{2}{|c|}{$\Delta s \approx  2\Delta x $}& \multicolumn{2}{|c|}{$\Delta s \approx  1.5 \Delta x$} &\multicolumn{2}{|c|}{$ \Delta s \approx  1 \Delta x$}& \multicolumn{2}{|c||}{$\Delta s \approx 0.75 \Delta x$}  \\
 \hline
 $\Delta x$ &\textcolor{blue}{ IBSL} & \textcolor{cyan}{IBDL} &\textcolor{blue}{  IBSL}& \textcolor{cyan}{IBDL}&\textcolor{blue}{ IBSL} &  \textcolor{cyan}{IBDL}&\textcolor{blue}{  IBSL} &  \textcolor{cyan}{IBDL} \\
 \hline
$2^{-6}$ &  \textcolor{blue}{  17} & \textcolor{cyan}{4} &    \textcolor{blue}{45}&\textcolor{cyan}{5}   &  \textcolor{blue}{249}&\textcolor{cyan}{5}  &   \textcolor{blue}{491}&\textcolor{cyan}{4}   \\
$2^{-7} $&   \textcolor{blue}{36}& \textcolor{cyan}{5}   &   \textcolor{blue}{57}&\textcolor{cyan}{4}    &  \textcolor{blue}{691}& \textcolor{cyan}{5} &  \textcolor{blue}{752}&\textcolor{cyan}{4}   \\
$2^{-8 }$&   \textcolor{blue}{49}& \textcolor{cyan}{5}  &    \textcolor{blue}{71}&\textcolor{cyan}{4}   & \textcolor{blue}{1233}& \textcolor{cyan}{4}   &  \textcolor{blue}{3057 }&\textcolor{cyan}{4}      \\
$2^{-9}$& \textcolor{blue}{ 61}&\textcolor{cyan}{4}    &     \textcolor{blue}{110}&\textcolor{cyan}{4}     &  \textcolor{blue}{1922}& \textcolor{cyan}{4}   &   \textcolor{blue}{5829}&\textcolor{cyan}{4}        \\
$ 2^{-10}$& \textcolor{blue}{68}&\textcolor{cyan}{4}   &     \textcolor{blue}{144}&\textcolor{cyan}{4}     &  \textcolor{blue}{1936}& \textcolor{cyan}{4}   &   \textcolor{blue}{ 8084 }&\textcolor{cyan}{4}    \\
$ 2^{-11}$ & \textcolor{blue}{95}&\textcolor{cyan}{4}  &   \textcolor{blue}{234}&\textcolor{cyan}{4}     &  \textcolor{blue}{ 4364}& \textcolor{cyan}{4}  &   \textcolor{blue}{9335 }&\textcolor{cyan}{4}    \\
$ 2^{-12}$  & \textcolor{blue}{142}&\textcolor{cyan}{4} & \textcolor{blue}{303}&\textcolor{cyan}{4}       &  \textcolor{blue}{ 4535}& \textcolor{cyan}{4}  &   \textcolor{blue}{10589  }&\textcolor{cyan}{4}   \\
  \hline
  \end{tabular}
\captionof{table}[Number of iterations of \texttt{minres} and \texttt{gmres} needed to solve the Helmholtz PDE in Equation \eqref{pde again again} using the IBSL and IBDL methods, respectively]{Number of iterations of \texttt{minres} and \texttt{gmres}, with tolerance $10^{-8}$, needed to solve Equation \eqref{pde again again} using the IBSL and IBDL methods, respectively. Both methods use finite differences.}
\label{iteration table 2}
\end{center}
\end{figure}

\begin{figure}
\centering
\begin{subfigure}{0.495\textwidth}
\centering
\includegraphics[width=\textwidth]{Figures/IBSLinteriorrefinement-eps-converted-to.pdf}
\caption{\normalsize IBSL Refinement}
\label{ibsl interior refinement}
\end{subfigure}
\begin{subfigure}{0.495\textwidth}
\centering
\includegraphics[width=\textwidth]{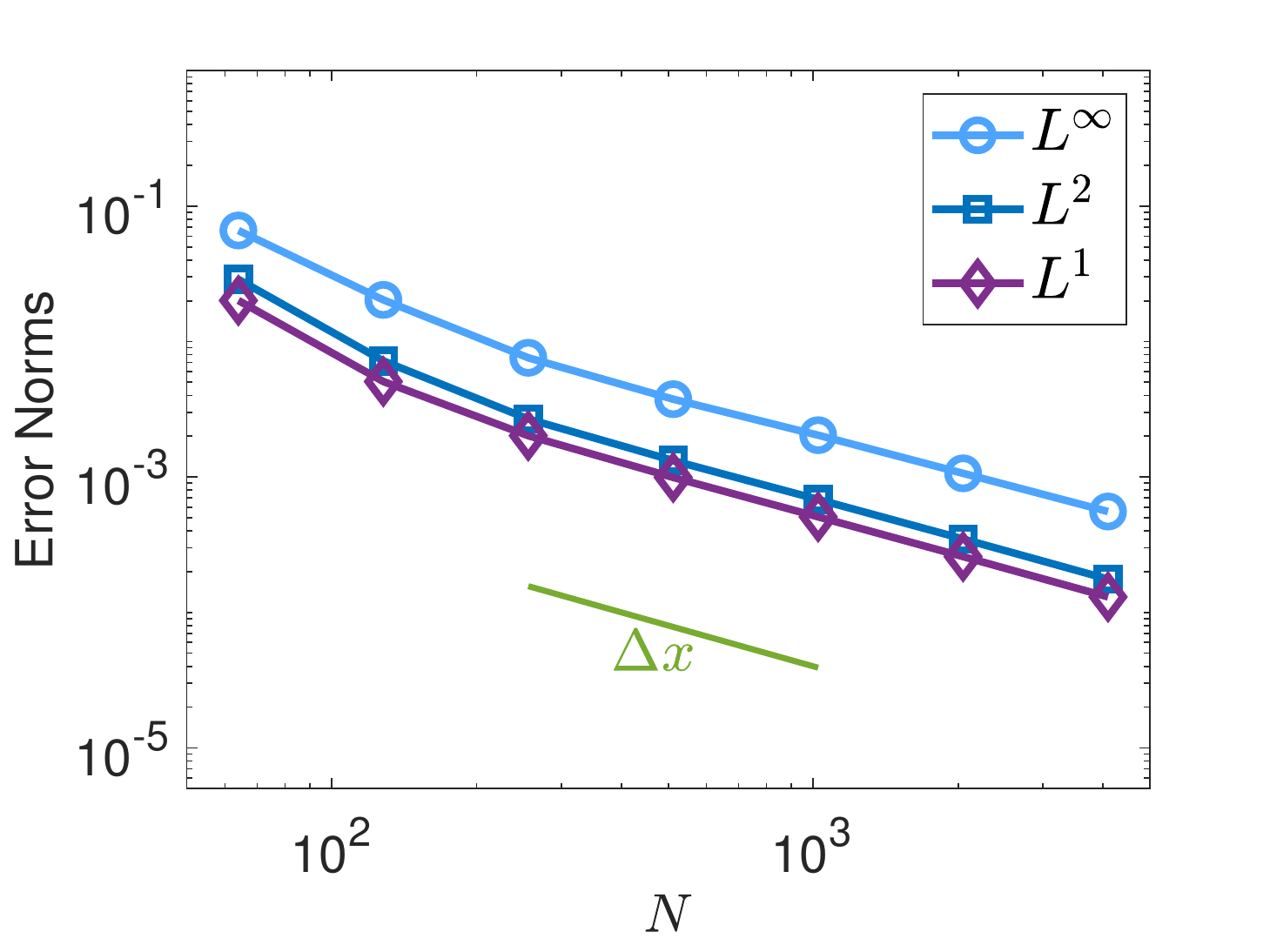}
\caption{\normalsize  IBDL Refinement}
\label{ibdl interior refinement}
\end{subfigure}
\begin{subfigure}{0.495\textwidth}
\centering
\includegraphics[width=\textwidth]{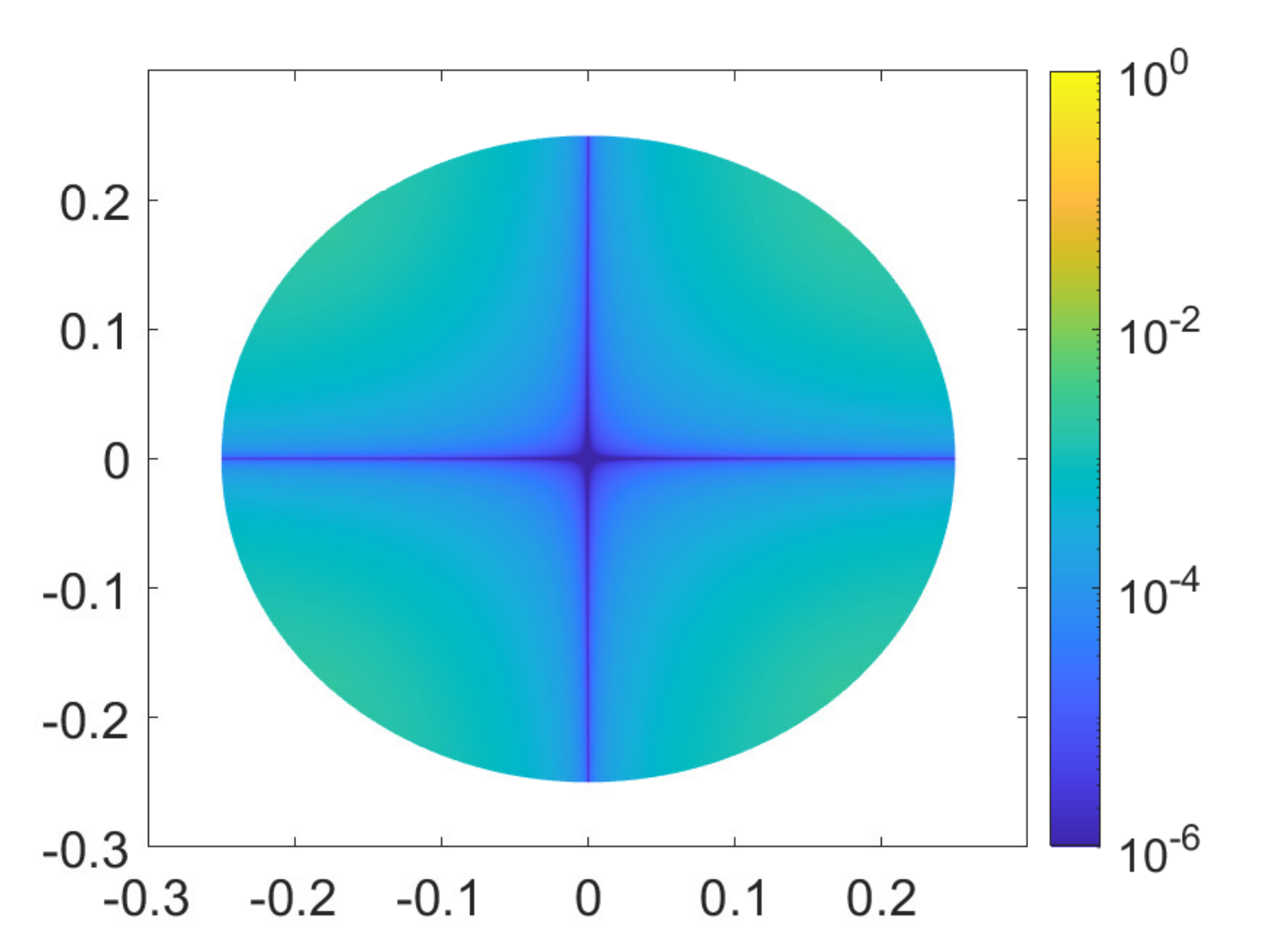}
\caption{\normalsize IBSL Error}
\label{ibsl interior error}
\end{subfigure}
\begin{subfigure}{0.495\textwidth}
\centering
\includegraphics[width=\textwidth]{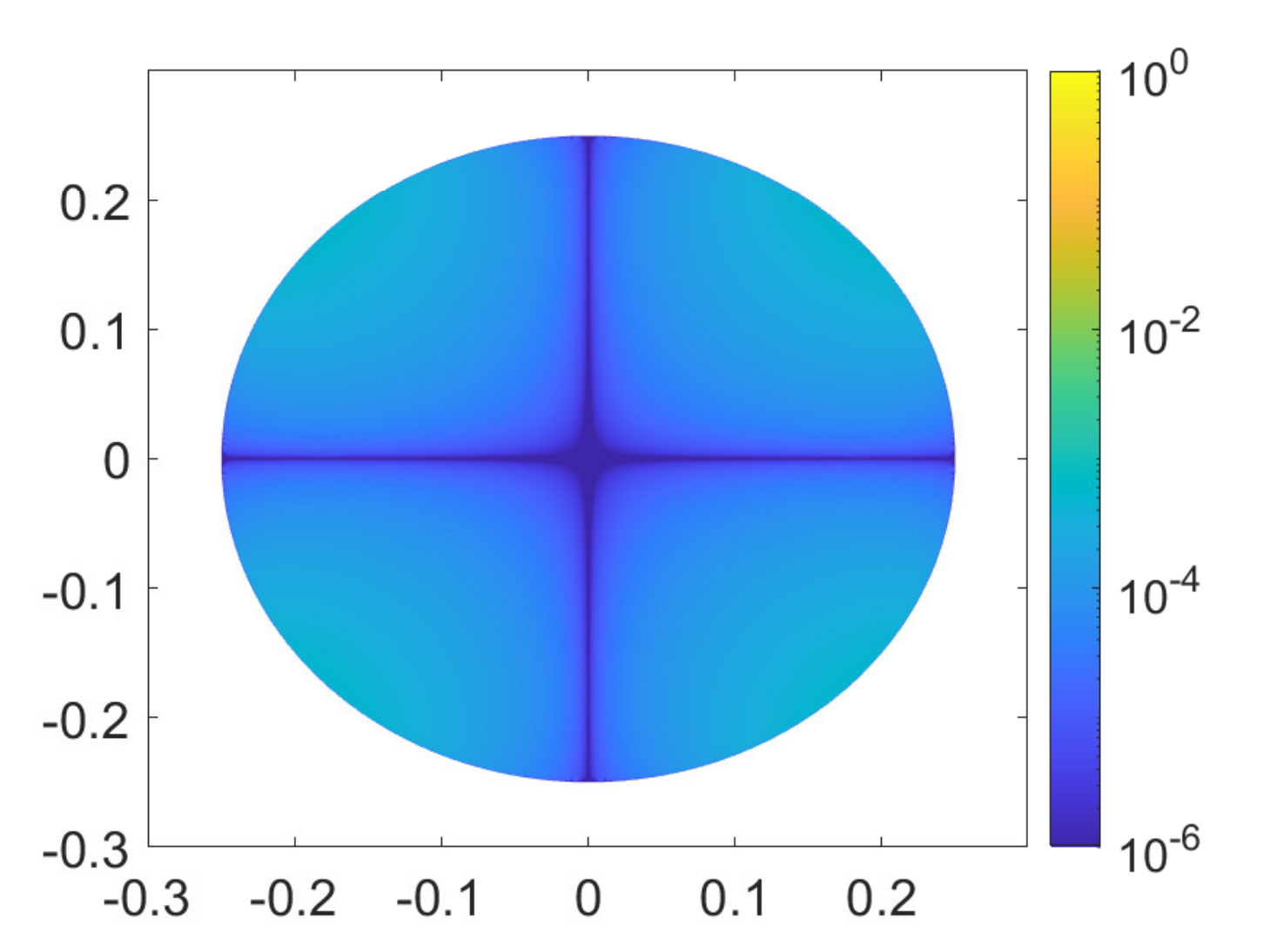}
\caption{\normalsize  IBDL Error}
\label{ibdl interior error}
\end{subfigure}
\caption[Error plots and refinement studies for solutions to the interior Helmholtz PDE in Equation \eqref{pde again again} found using the IBSL and IBDL methods]{Error plots and refinement studies for solutions to Equation \eqref{pde again again} found using the IBSL and IBDL methods, both with finite differences. The computational domain is the periodic box $[-0.5, 0.5]^2$, $\Omega$ is the interior of a circle of radius 0.25, and the prescribed boundary values are given by $U_b=\sin{2\theta}$. Figures \ref{ibsl interior error}-\ref{ibdl interior error} show the errors for the grid spacing $\Delta x = 2^{-12}$, and both methods use boundary point spacing of $ \Delta s  \approx 0.75 \Delta x$. The IBDL method replaces solution values within $m_1= 6$ meshwidths from the boundary using an interior interpolation point $m_2=8$ meshwidths away from the boundary. All plots use absolute errors.}
\end{figure}

Table \ref{iteration table 2} gives the iteration counts for the Krylov methods used to solve Equation \eqref{pde again again} for the IBSL and IBDL methods. This table illustrates the drastic improvement in efficiency that is a key advantage of the IBDL method. We see that only a small number of iterations are needed to solve for the strength of the IBDL  potential. Moreover, the number of iterations remains essentially constant as we refine the mesh or tighten the boundary points relative to the grid. The IBSL method, on the other hand, requires many more iterations, and this iteration count is greatly affected by the mesh size and the boundary point spacing. Tighter boundary points can result in lower quadrature error, while more widely spaced points allow for smoother force distributions and better conditioning in the case of the IBSL method \cite{GriffithDonev}. Making the IBSL method practical therefore often entails special handling, such as using more widely spaced points for discretization of the boundary, while using a denser set of points for the direct Eulerian-Lagrangian interaction \cite{Griffithpointspacing}. Moreover, to make the IBSL method efficient, one must use a preconditioner  \cite{Ceniceros, GriffithDonev, guyphilipgriffith, Stein}. On the other hand, the fact that the IBDL method presented here can be efficiently utilized for both coarsely or finely spaced boundary points means no such step is required for this method. 

In addition, the iteration counts for the IBSL method are affected by the tolerance for the Krylov method. For example, for $\Delta s \approx 2\Delta x$ and $\Delta x = 2^{-12}$, decreasing the tolerance from $10^{-8}$ to $10^{-10}$ increases the iteration count by more than fifty percent, from 142 to 233. For $\Delta s \approx 0.75 \Delta x$, the Krylov method cannot reach a tolerance of $10^{-10}$ for the IBSL method. For the IBDL method, on the other hand, the Krylov method reaches this smaller tolerance in only 5 iterations for either of these cases. 

Figures \ref{ibsl interior refinement}-\ref{ibdl interior refinement} compare the IBSL and IBDL refinement studies, and we see that we maintain first-order convergence of the solution. Figures \ref{ibsl interior error}-\ref{ibdl interior error} illustrate the solution errors for the two methods, and the errors from the IBDL method are actually lower for this problem. In general, the two methods typically give errors of comparable size. 
\subsection{Poisson equation on complex exterior domain}\label{5.4 poisson}
\begin{figure}
\begin{center}
 \begin{tabular}{||c | c | c | c | c ||} 
 \hline
 \multicolumn{5}{||c||}{Iteration Counts - Starfish boundary} \\
 \hline
 &\multicolumn{2}{|c|}{$\Delta s \approx  2\Delta x $}& \multicolumn{2}{|c||}{$\Delta s \approx  0.75 \Delta x$} \\ 
 \hline
  $\Delta x$ &\textcolor{blue}{ IBSL} & \textcolor{cyan}{IBDL} &\textcolor{blue}{  IBSL}& \textcolor{cyan}{IBDL}\\
 \hline 
$2^{-4}$ &  \textcolor{blue}{ 57 } & \textcolor{cyan}{13} &    \textcolor{blue}{1999}&\textcolor{cyan}{13}    \\
$2^{-5} $&   \textcolor{blue}{84}& \textcolor{cyan}{13}   &   \textcolor{blue}{5197}&\textcolor{cyan}{13}   \\
$2^{-6 }$&   \textcolor{blue}{87}& \textcolor{cyan}{13}  &    \textcolor{blue}{4289}&\textcolor{cyan}{13}     \\
$2^{-7}$& \textcolor{blue}{ 86}&\textcolor{cyan}{13}    &     \textcolor{blue}{3207}&\textcolor{cyan}{13}    \\
$ 2^{-8}$& \textcolor{blue}{74}&\textcolor{cyan}{14}   &     \textcolor{blue}{2869}&\textcolor{cyan}{14}   \\
$ 2^{-9}$ & \textcolor{blue}{65}&\textcolor{cyan}{14}  &   \textcolor{blue}{2403}&\textcolor{cyan}{14}   \\
$ 2^{-10}$  & \textcolor{blue}{55}&\textcolor{cyan}{14} & \textcolor{blue}{1554}&\textcolor{cyan}{14}     \\
\hline
  \end{tabular}
\captionof{table}[Number of iterations of \texttt{minres} and \texttt{gmres} to solve the Poisson PDE in Equation \eqref{pde again 3} using the IBSL and IBDL methods, respectively, on a PDE domain exterior to a starfish shape]{Number of iterations of \texttt{minres} and \texttt{gmres}, with tolerance $10^{-8}$, for the IBSL and IBDL methods, respectively, to solve Equation \eqref{pde again 3} on the periodic computational domain $[-2, 2]^2$. $\Omega$ is the region exterior to a ``starfish'' shape, and the prescribed boundary values are given by $U_b=\sin{(\pi x/2)} - \cos{(\pi x/2)}$. Both methods utilize finite differences.} \label{star iteration table 3}
\end{center}
\end{figure}

\begin{figure}
\centering
\begin{subfigure}{0.495\textwidth}
\centering
\includegraphics[width=\textwidth]{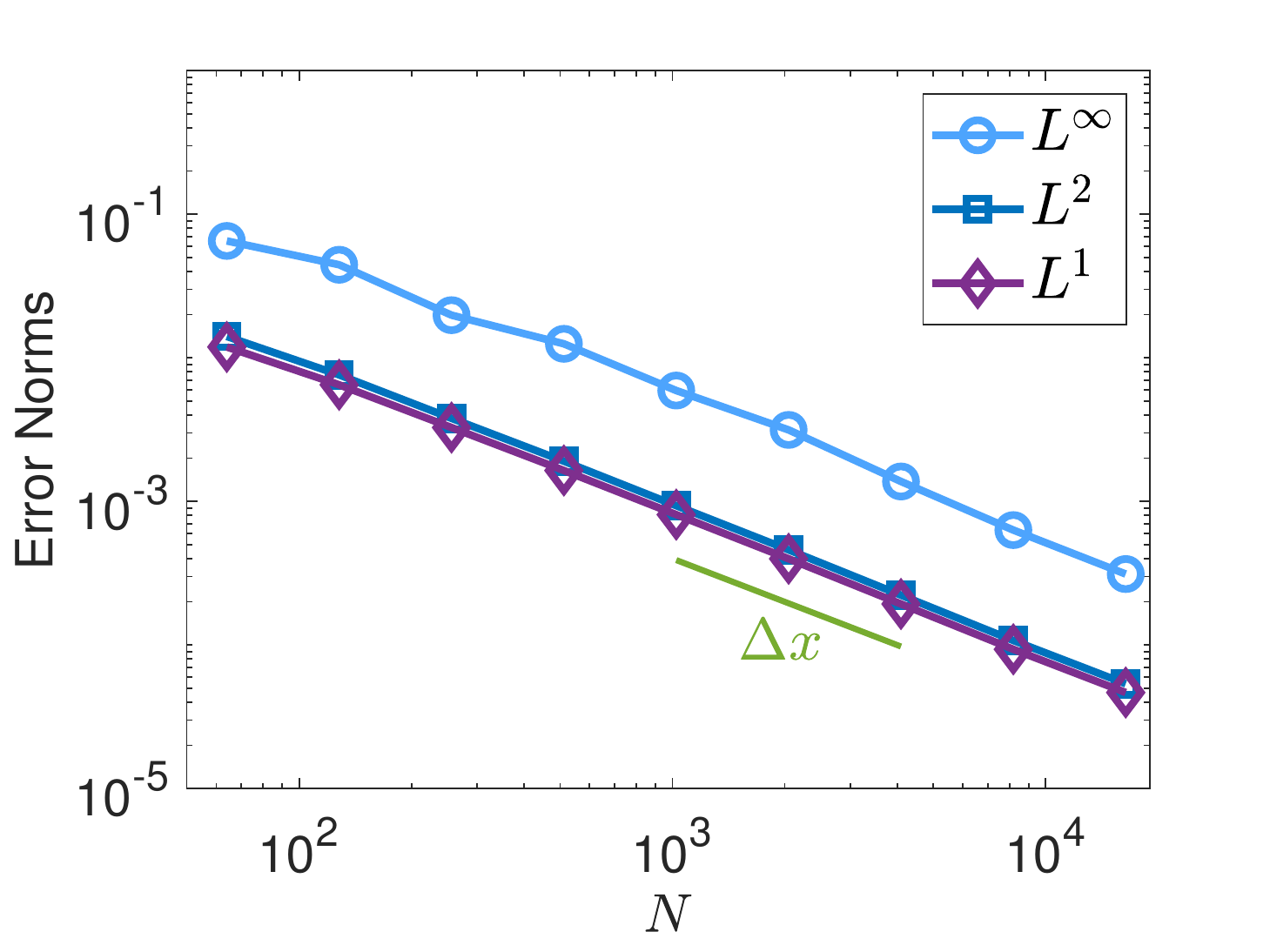}
\caption{\normalsize  IBSL Refinement}
\label{ibsl refinement ext}
\end{subfigure}
\begin{subfigure}{0.495\textwidth}
\centering
\includegraphics[width=\textwidth]{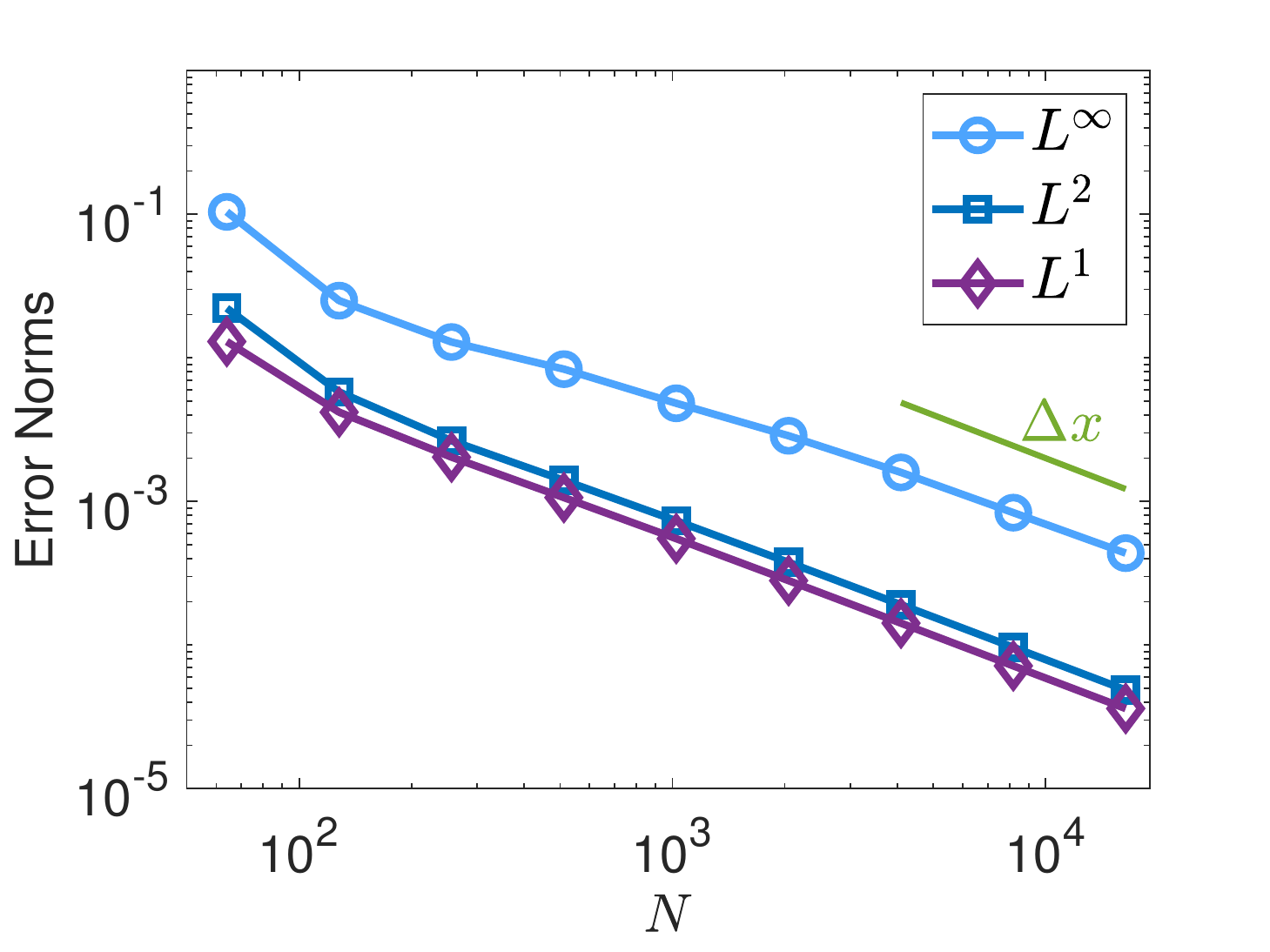}
\caption{\normalsize  IBDL Refinement}
\label{ibdl refinement ext}
\end{subfigure}
\begin{subfigure}{0.495\textwidth}
\centering
\includegraphics[width=\textwidth]{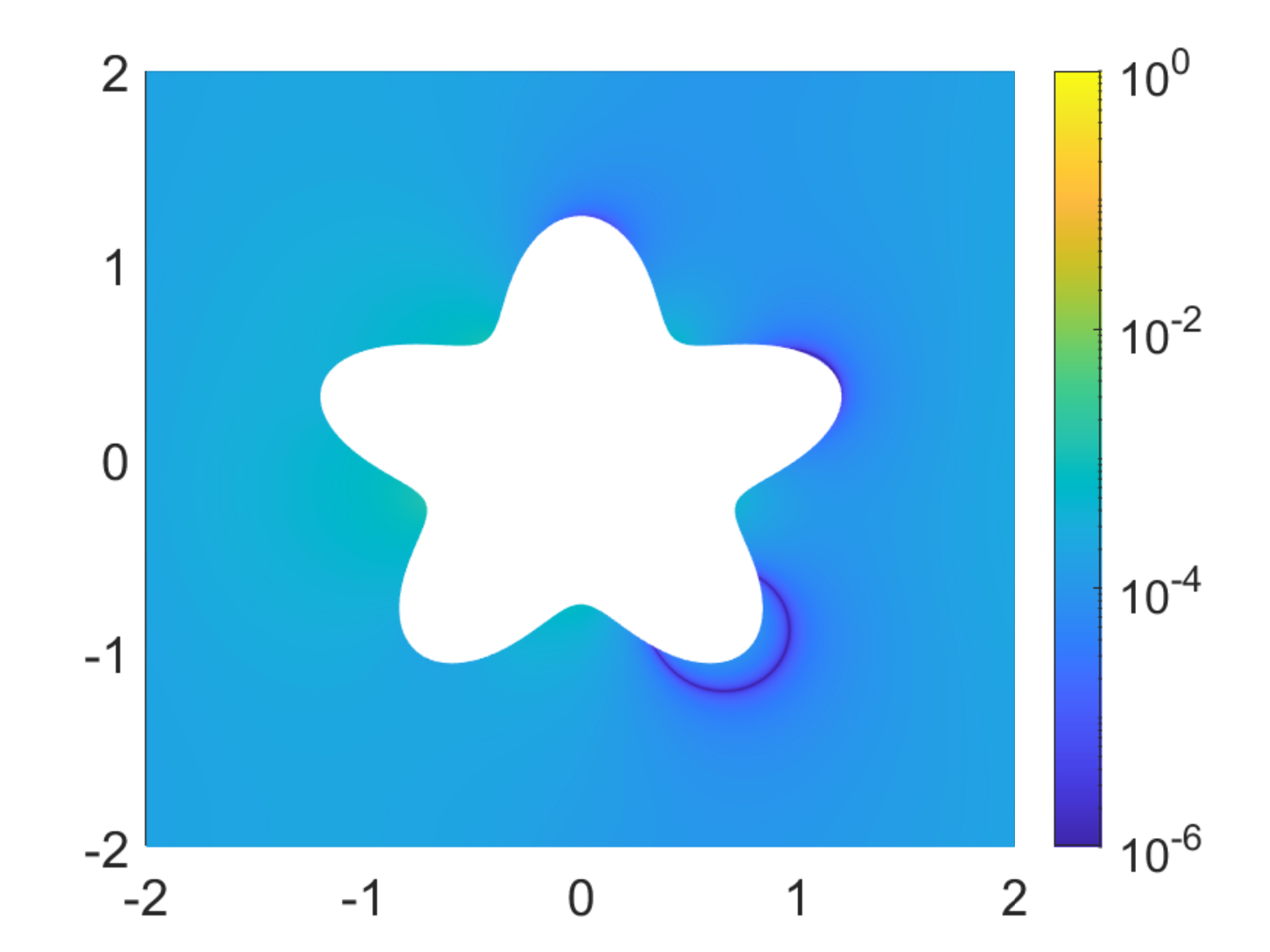}
\caption{\normalsize  IBSL Error}
\label{ibsl ext error}
\end{subfigure}
\begin{subfigure}{0.495\textwidth}
\centering
\includegraphics[width=\textwidth]{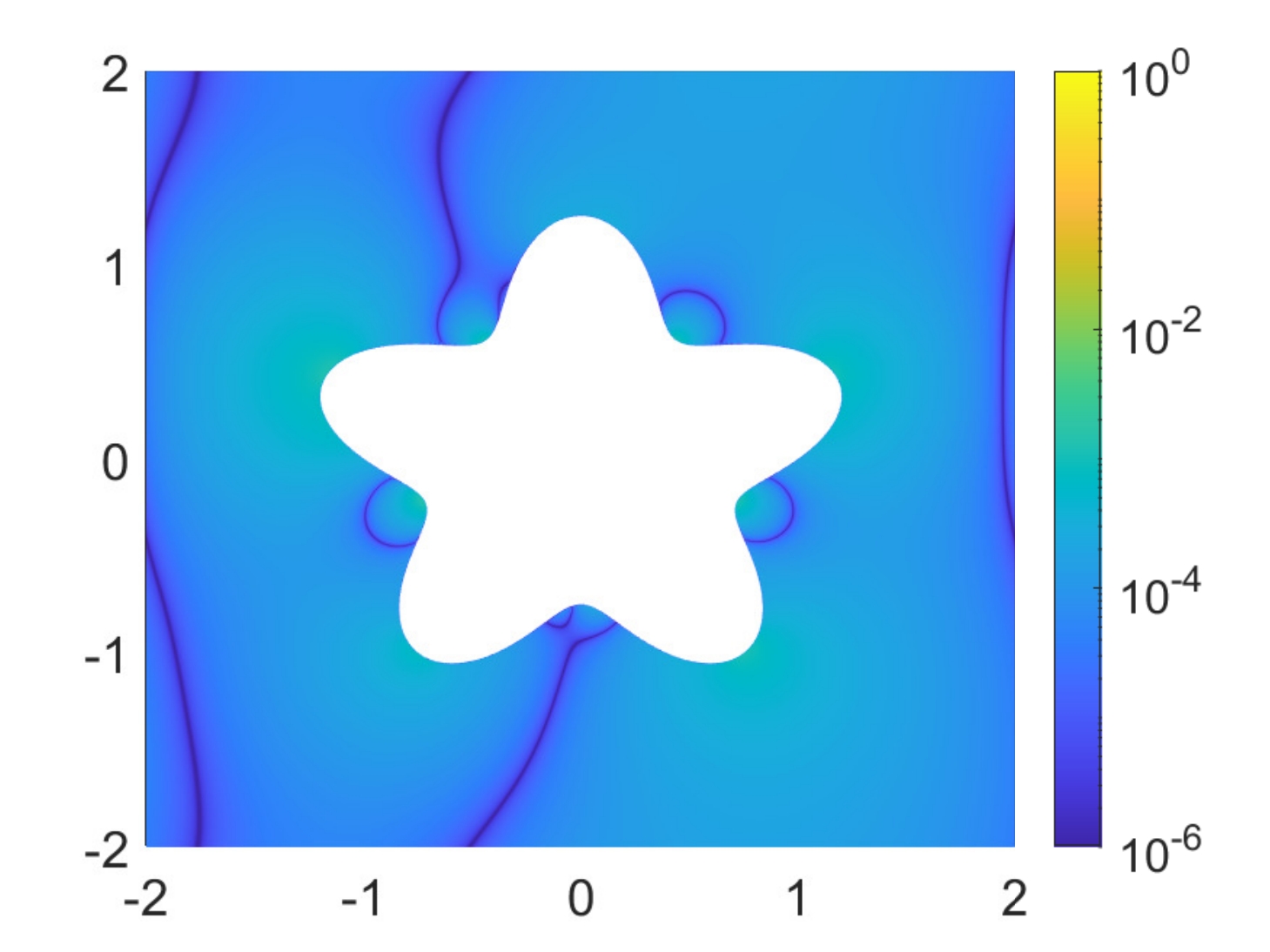}
\caption{\normalsize  IBDL Error}
\label{ibdl ext error}
\end{subfigure}
\caption[Refinement studies and error plots for solutions to the Poisson PDE in Equation \eqref{pde again 3} found using the IBSL and IBDL methods]{Refinement studies and error plots for solutions to Equation \eqref{pde again 3} found using the IBSL and IBDL methods, with finite differences. The computational domain is the periodic box $[-2, 2]^2$, $\Omega$ is the region exterior to a ``starfish'' shape, and the prescribed boundary values are given by $U_b=\sin{(\pi x/2)} - \cos{(\pi x/2)}$. A boundary point spacing of $ \Delta s \approx 0.75  \Delta x$ is used, and in the IBDL method, solution values within $m_1= 6$ meshwidths from the boundary are replaced using a second interpolation point $m_2=8$ meshwidths away from the boundary. Figures \ref{ibsl ext error}-\ref{ibdl ext error} show plots of the solution errors for a mesh size of $N=2^{12}$.}\label{poisson plots}
\end{figure}

We move now to the Poisson equation for which we use the methods described in Sections \ref{2.2 Laplacian L} and \ref{5.3 poisson} for the IBSL and IBDL methods, respectively, to handle the nullspace of the periodic Laplacian. Here we do not employ the completed double layer formulation, but use the original IBDL method as described in Equation \eqref{ibdl}. We also use this section to demonstrate the effectiveness of the IBDL method for a complex domain exterior to a non-convex boundary. We apply the IBSL and IBDL methods to the PDE given by 
\begin{subequations} \label{pde again 3}
\begin{alignat}{2}
& \Delta u  =\frac{\pi^2}{4}\Bigg( \cos{\Big(\frac{\pi y}{2}\Big)}-\sin{\Big(\frac{\pi x}{2}\Big)}\Bigg)\qquad && \text{in } \Omega  \label{pde1 again 3}\\
&u=\sin{\Big(\frac{\pi x}{2}\Big)} - \cos{\Big(\frac{\pi y}{2}\Big)} && \text{on } \Gamma,  \label{pde2 again 3}
\end{alignat}
\end{subequations}
where $\Omega$ is the region inside the periodic box $[-2,2]^2$ that is exterior to the ``starfish'' shape used in \cite{QBX}, which is given by
\begin{equation} \label{starfish eqn}
\begin{pmatrix}
x(\theta)\\
y(\theta)
\end{pmatrix}
= 
\Bigg( 1+\frac{\sin{(10\pi \theta)}}{4}\Bigg)\begin{pmatrix}
\cos{(2\pi \theta)}\\
\sin{(2\pi \theta)}
\end{pmatrix}, 
\end{equation}
for $0 \leq \theta \leq 1$. We use equally spaced boundary points, with $\Delta s \approx \alpha \Delta x$ for $\alpha=0.75$ and $2$. The solutions are computed for grid sizes ranging from $N=2^6$ to $2^{14}$, and the analytic solution is given by $ u=\sin{(\pi x/2)} - \cos{(\pi y/2)}$. For the IBDL method, we replace solution values within $m_1=6$ meshwidths of the boundary and use interpolation points located $m_2=8$ meshwidths from the boundary. Both methods utilize finite differences. Table \ref{star iteration table 3} gives the iteration counts for the Krylov methods used by the IBSL and IBDL methods. Although the IBSL iteration counts do not always increase with grid refinement, the iteration counts are still higher than the IBDL iteration counts, particularly for $\Delta s \approx 0.75 \Delta x$. Figure \ref{poisson plots} shows the refinement studies and plots of the solution errors for the IBSL and IBDL methods when $\alpha=0.75$. Again, we see that the convergence is first-order for both methods, with comparable error sizes. 

\subsection{Completed IBDL method}\label{5.4 poisson shit}

\begin{figure}
\centering
\begin{subfigure}{0.43\textwidth}
\centering
\includegraphics[width=\textwidth]{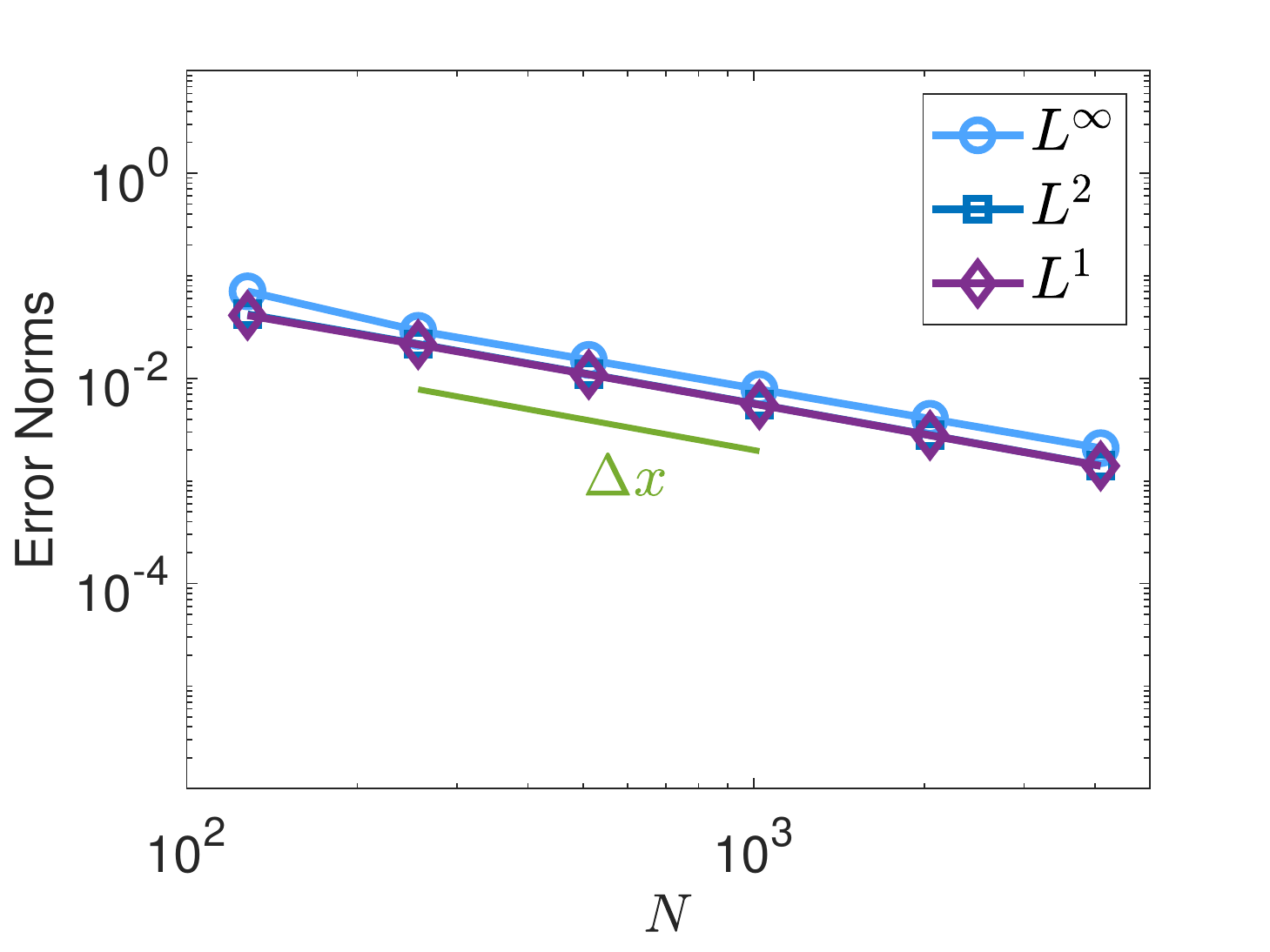}
\caption{\normalsize  Equation \eqref{poisson shit pde exp}, $L=1$}
\label{smallLexp}
\end{subfigure}
\begin{subfigure}{0.43\textwidth}
\centering
\includegraphics[width=\textwidth]{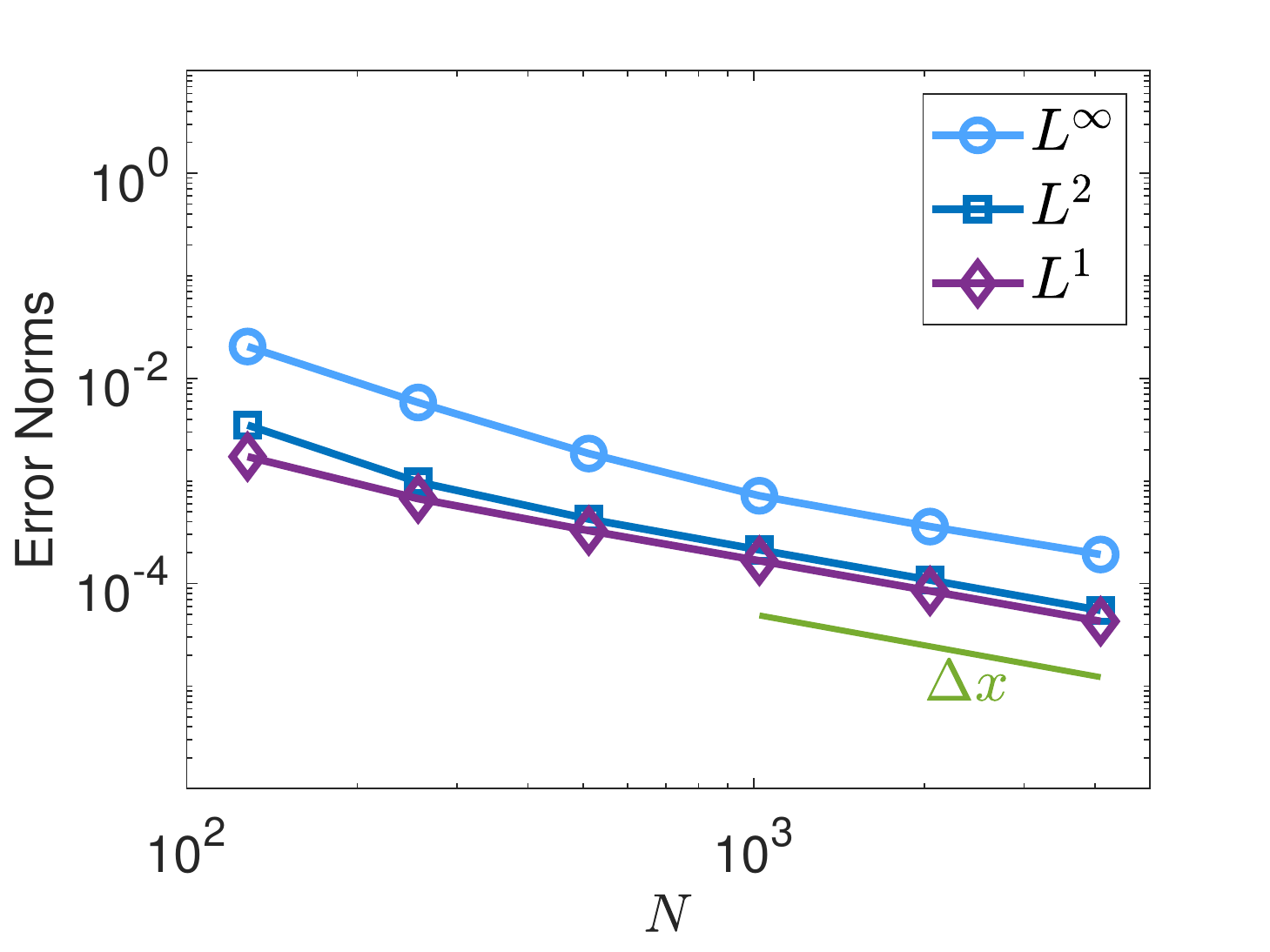}
\caption{\normalsize  Equation \eqref{poisson shit pde fine}, $L=1$}
\label{smallLfine}
\end{subfigure}
\begin{subfigure}{0.43\textwidth}
\centering
\includegraphics[width=\textwidth]{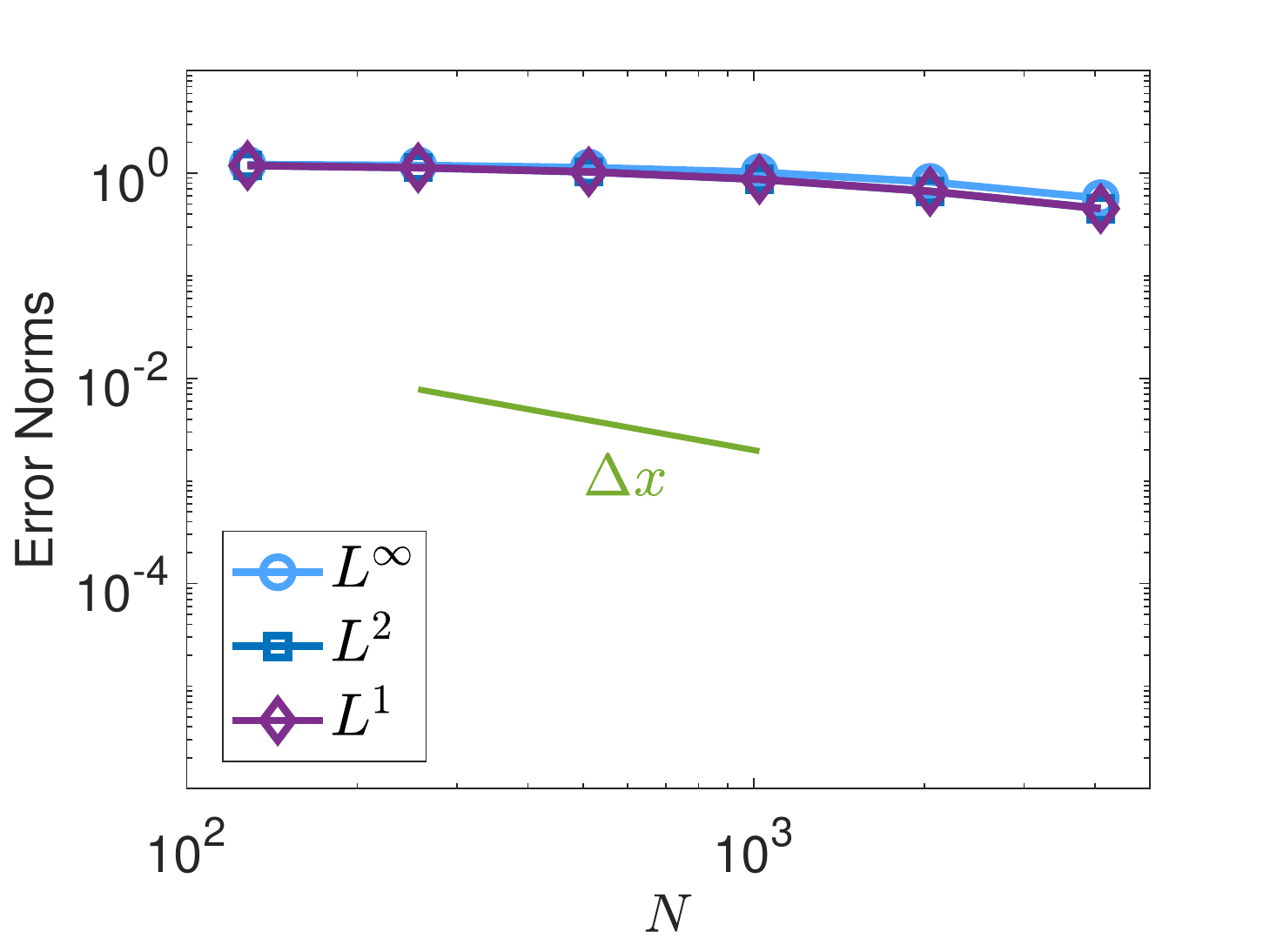}
\caption{\normalsize  Equation \eqref{poisson shit pde exp}, $L=8$}
\label{bigLexp}
\end{subfigure}
\begin{subfigure}{0.43\textwidth}
\centering
\includegraphics[width=\textwidth]{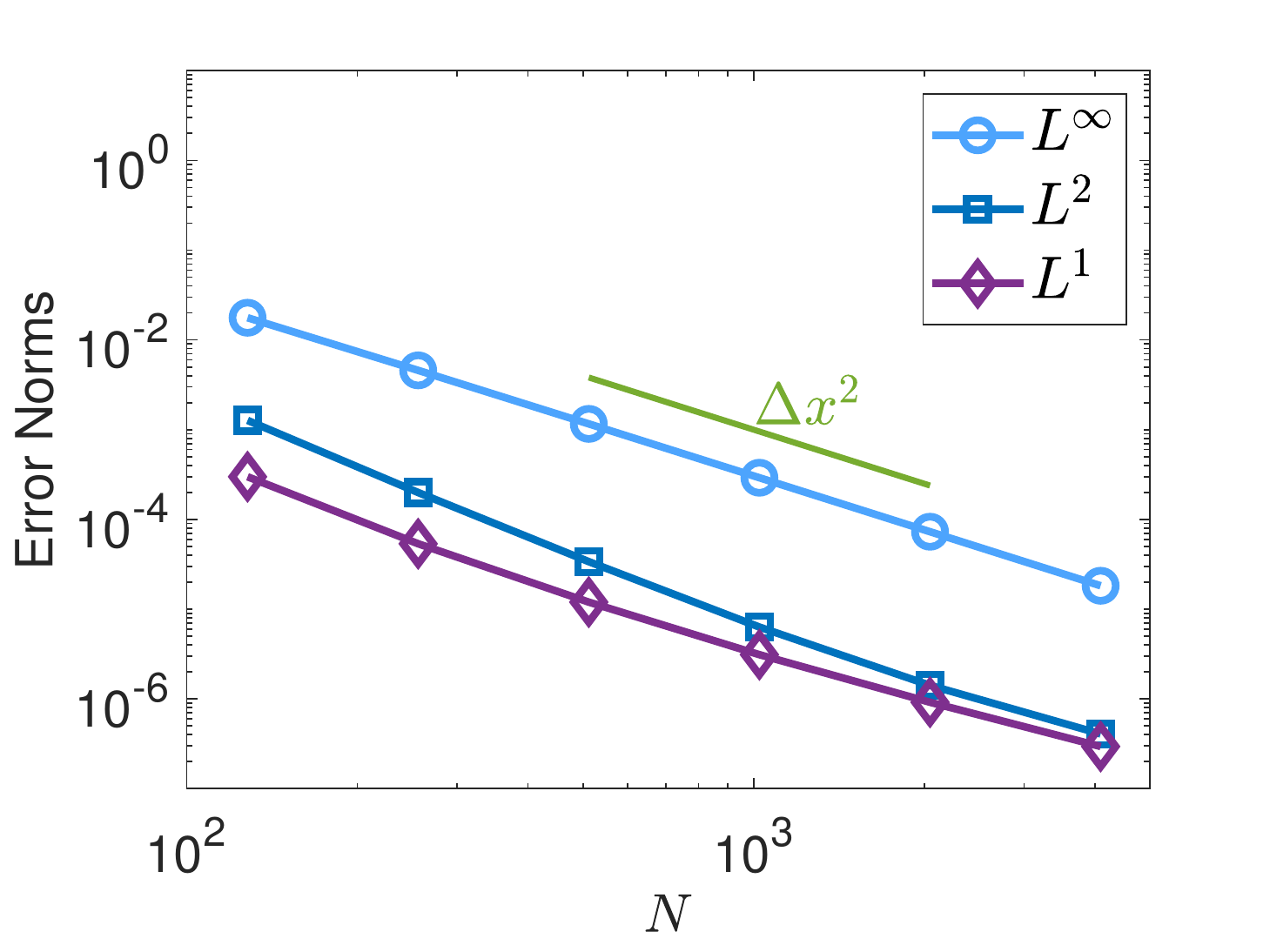}
\caption{\normalsize  Equation \eqref{poisson shit pde fine}, $L=8$}
\label{bigLfine}
\end{subfigure}
\begin{subfigure}{0.43\textwidth}
\includegraphics[width=\textwidth]{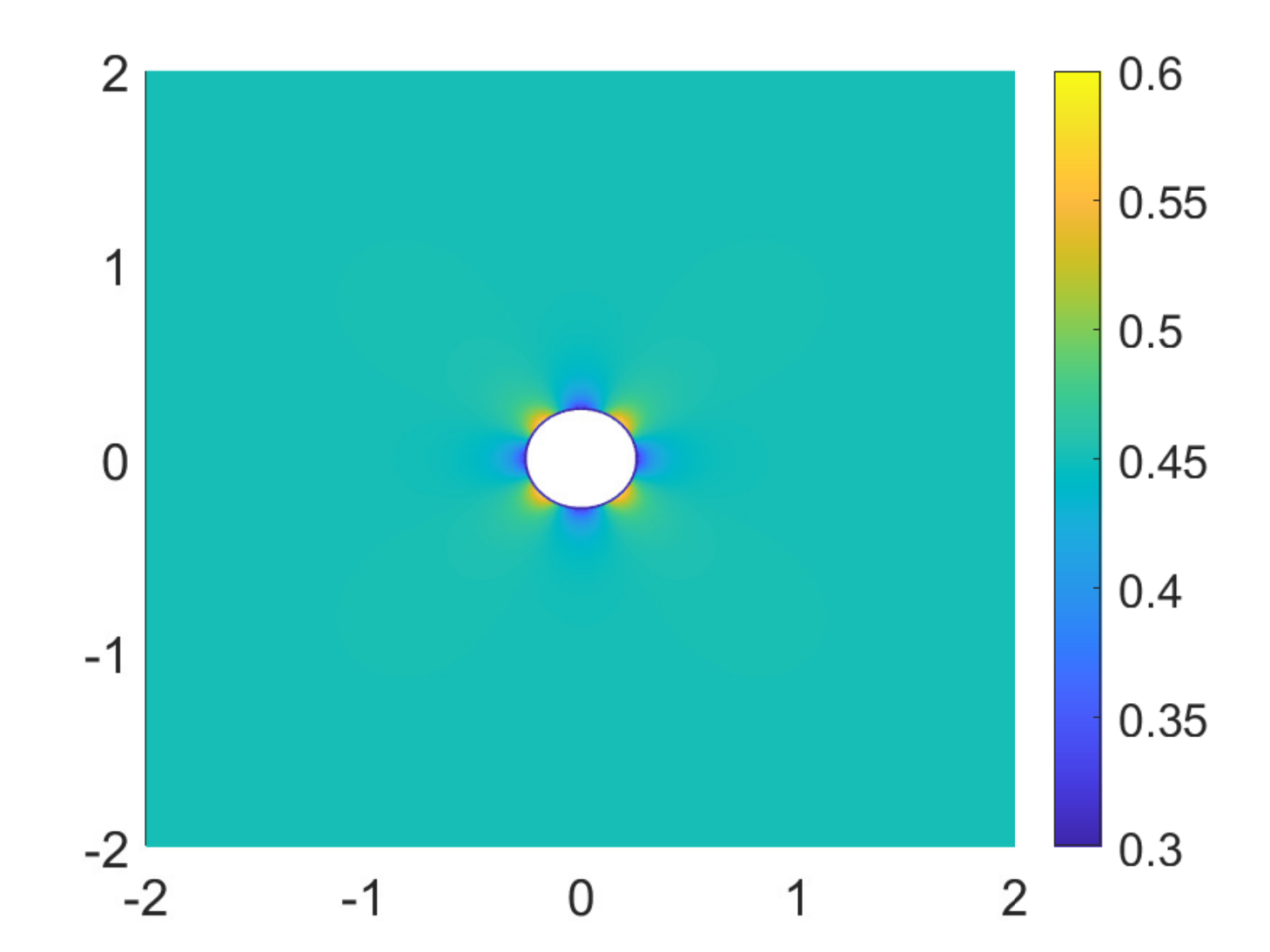}
\caption{\normalsize  Solution error, Equation \eqref{poisson shit pde exp}, $L=8$}
\label{bigLexperror}
\end{subfigure}
\caption[Refinement studies and error plot for solutions to the Poisson PDEs in Equations \eqref{poisson shit pde exp}-\eqref{poisson shit pde fine} found using the IBDL method on small and large exterior domains]{Refinement studies for solutions to Equations \eqref{poisson shit pde exp}-\eqref{poisson shit pde fine} found using the IBDL method, with finite differences. For Figures \ref{smallLexp}-\ref{smallLfine}, the computational domain is the periodic box $[-0.5, 0.5]^2$, and for Figures \ref{bigLexp}-\ref{bigLexperror}, the computational domain is the periodic box $[-4, 4]^2$. $\Omega$ is the region in the computational domain exterior to a circle of radius $0.25$, centered at the origin. A boundary point spacing of $ \Delta s \approx 0.75  \Delta x$ is used, and solution values within $m_1= 6$ meshwidths from the boundary are replaced using a second interpolation point $m_2=8$ meshwidths away from the boundary. All errors reported are absolute errors. Figure \ref{bigLexperror} gives a plot of the solution errors for Equation \eqref{poisson shit pde exp} for a mesh size of $N=2^{12}$.}\label{poisson shit plots}
\end{figure}

In this section, we report some observations regarding the need for the completed IBDL method. We examine two PDEs to demonstrate that this need depends on both the size of the exterior domain and on the given boundary values. For both PDEs, $\Omega$ is the region in the periodic box $[-L/2, L/2]^2$ that is exterior to a circle of radius $0.25$, centered at the origin. We will look at $L=1$ and $L=8$. The first PDE is
\begin{subequations} \label{poisson shit pde exp}
\begin{alignat}{2}
& \Delta u  =  \frac{4\pi^2}{L^2} e^{\sin{(2 \pi x/L )}}\Big(\cos^2{(2\pi x/L)}-\sin{(2\pi x/L)}\Big) \qquad && \text{in } \Omega \\
&u= e^{\sin{(2\pi x/L )}}&& \text{on } \Gamma, 
\end{alignat}
\end{subequations}
and the analytical solution is given by $u = e^{\sin{(2\pi x/L)}}$. The second PDE is 
\begin{subequations} \label{poisson shit pde fine}
\begin{alignat}{2}
& \Delta u  =  -\frac{4\pi^2}{L^2} \Big(\sin{(2\pi x/L)}-\cos{(2\pi y/L)}\Big) \qquad && \text{in } \Omega \\
&u=\sin{(2\pi x/L)}-\cos{(2\pi y/L)}&& \text{on } \Gamma, 
\end{alignat}
\end{subequations}
and the analytical solution is given by $u=\sin{(2\pi x/L)}-\cos{(2\pi y/L)}$. We use finite differences and interpolate values within $m_1=6$ meshwidths of the boundary. We first solve these problems using the IBDL method described in Equation \eqref{ibdl}, without the addition of the single layer potential. Figures \ref{smallLexp}-\ref{smallLfine} illustrate that on the small periodic domain of length $1$, both computed solutions have first-order convergence to the analytical solutions. However, Figures \ref{bigLexp}-\ref{bigLfine} demonstrate that for a larger periodic domain of length $8$, the two solutions behave very differently. The IBDL method is unable to achieve convergence to the analytical solution to Equation \eqref{poisson shit pde exp}. Figure \ref{bigLexperror} shows that the error is essentially constant in the region away from the boundary. We may therefore be seeing the addition of an element of the nullspace of the periodic Laplacian to our solution. On the other hand, we see second-order convergence to the analytical solution for Equation \eqref{poisson shit pde fine}. 

\begin{figure}
\centering
\begin{subfigure}{0.43\textwidth}
\centering
\includegraphics[width=\textwidth]{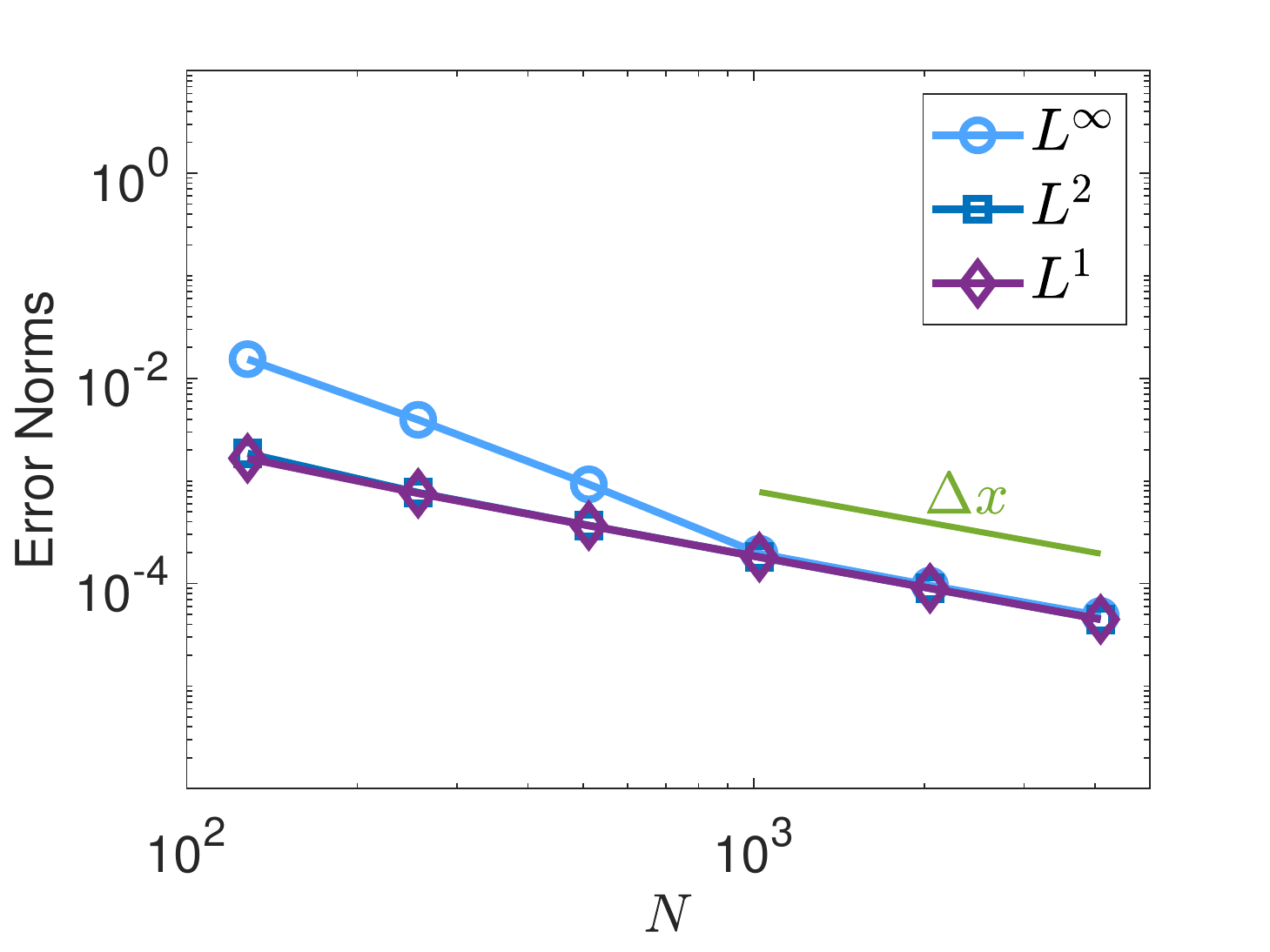}
\caption{\normalsize  Equation \eqref{poisson shit pde exp}, $L=8$}
\label{bigLexpcompleted}
\end{subfigure}
\begin{subfigure}{0.43\textwidth}
\centering
\includegraphics[width=\textwidth]{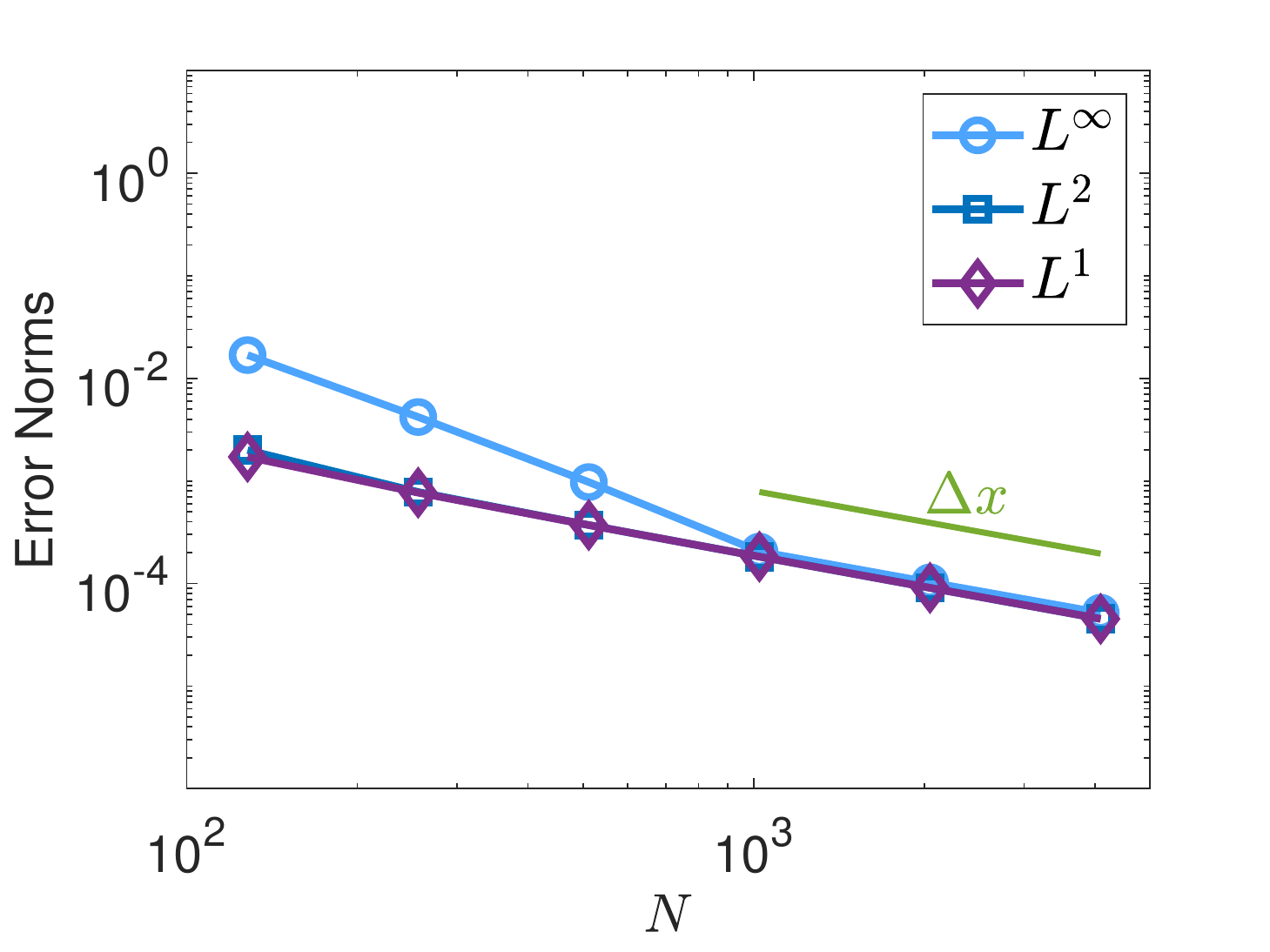}
\caption{\normalsize  Equation \eqref{poisson shit pde fine}, $L=8$}
\label{bigLfinecompleted}
\end{subfigure}
\caption[Refinement studies for solutions to the Poisson PDEs in Equations \eqref{poisson shit pde exp}-\eqref{poisson shit pde fine} found using the completed double layer method]{Refinement studies for solutions to Equations \eqref{poisson shit pde exp}-\eqref{poisson shit pde fine} found using the completed double layer method, with finite differences. The computational domain is the periodic box $[-4, 4]^2$. $\Omega$ is the region in the computational domain exterior to a circle of radius $0.25$, centered at the origin. A boundary point spacing of $ \Delta s \approx 0.75  \Delta x$ is used, and solution values within $m_1= 6$ meshwidths from the boundary are replaced using a second interpolation point $m_2=8$ meshwidths away from the boundary. All errors reported are absolute errors. }\label{poisson shit plots 2}
\end{figure}

We next apply the completed IBDL method, described in Equation \eqref{ibdl completed}, with $\eta=10$, to the two PDEs on the large computional domain of length $L=8$. Figure \ref{poisson shit plots 2} illustrates we get first order convergence for both PDEs with the completed method. However, we note that since we lose the second-order convergence for the PDE in Equation \eqref{poisson shit pde fine}, the errors for that PDE are larger. Further investigations need to be made into the different results from these two boundary value problems. 

\section{Near-boundary convergence}\label{ch 5 max norm}

As was discussed in Section \ref{5.1 discontinuity}, the discontinuity in the solution across the boundary prevents the method from achieving first-order convergence in the max norm near the boundary. As seen in Figure \ref{near boundary zoomed in}, the largest errors are confined to a region of about 2-3 meshwidths due to the support of the discrete delta function. However, there is also a wider region to which the numerical method spreads the errors. The physical length of the region on which the pointwise error fails to converge approaches $0$ as the grid is refined, so if one only needs the solution away from the boundary, one can proceed with the method as illustrated, or even omit the interpolation step altogether. However, in this section, we will examine this issue more closely and investigate choices that can be made to recover pointwise convergence for the entire PDE domain.

\subsection{Finite difference vs. Fourier spectral discretization }\label{5.5 FD vs spect}

\begin{figure}
\centering
\begin{subfigure}{0.495\textwidth}
\centering
\includegraphics[width=\textwidth]{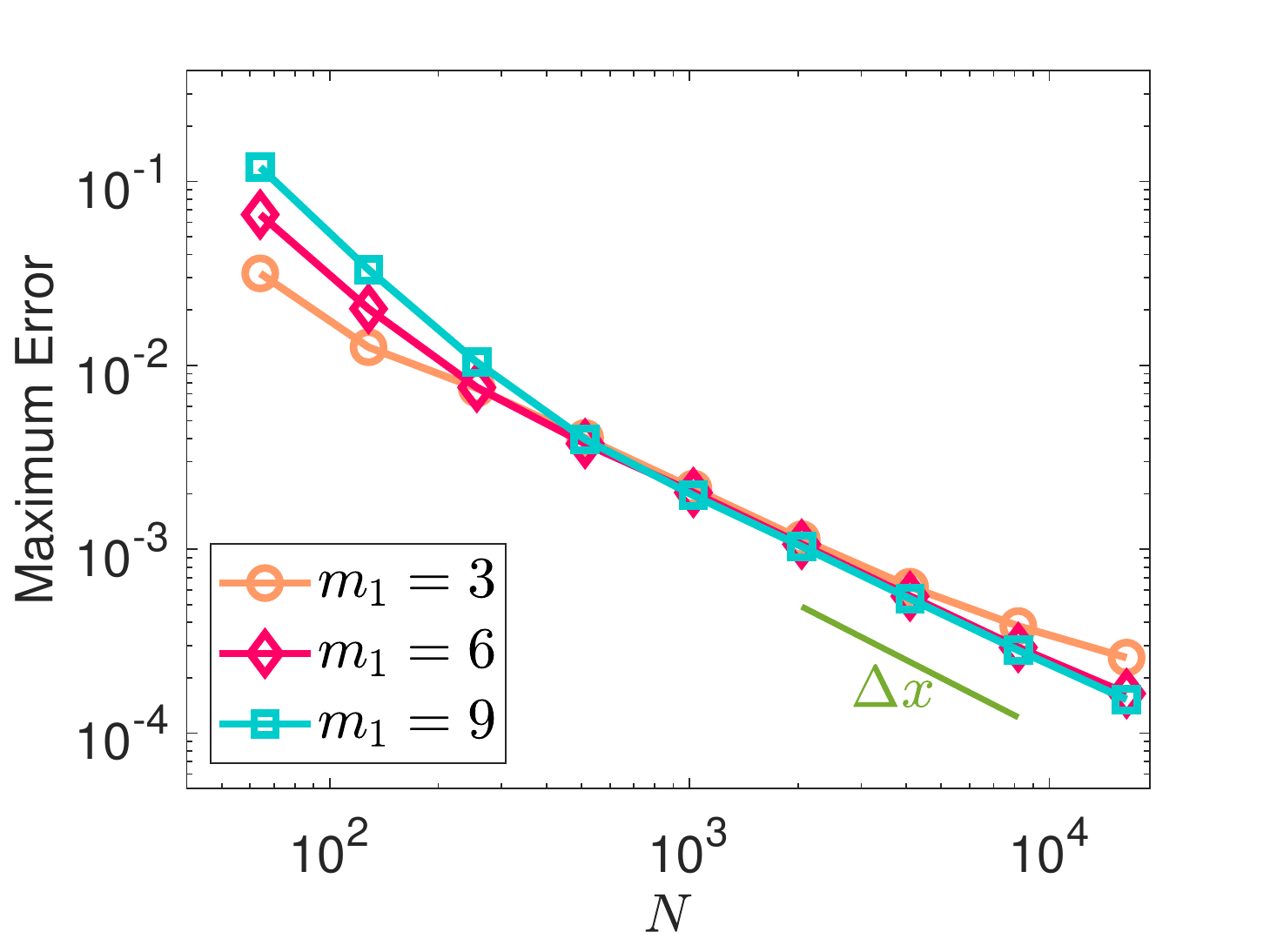}
\caption{\normalsize  Finite difference, Equation \eqref{pde again again}}
\label{sin FD}
\end{subfigure}
\begin{subfigure}{0.495\textwidth}
\centering
\includegraphics[width=\textwidth]{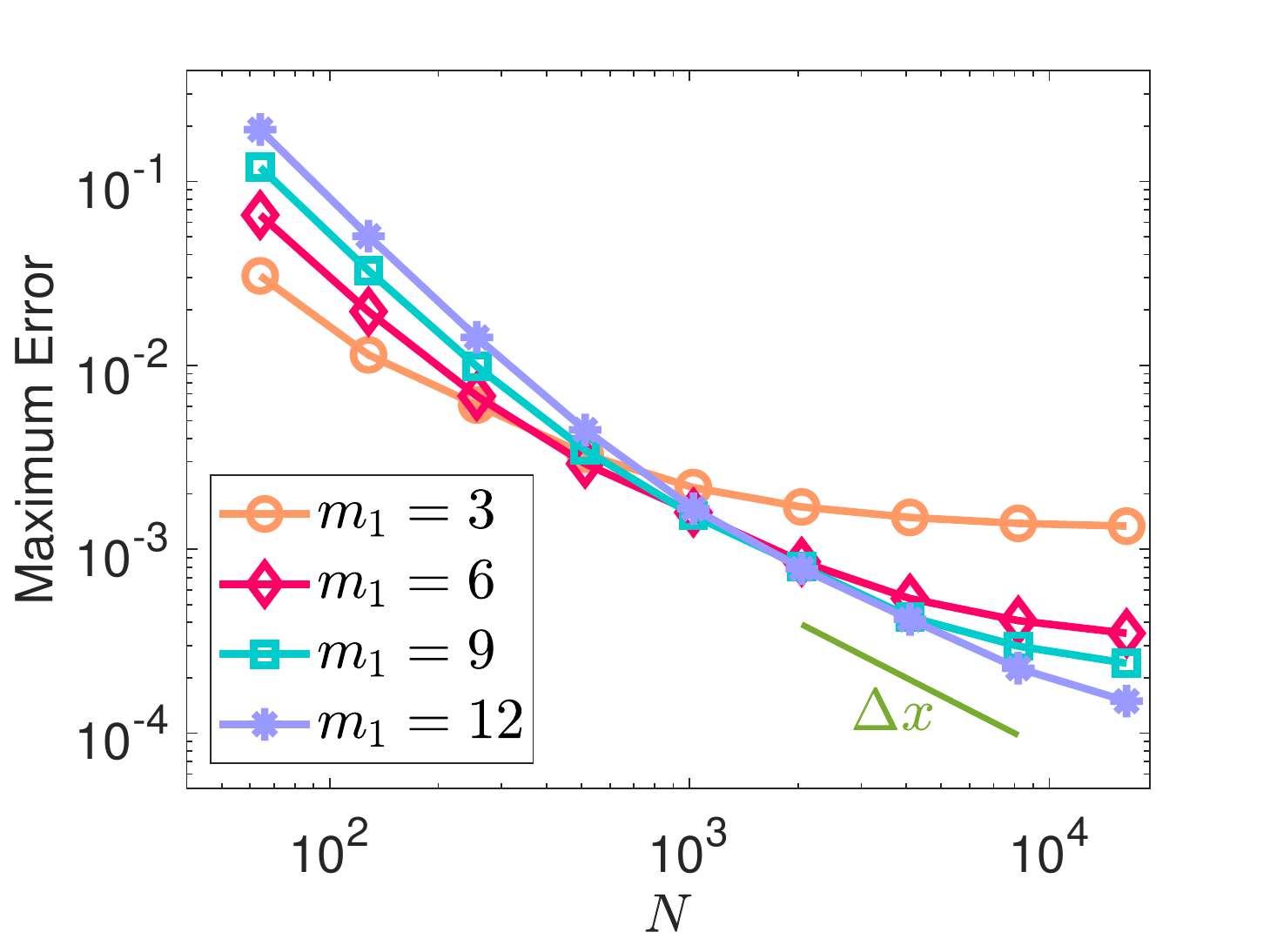}
\caption{\normalsize Fourier spectral, Equation \eqref{pde again again} }
\label{sin spectral}
\end{subfigure}
\caption[$L^{\infty}$ refinement studies for solving the Helmholtz PDE in Equation \eqref{pde again again} with the IBDL method using a range of interpolation widths and both finite difference and Fourier spectral discretization]{The $L^{\infty}$ refinement studies for solving Equation \eqref{pde again again} with the IBDL method. The computational domain is the periodic box $[-0.5, 0.5]^2$, $\Omega$ is the interior of a circle of radius 0.25, the boundary point spacing is $\Delta s \approx 0.75 \Delta x$, and $m_1$ is varied. Figure \ref{sin FD} is found using a finite difference method for discretization of the PDE, and Figure \ref{sin spectral} is found using a Fourier spectral method. }\label{FD vs spectral 1}
\end{figure}

We found that one key factor affecting the spread of error is the numerical method used to discretize the PDE operator. We therefore start by revisiting Equation \eqref{pde again again} and utilizing both finite difference and Fourier spectral methods for discretization. We extend our refinement study to an even finer mesh of $N=2^{14}$ and explore different values for $m_1$, the number of meshwidths from the boundary for which interpolation is used. For all of these, we use an interior interpolation point that is $m_2=m_1+2$ meshwidths from the boundary. Figure \ref{FD vs spectral 1} shows the corresponding $L^{\infty}$ refinement studies. In Figure \ref{sin FD}, we can see that $m_1=6$ is sufficient to maintain first-order pointwise convergence when using a finite difference discretization. Using a Fourier spectral method, on the other hand, requires a larger value of $m_1$. For low, fixed values of $m_1$, we see that once the error is small enough, the Fourier spectral method gives a diminishing rate of pointwise convergence and eventually a failure to converge.

\begin{figure}
\centering
\begin{subfigure}{0.495\textwidth}
\centering
\includegraphics[width=\textwidth]{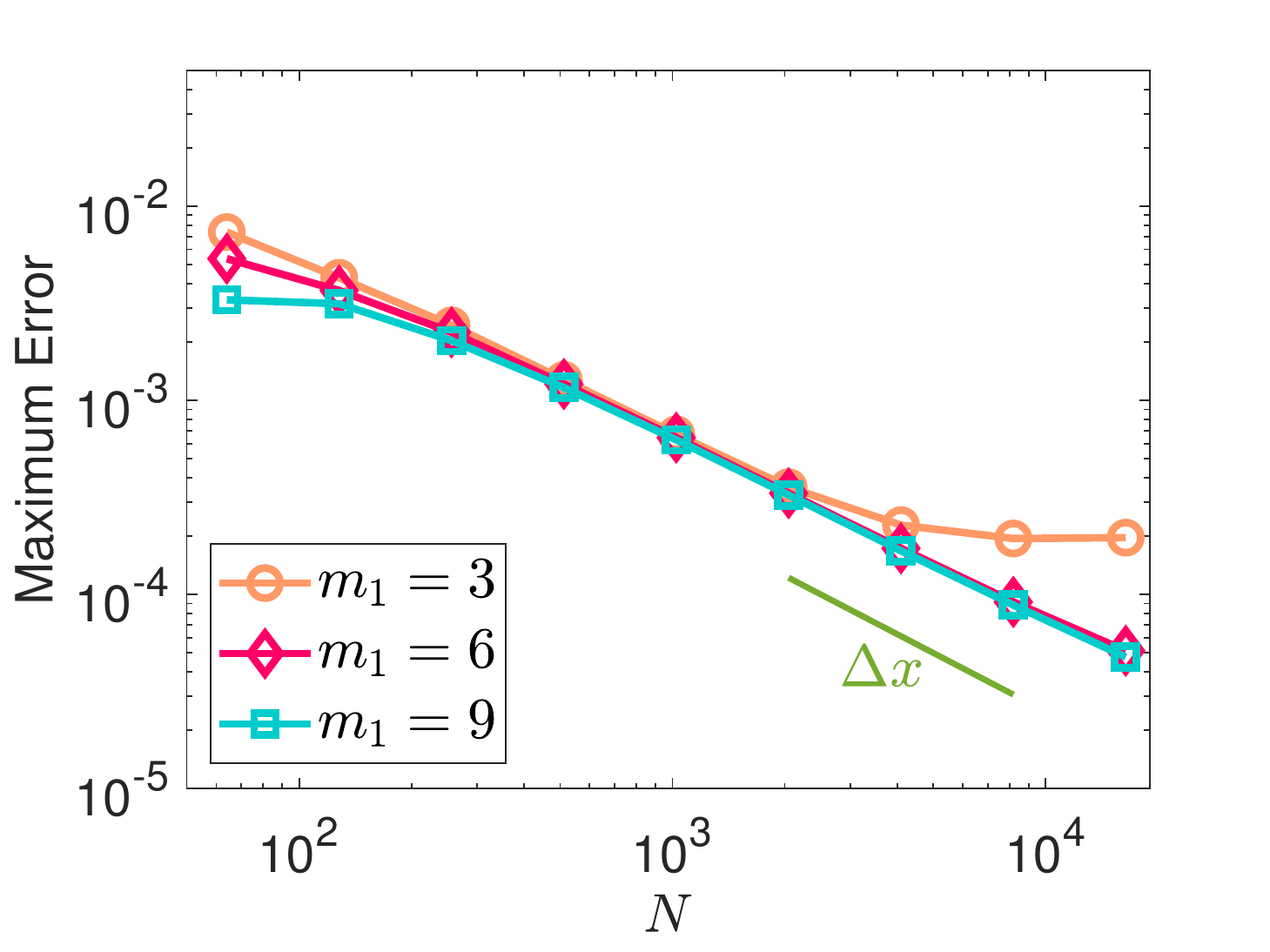}
\caption{\normalsize  Finite difference, Equation \eqref{pde again 2}}
\label{x FD}
\end{subfigure}
\begin{subfigure}{0.495\textwidth}
\centering
\includegraphics[width=\textwidth]{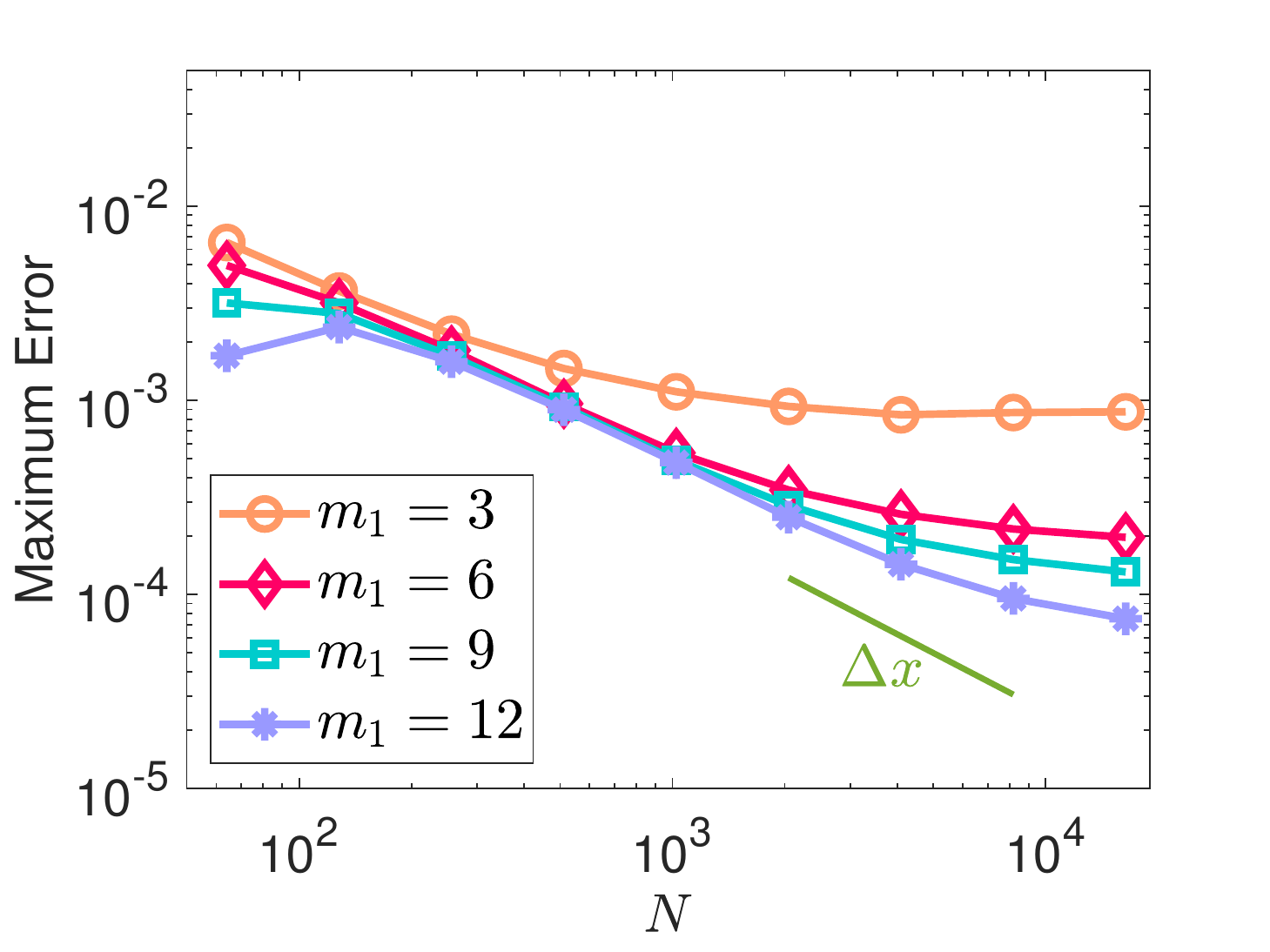}
\caption{\normalsize Fourier spectral, Equation \eqref{pde again 2} }
\label{x spectral}
\end{subfigure}
\begin{subfigure}{0.495\textwidth}
\centering
\includegraphics[width=\textwidth]{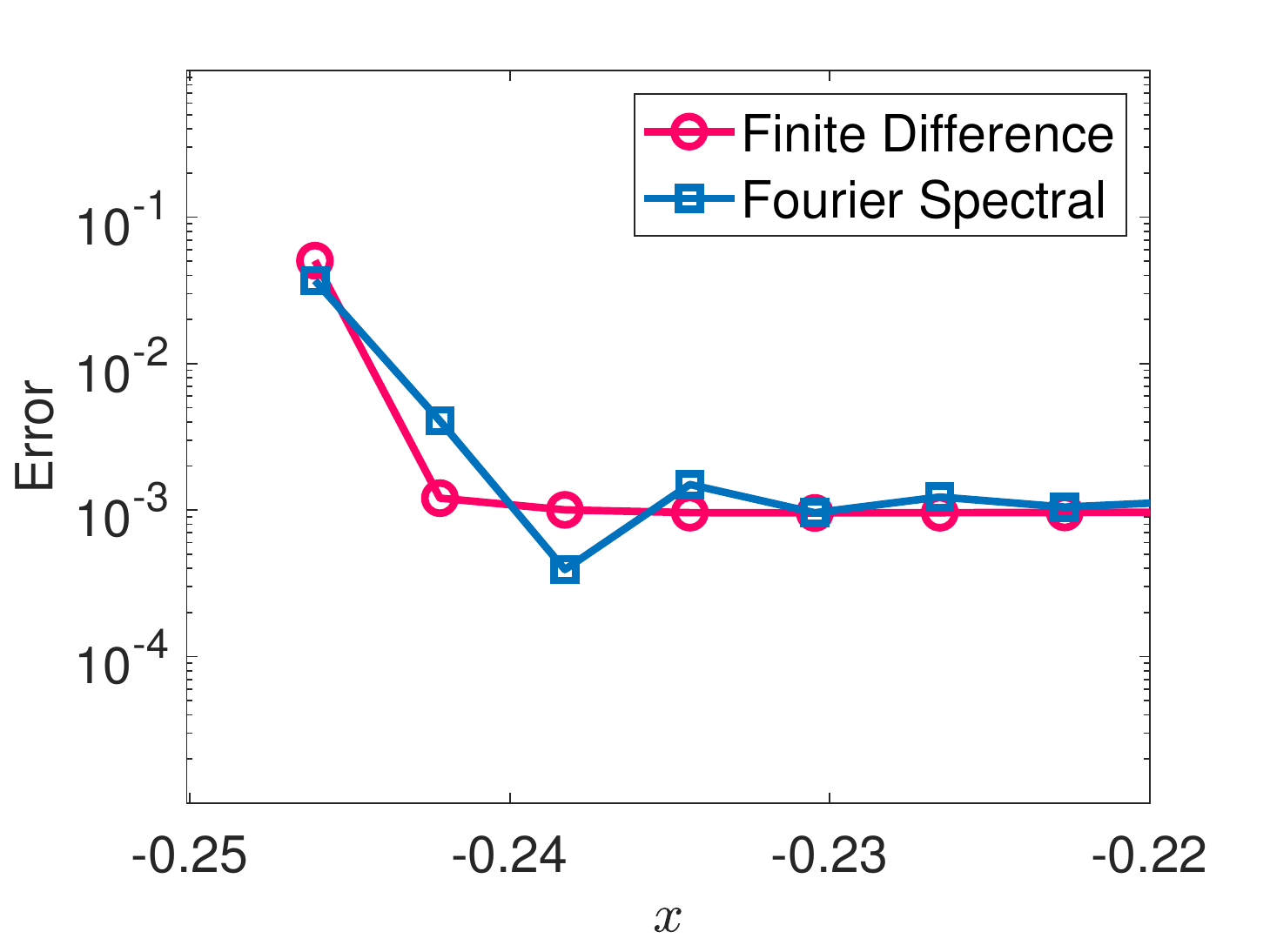}
\caption{\normalsize  Pointwise errors on coarse grid}
\label{coarse slice}
\end{subfigure}
\begin{subfigure}{0.495\textwidth}
\centering
\includegraphics[width=\textwidth]{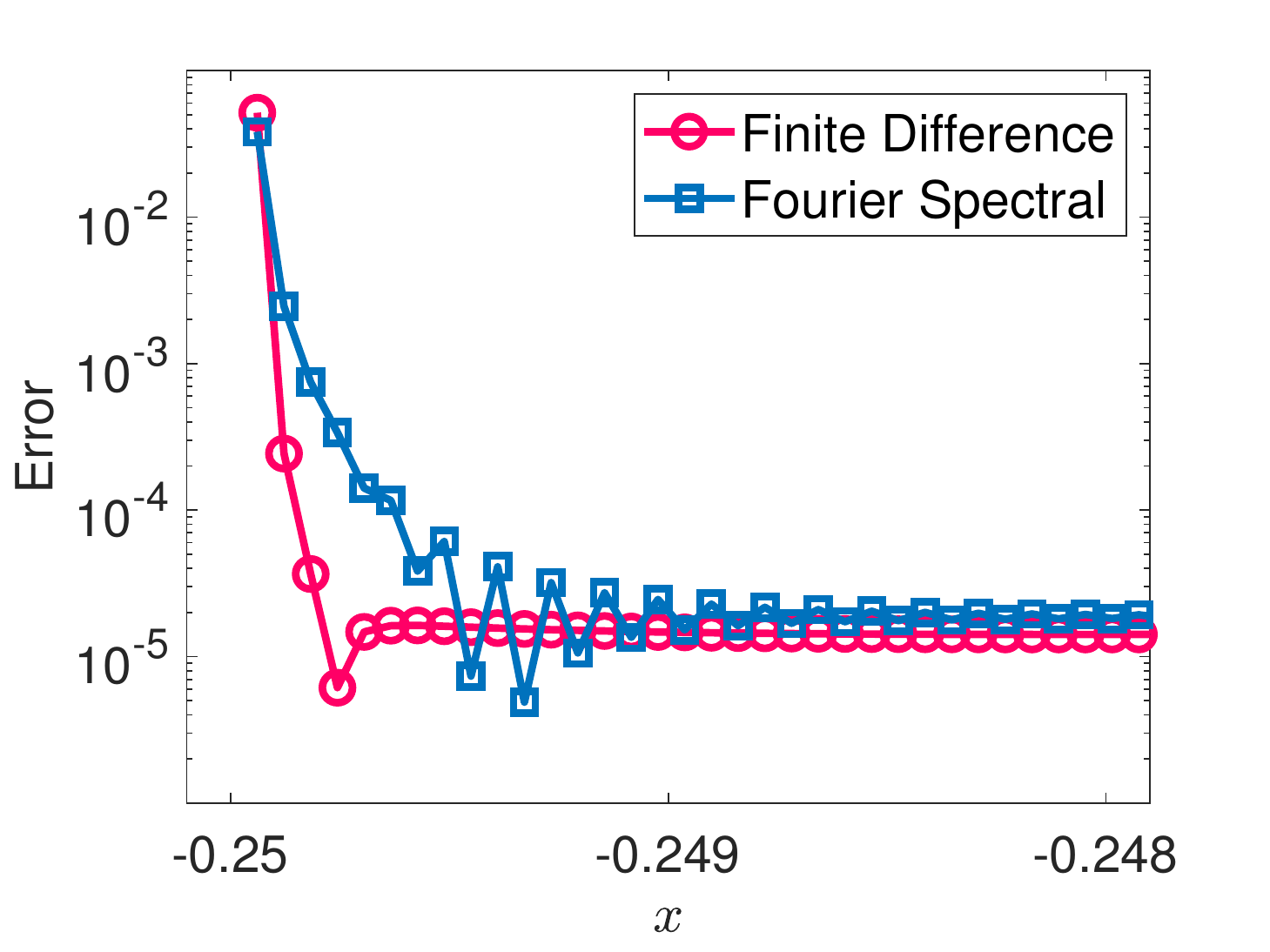}
\caption{\normalsize Pointwise errors on fine grid }
\label{fine slice}
\end{subfigure}
\caption[$L^{\infty}$ refinement studies using a range of interpolation widths and near-boundary pointwise error plots for the Helmholtz PDE in Equation \eqref{pde again 2} with a linear solution, solved with the IBDL method using both finite difference and Fourier spectral discretization]{Figures \ref{x FD}-\ref{x spectral} give the $L^{\infty}$ refinement studies for solving Equation \eqref{pde again 2} with the IBDL method for varying values of $m_1$. The computational domain is the periodic box $[-0.5, 0.5]^2$, $\Omega$ is the interior of a circle of radius 0.25, and the boundary point spacing is $\Delta s \approx 0.75 \Delta x$. Figure \ref{x FD} is found using a finite difference method for discretization of the PDE, and Figure \ref{x spectral} is found using a Fourier spectral method. Figures \ref{coarse slice}-\ref{fine slice} show the errors along the line $y=0$ when solving Equation \eqref{pde again 2} with the IBDL method with no interpolation and using both finite difference and Fourier spectral methods for discretization of the PDE. Figure \ref{coarse slice} uses $N=2^8$, and Figure \ref{fine slice} uses $N=2^{14}$. Each plot marker represents a gridpoint on the corresponding mesh.}\label{FD vs spectral 2}
\end{figure}

To further explore this phenomenon, we next look at the PDE
\begin{subequations} \label{pde again 2}
\begin{alignat}{2}
& \Delta u -  u = -(x+y) \qquad && \text{in } \Omega  \label{pde1 again 2}\\
&u=x+y \qquad && \text{on } \Gamma,  \label{pde2 again 2}
\end{alignat}
\end{subequations}
where $\Omega$ is the interior of a circle of radius 0.25, centered at the origin. The analytical solution is given by $u=x+y$. Figures \ref{x FD}-\ref{x spectral} show the $L^{\infty}$ refinement studies for this problem. Due to the linearity of the solution, the primary source of error is the discontinuity, and this allows us to observe the aforementioned behavior on coarser grids. We can see that the lack of convergence with the Fourier spectral method is worse for this problem, but we see that $m_1=6$ is still sufficient to maintain convergence using a finite difference method. 

To explore this phenomenon a little more closely, we also omit the interpolation step of the method and examine cross-sections of the error, located at $y=0$, for which the boundary point is located at $x=-0.25$. Figures \ref{coarse slice}-\ref{fine slice} show these error cross-sections on a relatively coarse grid of $N=2^8$ and a fine grid of $N=2^{14}$. Each marker on the plot represents a grid point on the mesh. For the coarse grid, we see in Figure \ref{coarse slice} that the pointwise errors for both methods decrease to the size of the interior error within a couple meshwidths. For the fine grid, however, the Fourier spectral method takes more meshwidths to decrease to the level of the interior error than the finite difference method. This illustrates the further spreading out of the error that takes place when using the Fourier spectral method due to the slow convergence of errors in a truncated Fourier series for a discontinuity. It is important to note, however, that while the number of meshwidths affected increases as the grid is refined, the physical distance affected still decreases. This is clear when observing the different scaling on the $x$-axes between Figures \ref{coarse slice} and \ref{fine slice}.

\begin{figure}
\centering
\begin{subfigure}{0.495\textwidth}
\centering
\includegraphics[width=\textwidth]{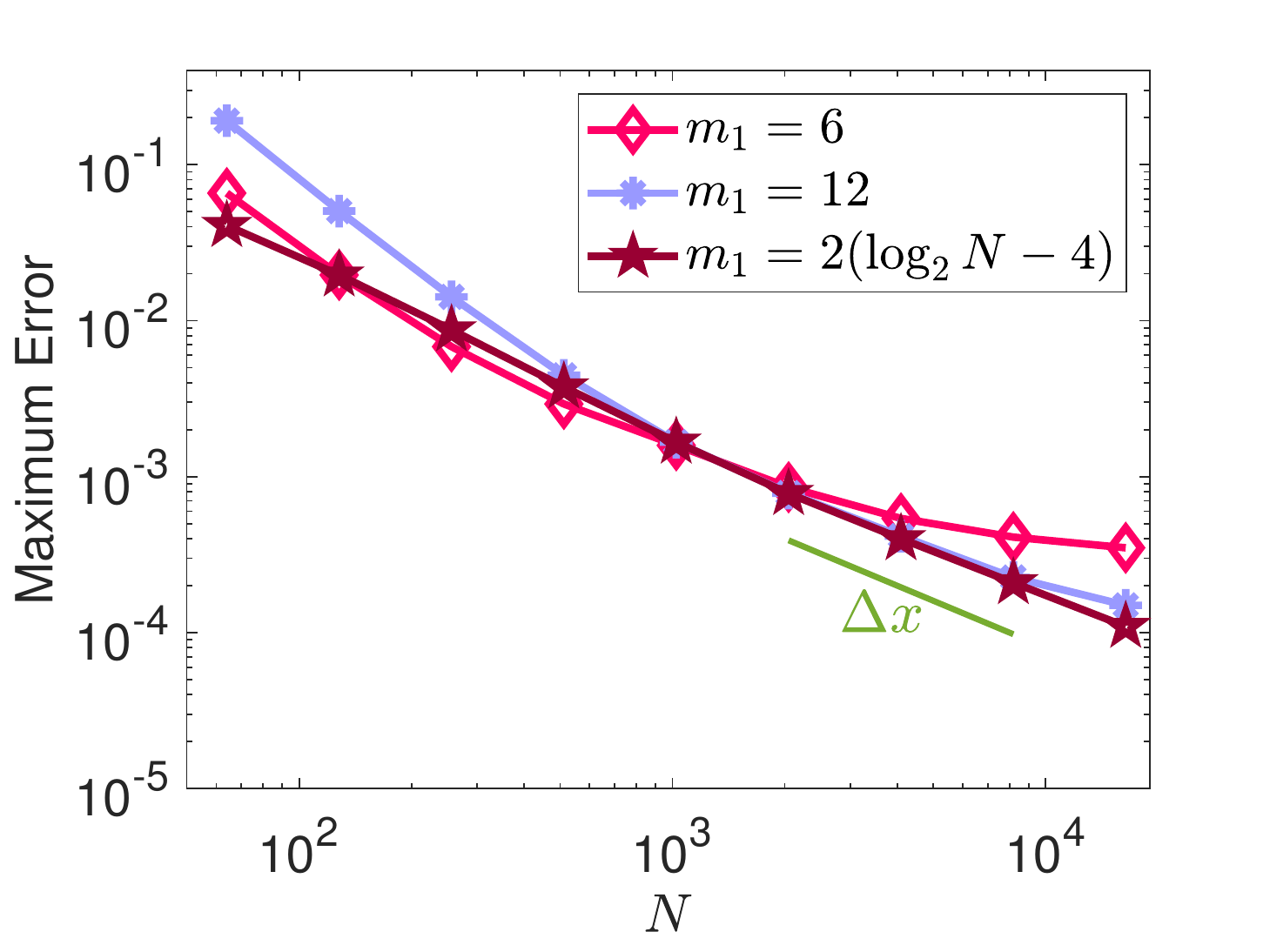}
\caption{\normalsize Fourier spectral, Equation \eqref{pde again again}}
\label{sin spectral 2}
\end{subfigure}
\begin{subfigure}{0.495\textwidth}
\centering
\includegraphics[width=\textwidth]{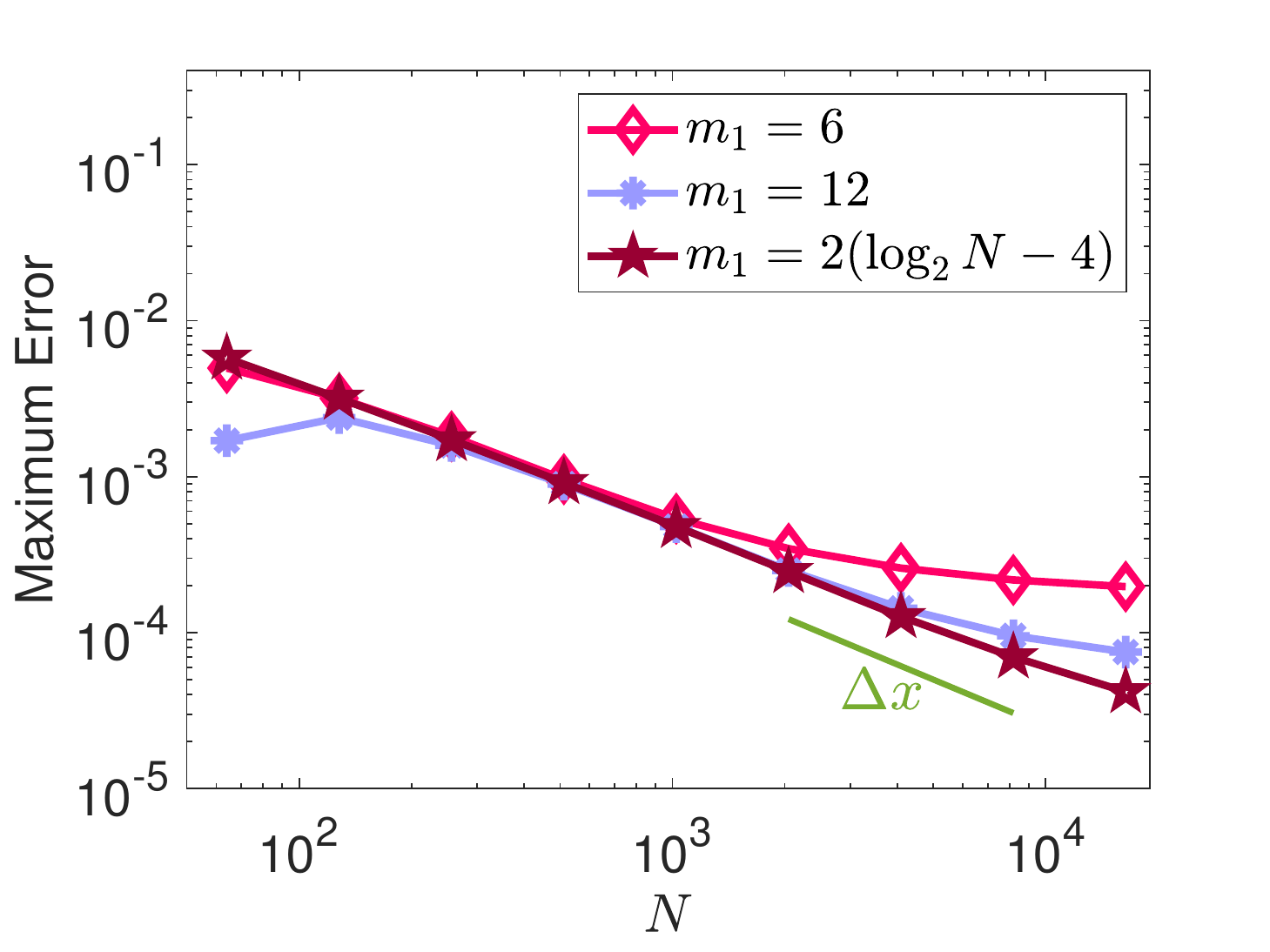}
\caption{\normalsize Fourier spectral, Equation \eqref{pde again 2}}
\label{x spectral 2}
\end{subfigure}
\caption[$L^{\infty}$ refinement studies for solving the Helmholtz PDEs in Equations \eqref{pde again again} and \eqref{pde again 2} with the IBDL method using a Fourier spectral discretization and an increasing interpolation width]{The $L^{\infty}$ refinement studies for solving Equations \eqref{pde again again} and \eqref{pde again 2}  using a Fourier spectral method for discretization of the PDE. The value of $m_1$ is again varied, but also included is an increasing function for $m_1$ given by $m_1=2(\log_2{(N)}-4)$. }\label{FD vs spectral 3}
\end{figure}

In order to recover pointwise convergence when utilizing a Fourier spectral method, we can choose a function for $m_1$ that increases as the mesh is refined. We have found that increasing logarithmically seems sufficient for these problems. For example, one could simply increase $m_1$ by 2 each time the number of grid points, $N$, is doubled. By doing this, we use interpolation at a greater number of meshwidths as we refine the grid, but, again, the physical distance on which we use interpolation is still decreasing relatively quickly towards 0. Figures \ref{sin spectral 2} and \ref{x spectral 2} illustrate using such a function for $m_1$ for Equations \eqref{pde again again} and \eqref{pde again 2}, respectively. The function we use is $m_1=2(\log_2{N}-4)$, which gives $m_1$ values ranging from $ 4$ to $20$ on these grid sizes. We see that this option recovers pointwise convergence.

\subsection{Boundary point spacing }\label{5.5 ds}
  
\begin{figure}
\centering
\begin{subfigure}{0.495\textwidth}
\centering
\includegraphics[width=\textwidth]{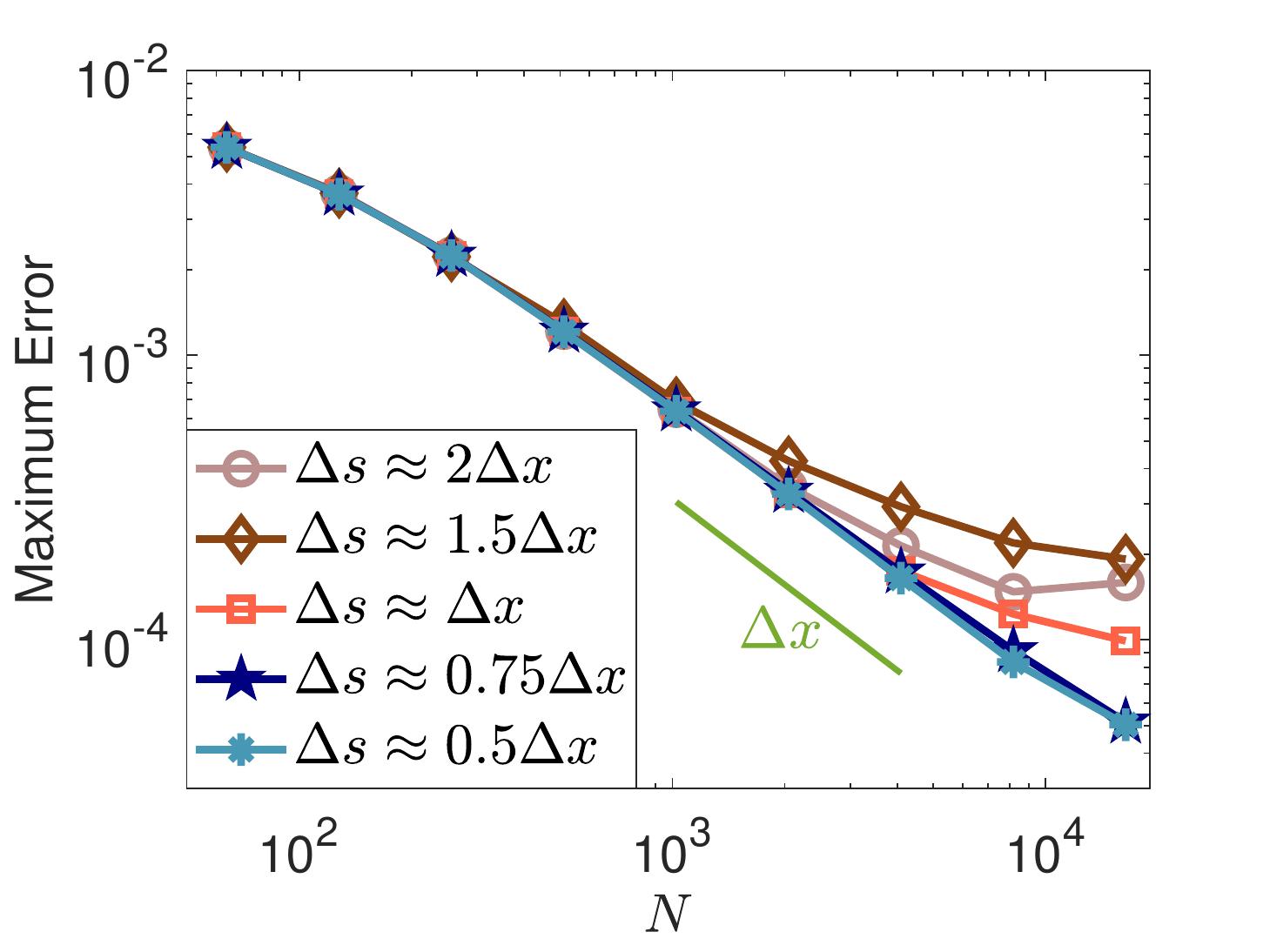}
\caption{\normalsize  Fixed $m_1=6$}
\label{x FD fixed m}
\end{subfigure}
\begin{subfigure}{0.495\textwidth}
\centering
\includegraphics[width=\textwidth]{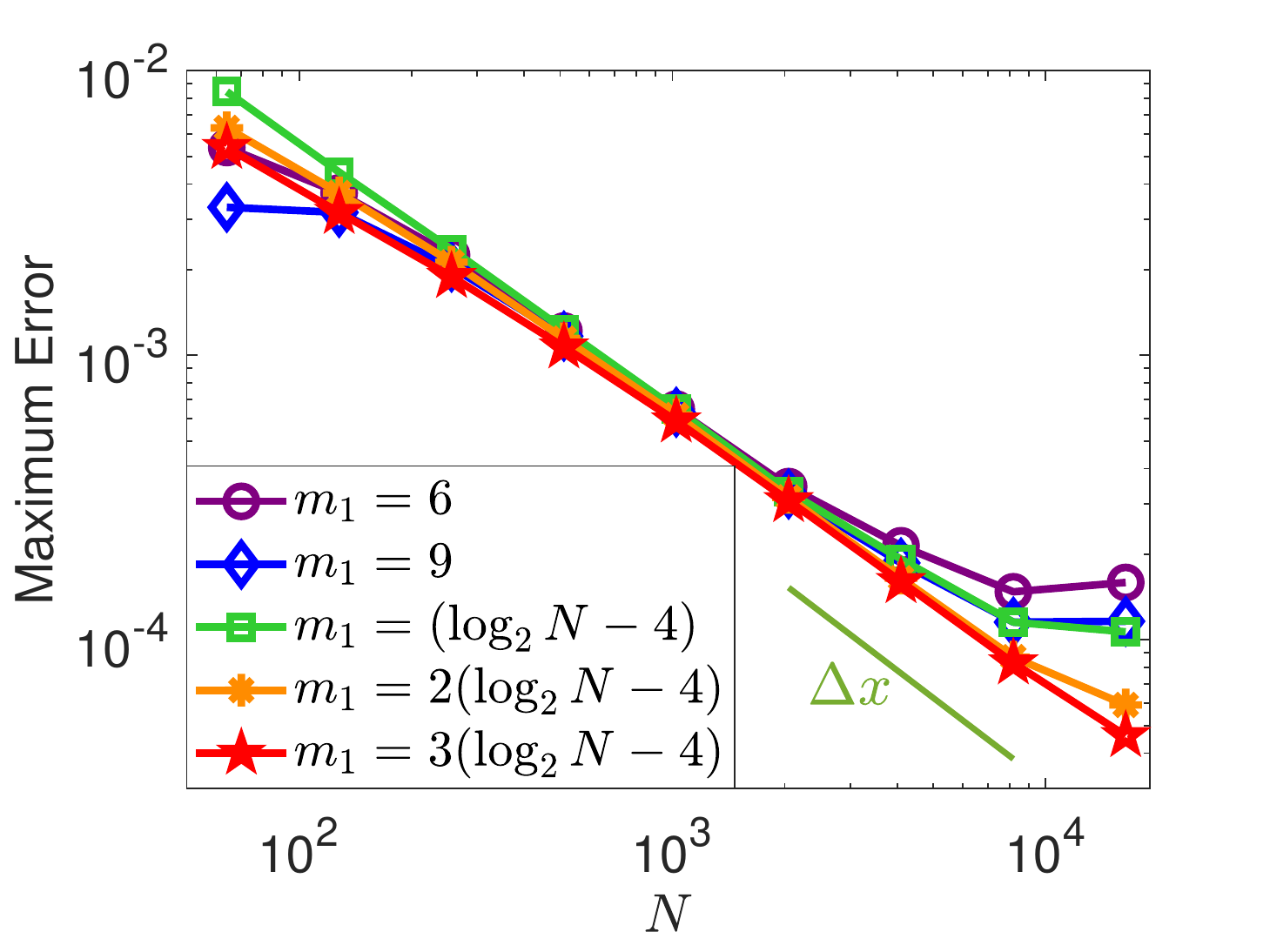}
\caption{\normalsize Fixed $\Delta s \approx 2 \Delta x$ }
\label{x FD fixed 2dx}
\end{subfigure}
\caption[$L^{\infty}$ refinement studies for solving Helmholtz Equation\eqref{pde again 2} with the IBDL method using a finite difference discretization, first for a range of boundary point spacings and then for a wide spacing and varied interpolation widths] {The $L^{\infty}$ refinement studies for solving Equation \eqref{pde again 2} with the IBDL method using a finite difference discretization of the PDE. Figure \ref{x FD fixed m} uses a fixed interpolation width of $m_1=6$ meshwidths and a second interpolation point $m_2=8$ meshwidths from the boundary, and the boundary point spacing is varied. Figure \ref{x FD fixed 2dx} has a wide boundary point spacing of  $\Delta s \approx 2 \Delta x$, and the interpolation width $m_1$ is varied.}\label{FD point spacing}
\end{figure}

  In addition to the PDE discretization and interpolation width, one other factor affecting the error spread is the boundary point spacing. We have observed that in general, using tightly spaced boundary points gives better pointwise convergence for a fixed interpolation width, especially in the case of the finite difference method. Specifically, we have seen that keeping $\Delta s \leq 0.75 \Delta x$ will help to avoid a slowing pointwise convergence rate. This is illustrated in Figure \ref{x FD fixed m}, which gives the $L^{\infty}$ refinement studies for Equation \eqref{pde again 2} using a finite difference discretization, $m_1=6$, and various the boundary point spacings. Furthermore, since we saw in Sections \ref{5.4 hh} and \ref{5.4 poisson} that the iteration count remains essentially unchanged by tightening the boundary points, the only numerical drawback to using this tighter point spacing is that each iteration is a little more expensive.

If one does want to use more widely spaced points, one can simply increase the interpolation widths suggested in the previous section. Figure \ref{x FD fixed 2dx} shows the $L^{\infty}$ refinement studies for Equation \eqref{pde again 2} when the boundary point spacing is given by $\Delta s \approx 2 \Delta x$ and the interpolation width is varied. In this case of more widely spaced boundary points, increasing $m_1$ logarithmically as the mesh is refined again recovers pointwise convergence.

In accordance with the observations seen in Section \ref{ch 5 max norm}, for the majority of Chapter \ref{chapter 5}, we have elected to use a finite difference PDE discretization, a boundary point spacing of $\Delta s \approx 0.75 \Delta x$, an interpolation width of $m_1=6$ meshwidths, and a second interpolation point $m_2=8$ meshwidths from the boundary. Similar results can be found by using a Fourier spectral method, a boundary point spacing of $\Delta s \approx 0.75 \Delta x$, an interpolation width given by $m_1=2(\log_2{(N)}-4)$, and a second point $m_2=m_1+2$ meshwidths from the boundary. We reiterate that using these options is not essential if one is only concerned with the solution away from the boundary because the width of the region on which the error does not converge pointwise does approach $0$ as the grid is refined. For reasons discussed in Sections \ref{6.4 space} and \ref{6.5 fd}, we primarily make use of Fourier spectral methods in Chapter \ref{chapter 6}.

\section{Results: convergence of potential strength}\label{ch 5 potential}

In addition to the increased speed with which we solve the system for the potential strength, $Q$, another benefit of the better conditioning of the IBDL method is that we are more easily able get convergence in this potential strength. It is a well-known problem that the poor conditioning of the Schur complement results in a noisy force distribution in the IBSL method \cite{Goza}, and this noise is larger for tighter boundary point spacing. Here, we illustrate that in the IBDL method, we obtain much smoother distributions with little need for filtering. 

To demonstrate this, we compute the IBSL and IBDL force distributions, $F$ and $Q$, obtained by solving the interior problem from Section \ref{5.4 hh}, where the computational domain is instead the periodic box $[-10,10]^2$. For estimating convergence, we approximate the exact distributions using a boundary element method, and this larger computational domain allows us to approximate the periodic Green's function with the free-space Green's function. 

\begin{figure}
\centering
\begin{subfigure}{0.33\textwidth}
\centering
\includegraphics[width=\textwidth]{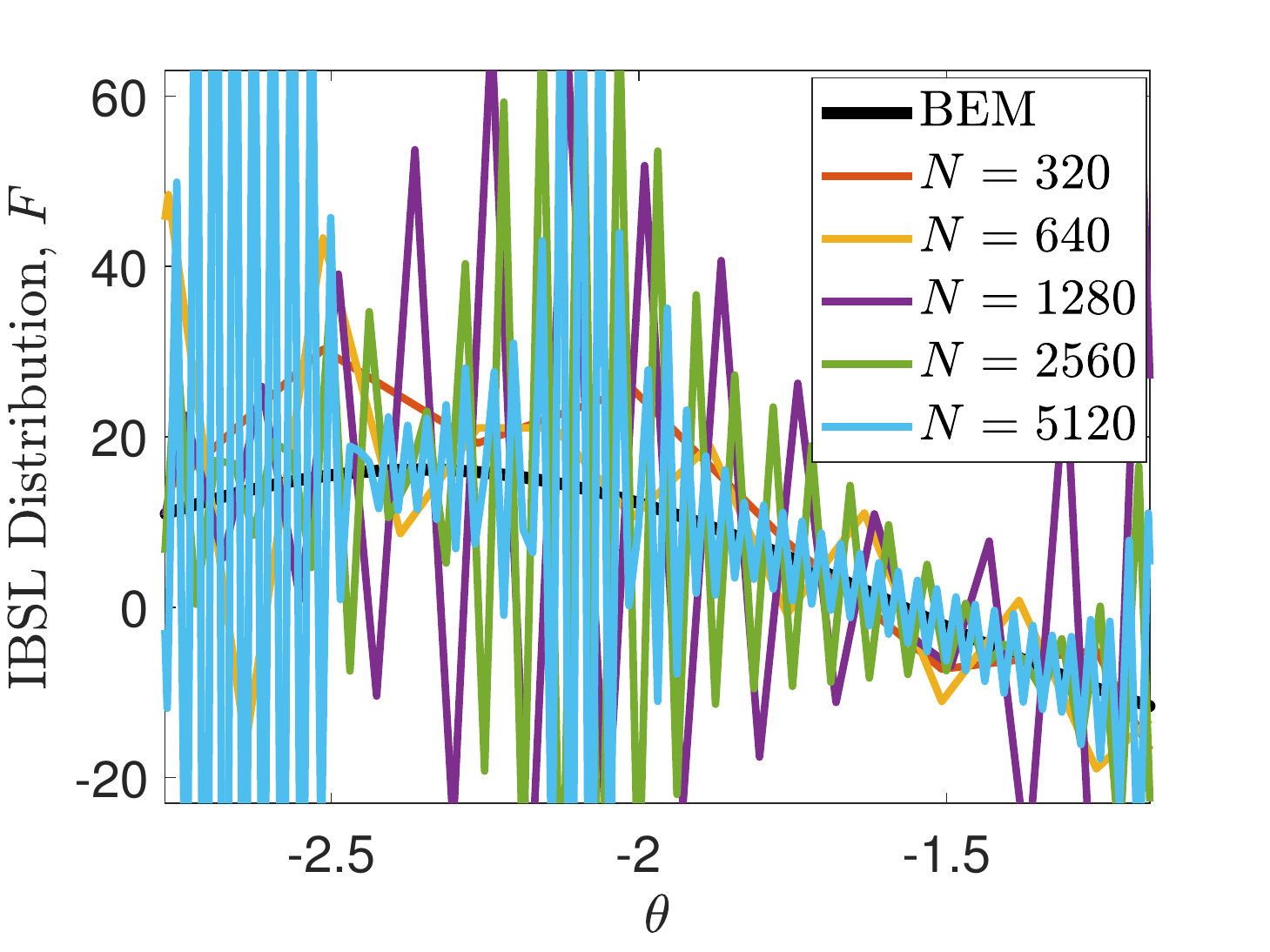}
\caption{ IBSL, $\Delta s \approx \Delta x$}
\label{ibsl f plot ds1}
\end{subfigure}
\begin{subfigure}{0.33\textwidth}
\centering
\includegraphics[width=\textwidth]{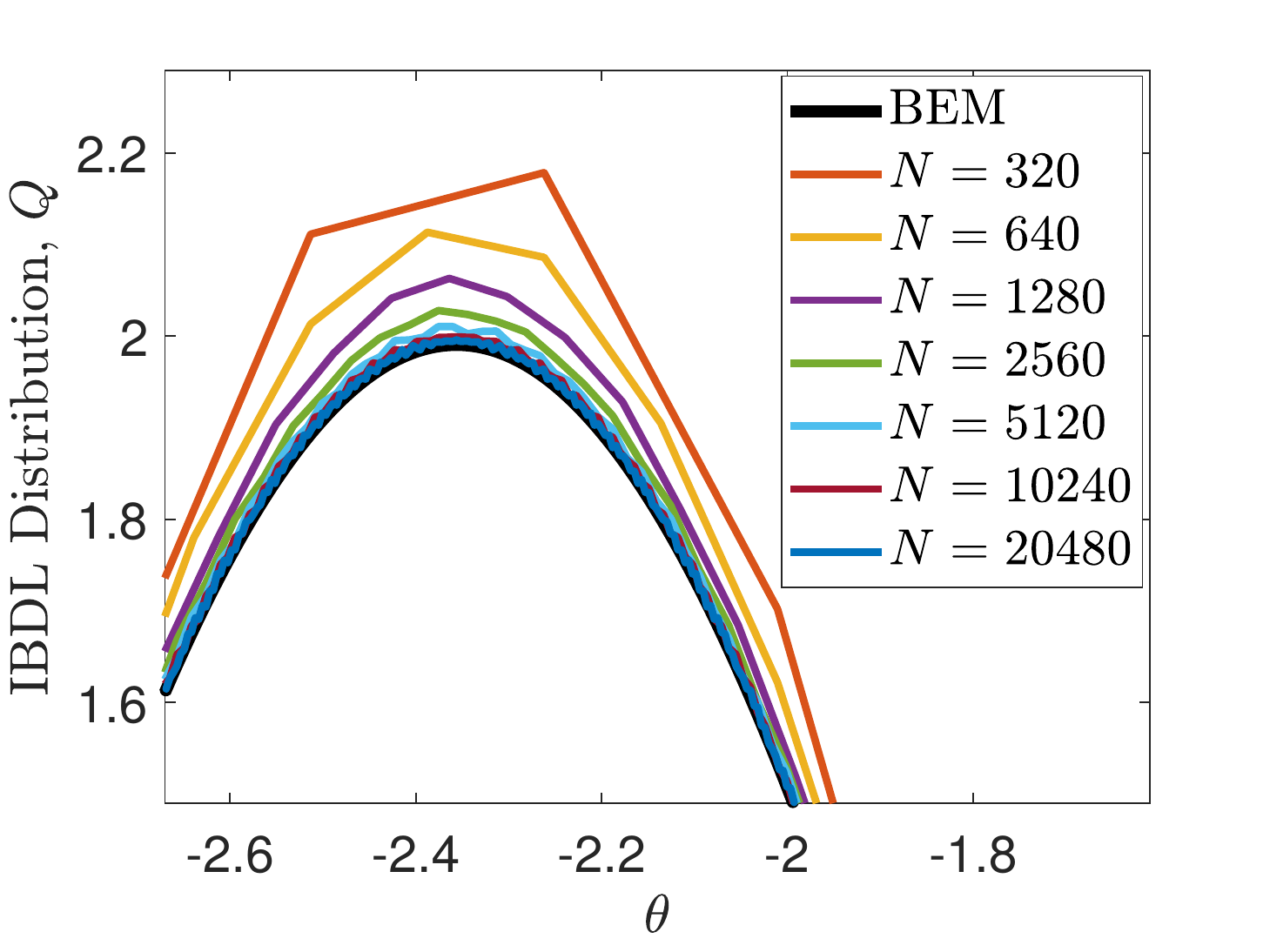}
\caption{ IBDL, $\Delta s \approx \Delta x$}
\label{ibdl q plot ds1}
\end{subfigure}
\begin{subfigure}{0.33\textwidth}
\centering
\includegraphics[width=\textwidth]{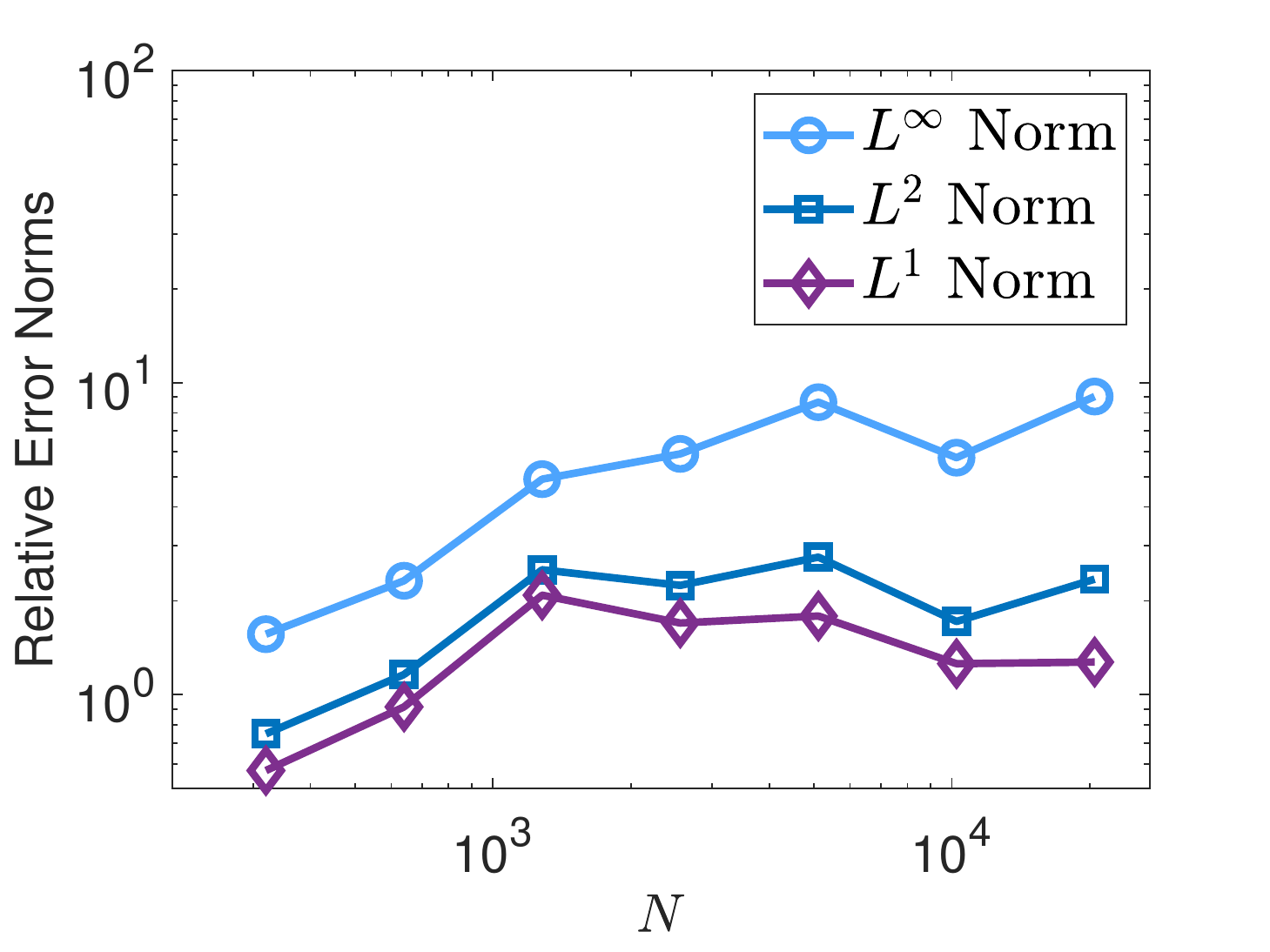}
\caption{ IBSL, $\Delta s \approx \Delta x$}
\label{ibsl f refinement}
\end{subfigure}
\begin{subfigure}{0.33\textwidth}
\centering
\includegraphics[width=\textwidth]{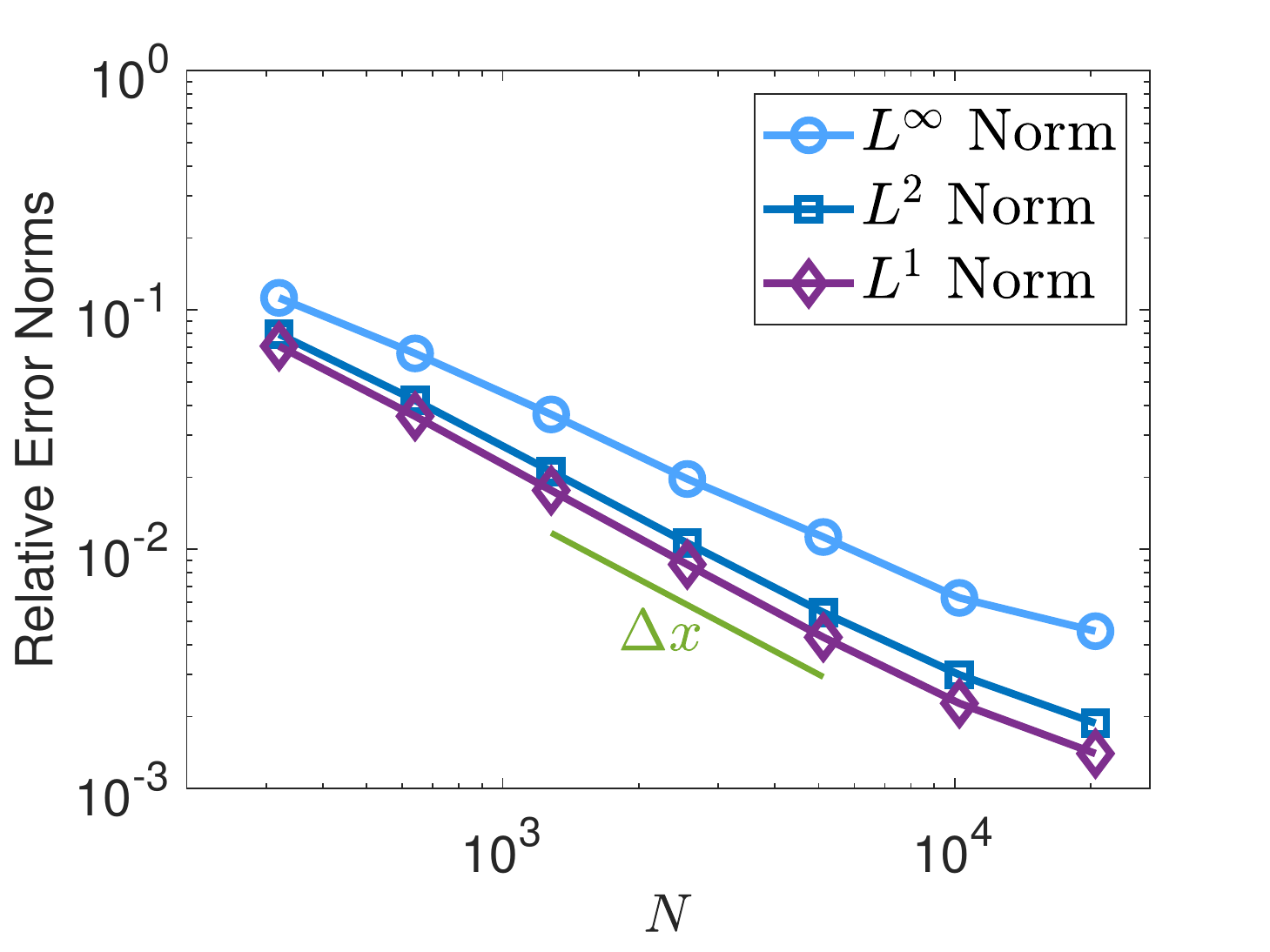}
\caption{ IBDL, $\Delta s \approx \Delta x$}
\label{ibdl q refinement}
\end{subfigure}
\begin{subfigure}{0.33\textwidth}
\centering
\includegraphics[width=\textwidth]{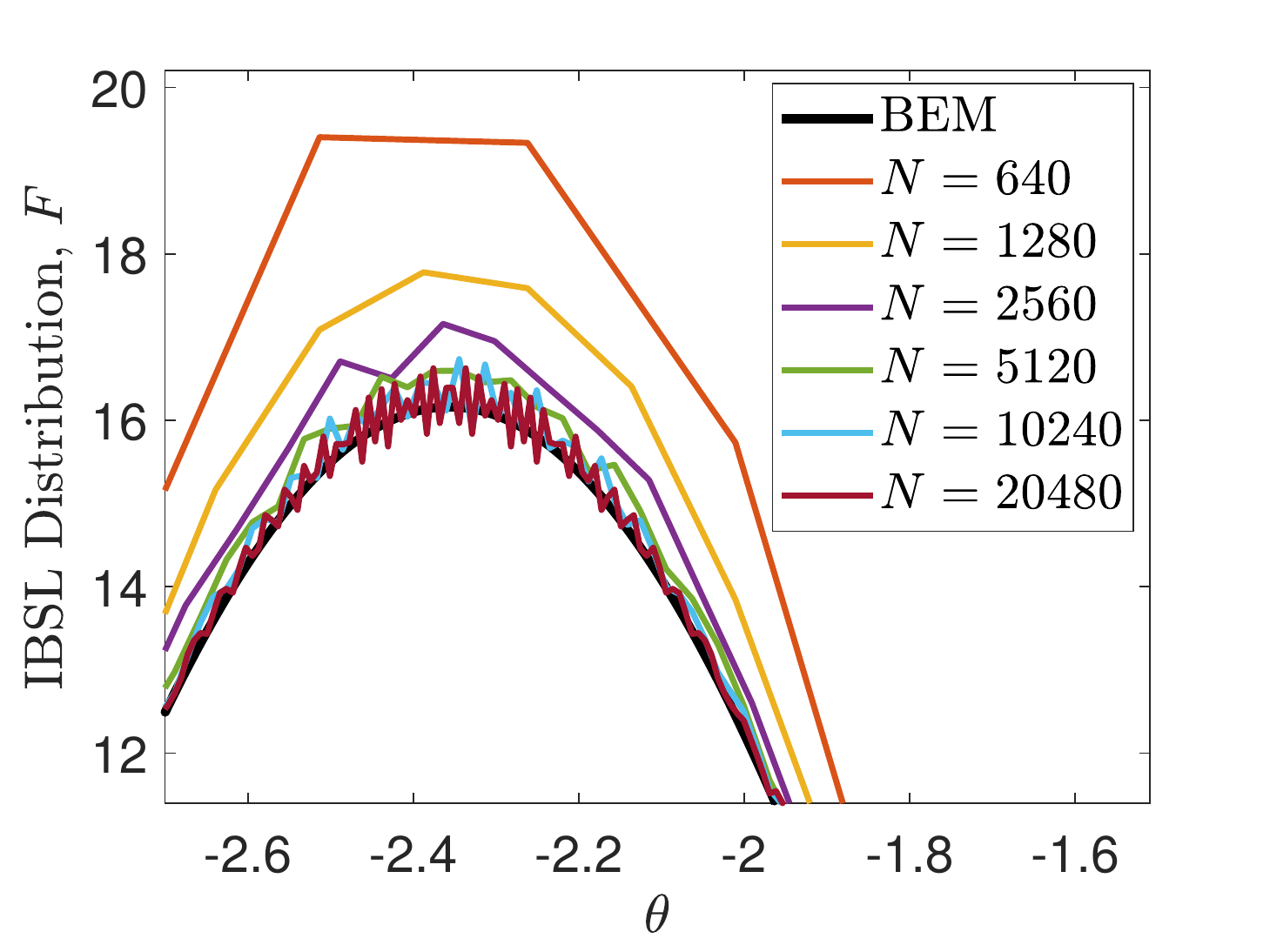}
\caption{ IBSL, $\Delta s \approx 2\Delta x$}
\label{ibsl f plot ds2}
\end{subfigure}
\begin{subfigure}{0.33\textwidth}
\centering
\includegraphics[width=\textwidth]{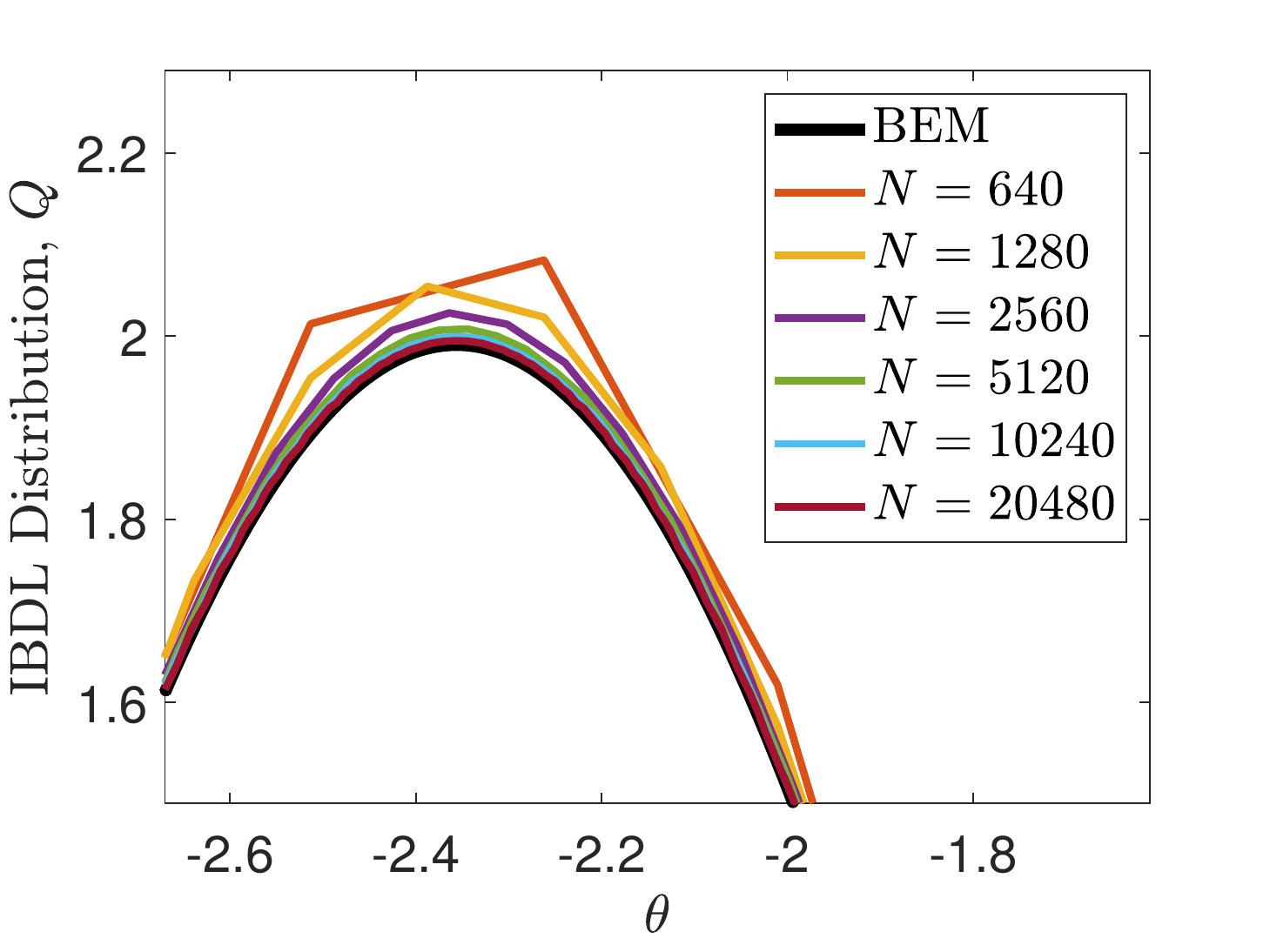}
\caption{ IBDL, $\Delta s \approx2 \Delta x$}
\label{ibdl q plot ds2}
\end{subfigure}
\begin{subfigure}{0.33\textwidth}
\centering
\includegraphics[width=\textwidth]{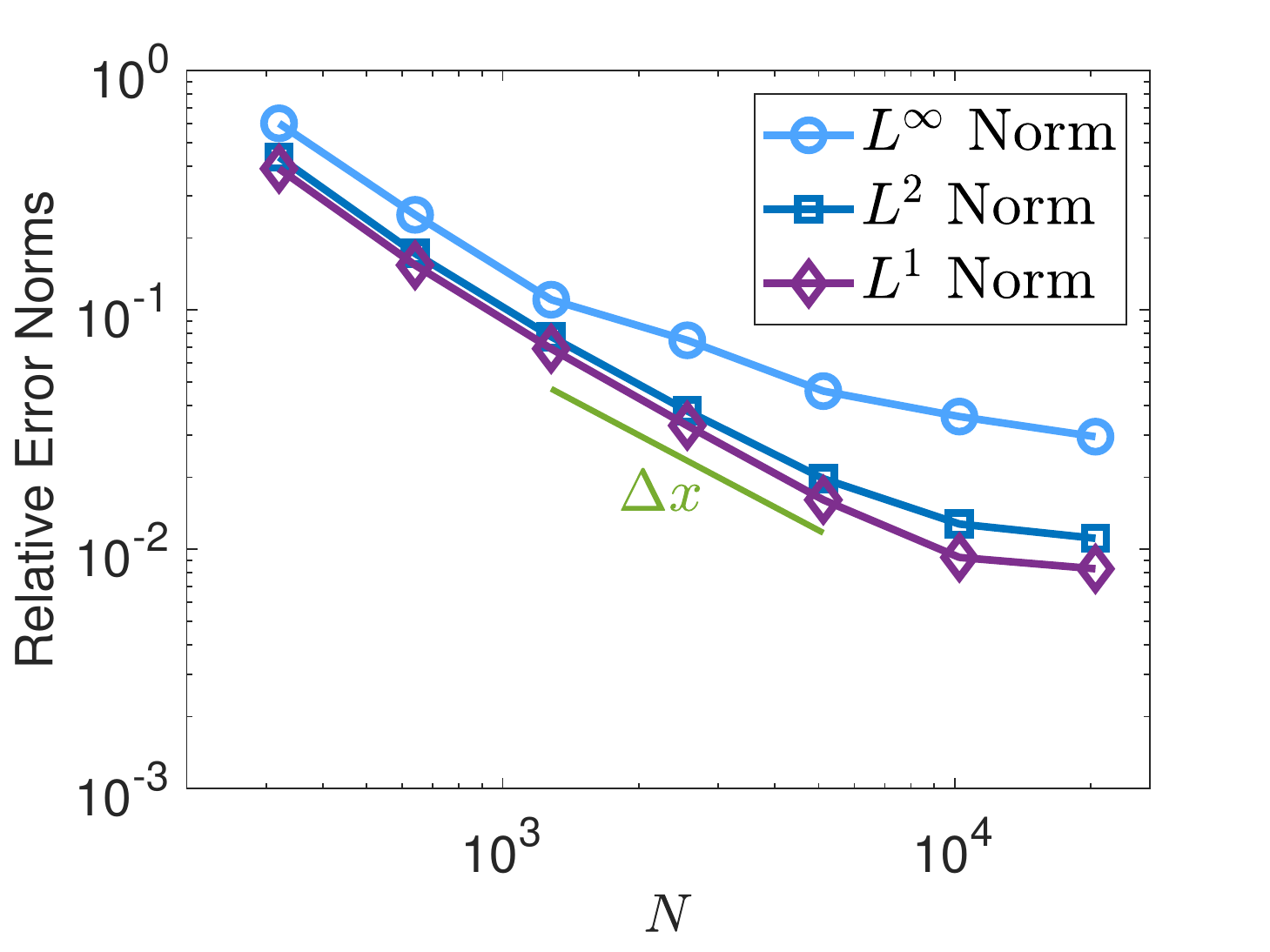}
\caption{ IBSL, $\Delta s \approx 2\Delta x$}
\label{ibsl f refinement ds2}
\end{subfigure}
\begin{subfigure}{0.33\textwidth}
\centering
\includegraphics[width=\textwidth]{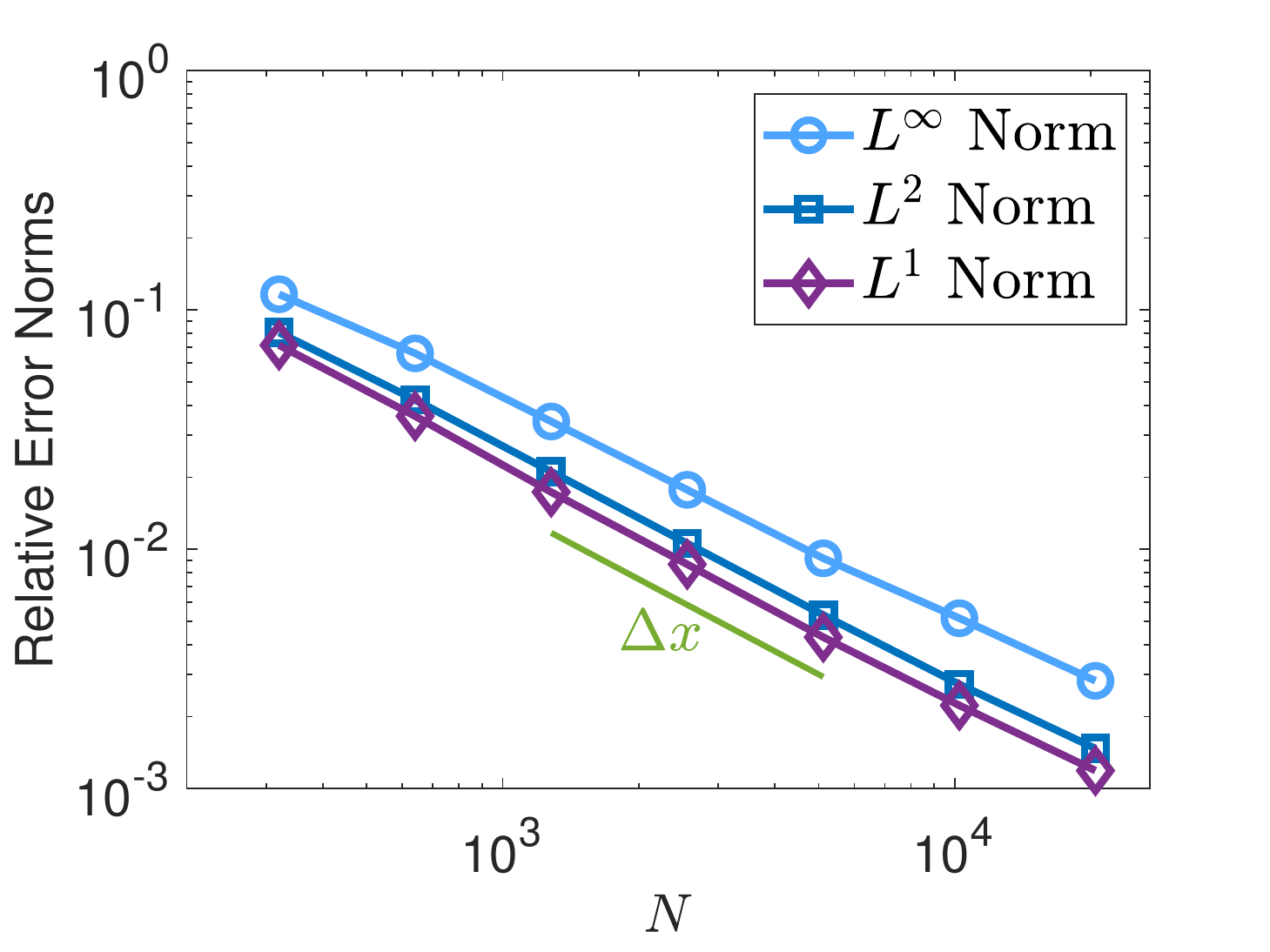}
\caption{ IBDL, $\Delta s \approx2 \Delta x$}
\label{ibdl q refinement ds2}
\end{subfigure}
\caption[Plots of IBSL and IBDL potential strengths for the Helmholtz PDE in Equation \eqref{pde again again} and refinement studies, approximating error using a boundary element potential strength]{Figures \ref{ibsl f plot ds1}-\ref{ibdl q plot ds1} and  \ref{ibsl f plot ds2}-\ref{ibdl q plot ds2} show portions of the IBSL and IBDL potential strengths for Equation \eqref{pde again again}, plotted against the angle $\theta$. Boundary point spacing is given by $\Delta s \approx \Delta x$ for Figures \ref{ibsl f plot ds1}-\ref{ibdl q plot ds1} and $\Delta s \approx 2\Delta x$ for Figures \ref{ibsl f plot ds2}-\ref{ibdl q plot ds2}. Black curves correspond to distributions found using a boundary element method with uniform, straight line elements and Gauss-Legendre quadrature. The BEM distributions are used to approximate the error, normalized by the maximum value. Figures \ref{ibsl f refinement}-\ref{ibdl q refinement} and \ref{ibsl f refinement ds2}-\ref{ibdl q refinement ds2} show the corresponding refinement studies. }\label{Q and F stuff }
\end{figure}

In Figures \ref{ibsl f plot ds1}-\ref{ibdl q refinement}, for which we use $\Delta s \approx \Delta x$, we see approximately first-order convergence for the IBDL distribution $Q$ and no convergence for the IBSL distribution $F$. Figure \ref{ibsl f plot ds1} shows that $F$ quickly becomes very noisy as we refine the grid. As seen in Figure \ref{ibdl q plot ds1}, there is much less noise in $Q$.  Note that the plots only show a portion of the distributions in order to more easily view this noise. 

For Figures \ref{ibsl f plot ds2}-\ref{ibdl q refinement ds2}, we increase the spacing for the boundary points to $\Delta s \approx 2 \Delta x$. Then we see some convergence in the IBSL distribution $F$. However, the IBDL potential strength $Q$ achieves greater convergence and a much lower relative error than $F$. Therefore, we see that the IBDL method in general allows us to achieve a much more accurate distribution, and the next section demonstrates an application for which this convergence is vital.

\section{Results: Neumann boundary conditions}\label{ch 5 neumann results}

As discussed in Section \ref{ch 5 neumann}, we are able to use the connection to integral equations to solve a PDE with Neumann boundary conditions by solving 
\begin{equation}
(-S^* \L^{-1}\widetilde S)U_b-\frac12 U_b=S^*\L^{-1}SV_b-S^*\L^{-1}\tilde g \label{Neumann saddle point again}
\end{equation} 
for the unknown boundary values $U_b$. By doing this, we are again inverting an operator corresponding to a second-kind integral operator with a small condition number. Therefore, we are able to get convergence in this boundary distribution, as discussed in the previous section. When we have Neumann boundary conditions, the convergence of the distribution, $U_b$ is especially important because $U_b$ is used when performing the interpolation step for grid points near the boundary. In the Dirichlet case, $U_b$ is known, providing us with more accurate interpolation data, but in the case of Neumann, the values both in the interior and on the boundary are approximate. 

\begin{figure}
\centering
\begin{subfigure}{0.495\textwidth}
\centering
\captionsetup{justification=centering}
\includegraphics[width=\textwidth]{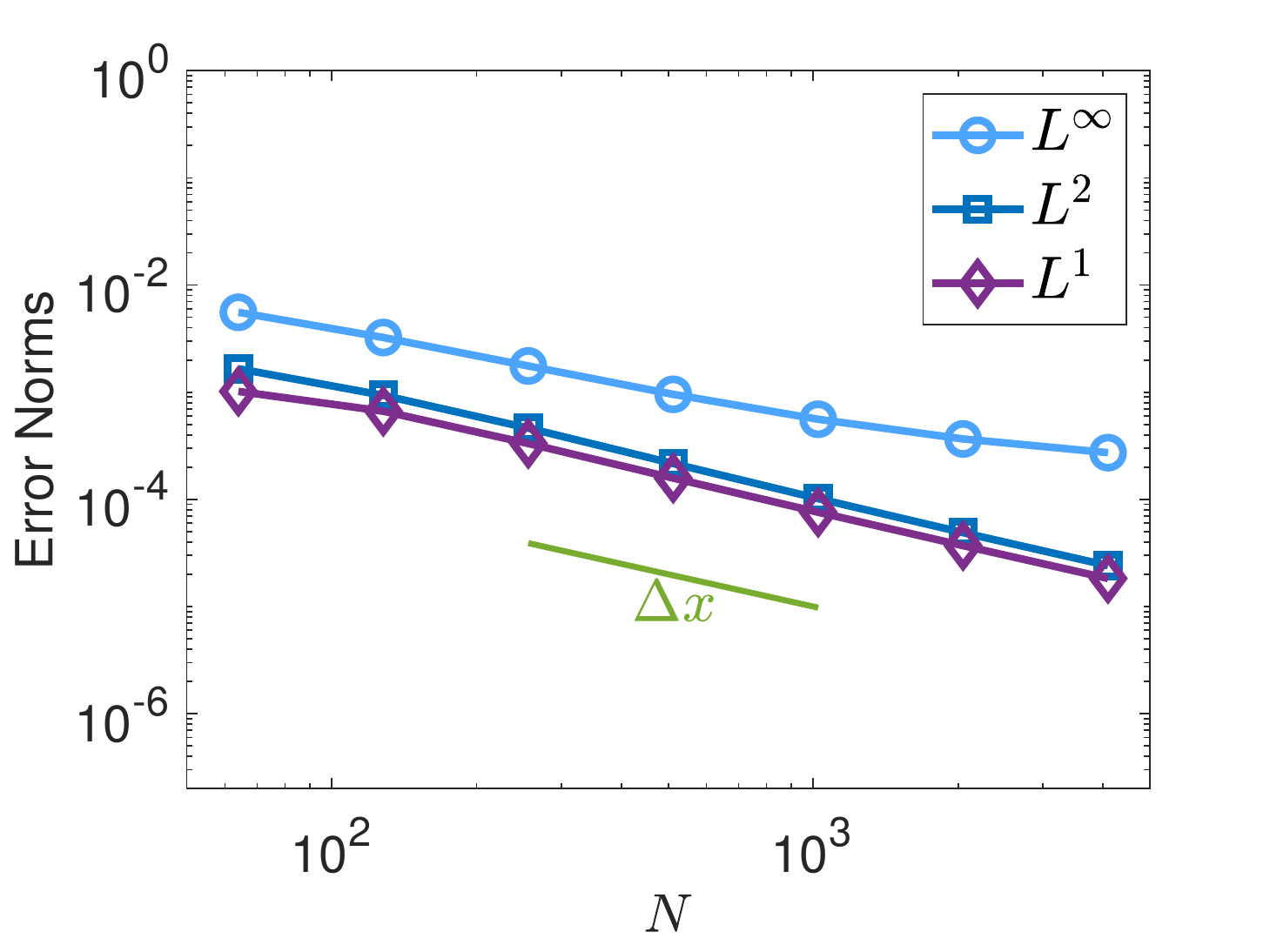}
\caption{\normalsize  Solution refinement study}
\label{neumann no filter 075}
\end{subfigure}
\begin{subfigure}{0.495\textwidth}
\centering
\captionsetup{justification=centering}
\includegraphics[width=\textwidth]{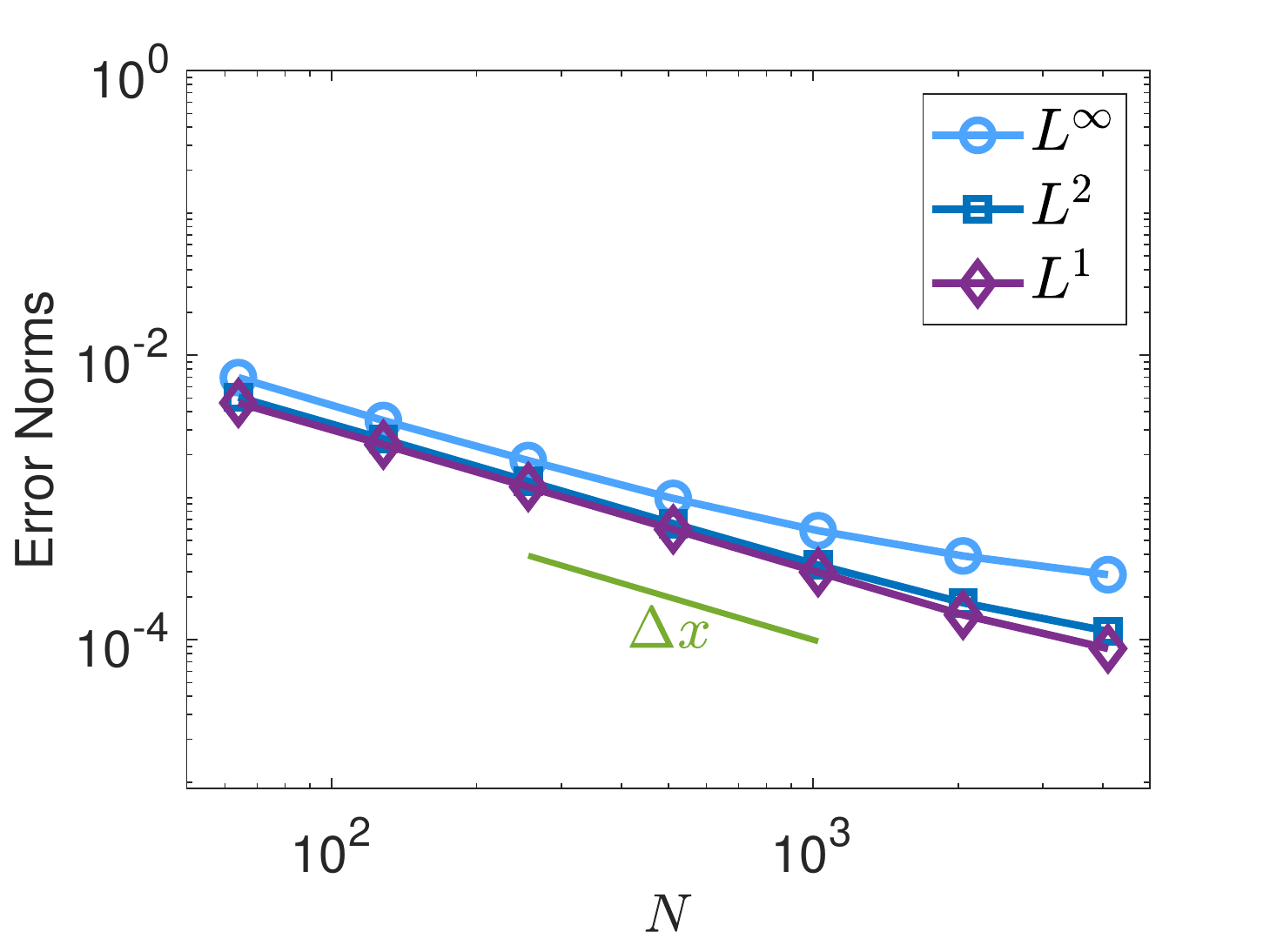}
\caption{\normalsize  $U_b$ refinement study}
\label{neumann ub no filter 075}
\end{subfigure}
\caption[Refinement studies for $U_b$ and the solution from solving the Helmholtz PDE in Equation \eqref{pde again 4} with Neumann boundary conditions]{Refinement studies for solving Equation \eqref{pde again 4} with the method discussed in Section \ref{ch 5 neumann}. The computational domain is the periodic box $[-0.5, 0.5]$, $\Omega$ is the interior of a circle of radius 0.25, and the prescribed normal derivatives on the boundary are given by $\partial u/\partial n|_{\Gamma}= 8(x^2-y^2)$. We use a boundary point spacing of $\Delta s \approx 0.75 \Delta x$ and replace solution values within $m_1 = 6$ meshwidths from the boundary using interior interpolation points $m_2=8$ meshwidths away from the boundary. Figures \ref{neumann no filter 075} and \ref{neumann ub no filter 075} show the refinement studies for the absolute solution errors and the boundary values, $U_b$, respectively.}\label{Neumann plots}
\end{figure}

We consider the PDE given by 
\begin{subequations} \label{pde again 4}
\begin{alignat}{2}
& \Delta u -  u =-(x^2-y^2)\hspace{0.2cm} && \text{in } \Omega  \label{pde1 again 4}\\
&\frac{\partial u}{\partial n}= 8(x^2-y^2)&& \text{on } \Gamma,  \label{pde2 again 4}
\end{alignat}
\end{subequations}
where $\Omega$ is the interior of a circle of radius 0.25, centered at the origin. The analytical solution is given by $u=x^2-y^2$. Our computational domain here is the periodic box $[-0.5,0.5]^2$, and we use a finite difference method. Our boundary point spacing is given by $\Delta s \approx0.75 \Delta x$, and we replace solution values within $m_1=6$ meshwidths of the boundary using interior interpolation points located $m_2=8$ meshwidths from the boundary. For this problem, our method only requires $4-6$ iterations of \texttt{gmres}. Figure \ref{Neumann plots} illustrates the refinement studies for the solution values and $U_b$, and we see approximately first-order convergence in both of these.

   \chapter{The Immersed Boundary Double Layer (IBDL) method: Brinkman, Stokes and Navier-Stokes equations}
   \label{chapter 6}
   
In this chapter, we present the Immersed Boundary Double Layer method for problems governing fluid flow, including Stokes, Brinkman, and Navier-Stokes equations. As seen in Section \ref{ch 2 stokes}, the IBSL method can be adapted to the Stokes equation quite simply. The term that introduces the required boundary density, $S\F$, remains unchanged. However, adapting the IBDL method to Stokes equation is nontrivial, and we must use the form of the Green's function stress tensor, $\T$, to establish the required form.

This chapter is organized as follows. In Section \ref{ch 6 IBDL}, we introduce the method and demonstrate the connection to a regularized double layer integral equation for Stokes flow. We also modify the method for an exterior domain in order to represent flows with net force and torque. In Section \ref{ch 6 neumann}, we present the method as applied to a fluid equation with Neumann boundary conditions. In Section \ref{ch 6 numerical}, we discuss the numerical implementation, focusing on aspects that differ from the previous IBDL method for scalar elliptic problems. In Section \ref{ch 6 results}, we present numerical results for applications to Stokes and Brinkman equations, and in Section \ref{ch 6 Navier Stokes}, we present numerical results for the Navier-Stokes equation.


\section{Mathematical description of the method}\label{ch 6 IBDL}
We now look at the boundary value problem given by 
\begin{subequations} \label{stokes pde}
\begin{alignat}{2}
& \L \u  -\grad p= \g \qquad && \text{in } \Omega  \label{s pde 1}\\
& \grad \dotp \u = 0 \qquad && \text{in } \Omega  \label{s pde 2}\\
&\u=\U_b \qquad && \text{on } \Gamma,   \label{s pde 3}
\end{alignat}
\end{subequations}
where $\L=\mu \Delta - k^2$, $\u$ is the two-dimensional velocity, and $p$ is the pressure. We again assume that $\Omega$ is embedded in a larger computational domain $\C$ and that the boundary of $\Omega$ is given by $\Gamma$, and it is a smooth one-dimensional curve parametrized by $\X(s)$. 

\subsection{Formulation and motivation}\label{6.1 equation}

The Immersed Boundary Double Layer formulation for Equation \eqref{stokes pde} is given by 
\begin{subequations} \label{ibdl stokes}
\begin{alignat}{2}
& \L \u -\grad p +\mu\Div (S \A) = \gt \qquad && \text{in } \mathcal{C}  \label{ibdl stokes 1}\\
& \Div \u+S(\Q\dotp\n)= 0\qquad && \text{in } \mathcal{C} \label{ibdl stokes 2}\\
&S^* \u + \frac12 \Q = \U_b \qquad && \text{on } \Gamma , \label{ibdl stokes 3S}
\end{alignat}
\end{subequations}
where we define $\A$ as 
\begin{equation}
\A\equiv  \Q\otimes \n + \n \otimes\Q,
\end{equation}
or in Einstein notation,
\begin{equation}
A_{ij}=Q_in_j+Q_jn_i.
\end{equation}
We again assume an arclength parametrization of $\Gamma$, and $\n$ is the unit normal pointing out of $\Omega$. Here, $\Q$ is a vector-valued distribution supported on $\Gamma$, and it plays the role of the Lagrange multiplier enforcing the boundary condition. The term given by $\Div(S \A)$ can be seen as a generalization of the term $\widetilde S Q \equiv \grad\dotp (SQ\n)$ in the IBDL method for the Poisson and Helmholtz equations. Also note the presence of an additional term in Equation \eqref{ibdl stokes 2} given by $ S (\Q\dotp \n)$. Lastly, recall that $\gt$ represents an extension of the function $\g$ from $\Omega$ to the computational domain $\C$, and we define it as $\gt \equiv \g\chi_{\scaleto{\Omega}{4.5pt}} + \g_e \chi_{\scaleto{\C \setminus \Omega}{6pt}}$. 

We now present motivation for the form of this system of equations, and in Section \ref{6.1 connection}, we explicitly verify the connection between the IBDL method and a regularized double layer integral equation. To motivate the structure of this system, we look again to the double layer integral representation presented in Section \ref{3.5 DL}, given by 
\begin{equation}
u_j(\x_0)=\int_{\Gamma} \Psi_i(\x) T_{ijk}(\x,\x_0)n_k(\x) dl(\x). \label{T eqn}
\end{equation}
If we decompose $\T$ into its constituents, we get
\begin{multline}
u_j(\x_0)=\int_{\Gamma} \Psi_i(\x) \Big(\frac{\partial G_{ij}}{\partial \x_k}(\x,\x_0)+\frac{\partial G_{kj}}{\partial \x_i}(\x,\x_0)\Big)n_k(\x) dl(\x)\\-\int_{\Gamma} \Psi_i(\x) \delta_{ik} p_j(\x,\x_0) n_k(\x) dl(\x) , 
\end{multline}
where $\G$ is the Stokes Green's function and $\p$ is the corresponding pressure vector. This can also be rewritten as 
\begin{multline}
u_j(\x_0)=\int_{\Gamma} \frac{\partial G_{ij}}{\partial \x_k}(\x,\x_0)\Big( \Psi_i(\x)n_k(\x)+\Psi_k(\x)n_i(\x) \Big) dl(\x)\\-\int_{\Gamma} p_j(\x,\x_0)  \Psi_k n_k(\x) dl(\x) . \label{T eqn 3}
\end{multline}
The first term represents a symmetric distribution of point force dipoles, and this corresponds to the IBDL term in Equation \eqref{ibdl stokes 1}. Additionally, as in Chapter \ref{chapter 5}, we can associate the potential strength $\Psi$ with the IBDL distribution $\Q$. By doing this, we can see that the strength of the point force dipoles matches our definition of $\A$. The second term corresponds to a distribution of point sources, and this explains the need for the extra source term in Equation \eqref{ibdl stokes 2}, and again the strength of this source matches our use of $\Q\dotp \n$. 

\subsection{Explicit connection between Stokes IBDL and boundary integral equation}\label{6.1 connection}

We will now demonstrate the connection by focusing on the homogeneous Immersed Boundary Double Layer system for Stokes equation, and for simplicity, we assume the invertibility of $\Delta$. The homogeneous IBDL system for Stokes is 
\begin{subequations} \label{ibdl stokes again}
\begin{alignat}{2}
&\mu  \Delta \u -\grad p +\mu\Div (S \A) = 0 \qquad && \text{in } \mathcal{C}  \label{ibdl stokes 1 again}\\
& \Div \u+S(\Q\dotp\n)= 0\qquad && \text{in } \mathcal{C} \label{ibdl stokes 2 again }\\
&S^* \u + \frac{1}{2} \Q = \U_b \qquad && \text{on } \Gamma .\label{ibdl stokes 3S again}
\end{alignat}
\end{subequations}
Taking the divergence of Equation \eqref{ibdl stokes 1 again} and using Equation \eqref{ibdl stokes 2 again }, we get 
\begin{equation}
\Delta p =  \mu \Div \Div (S\A)+ \mu \Delta (\Div \u)  = \mu  \Div \Div (S\A) - \mu\Delta (S(\Q\dotp\n)).
\end{equation}
Using the definition of the spread operator, this becomes 
\begin{equation}
\Delta p =\mu \Div\Div \int_{\Gamma} \A (s) \delta_h(\x-\X(s)) ds - \mu \Delta  \int_{\Gamma} \Q(s) \dotp\n(s)\delta_h(\x-\X(s)) ds .
\end{equation}
We can then invert the operator and use the property of the regularized Green's function for Laplace's equation, $\Delta G_h^L(\x,\x_0)= -\delta_h(\x-\x_0)$. We then have
\begin{equation}
 p(\x) = -\mu \Div\Div \int_{\Gamma} \A(s) G_h^L(\x,\X(s)) ds + \mu \Delta  \int_{\Gamma} \Q(s) \dotp\n(s)G_h^L(\x,\X(s)) ds .
\end{equation}
For clarity, we will now switch to Einstein notation. Bringing the operators into the integrals and manipulating the expressions, we get 
\begin{equation}
 p(\x) = -\mu \int_{\Gamma} \frac{\partial^2 G_h^L}{\partial x_i\partial x_k}(\x,\X(s))A_{ik}(s) ds + \mu \int_{\Gamma} \frac{\partial^2 G_h^L}{\partial x_m\partial x_m}(\x,\X(s)) Q_k(s) n_k (s) ds .\label{pressure ibdl}
\end{equation}
Using this expression for pressure in Equation \eqref{ibdl stokes 1 again} and canceling the $\mu$ on all terms, we get
\begin{multline}
 \Delta u_j = -\int_{\Gamma} \frac{\partial^3 G_h^L}{\partial x_i\partial x_j\partial x_k}(\x,\X(s))A_{ik}(s) ds + \int_{\Gamma} \frac{\partial^3 G_h^L}{\partial x_m\partial x_m\partial x_j}(\x,\X(s)) Q_k(s) n_k(s) ds \\ -\frac{\partial}{\partial x_k}\int_{\Gamma} A_{jk}(s)\delta_h(\x-\X(s)) ds .
\end{multline}
We next invert the operator. For the first term, we use that $\Delta H_h=G_h^L$ from Equation \eqref{G and H relationship}. For the second term, we simply cancel the $\Delta$ in the integral, and for the third term, we use $\Delta^2H_h=-\delta_h$ from Equation \eqref{Hh2}. We then get 
\begin{multline}
 u_j(\x) = -\int_{\Gamma} \frac{\partial^3 H_h}{\partial x_i\partial x_j\partial x_k}(\x,\X(s))A_{ik}(s) ds + \int_{\Gamma} \frac{\partial G_h^L}{\partial x_j}(\x,\X(s)) Q_k(s) n_k(s) ds  \\+\int_{\Gamma} \frac{\partial^3 H_h}{\partial x_m\partial x_m\partial x_k}(\x,\X(s))A_{jk}(s) ds .
\end{multline}
Switching the last two integrals and rearranging, we get
\begin{multline}
 u_j(\x) = \int_{\Gamma} \frac{\partial}{\partial x_k}\Bigg(-\frac{\partial^2 H_h}{\partial x_i\partial x_j}(\x,\X(s))A_{ik}(s)  + \frac{\partial^2 H_h}{\partial x_m\partial x_m}(\x,\X(s))A_{jk}(s) \Bigg)ds  \\+ \int_{\Gamma} \frac{\partial G_h^L}{\partial x_j}(\x,\X(s)) Q_k(s) n_k(s) ds.
\end{multline}
We can also write this as
\begin{multline}
  u_j(\x) = \int_{\Gamma} \frac{\partial}{\partial x_k}\Bigg(-\frac{\partial^2 H_h}{\partial x_i\partial x_j}(\x,\X(s))  + \delta_{ij}\frac{\partial^2 H_h}{\partial x_m\partial x_m}(\x,\X(s)) \Bigg)A_{ik}(s)ds  \\+ \int_{\Gamma} \frac{\partial G_h^L}{\partial x_j}(\x,\X(s)) Q_k(s) n_k(s) ds.
\end{multline}
Recognizing the form of the regularized Stokes Green's function from Equation \eqref{stokes g to other ones h}, we have 
\begin{equation}
 u_j(\x) = \int_{\Gamma} \frac{\partial G_{ij}^h}{\partial x_k}(\x,\X(s)) A_{ik}(s)ds  + \int_{\Gamma} \frac{\partial G_h^L}{\partial x_j}(\x,\X(s)) Q_k(s) n_k(s) ds.
\end{equation}
Lastly, we can use our definition of $\A$ to write this as 
\begin{equation}
u_j(\x) = \int_{\Gamma} \frac{\partial G_{ij}^h}{\partial x_k}(\x,\X(s)) \Big(Q_i(s)n_k(s)+Q_k(s)n_i(s)\Big)ds  + \int_{\Gamma} \frac{\partial G_h^L}{\partial x_j}(\x,\X(s)) Q_k(s) n_k(s) ds.
\end{equation}
Using the the identity for the pressure vector given in Equation \eqref{p to GL}, this is 
\begin{equation}
  u_j(\x) = \int_{\Gamma} \frac{\partial G_{ij}^h}{\partial x_k}(\x,\X(s)) \Big(Q_i(s)n_k(s)+Q_k(s)n_i(s)\Big)ds  - \int_{\Gamma}p_j^h(\x,\X(s)) Q_k(s) n_k(s) ds.
\end{equation}
Then, we recognize that by linearity, the decomposition of the stress tensor used in Equation \eqref{T eqn} also holds using the regularized functions. Then by reversing the process we took to decompose $\T$ in Equations \eqref{T eqn}-\eqref{T eqn 3}, we get 
\begin{equation}
u_j(\x)=\int_{\Gamma} Q_i(s) T_{ijk}^h(\x, \X(s)) n_k(s)ds, \label{almost there 2}
\end{equation}
for $\x\in\Omega\setminus \Gamma$. Next, we use Equation \eqref{ibdl stokes 3S again} to get 
\begin{equation}
\U_b(s') =\int_{\C} \u(\x)\delta_h(\x-\X(s')) d\x +\frac{1}{2} Q_j(s'). 
\end{equation}
We then follow the same steps seen in Section \ref{4.1 connect}. We combine this with Equation \eqref{almost there 2}, change the order of integration, recognize the presence of $T_{ijk}^h * \delta_h$, and denote the twice-regularized stress tensor by $T_{ijk}^{hh}$. This leaves us with 
\begin{equation}
U_{b_j}(s') =\int_{\Gamma} Q_i(s) T_{ijk}^{hh}(\X(s'), \X(s)) n_k(s)ds +\frac{1}{2} Q_j(s')  . 
\end{equation}
The odd symmetry of the stress tensor, which is preserved through convolutions with the regularized delta function, gives us that $T_{ijk}^{hh}(\X(s'),\X(s)) = -T_{ijk}^{hh}(\X(s),\X(s'))$ \cite{Pozred}. Using this to switch the arguments of $T_{ijk}$ and appropriately redefining $\Q$ and $\U_b$ as functions of $\x$, we get 
\begin{equation}
U_{b_j}(\X(s')) =-\int_{\Gamma} Q_i(\X(s)) T_{ijk}^{hh}(\X(s), \X(s')) n_k(\X(s))ds +\frac{1}{2} Q_j(\X(s')). \label{IBDL stokes one} 
\end{equation}

Recall that the double layer integral equation for $\X_0$ on the boundary from Section \ref{ch 3 stokes integ eqns} is given by 
\begin{equation}
U_{b_j} (\X_0) = \int_{\Gamma}^{PV} \Psi_i(\x) T_{ijk}(\x,\X_0)n_k(\x) dl(\x)   - \frac12 \Psi_j(\X_0). \label{stokes double layer eqn}
\end{equation}
Using our arclength parametrization, we can write this as 
\begin{equation}
U_{b_j} (\X(s')) = \int_{\Gamma}^{PV} \Psi_i(\X(s)) T_{ijk}(\X(s),\X(s'))n_k(\X(s)) ds   - \frac12 \Psi_j(\X(s')). \label{stokes double layer eqn 2}
\end{equation}
By comparing Equations \eqref{IBDL stokes one} and \eqref{stokes double layer eqn 2}, we can see that in the limit that the regularization width, $h$, approaches $0$, we have
\begin{equation}
\Q(\X)=-\bm{\Psi}(\X). \label{ stokesF eqn}
\end{equation}
We then see from Section \ref{3.5 DL} that $\Q$ gives the jump in solution values across the boundary. We have now established that the IBDL method for Stokes is equivalent to a \emph{regularized} double layer integral equation.

\subsection{Exterior domains}\label{6.1 Hebeker}

As discussed in Section \ref{3.5 DL}, the double layer potential is only capable of representing a flow with no net force or net torque. Therefore, for exterior domains, we will use the completed IBDL method, which was presented for the Poisson case in Section \ref{5.3 poisson}. Specifically, we add a single layer potential whose strength is given by a constant multiple of the double layer potential strength. For Stokes and Brinkman equations, the formulation is 
\begin{subequations} \label{ibdl stokes hebeker}
\begin{alignat}{2}
& \L \u -\grad p +\eta S \Q + \mu \Div (S \A) = \gt \qquad && \text{in } \mathcal{C}  \label{ibdl stokes hebeker 1}\\
& \Div \u+S(\Q\dotp\n)= 0\qquad && \text{in } \mathcal{C} \label{ibdl stokes hebeker 2}\\
&S^* \u + \frac{1}{2} \Q = \U_b \qquad && \text{on } \Gamma \label{ibdl stokes hebeker 3}
\end{alignat}
\end{subequations}
where $\L = \mu\Delta -k^2$ and $\A$ is defined by 
\begin{equation}
A_{ij}=Q_in_j+Q_jn_i.
\end{equation}
Additionally, $\eta$ is a positive constant, and in Section \ref{6.5 eta}, we explore the affects of different choices for $\eta$. This supplementary flow is only absolutely required in the case of the Stokes equation, or $k=0$. However, for the Brinkman equation, as $k$ approaches $0$, including this term becomes important for maintaining low errors and first-order convergence. However, as is discussed further in Section \ref{6.5 eta}, including the term can lower errors even for $k$ away from $0$.

\textbf{Force and Torque. } We now derive the equations for net force and net torque corresponding to the completed IBDL method. If we start from Equation \eqref{ibdl stokes hebeker} and use the derivation in the previous section, the expression for pressure given in Equation \eqref{pressure ibdl} would have an additional term resulting from the added single layer potential. We would then have
\begin{multline}
 p(\x) = -\mu \int_{\Gamma} \frac{\partial^2 G_h^L}{\partial x_i\partial x_k}(\x,\X(s))A_{ik}(s) ds + \mu \int_{\Gamma} \frac{\partial^2 G_h^L}{\partial x_m\partial x_m}(\x,\X(s)) Q_k(s) n_k (s) ds \\
 - \eta\int_{\Gamma}\frac{\partial G_h^L}{\partial x_i} (\x,\X(s)) Q_i(s)ds.
\end{multline}
Using this expression for pressure in Equation \eqref{ibdl stokes hebeker 1} and dividing by $\mu$, we get
\begin{alignat}{2}
 &&\Delta u_j =& -\int_{\Gamma} \frac{\partial^3 G_h^L}{\partial x_i\partial x_j\partial x_k}(\x,\X(s))A_{ik}(s) ds + \int_{\Gamma} \frac{\partial^3 G_h^L}{\partial x_m\partial x_m\partial x_j}(\x,\X(s)) Q_k(s) n_k(s) ds \nonumber \\  
 && &-\frac{\partial}{\partial x_k}\int_{\Gamma} A_{jk}(s)\delta_h(\x-\X(s)) ds \nonumber\\
 && &- \frac{\eta}{\mu} \int_{\Gamma}\frac{\partial^2 G_h^L}{\partial x_i \partial x_j} (\x,\X(s)) Q_i(s)ds   - \frac{\eta}{\mu}\int_{\Gamma} Q_j(s)\delta_h(\x,\X(s))ds,
\end{alignat}
where the last two integrals are the contribution from the single layer term. We next invert the operator and use the same identities used in the previous section for the first three terms. For the last two terms, we use the identities used in Section \ref{4.2 connect} to go from Equation \eqref{sl part 1} to Equation \eqref{almost there}. Then we get
\begin{alignat}{2}
&& u_j(\x) =& -\int_{\Gamma} \frac{\partial^3 H_h}{\partial x_i\partial x_j\partial x_k}(\x,\X(s))A_{ik}(s) ds + \int_{\Gamma} \frac{\partial G_h^L}{\partial x_j}(\x,\X(s)) Q_k(s) n_k(s) ds  \nonumber \\
&& &+\int_{\Gamma} \frac{\partial^3 H_h}{\partial x_m\partial x_m\partial x_k}(\x,\X(s))A_{jk}(s) ds \nonumber \\
&& &- \frac{\eta}{\mu}\int_{\Gamma}  \frac{\partial^2 H_h}{\partial x_i\partial x_j }(\x,\X(s)) Q_i(s) ds +\frac{\eta}{\mu}\int_{\Gamma} Q_j(s)\frac{\partial^2 H_h}{\partial x_k \partial x_k}(\x,\X(s)) ds .
\end{alignat}
From the previous section, we know that the first three integrals simplify to the double layer potential. Then, from Section \ref{4.2 connect}, we know that the last two integrals simplify to a single layer potential. We therefore have 
\begin{equation}
u_j(\x)=\int_{\Gamma} Q_i(s) T_{ijk}^h(\x, \X(s)) n_k(s)ds+\frac{\eta}{\mu}\int_{\Gamma} G_{ij}^h(\x,\X(s))Q_i(s)ds. 
\end{equation}
Specifically, we see that with the formulation in Equation \eqref{ibdl stokes hebeker}, we retrieve the integral representation from Equation \eqref{stokes completed double layer rep}, with the correct constant on the second term. Using that $\Q=-\bm{\Psi}$, we therefore know from Section \ref{3.5 DL} that the net force and net torque through a curve enclosing $\Gamma$ is given by 
\begin{subequations}
\begin{alignat}{2}
&\bm{B}=-\int_{ \Gamma}\bm{\Q}(s) ds\\
&L=- \int_{ \Gamma} \X(s) \times \bm{\Q} (s) ds,
\end{alignat}
\end{subequations}
and we will use this for calculating net force in Section \ref{ch 6 results}.

\section{Neumann boundary conditions}\label{ch 6 neumann}

We now look at the PDE with Neumann boundary conditions given by 
\begin{subequations} 
\begin{alignat}{2}
& \L \u  -\grad p= \g \qquad && \text{in } \Omega  \\
& \grad \dotp \u = 0 \qquad && \text{in } \Omega \\
&\sig\dotp \n=\bm{F}_b \qquad && \text{on } \Gamma, 
\end{alignat}
\end{subequations}
where $\L=\mu \Delta - k^2$. Here, we present the IBDL formulation for this PDE. For brevity, we omit the explicit connection to boundary integral equations. The steps follow closely with those from Sections \ref{ch 5 neumann}, \ref{6.1 connection}, and \ref{6.1 Hebeker}. The formulation is 
\begin{subequations} \label{ibdl stokes neumann}
\begin{alignat}{2}
& \L \u -\grad p +S(\bm{F}_b)+\mu\Div (S \A) = \gt \qquad && \text{in } \mathcal{C}  \\
& \Div \u+S(\U_b\dotp\n)= 0\qquad && \text{in } \mathcal{C} \\
&S^* \u  =\frac12 \U_b \qquad && \text{on } \Gamma , 
\end{alignat}
\end{subequations}
where $\A$ is defined by 
\begin{equation}
A_{ij}=U_{b_i}n_j+U_{b_j}n_i.
\end{equation}
Here, the unknown potential strength that we solve for is the distribution of boundary values, $\U_b$.

\section{Numerical implementation}\label{ch 6 numerical}

There is no change in the discretization of the spread and interpolation operators discussed in Section \ref{ch 5 numerical}. We also continue to use the same linear interpolation for solution values near the grid. In this section, we discuss the discretization of space and then the methods for solving the discrete Brinkman and Stokes equations. 

\subsection{Discretization of differential operators }\label{6.4 space}

As discussed in Section \ref{2.3 space and op}, taking a discrete divergence of a discrete gradient results in a wide Laplacian. In the case of the IBSL method, we have the option to avoid this issue by using a staggered grid, in which vector quantities are defined on the centers of the cell edges and scalar quantities are defined in the centers of the cells. In the IBSL method, this requires that the spread action, $S\F$, consists of spreading the $F_1$ components to vertical edge centers and the $F_2$ components to horizontal edge centers. The discrete divergence of $S\F$, which is used in the solution of pressure, then maps the $\Div ( S\F)$ quantities to the cell centers. In the IBDL method, however, one must spread a \textit{matrix}, $\A$. Using the definition of $\A$, we could also say we need to spread the vectors $\Q$ and $\n$, but we also need to combine elements on different edges, for instance, when we calculate $Q_1n_2$. Therefore, using a staggered grid for the IBDL method requires a nontrivial adaptation, and this is an area for future research. 

In the IBSL method, we can also use the wide Laplacian directly to solve for pressure, but in the case of the IBDL method, we fail to get first-order convergence with this approach. We conjecture that the interaction of two different differential operators results in inconsistent numerical Green's functions, and in the case of the IBDL method, this may lead to larger numerical errors due to the high frequency oscillations that would result from taking the derivative of the regularized delta function. This is again another area for future research. We illustrate this issue in Section \ref{6.5 fd}, and we demonstrate that we are able to recover first order convergence by approximating the wide Laplacian with the usual 5-point second-order Laplacian and utilizing a smoother delta function. When we use finite differences for the IBDL method for Stokes, we therefore use the usual 5-point Laplacian and the 6-point B-spline delta function utilized in \cite{deltapaper}. It is given by 
\begin{equation}
\delta_h=\frac{1}{h^2} \phi\bigg(\frac{x}{h}\bigg)\phi\bigg(\frac{y}{h}\bigg),
\end{equation}
for $h=\Delta x = \Delta y$ and 
\begin{equation}
\phi (r)=\begin{cases}
\frac{11}{20} - \frac12 r^2 +\frac14 r^4 -\frac{1}{12} r^5 & 0\leq |r|\leq 1 \\ 
 \frac{17}{40} +\frac58r-\frac74r^2+\frac54r^3-\frac38r^4+\frac{1}{24}r^5 & 1\leq |r|\leq 2 \\
  \frac{81}{40} -\frac{27}{8} r+\frac94r^2-\frac34 r^3 +\frac18r^4-\frac{1}{120}r^5 & 2\leq |r|\leq 3 \\
	0 & |r|\geq 3  . \end{cases}
\end{equation}
However, we primarily utilize a Fourier spectral discretization for the IBDL method applied to Stokes, and this discretization is described in Section \ref{2.3 space and op}.

\subsection{Solution to discrete system for Brinkman equation  }\label{6.4 Brinkman}

Recall that the IBDL formulation of the PDE in Equation \eqref{stokes pde} is given by
\begin{subequations} \label{ibdl stokes again again}
\begin{alignat}{2}
& \L \u -\grad p +\eta S\Q + \mu \Div (S \A) = \gt \qquad && \text{in } \mathcal{C}  \label{ibdl stokes 1 again again}\\
& \Div \u+S(\Q\dotp\n)= 0\qquad && \text{in } \mathcal{C} \label{ibdl stokes 2 again again}\\
&S^* \u + \frac{1}{2} \Q = \U_b \qquad && \text{on } \Gamma , \label{ibdl stokes 3S again again}
\end{alignat}
\end{subequations}
where $\eta \geq 0$, and for this section we focus on $k>0$. We may take $\eta=0$ in the case of an interior domain or large $k$. 

Taking the divergence of Equation \eqref{ibdl stokes 1 again again} and using Equation \eqref{ibdl stokes 2 again again}, we get
\begin{equation}
\Delta p=-\L S(\Q\dotp\n)+\eta \Div (S\Q)+\mu \Div\Div(S\A)-\Div\gt . \label{pressure from dl}
\end{equation}
Equation \eqref{pressure from dl} is solvable on the periodic domain, and we invert the periodic Laplacian with the operator $\Delta_0^{-1}$, which was described in Section \ref{2.2 Laplacian L}. We therefore fix its mean value on $\C$ to be $0$. Completing this inversion, we get
\begin{equation}
p=-\Delta_0^{-1}\L  S(\Q\dotp\n)  - \Delta_0^{-1}\ \Div \Big(\gt-\eta S \Q-\mu\Div (S\A)\Big). \label{pressure again dl}
\end{equation}
Let us again denote the operator that projects onto divergence-free fields as $\mathbb{P}=\mathds{I} - \grad \Delta_0^{-1} \Div$. Using this operator and the expression for $p$ from Equation \eqref{pressure again dl},  Equation \eqref{ibdl stokes 1 again again} gives us
\begin{equation}
\L \u = -\grad \Delta_0^{-1} \L  S(\Q\dotp\n) + \mathbb{P} \Big(\gt - \eta S\Q-\mu \Div (S\A)\Big) .\label{onC}
\end{equation}
Inverting $\L$ and applying the interpolation operator $S^*$, we get the following equation for $\Q$. 
\begin{equation}\label{brink ibdl eqn}
-S^*\grad \Delta_0^{-1}S(\Q\dotp\n) - S^*\L^{-1}\mathbb{P}\Big(\eta S\Q +\mu \Div (S \A)\Big) +\frac12\Q= \U_b - S^*\L^{-1} \mathbb{P} \gt,
\end{equation}
where we recall that $\A$ is a function of $\Q$, given by $\A\equiv \Q \otimes \n + \n \otimes \Q$. As we did for the elliptic PDE, we solve this equation for $\Q$ using \texttt{gmres} and then solve Equations \eqref{ibdl stokes 2 again again} and \eqref{ibdl stokes 2 again again} for the velocity $\u$ and pressure $p$.

\subsection{Solution to discrete system for Stokes equation  }\label{6.4 Stokes}
We now look at the IBDL formulation for the Stokes equation, given by 
\begin{subequations} \label{ibdl stokes again 2}
\begin{alignat}{2}
&\mu \Delta \u -\grad p +\eta S\Q + \mu\Div (S \A) = \gt \qquad && \text{in } \mathcal{C}  \label{ibdl stokes 1 again 2}\\
& \Div \u+S(\Q\dotp\n)= 0\qquad && \text{in } \mathcal{C} \label{ibdl stokes 2 again 2}\\
&S^* \u + \frac{1}{2} \Q = \U_b \qquad && \text{on } \Gamma , \label{ibdl stokes 3S 2}
\end{alignat}
\end{subequations}
where $\eta> 0$ for an exterior PDE domain $\Omega$, and $\eta=0$ for an interior $\Omega$. In this case, we must again adjust the method of solution since the periodic Laplacian is not invertible. Additionally, note that as discussed for the Poisson equation in Section \ref{5.3 poisson}, for the case of $\eta=0$, the constraint matrix has a one-dimensional nullspace, but since the eigenfunction spanning the nullspace again gives values of $0$ on the physical domain, we need not alter the method for an interior domain. For an exterior domain, the addition of the single layer term results in a matrix that is invertible. Therefore, we need only to handle the nullspace of the periodic Laplacian. 

We begin by decomposing the solution $\u$ as 
\begin{equation}
\u=\u_0+\ub,
\end{equation}
where $\u_0$ has mean $0$ on $\C$ and $\ub=(\bar u, \bar v)$ gives the mean values of $u$ and $v$, respectively.
Following the same process as the previous section and noting that $\Delta \ub=0$, we get that the equation that must be solved for $\u_0$ on $\C$ is given by 
\begin{equation}
\mu \Delta \u_0 =-\mu\grad \Delta_0^{-1} \Delta  S(\Q\dotp\n) +  \mathbb{P} \Big(\gt - \eta S\Q-\mu \Div (S\A)\Big).\label{onCagain}
\end{equation}
By integrating this over a general computational domain $\C$, we get
\begin{multline}
\mu \int_{\C} \Delta \u_0 d\x = -\mu\grad \Delta_0^{-1} \Delta\int_{\C}  S(\Q\dotp\n) d\x -\grad \Delta_0^{-1}\int_{\C} \Div \Big(\gt-\eta S \Q-\mu \Div (S\A)\Big)d\x  \\-\int_{\C}\Big( \eta S\Q+\mu \Div (S\A)-\gt\Big)d \x.
\end{multline}
Using the divergence theorem where appropriate, we get 
\begin{multline}
\mu \int_{\partial \C} \grad \u_0 \dotp \n_{\C} dl(\x) = - \mu  \Delta_0^{-1} \Delta \int_{\partial \C} S(\Q\dotp\n)\n_{\C} dl(\x) -  \grad \Delta_0^{-1} \int_{\partial \C}\gt \dotp \n_{\C}  dl(\x)  \\
\hspace{2cm}+  \grad \Delta_0^{-1}\int_{\partial \C} \Big( \eta S\Q   + \mu S\A\Big) \dotp\n_{\C}dl(\x) 
- \eta \int_{\C} S\Q d\x -\mu  \int_{\partial \C}  (S\A)\dotp\n_{\C}dl(\x)+ \int_{\C} \gt d\x\label{big integral}
\end{multline}
If we assume $\Gamma$ is away from the boundary of $\C$, the first, third, and fifth integrals on the right-hand-side vanish.  If we take $\C$ to be a periodic box, the left-hand-side also vanishes, leaving us with 
\begin{equation}
0 = -  \grad \Delta_0^{-1} \int_{\partial \C}\gt \dotp \n_{\C}  dl(\x)    - \eta \int_{\C} S\Q d\x+ \int_{\C} \gt d\x\label{big integral 2}
\end{equation}

\textbf{Interior PDE domain.  }  Let us now first consider the case in which the PDE domain $\Omega$ is an interior domain. In this case, $\eta=0$, and the second integral in Equation \eqref{big integral 2} vanishes. We are left with 
\begin{equation}
 \grad \Delta_0^{-1} \int_{\partial \C} \gt \dotp \n_{\C}  dl(\x)  = \int_{\C} \gt d\x.
\end{equation}
If we make the same choice for $g_e$ as we did for the Poisson equation, 
\begin{equation}
g_e =-\frac{1}{|\C\setminus \Omega|} \int_{\Omega} gd\x, 
\end{equation} 
then the right-hand-side is $0$, and since on $\partial \C$, $\gt$ is a constant, the periodicity of $\C$ gives us $0$ on the left-hand-side. With this solvability constraint satisfied, the solution to Equation \eqref{laplaces with ib} can be found, and it is unique up to an additive constant. We therefore find the solution with mean $0$ on $\C$, and the potential strength, $\Q$, enforces the boundary condition, giving the correct solution on the PDE domain $\Omega$. 
Then by inverting $\Delta$ in Equation \eqref{onCagain} and applying the interpolation operator $S^*$, we get the following equation for $\Q$. 
\begin{equation}
-S^*\Delta_0^{-1}\grad \Delta_0^{-1}\Delta S(\Q\dotp\n) - S^*\Delta_0^{-1}\mathbb{P}\Div (S \A)+\frac12\Q= \U_b - \frac{1}{\mu}S^*\Delta_0^{-1} \mathbb{P} \gt. 
\end{equation}
After solving this equation using \texttt{gmres}, we then solve Equations \eqref{ibdl stokes 1 again} and \eqref{ibdl stokes 2 again } for $p$ and $\u$. 

\textbf{Exterior PDE domain.  }  Now let us consider the case in which $\Omega$ is an exterior domain. Returning to Equation \eqref{big integral 2}, we then have that $\gt$ on $\partial \C$ is $\g$, the periodic forcing function in the PDE. Therefore, the first integral vanishes. Using that $\int_{\C} S\Q d\x=\int_{\Gamma} \Q ds$, we then get the solvability constraint
\begin{equation}
\eta \int_{\Gamma} \Q ds = \int_{\C} \gt d\x.
\end{equation}
We therefore use this constraint to solve for the unknown $\ub$, as we did for the IBSL method in Sections \ref{2.2 Laplacian L} and \ref{2.3 stokes}. Discretizing this constraint, we get the following system of equations
\begin{subequations}
\begin{alignat}{1}
& \mu \Delta \u_0 -\grad p +\eta S\Q + \mu \Div (S \A) = \gt \qquad  \\
& \Div \u_0+S(\Q\dotp\n)= 0\\
&S^* \u_0+\ub \mathds{1}_{N_{IB}}+\frac{1}{2}\Q = \U_b  \\
&\eta (\Delta s) \mathds{1}_{N_{IB}}^\intercal \Q = (\Delta x \Delta y)  \mathds{1}_{N^2}^\intercal \gt
\end{alignat}
\end{subequations}
We therefore have the following discrete system of equations, 
\begin{subequations}\label{discrete operator for Q}
\begin{multline}
-S^*\Delta_0^{-1}\grad \Delta_0^{-1}\Delta S(\Q\dotp\n) - S^*\Delta_0^{-1}\mathbb{P}\Big(\frac{\eta}{\mu} S\Q + \Div (S \A)\Big) +\ub\mathds{1}_{N_{IB}}  +\frac{1}{2}\Q\\= \U_b -\frac{1}{\mu} S*\Delta_0^{-1} \mathbb{P} \gt
\end{multline}
\begin{equation}
\eta (\Delta s) \mathds{1}_{N_{IB}}^\intercal \Q = (\Delta x\Delta y)  \mathds{1}_{N^2}^\intercal \gt, 
\end{equation}
\end{subequations}
where the unknowns are $\Q$ and $\ub$. We solve this equation using \texttt{gmres}, and then solve Equations \eqref {ibdl stokes 1 again again} and \eqref{ibdl stokes 2 again again} for $p$ and $\u_0$. Finally, we find $\u$ by adding $\ub$.

\section{Results: Stokes and Brinkman equations}\label{ch 6 results}

In this section, we apply the Immersed Boundary Double Layer method to Stokes and Brinkman equations with Dirichlet boundary conditions and discuss numerical results. In Section \ref{6.5 stokes 1}, we revisit the Brinkman equation seen in Section \ref{2.4 Stokes} and demonstrate the improved efficiency of the IBDL method. In Section \ref{6.5 eta}, we examine the affect of different choices of $\eta$ in the completed IBDL formulation. In Section \ref{6.5 fd}, we observe the affects of using two finite difference forms for the Laplacian. In Section \ref{6.5 cylinder}, we apply the method to the problem of Stokes flow past a periodic array of cylinders and compare to previous numerical and asymptotic results. Lastly, in Section \ref{6.5 many}, we apply the IBDL and IBSL methods to flow past $9$ objects of varying size to demonstrate that as the domain increases in complexity, the efficiency of the IBDL method over the IBSL method can be even more extreme.

\subsection{Brinkman equation with analytical solution}\label{6.5 stokes 1}

We first revisit the Brinkman PDE from Section \ref{2.4 Stokes}, given by 
\begin{subequations} \label{brink pde again}
\begin{alignat}{2}
& \Delta \u -  \u -\grad p = \g \qquad && \text{in } \Omega  \label{brink pde 1 again}\\
& \Div \u=0 \qquad && \text{in } \Omega  \label{brink pde 2 again}\\
&\u=\U_b \qquad && \text{on } \Gamma,  \label{brink pde 3 again}
\end{alignat}
\end{subequations}
where we use the analytical solution given by
\begin{subequations} \label{brink soln again}
\begin{alignat}{2}
& u=e^{\sin{x}}\cos{y} \\
& v=-\cos{x}\hspace{0.1cm}  e^{\sin{x}}\sin{y}\\
& p= e^{\cos{y}}
\end{alignat}
\end{subequations}
to find the boundary values, $\U_b=\u|_{\Gamma}$, and the forcing function,  
\begin{equation} \label{brink g again}
\g=\begin{pmatrix}
e^{\sin{x}}\cos{y}\hspace{0.1cm} (\cos^2{x}-\sin{x}-2)\\
-\cos{x} \hspace{0.1cm} e^{\sin{x}}\sin{y}\hspace{0.1cm} (\cos^2{x}-3\sin{x}-3)+\sin{y}\hspace{0.1cm} e^{\cos{y}}\end{pmatrix}.
\end{equation}
Here $\Omega$ is the interior of a circle of radius 0.75, centered at the origin, and the computational domain is the periodic box, $\C=[-1, 1]^2$. We use a Fourier spectral method and equally spaced boundary points with $\Delta s \approx \alpha \Delta x$ for various values of $\alpha$, and the solutions are computed for grid sizes ranging from $N=2^6$ to $2^{12}$.  For the IBDL method, we use interpolation for grid points within $m_1=2(\log_2{N}-4)$ meshwidths of $\Gamma$.

\begin{figure}
\begin{center}
 \begin{tabular}{||c | c | c | c | c | c | c | c | c||} 
 \hline
 \multicolumn{9}{||c||}{Iteration Counts - Equation \eqref{brink pde again}} \\
 \hline
  &\multicolumn{2}{|c|}{$\Delta s \approx  2\Delta x $}& \multicolumn{2}{|c|}{$\Delta s \approx  1.5 \Delta x$} &\multicolumn{2}{|c|}{$ \Delta s \approx  1 \Delta x$}& \multicolumn{2}{|c||}{$\Delta s \approx 0.75 \Delta x$}  \\ 
 \hline
 $\Delta x$ &\textcolor{blue}{ IBSL} & \textcolor{cyan}{IBDL} &\textcolor{blue}{  IBSL}& \textcolor{cyan}{IBDL}&\textcolor{blue}{ IBSL} &  \textcolor{cyan}{IBDL}&\textcolor{blue}{  IBSL} &  \textcolor{cyan}{IBDL} \\
 \hline
$2^{-5}$ &  \textcolor{blue}{ 123 } & \textcolor{cyan}{10} &    \textcolor{blue}{296}&\textcolor{cyan}{10}   &  \textcolor{blue}{2558}&\textcolor{cyan}{10}  &   \textcolor{blue}{16190 stag. }&\textcolor{cyan}{10}   \\
$2^{-6} $&   \textcolor{blue}{185}& \textcolor{cyan}{10}   &   \textcolor{blue}{472}&\textcolor{cyan}{10}    &  \textcolor{blue}{5578}& \textcolor{cyan}{10} &  \textcolor{blue}{47863 stag.}&\textcolor{cyan}{10}   \\
$2^{-7 }$&   \textcolor{blue}{257}& \textcolor{cyan}{10}  &    \textcolor{blue}{678}&\textcolor{cyan}{10}   & \textcolor{blue}{14438}& \textcolor{cyan}{10}   &  \textcolor{blue}{95009 stag. }&\textcolor{cyan}{9}      \\
$2^{-8}$& \textcolor{blue}{ 351}&\textcolor{cyan}{10}    &     \textcolor{blue}{825}&\textcolor{cyan}{10}     &  \textcolor{blue}{17475}& \textcolor{cyan}{10}   &   \textcolor{blue}{68861 stag.}&\textcolor{cyan}{9}        \\
$ 2^{-9}$& \textcolor{blue}{478}&\textcolor{cyan}{10}   &     \textcolor{blue}{1137}&\textcolor{cyan}{10}     &  \textcolor{blue}{17789}& \textcolor{cyan}{10}   &   \textcolor{blue}{62697 stag.}&\textcolor{cyan}{9}    \\
$ 2^{-10}$ & \textcolor{blue}{622}&\textcolor{cyan}{10}  &   \textcolor{blue}{1620}&\textcolor{cyan}{10}     &  \textcolor{blue}{ 20333}& \textcolor{cyan}{10}  &   \textcolor{blue}{ 75910 stag.}&\textcolor{cyan}{9}    \\
$ 2^{-11}$  & \textcolor{blue}{821}&\textcolor{cyan}{10} & \textcolor{blue}{2180}&\textcolor{cyan}{10}       &  \textcolor{blue}{27467 }& \textcolor{cyan}{10}  &   \textcolor{blue}{ 79699 stag. }&\textcolor{cyan}{9}   \\
  \hline
  \end{tabular}
\captionof{table}[Number of iterations of \texttt{minres} and \texttt{gmres} to solve the Brinkman PDE in Equation \eqref{brink pde again} using the IBSL and IBDL methods]{Number of iterations of \texttt{minres} and \texttt{gmres}, with tolerance $10^{-8}$, needed to solve Equation \eqref{brink pde again} using the IBSL and IBDL methods, respectively, on the periodic computational domain $[-1, 1]^2$.  The word \textit{stag.} is used to indicate that the method stagnated at this number of iterations without converging to the desired tolerance.  }
\label{iteration table 6.1}
\end{center}
\vspace{-1cm}
\end{figure}

\begin{figure}
\begin{subfigure}{0.495\textwidth}
\centering
\includegraphics[width=\textwidth]{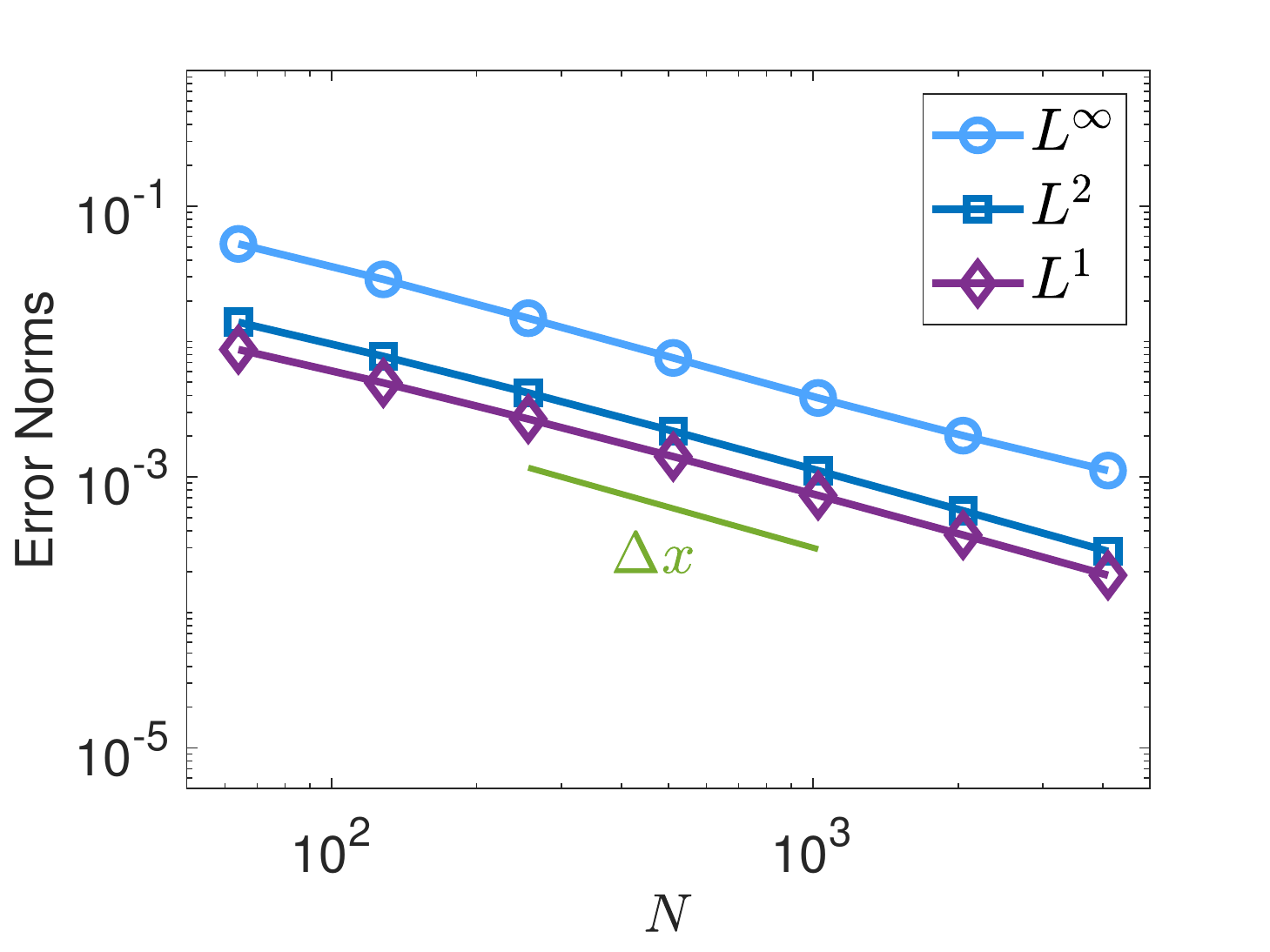}
\caption{\normalsize  IBDL, $u$}
\label{IBDL brink refinement u}
\end{subfigure}
\begin{subfigure}{0.495\textwidth}
\centering
\includegraphics[width=\textwidth]{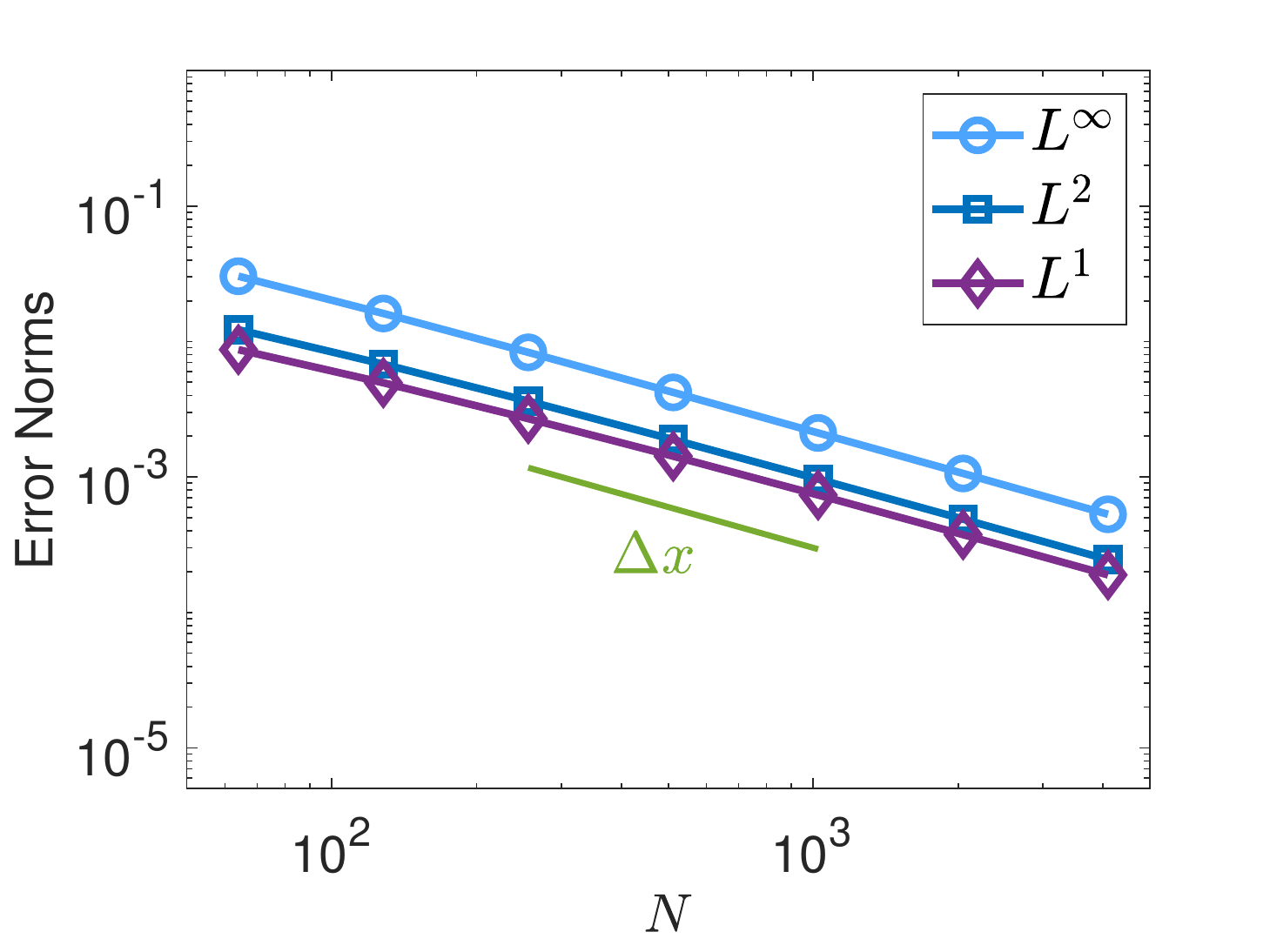}
\caption{\normalsize IBDL, $v$}
\label{IBDL brink refinement v}
\end{subfigure}
\begin{subfigure}{0.495\textwidth}
\centering
\includegraphics[width=\textwidth]{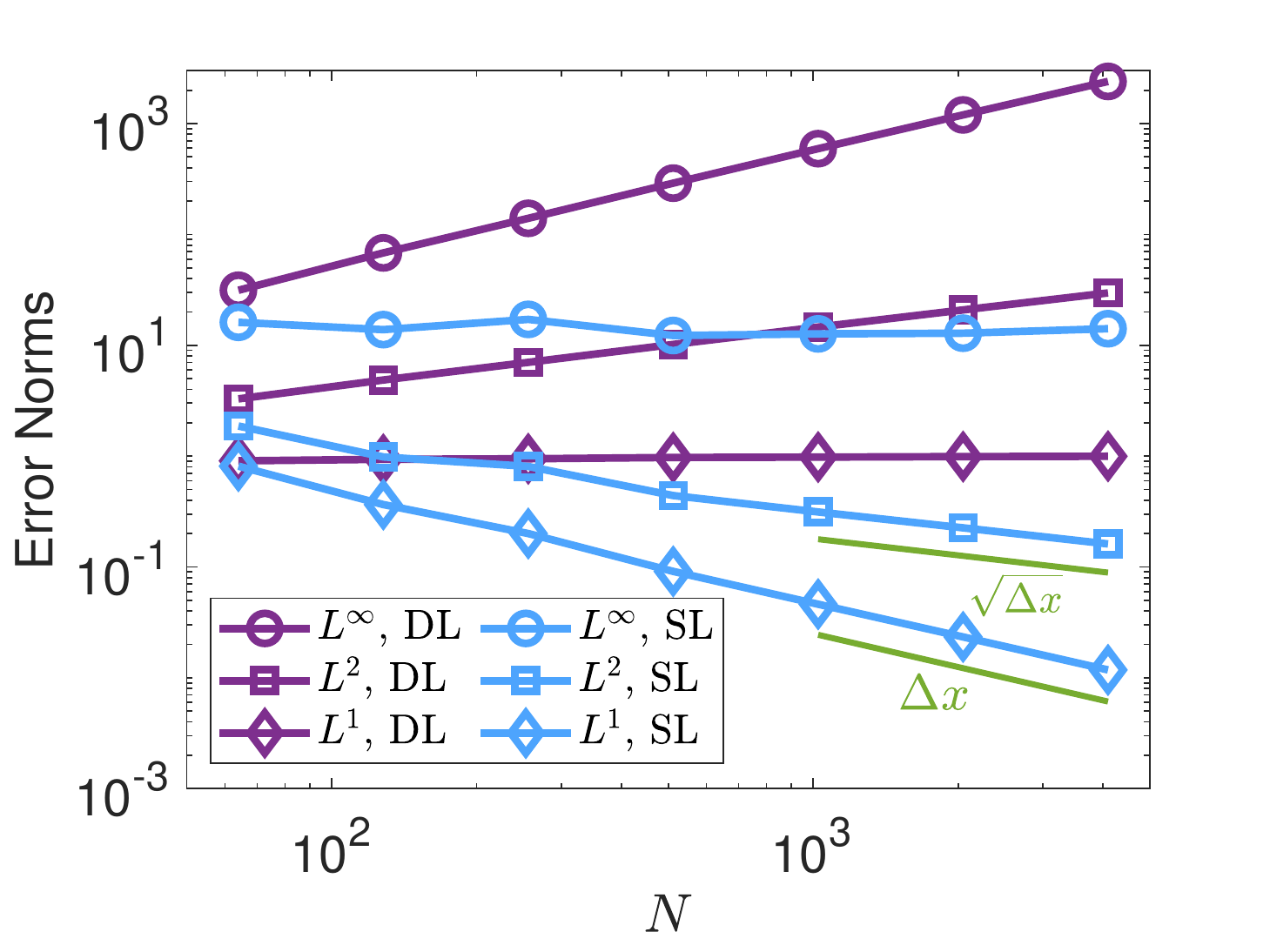}
\caption{\normalsize pressure, $\Omega$}
\label{pressure omega}
\end{subfigure}
\begin{subfigure}{0.495\textwidth}
\centering
\includegraphics[width=\textwidth]{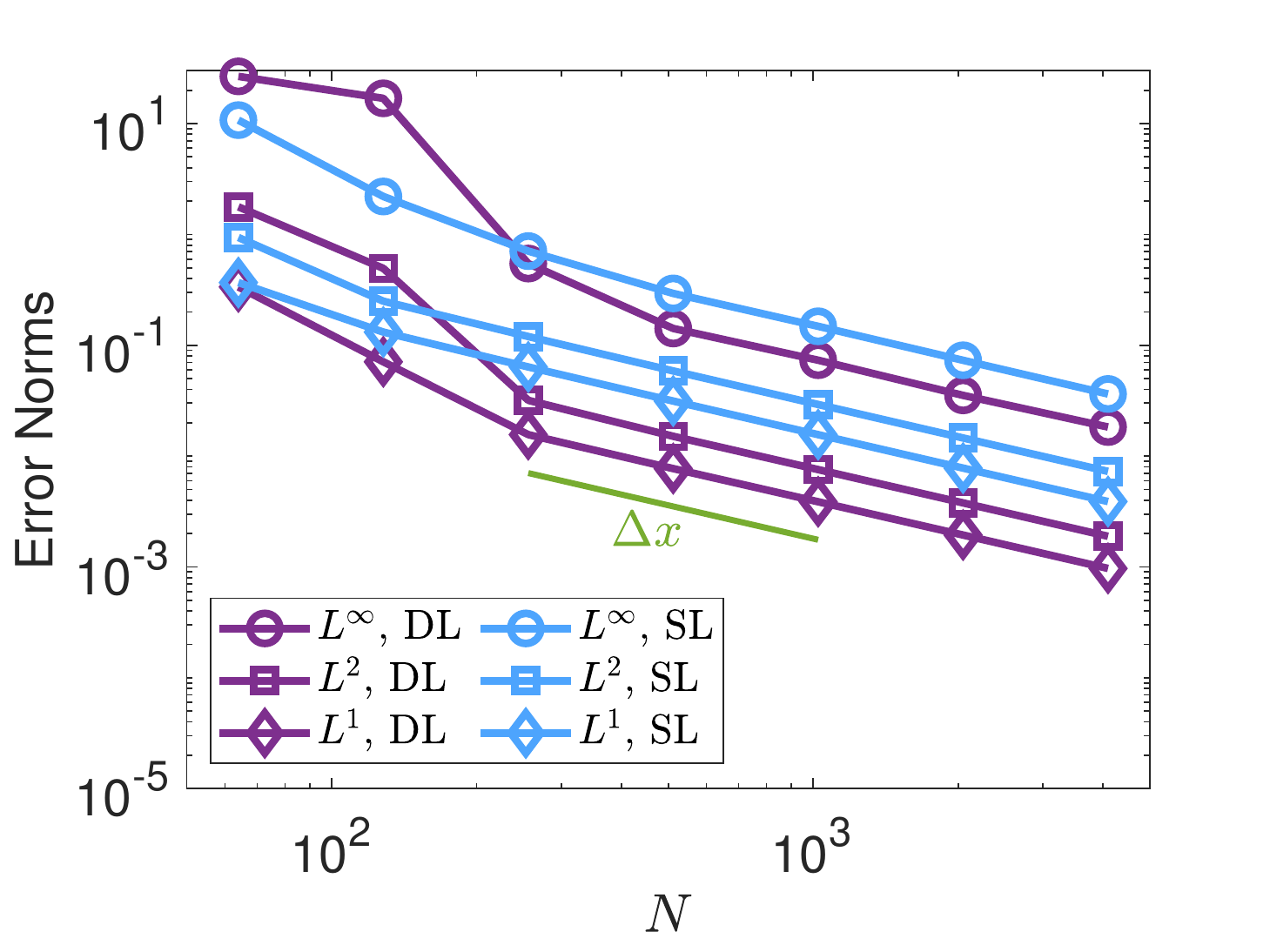}
\caption{\normalsize pressure, away from $\Gamma$}
\label{pressure away}
\end{subfigure}
\caption[Refinement studies for velocity solutions to the Brinkman equation \eqref{brink pde again} using the IBDL method and refinement studies for pressure solutions using the IBDL and IBSL methods]{Refinement studies for solutions to Equation \eqref{brink pde again} found using the IBSL and IBDL methods. The computational domain is the periodic box $[-1, 1]^2$, and the boundary point spacing is $ \Delta s  \approx 1.5 \Delta x$. We use a Fourier spectral method, and for IBDL, we use interpolation for grid points within $m_1=2(\log_2{N}-4)$. Figures \ref{IBDL brink refinement u}-\ref{IBDL brink refinement v} give the IBDL refinement studies for velocity. Figure \ref{pressure omega} gives the refinement study for pressure on $\Omega$ for both methods, and Figure \ref{pressure away} gives the refinement study for pressure on the portion of $\Omega$ that is away from $\Gamma$ by a distance of $0.02$. All errors are absolute.}\label{brink interior refinement dl}
\end{figure}

Table \ref{iteration table 6.1} gives the iteration counts for using \texttt{minres} and \texttt{gmres} to solve Equation \eqref{brink pde again} with the IBSL and IBDL methods, respectively. As expected from the formulation of the IBDL method, we see that the IBDL iteration counts are small and essentially constant as we refine the mesh or boundary point spacing. Furthermore, while the IBSL method fails to reach the desired tolerance for the $\Delta s \approx 0.75 \Delta x$, we have no such problem with the IBDL method. Figures \ref{IBDL brink refinement u}-\ref{IBDL brink refinement v} demonstrate the first-order convergence of the velocity solutions, and again, the error sizes are comparable between the two methods. Through this example, we once again see the greater efficiency of the IBDL method. 

We also include refinement studies for the pressure for the IBSL and IBDL methods. Figure \ref{pressure omega} shows the refinement studies for both methods on the entire PDE domain $\Omega$. Since the pressure is discontinuous with the IBSL method, it fails to converge pointwise and gets lower order convergence for the $L^2$ norm. However, in the IBDL method, the pressure actually blows up near the boundary, giving no convergence. Figure \ref{pressure away}, on the other hand, uses the errors for the portion of $\Omega$ that is away from $\Gamma$ by a distance of $0.02$. We can see that both solutions give first-order convergence of the pressure away from $\Gamma$ by this distance. 

\subsection{Choice of $\eta$}\label{6.5 eta}

For this section, we use the following exterior Brinkman equation to explore the choice of $\eta$ in the completed IBDL method. The PDE is  
\begin{subequations} \label{brink pde ext}
\begin{alignat}{2}
& \Delta \u - k^2  \u -\grad p = \g \qquad && \text{in } \Omega \\
& \Div \u=0 \qquad && \text{in } \Omega  \\
&\u=\U_b \qquad && \text{on } \Gamma, 
\end{alignat}
\end{subequations}
where we use the analytical solution given by
\begin{subequations} 
\begin{alignat}{2}
& u=e^{\sin{(\pi x)}}\cos{(\pi y)} \\
& v=-\cos{(\pi x)}\hspace{0.1cm}  e^{\sin{(\pi x)}}\sin{(\pi y)}\\
& p= e^{\cos{(\pi y)}}
\end{alignat}
\end{subequations}
to determine the boundary values, $\U_b=\u|_{\Gamma}$, and the forcing function,
\begin{equation}
\g=\begin{pmatrix}
\pi^2e^{\sin{(\pi x)}}\cos{(\pi y)}\hspace{0.1cm} \Big(\cos^2{(\pi x)}-\sin{(\pi x)}-1-k^2/\pi^2\Big)\\
-\pi^2\cos{(\pi x)} \hspace{0.1cm} e^{\sin{(\pi x)}}\sin{(\pi y)}\hspace{0.1cm} \Big(\cos^2{(\pi x)}-3\sin{(\pi x)}-2-k^2/\pi^2\Big)+\pi \sin{(\pi y)}\hspace{0.1cm}  e^{\cos{(\pi y)}}\end{pmatrix}.
\end{equation}
Here, we will use $k=0,0. 1,$ and $3$. $\Omega$ is the region of the periodic box, $\C=[-1, 1]^2$ that is exterior to a circle of radius 0.75, centered at the origin. We use a Fourier spectral method with equally spaced boundary points with $\Delta s \approx  \Delta x$, and we use interpolation for grid points within $m_1=2(\log_2{N}-4)$ of $\Gamma$.

\begin{figure}
\begin{subfigure}{0.495\textwidth}
\centering
\includegraphics[width=\textwidth]{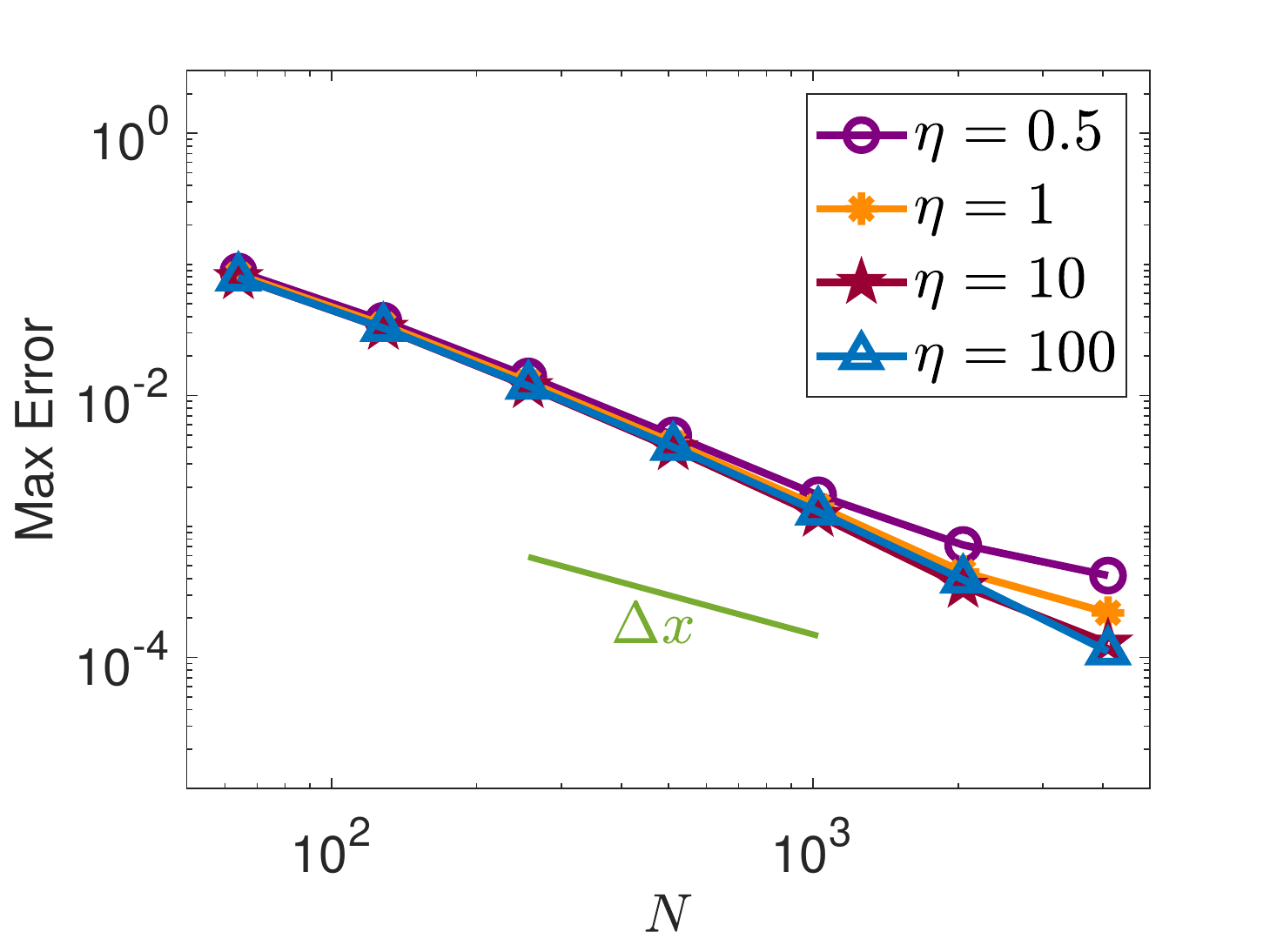}
\caption{\normalsize  Stokes equation}
\label{stokes etas}
\end{subfigure}
\begin{subfigure}{0.495\textwidth}
\centering
\includegraphics[width=\textwidth]{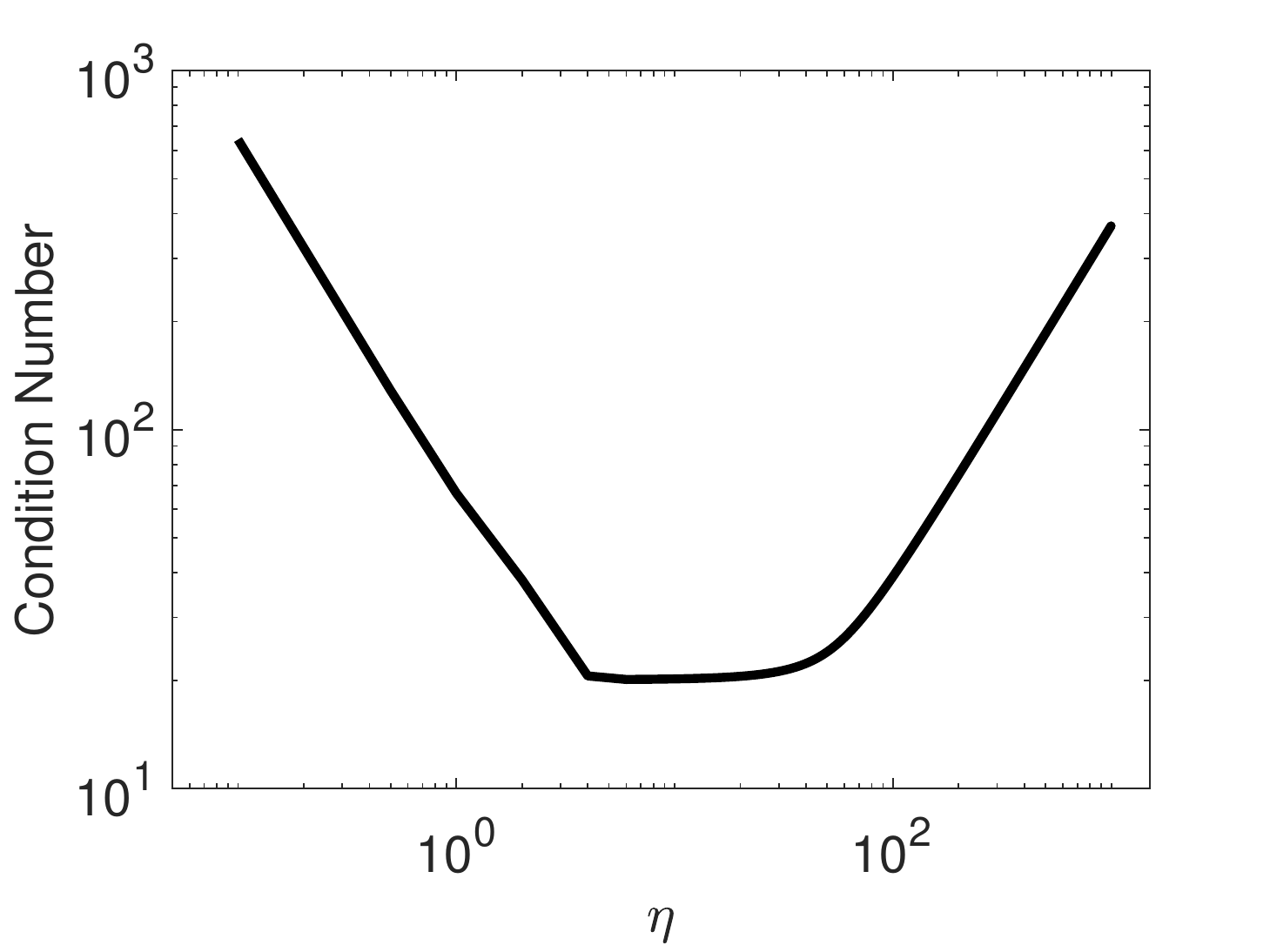}
\caption{\normalsize Stokes equation}
\label{condplot}
\end{subfigure}
\begin{subfigure}{0.495\textwidth}
\centering
\includegraphics[width=\textwidth]{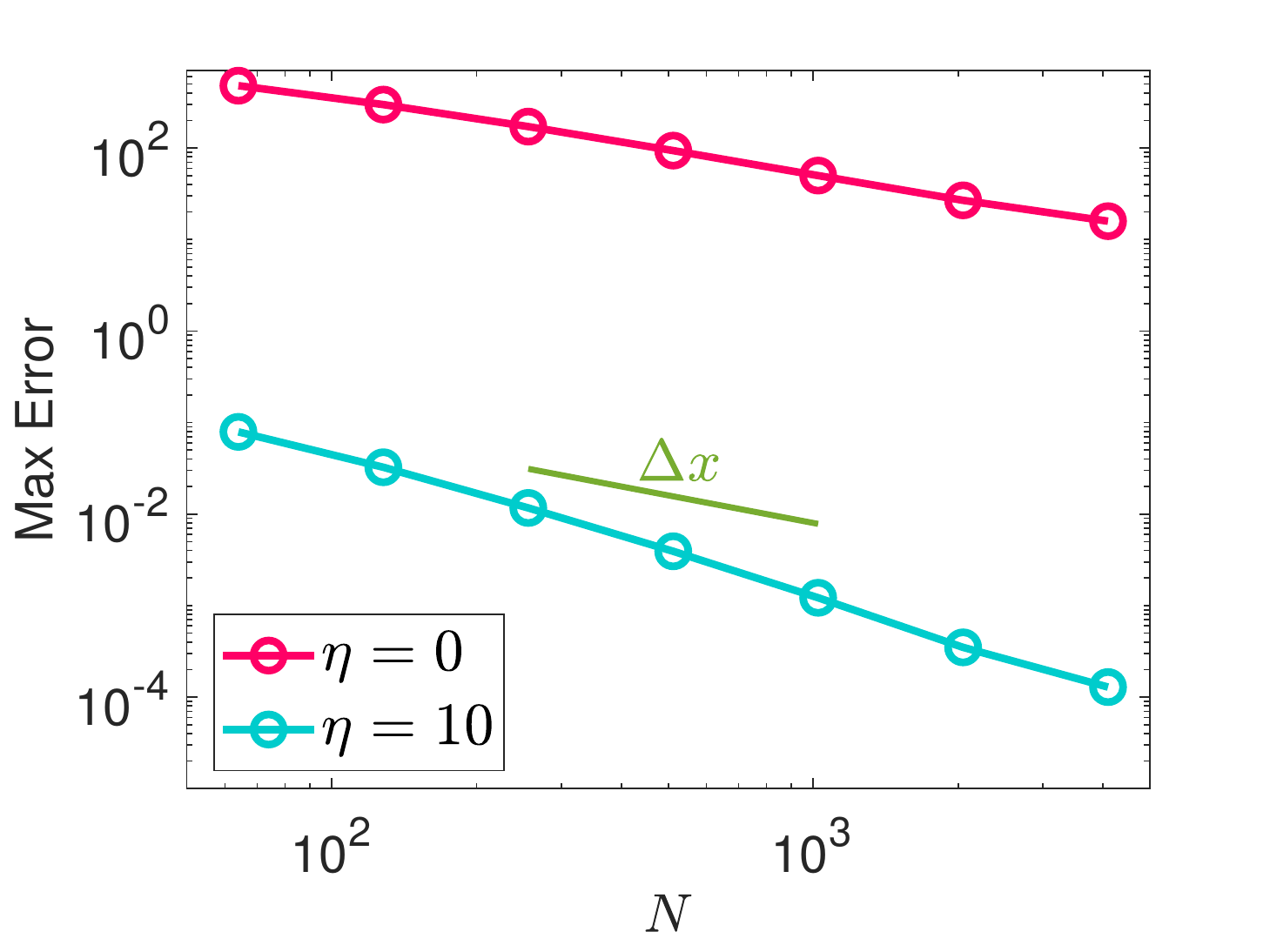}
\caption{\normalsize Brinkman, $k=0.01$}
\label{brink01}
\end{subfigure}
\begin{subfigure}{0.495\textwidth}
\centering
\includegraphics[width=\textwidth]{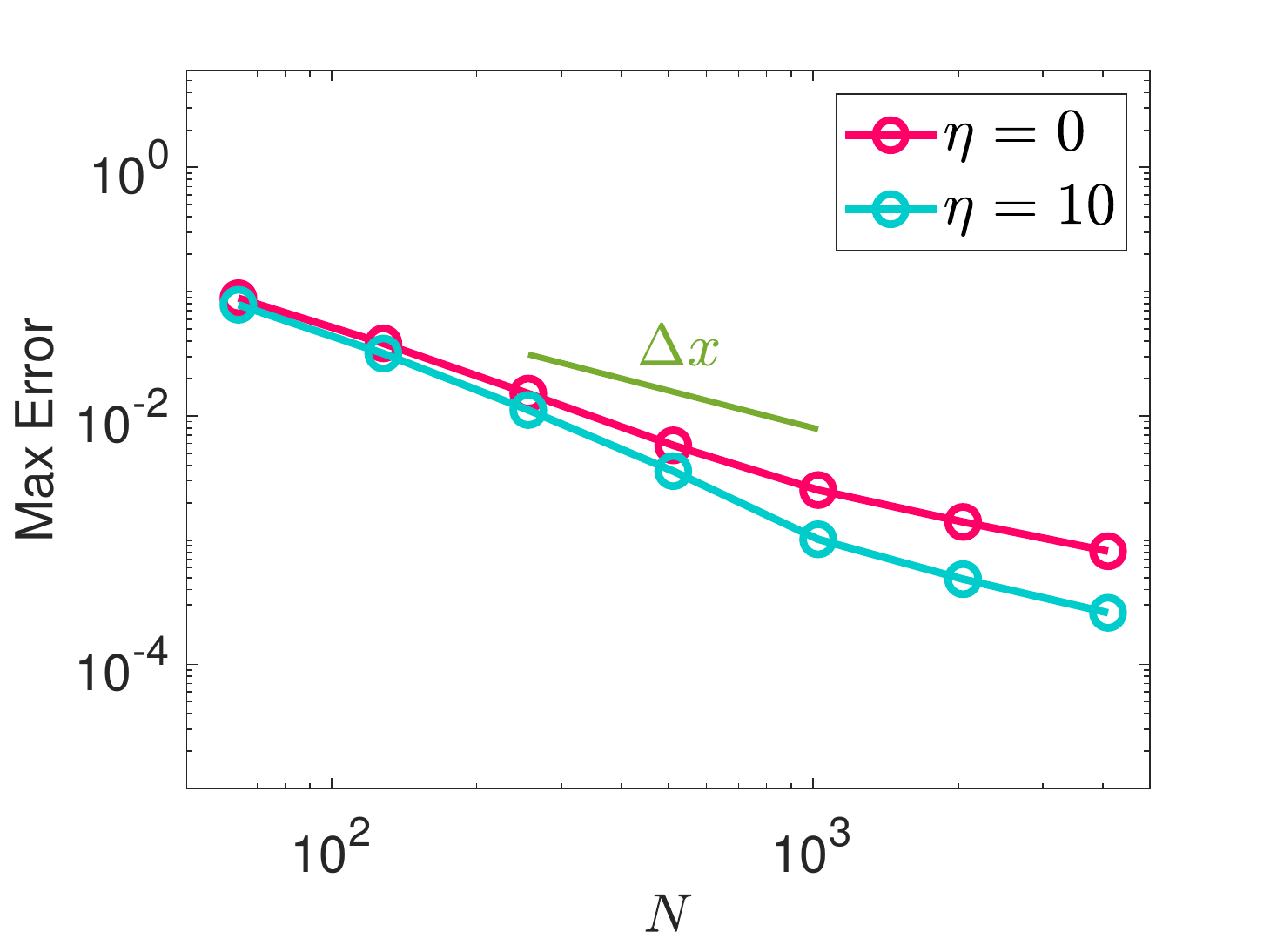}
\caption{\normalsize Brinkman, $k=3$}
\label{brink3}
\end{subfigure}
\caption[Plot of condition number versus $\eta$ and $L^{\infty}$ refinement studies for horizontal velocity solutions to the exterior Stokes and Brinkman equations \eqref{brink pde ext} using the completed IBDL method and various values of $\eta$]{Figures \ref{stokes etas}, \ref{brink01}, and \ref{brink3} give $L^{\infty}$ refinement studies for horizontal velocity solutions to the exterior Stokes and Brinkman equations \eqref{brink pde ext} using the completed IBDL method and various values of $\eta$. Figure \ref{condplot} gives a plot of the condition number versus $\eta$ for operator that must be inverted for Stokes equation. $\Omega$ is the region of periodic box, $\C=[-1, 1]^2$ that is exterior to a circle of radius 0.75, centered at the origin. The method uses a Fourier spectral method with equally spaced boundary points with $\Delta s \approx  \Delta x$, and interpolation is used for grid points within $m_1=2(\log_2{N}-4)$ of $\Gamma$. All errors are absolute.}\label{eta plots}
\end{figure}

We begin with Stokes equation, or $k=0$. In this case, in order to represent flows with net force and net torque, we must use the completed IBDL method, with $\eta>0$. Figure \ref{stokes etas} shows the $L^{\infty}$ refinement studies for this problem for four different values of $\eta$. We can see that using a larger value of $\eta$ can result in lower pointwise error on fine grids. For this problem, there is no meaningful benefit to increasing $\eta$ above $10$. We have observed this to be true for other problems as well. For a range of $\eta$ values, Figure \ref{condplot} plots the condition number of the operator that must be inverted to solve for $\Q$. This operator is seen in Equation \eqref{discrete operator for Q}. Selecting a very small $\eta$ value will result in a large condition number, as the operator approaches that for the IBDL without the completion. The smallest condition number occurs at $\eta\approx 6$. After this, the condition number grows as we increase $\eta$. However, for moderate values of $\eta$ there is not a large impact on the iteration counts. For this problem, the iteration count for $\eta=10$ is about 8, and the iteration count for $\eta=100$ is still only about 13. To compare, the condition number for the IBSL operator for $N=2^{7}$ and $\Delta s \approx \Delta x$ is approximately $2.6(10)^9$. For the wider boundary point pacing of $\Delta s\approx 2\Delta x$, the condition number is approximately $1.7(10)^{7}$. These condition numbers are much higher than any of the condition numbers we see for $\eta$ values ranging from $1$ to $1000$. Additionally, finer grids would result in even larger condition numbers for the IBSL method. Therefore, we see that there is flexibility in the choice of $\eta$, and we need not fear approaching the poor conditioning of the IBSL method. For the majority of this dissertation, we use $\eta=10$. 

We next look at the Brinkman equation for $k=0.01$ and $k=3$. Figures \ref{brink01}-\ref{brink3} show the refinement studies for utilizing the IBDL method without the completion and with the completion for $\eta=10$. We can see that for a small $k$ value, we may still obtain convergence without the completion, but including the single layer term is important to obtain small errors. Using $\eta=10$ in this case drastically lowers the errors. For larger $k$ values, as the Brinkman equation diverges from Stokes equation, the need to include the single layer term diminishes, but including it often lowers the error, as we see for $k=3$. 

\subsection{Finite difference method}\label{6.5 fd}

For this section, we use the following PDE to examine the results from the IBDL method with two different finite difference discretizations for $\Div \grad$.  The PDE is  
\begin{subequations} \label{stokes pde another}
\begin{alignat}{2}
& \Delta \u  -\grad p = \g \qquad && \text{in } \Omega \\
& \Div \u=0 \qquad && \text{in } \Omega  \\
&\u=\U_b \qquad && \text{on } \Gamma, 
\end{alignat}
\end{subequations}
where we use the analytical solution given by
\begin{subequations} 
\begin{alignat}{2}
& u=\sin{y}-x e^{xy} \\
& v=\cos{x}+ye^{xy}\\
& p= e^{x+y}
\end{alignat}
\end{subequations}
to determine the boundary values, $\U_b=\u|_{\Gamma}$, and the forcing function,
\begin{equation}
\g=\begin{pmatrix}
-e^{xy} \big(2y+xy^2+x^3\big) - \sin{y} - e^{x+y} \\
e^{xy} \big(2x+x^2y+y^3\big) - \cos{x} - e^{x+y}\end{pmatrix}.
\end{equation}
Here $\Omega$ is the interior of a circle of radius $0.25$, centered at the origin. The computational domain is $\C=[-0.5, 0.5]^2$. We use equally spaced boundary points with $\Delta s \approx 0.75 \Delta x$, and we use interpolation for grid points within $m_1=6$ meshwidths of $\Gamma$. We also use the smoother 6-point B-spline delta function defined in Section \ref{6.4 space}. 

\begin{figure}
\begin{subfigure}{0.43\textwidth}
\centering
\includegraphics[width=\textwidth]{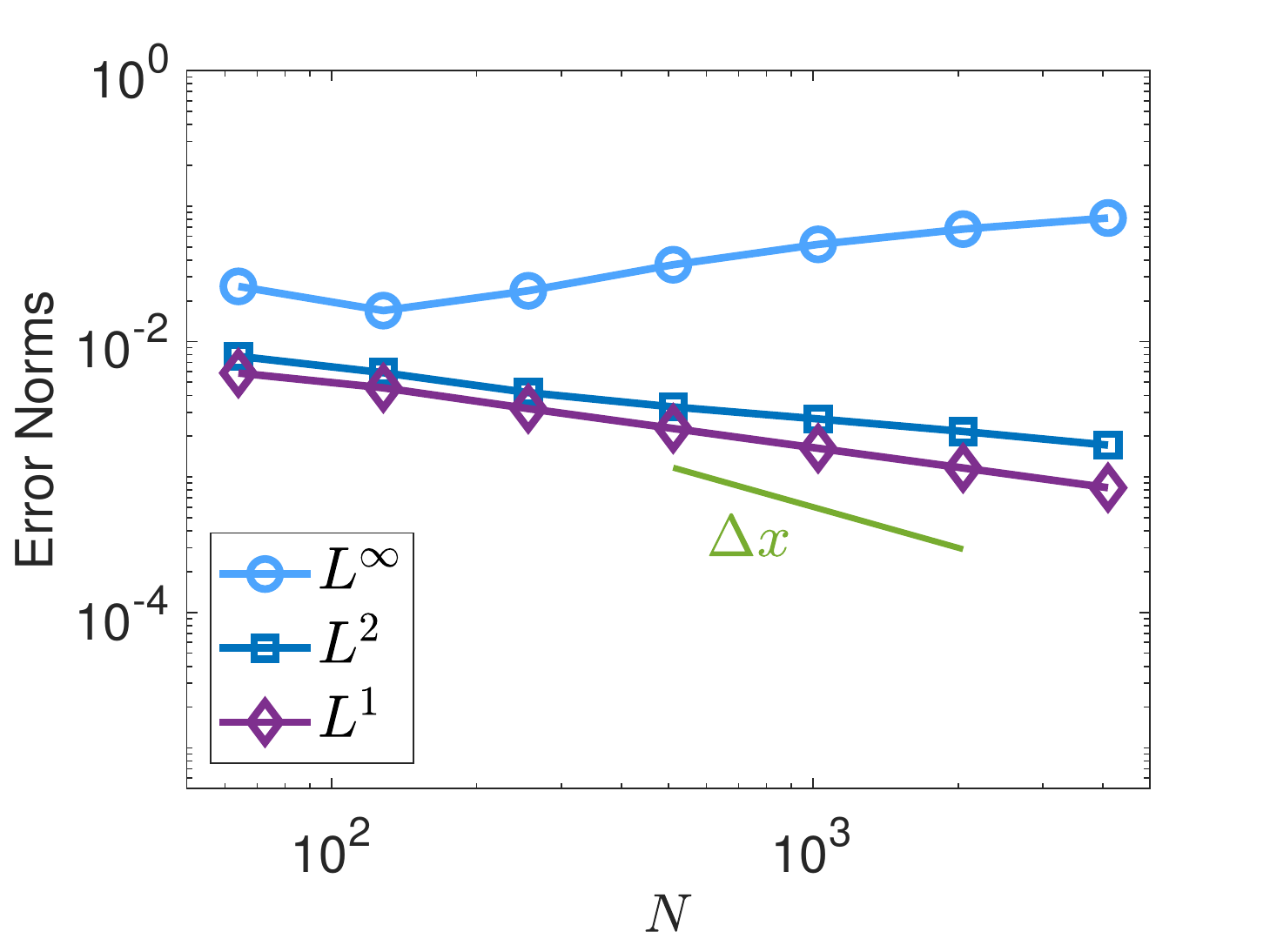}
\caption{\normalsize  refinement, method 1}
\label{wide refinement}
\end{subfigure}
\begin{subfigure}{0.43\textwidth}
\centering
\includegraphics[width=\textwidth]{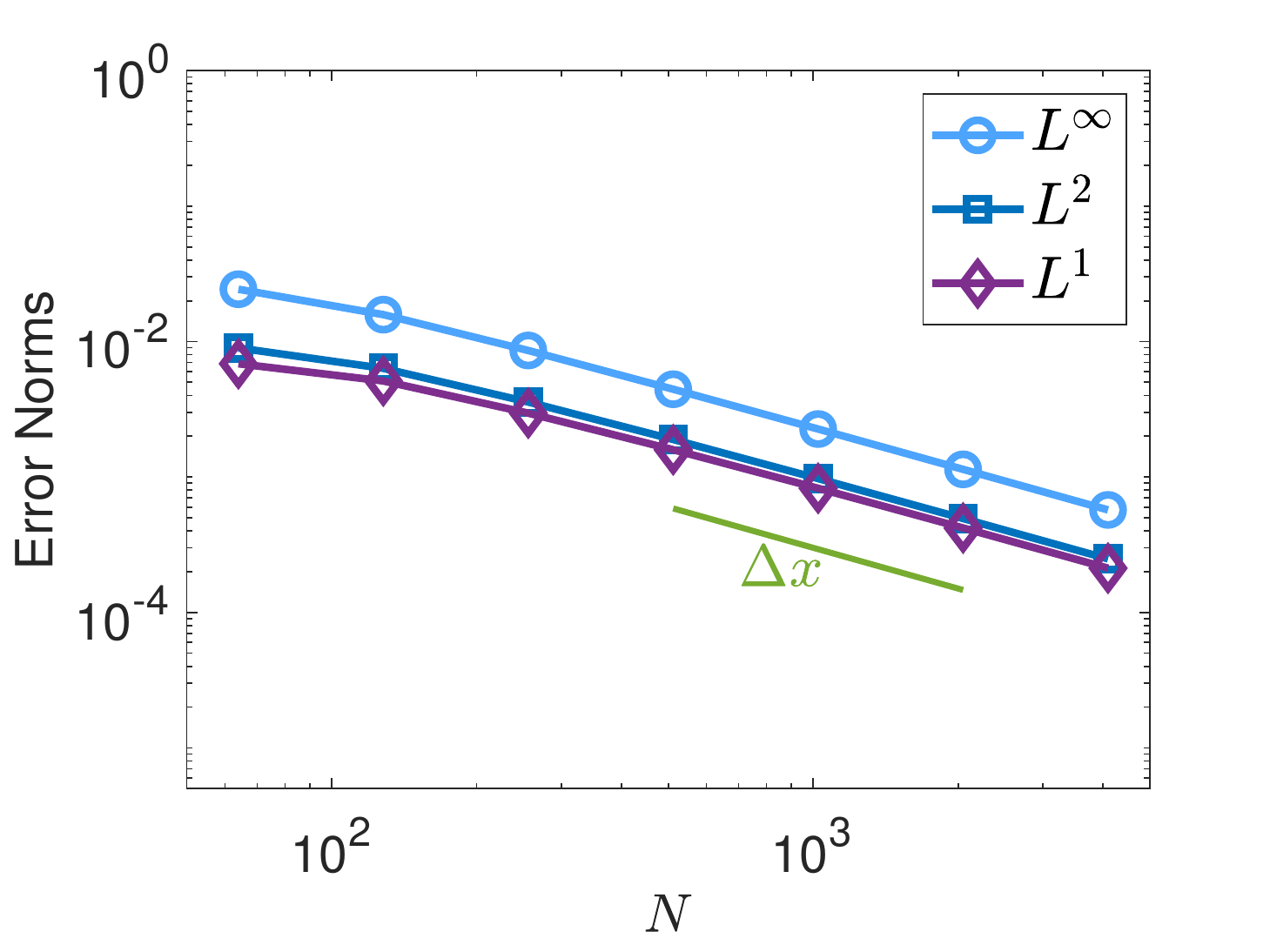}
\caption{\normalsize refinement, method 2}
\label{regular refinement}
\end{subfigure}
\begin{subfigure}{0.43\textwidth}
\centering
\includegraphics[width=\textwidth]{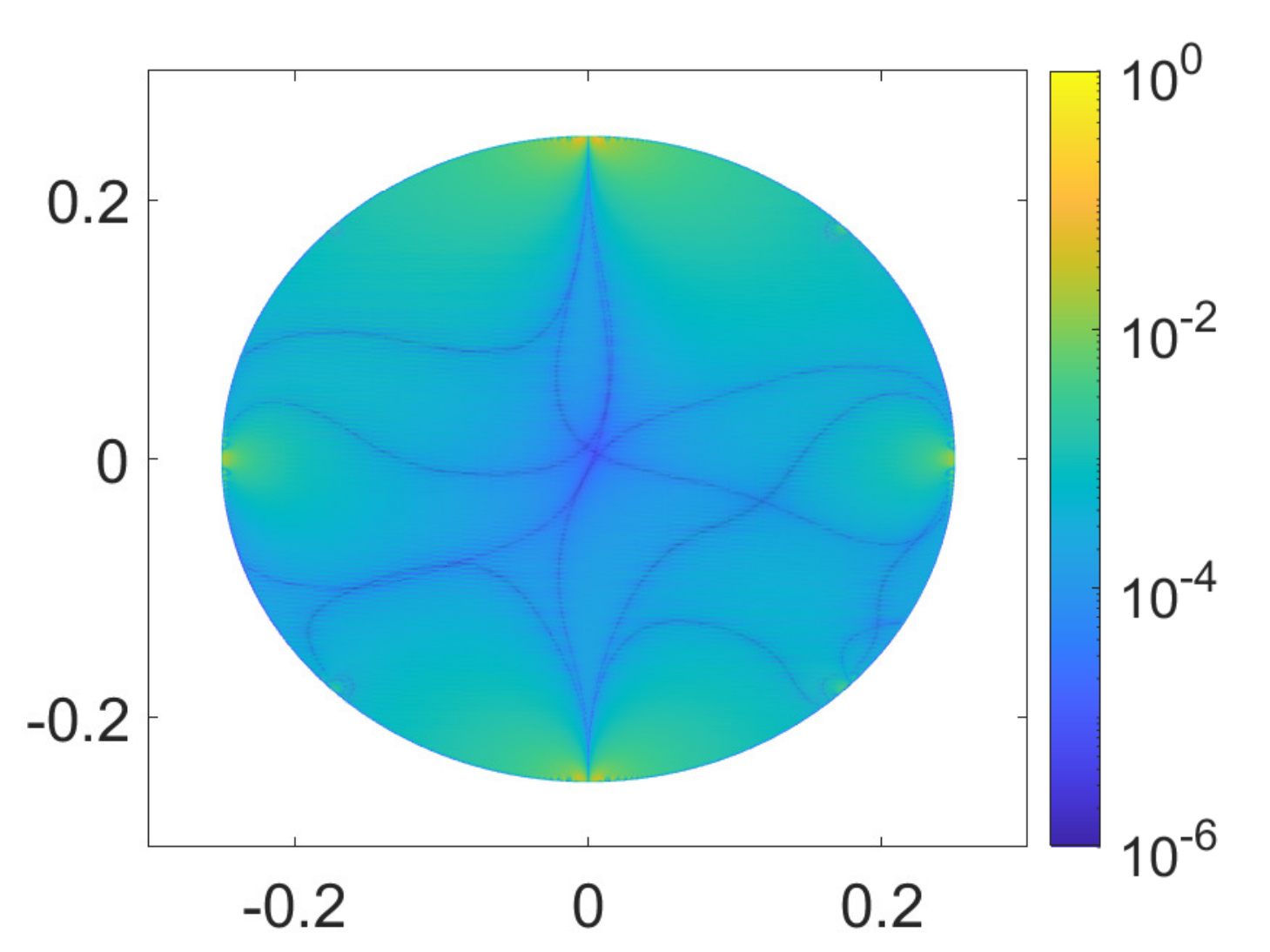}
\caption{\normalsize error, method 1 }
\label{wide error}
\end{subfigure}
\begin{subfigure}{0.43\textwidth}
\centering
\includegraphics[width=\textwidth]{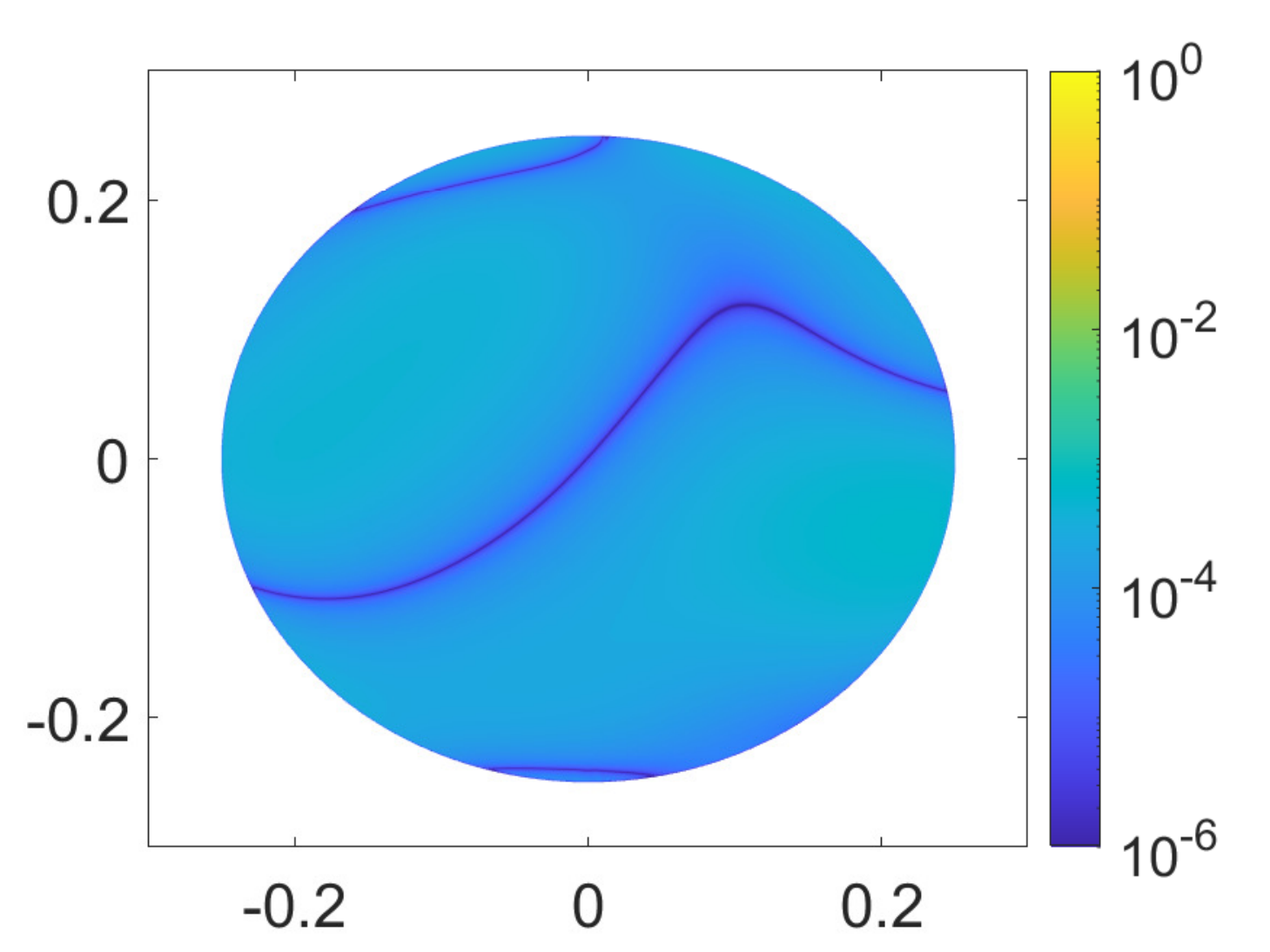}
\caption{\normalsize error, method 2}
\label{regular error}
\end{subfigure}
\begin{subfigure}{0.43\textwidth}
\centering
\includegraphics[width=\textwidth]{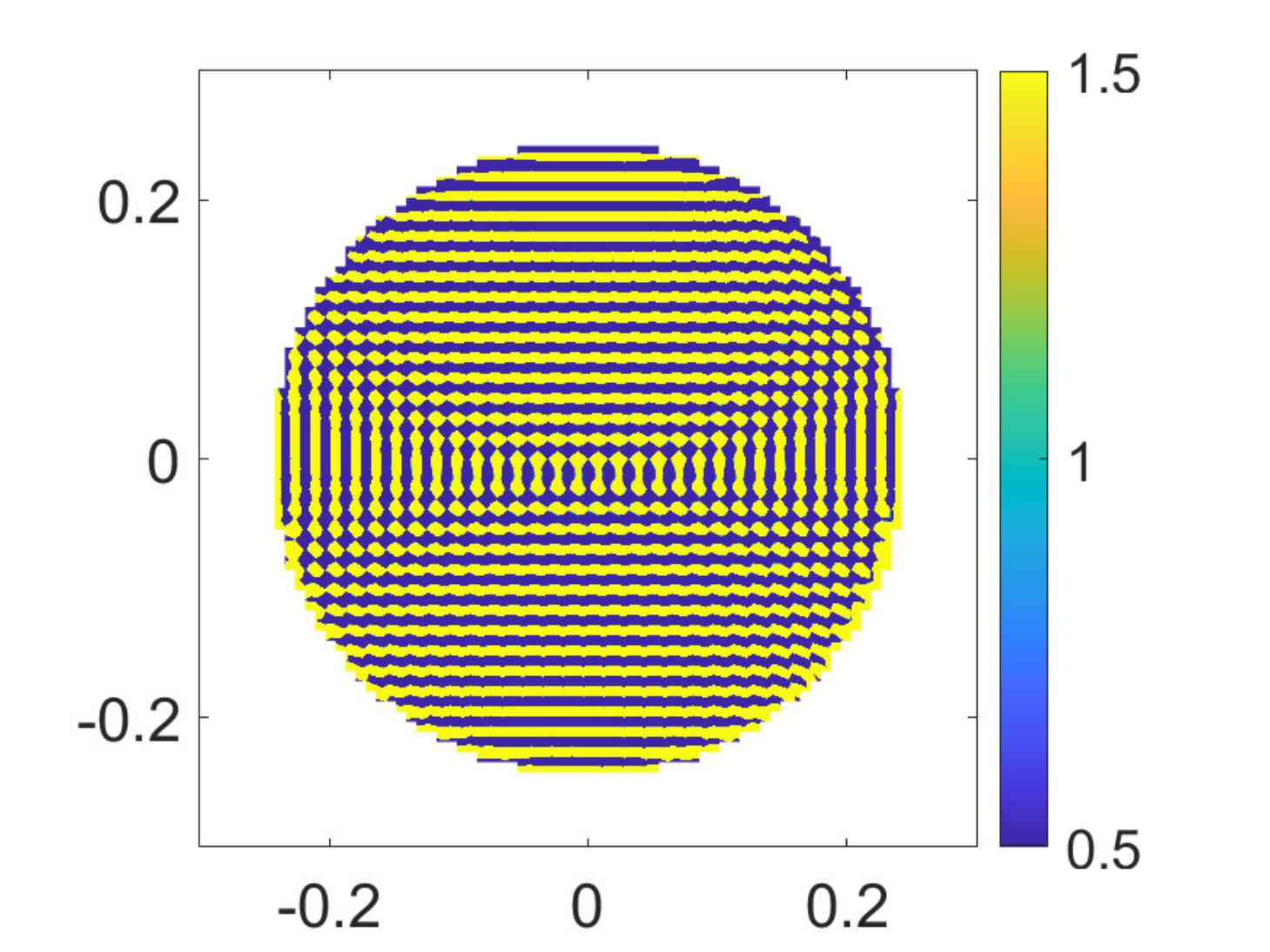}
\caption{\normalsize pressure, method 1 }
\label{wide pressure}
\end{subfigure}
\begin{subfigure}{0.43\textwidth}
\centering
\includegraphics[width=\textwidth]{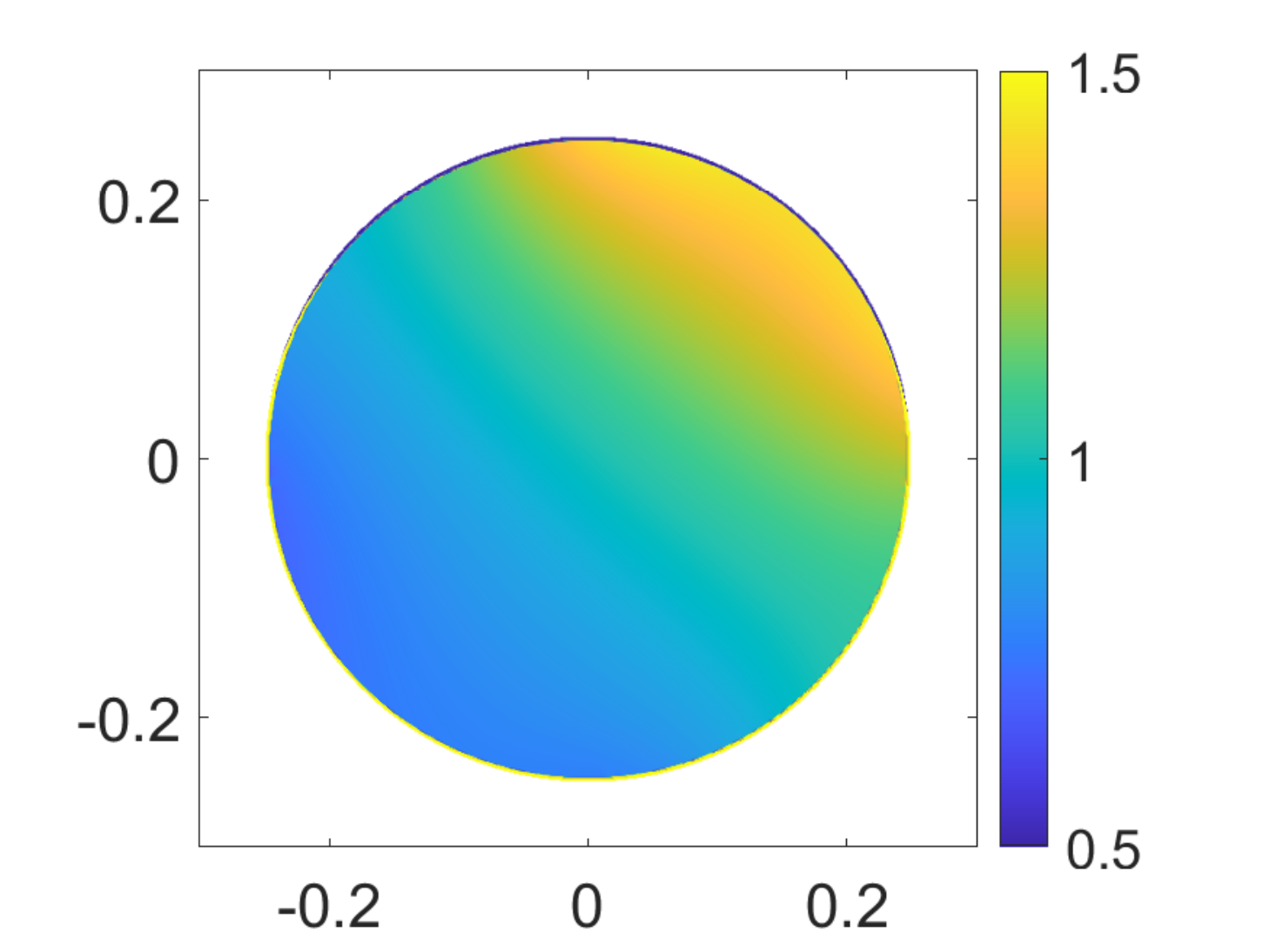}
\caption{\normalsize pressure, method 2}
\label{regular pressure}
\end{subfigure}
\caption[Refinement studies and error plots for horizontal velocity solutions and pressure solution plots for the Stokes PDE in Equation \eqref{stokes pde another}, found using the IBDL method with two forms of finite difference discretizations]{Refinement studies and plots for Equation \eqref{stokes pde another}, found using the IBDL method with the two forms of finite difference discretizations described in Section \ref{6.5 fd}. Figures \ref{wide refinement}-\ref{regular refinement} give the refinement studies for horizontal velocity, and Figures \ref{wide error}-\ref{regular err} show the horizontal velocity errors for $N=2^{13}$. Figures \ref{wide pressure}-\ref{regular pressure} give the pressure solutions for $N=2^7$. $\Omega$ is the interior of a circle of radius $0.25$, centered at the origin, and the computational domain is $\C=[-0.5, 0.5]^2$. The methods use the 6-point B-spline delta function defined in Section \ref{6.4 space} and equally spaced boundary points with $\Delta s \approx 0.75 \Delta x$. Interpolation is used for grid points within $m_1=6$ meshwidths of $\Gamma$. The errors used in the refinement studies are absolute, and Figures \ref{wide error}-\ref{regular pressure} use a grid size of $N=2^{12}$.}\label{eta plots}
\end{figure}

As discussed in Section \ref{2.3 space and op}, the discrete divergence of the discrete gradient results in a Laplacian with a wider stencil. In Method 1, we invert this wide Laplacian when solving for pressure, as we do in the IBSL method. Figure \ref{wide refinement} shows the refinement study for horizontal velocity, and Figure \ref{wide error} shows the horizontal velocity error for the grid size $N=2^{12}$. Figure \ref{wide pressure} shows the pressure solution for grid size $N=2^7$. To plot the pressure, we again find the difference between the computed pressure and analytical solution at the origin and add this constant to the computed pressure. We do not achieve first-order convergence, and Figure \ref{wide error} shows larger pointwise errors, especially as we near the boundary, and we note that increasing the interpolation width is not sufficient to recover pointwise convergence.  Furthermore, Figure \ref{wide pressure} demonstrates that the interaction between the two discrete Laplacian operators case high frequency oscillations in the pressure solution. 

In Method 2, we approximate the discrete divergence of the discrete gradient with the usual 5-point, second-order accurate discrete Laplacian. Figures \ref{regular refinement}, \ref{regular error}, and \ref{regular pressure} give the corresponding refinement study, velocity error, and pressure solution. We can see that this approximation, combined with the smoother delta function, gives us first-order convergence of the solution and avoids the high frequency oscillations in the pressure. More investigation into other possible finite difference discretizations is needed, but Method 2 provides a simple option.

\subsection{Stokes flow around a periodic array of cylinders}\label{6.5 cylinder}

We next look at steady Stokes flow around a doubly periodic array of cylinders. This problem has been studied at length, and we can compare our numerical results for the drag force on the cylinder to literature and asymptotic approximations \cite{SandA1982a, Hasimoto,GreengardKropinski, GriffithDonev}. To do this, we look at the PDE 
\begin{subequations} \label{stokes pde flow}
\begin{alignat}{2}
& \Delta \u  -\grad p =\g \qquad && \text{in } \Omega  \label{stokes pde flow 1}\\
& \Div \u=0 \qquad && \text{in } \Omega  \\
&\u=0 \qquad && \text{on } \Gamma, 
\end{alignat}
\end{subequations}
where $\Omega$ is the region in the periodic box $[-0.5, 0.5]^2$ that is exterior to a disk centered at the origin, and $\Gamma$ is the boundary of this disk. Let us define $c$ as the fraction of the two-dimensional area occupied by the disk. Here, $c=\pi r^2$, and the maximum value is $c_{max}=\pi/4\approx 0.785$. We explore the dimensionless drag force per unit length of a cylinder obtained from the IBDL method for different values of $c$.  Let the fluid viscosity be $\mu$, which is chosen as $1$ in Equation \eqref{stokes pde flow 1}. Let $U$ be the average velocity in the $x$-direction on the left side of the periodic box. In other words, 
\begin{equation}
U=\int_{-0.5}^{0.5} u(-0.5, y)dy.
\end{equation}
We drive the flow with the forcing function
\begin{equation}
\g=\begin{pmatrix} -1 \\0\end{pmatrix},
\end{equation} 
and numerically estimate the integral for $U$. Recall that in Chapter \ref{chapter 3}, we denoted the net force on an object as $\bm{B}$. For the IBDL method, we solve Equation \eqref{stokes pde flow} using the completed formulation described in Section \ref{6.1 Hebeker}. We therefore add a supplementary flow provided by a single layer potential with strength $\eta \Q$, with $\eta=10$. Then, as discussed in Section \ref{6.1 Hebeker}, we calculate the net force by numerically approximating 
\begin{equation}
\bm{B}=-\eta\int_{\Gamma} \Q(s)ds. 
\end{equation}
Similarly, as discussed in Section \ref{3.5 SL}, we calculate the net force in the IBSL case by numerically approximating 
\begin{equation}
\bm{B}=-\int_{\Gamma} \F(s)ds. 
\end{equation}
Then we calculate the dimensionless drag force, $D=B_1/(\mu U)$. We compare the values of $D$ from the IBDL method with those obtained from the IBSL method and the integral method of Greengard and Kropinski \cite{GreengardKropinski}. We also compare the results to the asymptotic expressions for the dilute ($c\ll 1$) and dense ($c_{max}-c\ll1$) cases \cite{SandA1982a, Hasimoto}. These expressions are given by
\begin{equation}
D_{den}\approx \frac{9\pi}{2\sqrt{2}}\Bigg(1-\sqrt{\frac{c}{c_{max}}}\Bigg)^{-5/2}
\end{equation}
\begin{equation}
D_{dil}\approx \frac{8\pi}{\log{1/c} - 1.47633597+2c-1.77428264c^2+4.07770444c^3-4.84227402 c^4}.
\end{equation}

\begin{figure}
\begin{center}
 \begin{tabular}{||c | c | c | c | c | c ||} 
 \hline
 \multicolumn{6}{||c||}{Dimensionless Drag, $D$} \\
 \hline
  $c$ & IBDL & IBSL&G\&K & $D_{dil}$ & $D_{den}$   \\ 
 \hline
0.05 & 15.5591     & 15.5648       & 15.5578&15.5578&--     \\
0.1&   24.8368     &  24.8431      &24.8317&24.8323&--\\
0.2&    51.5459    & 51.5553       &51.5269&51.6068&--\\
0.3& 102.940       &  102.953      &102.881&105.234&--\\
0.4&  218.094      & 218.093       &217.894&270.011&--\\
0.5&  533.420      &   533.226     &532.548&  --       &544.313\\
0.6&   1770.23     & 1767.14       &1763.57&--    &1775.23\\
0.7&   13775.8     &  13580.9      &13519.3&--&13512.3\\
0.71& 19087.8&    18706.3           &--     & --   &18607.8\\
0.72& 27734.8&       26961.8       &--   &--& 26785.0\\
0.75& 147482     & 127572       &127543&--&127424\\
  \hline
  \end{tabular}
\captionof{table}[Dimensionless drag forces for IBDL and other methods for Stokes flow past a periodic array of cylinders]{Dimensionless drag, $D$, for various values of area concentration, $c$, for Stokes flow past a periodic array of cylinders. Included are results from the IBDL method, the IBSL method, and the boundary integral method of Greengard and Kropinski \cite{GreengardKropinski}. The last two columns give dilute and dense asymptotic approximations \cite{SandA1982a, Hasimoto}. }
\label{drag table}
\end{center}
\end{figure}
\begin{figure}
\begin{subfigure}{0.495\textwidth}
\centering
\includegraphics[width=\textwidth]{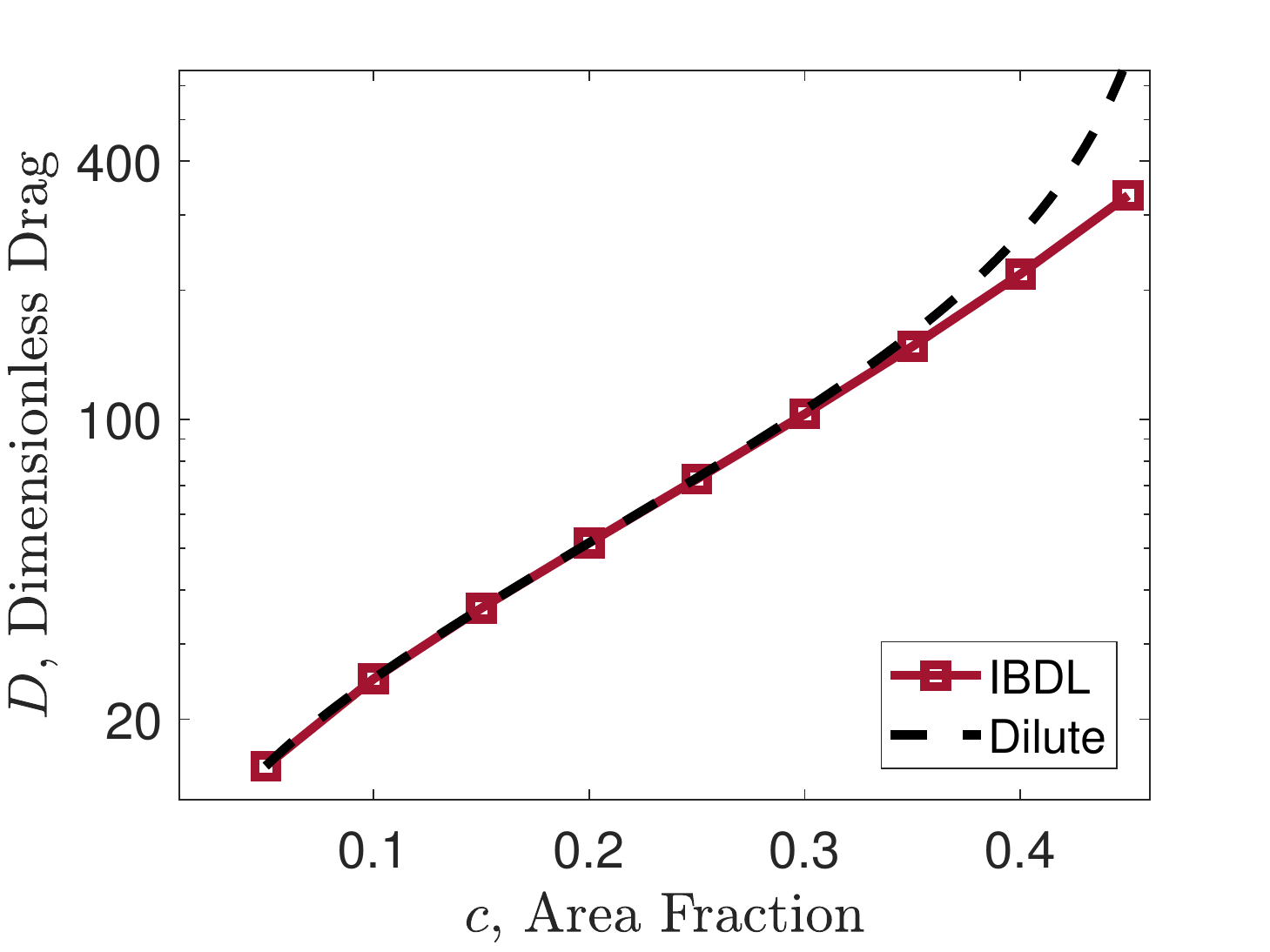}
\caption{\normalsize  small $c$}
\label{drag plot dilute}
\end{subfigure}
\begin{subfigure}{0.495\textwidth}
\centering
\includegraphics[width=\textwidth]{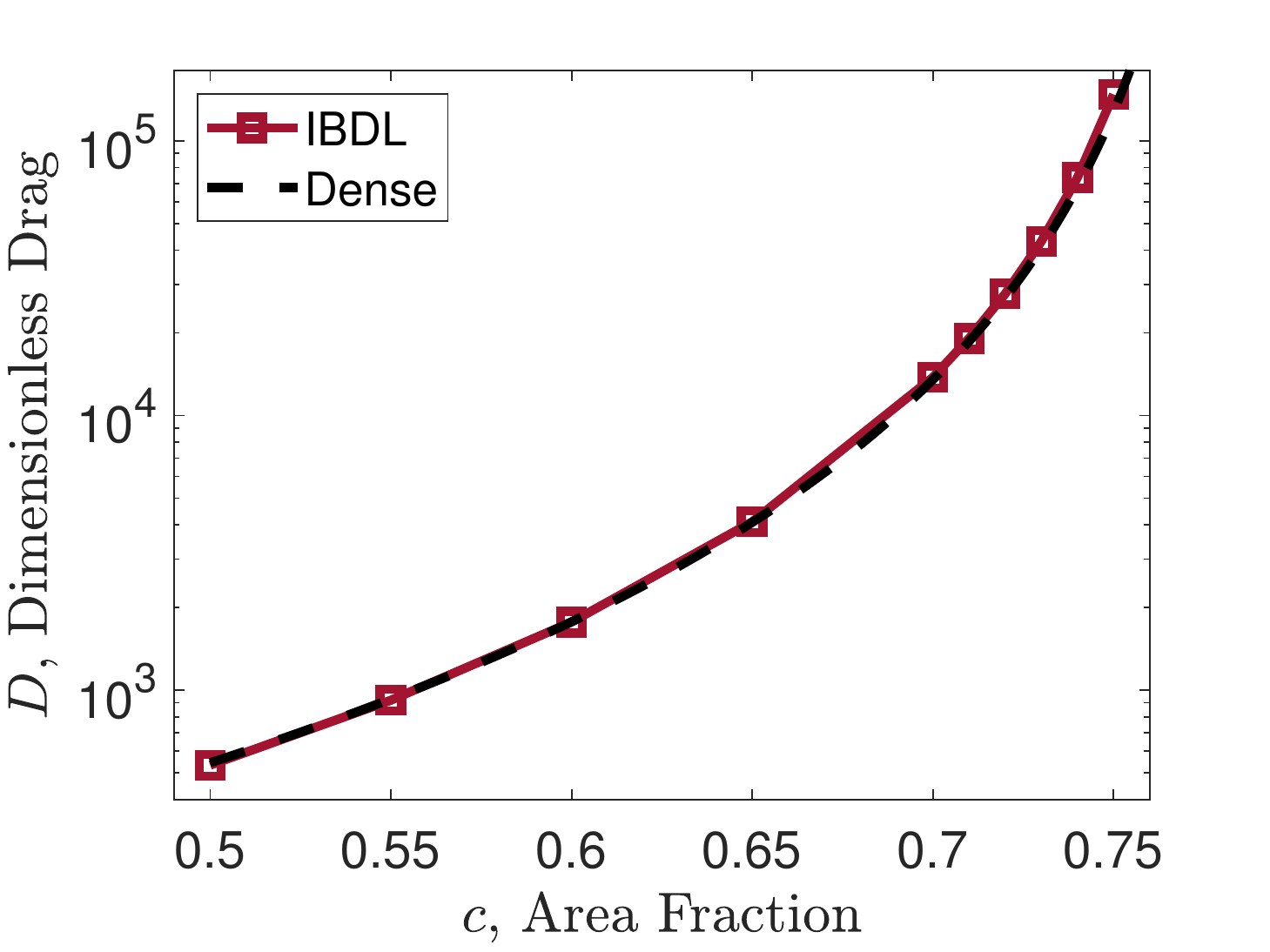}
\caption{\normalsize  large $c$}
\label{drag plot dense}
\end{subfigure}
\caption[Plots of dimensionless drag versus area concentration for Stokes flow past a periodic array of cylinders]{Dimensionless drag, $D$ for different area concentrations for Stokes flow past a periodic array of cylinders. Included are results from the IBDL method and the dilute and dense asymptotic approximations. For the IBDL method, we use a grid size of $N=2^{13}$, boundary point spacing of $\Delta s \approx \Delta x$, and interpolation width $m_1=2(\log_2{N}-4 )$.}\label{drag plot}
\end{figure}

\begin{figure}
\begin{subfigure}{0.495\textwidth}
\centering
\includegraphics[width=\textwidth]{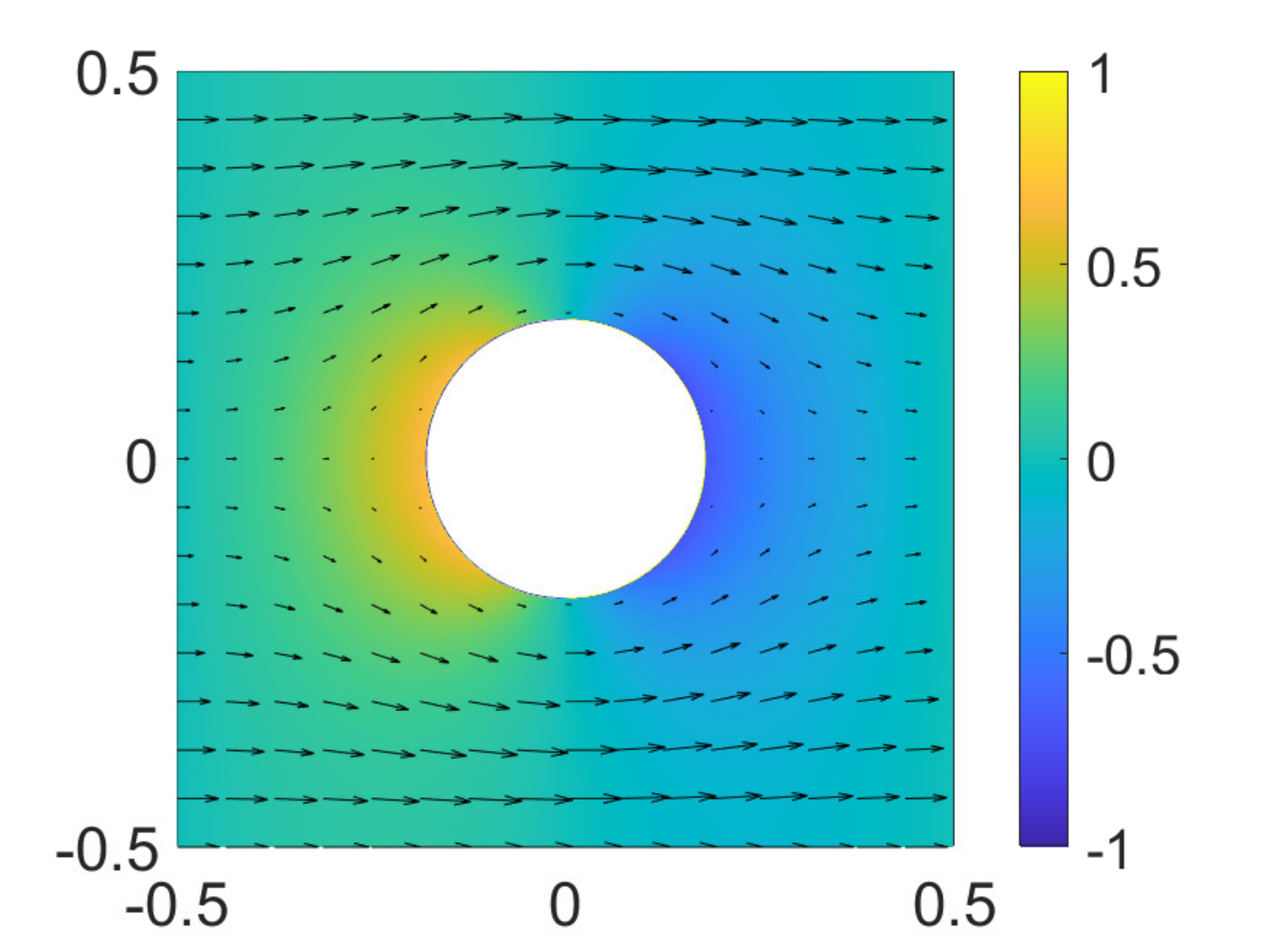}
\caption{\normalsize  Solution plot}
\label{flow past cylinder solution}
\end{subfigure}
\begin{subfigure}{0.495\textwidth}
\centering
\includegraphics[width=\textwidth]{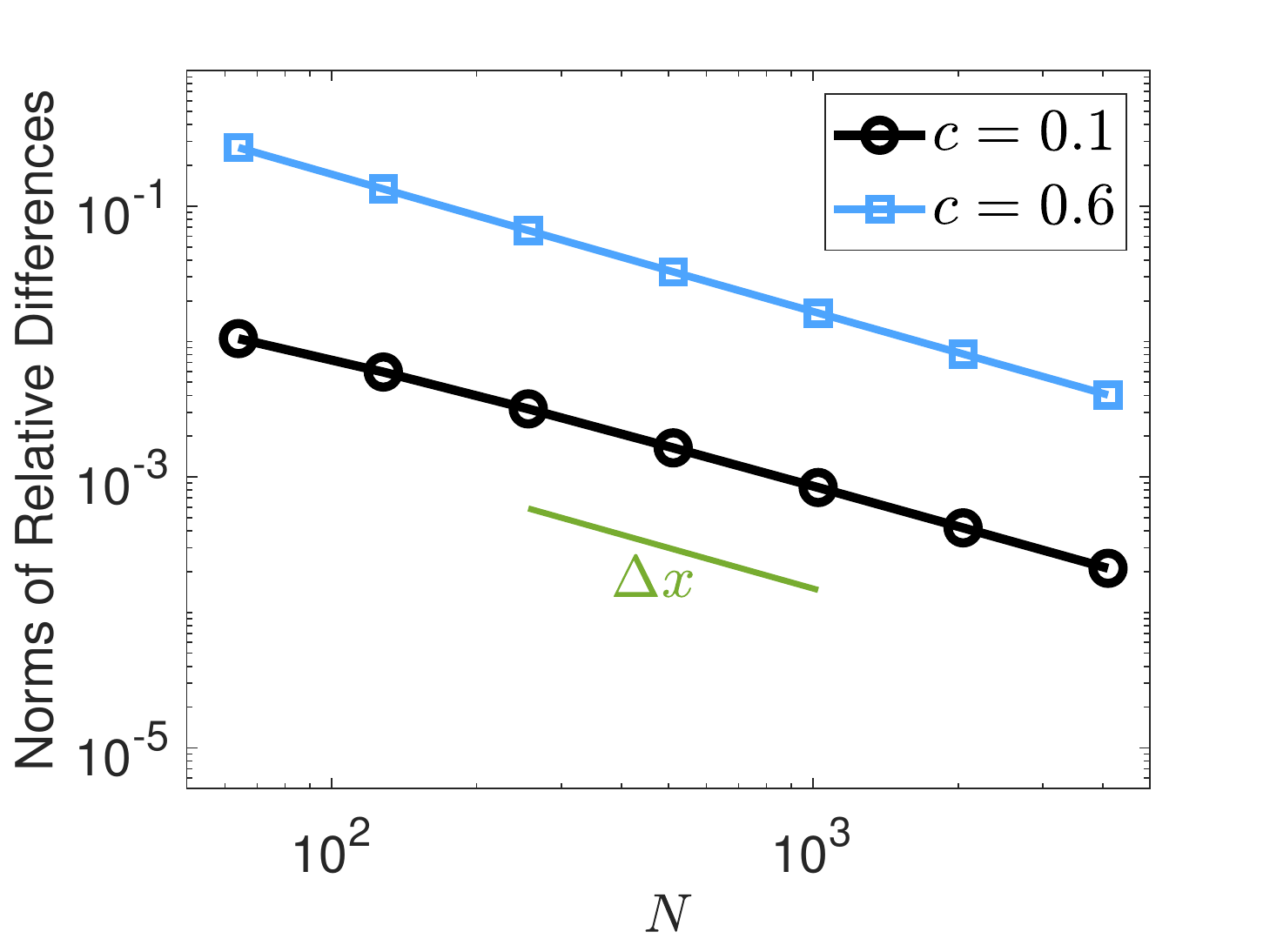}
\caption{\normalsize  Drag}
\label{flow past cylinder drag ref}
\end{subfigure}
\begin{subfigure}{0.495\textwidth}
\centering
\includegraphics[width=\textwidth]{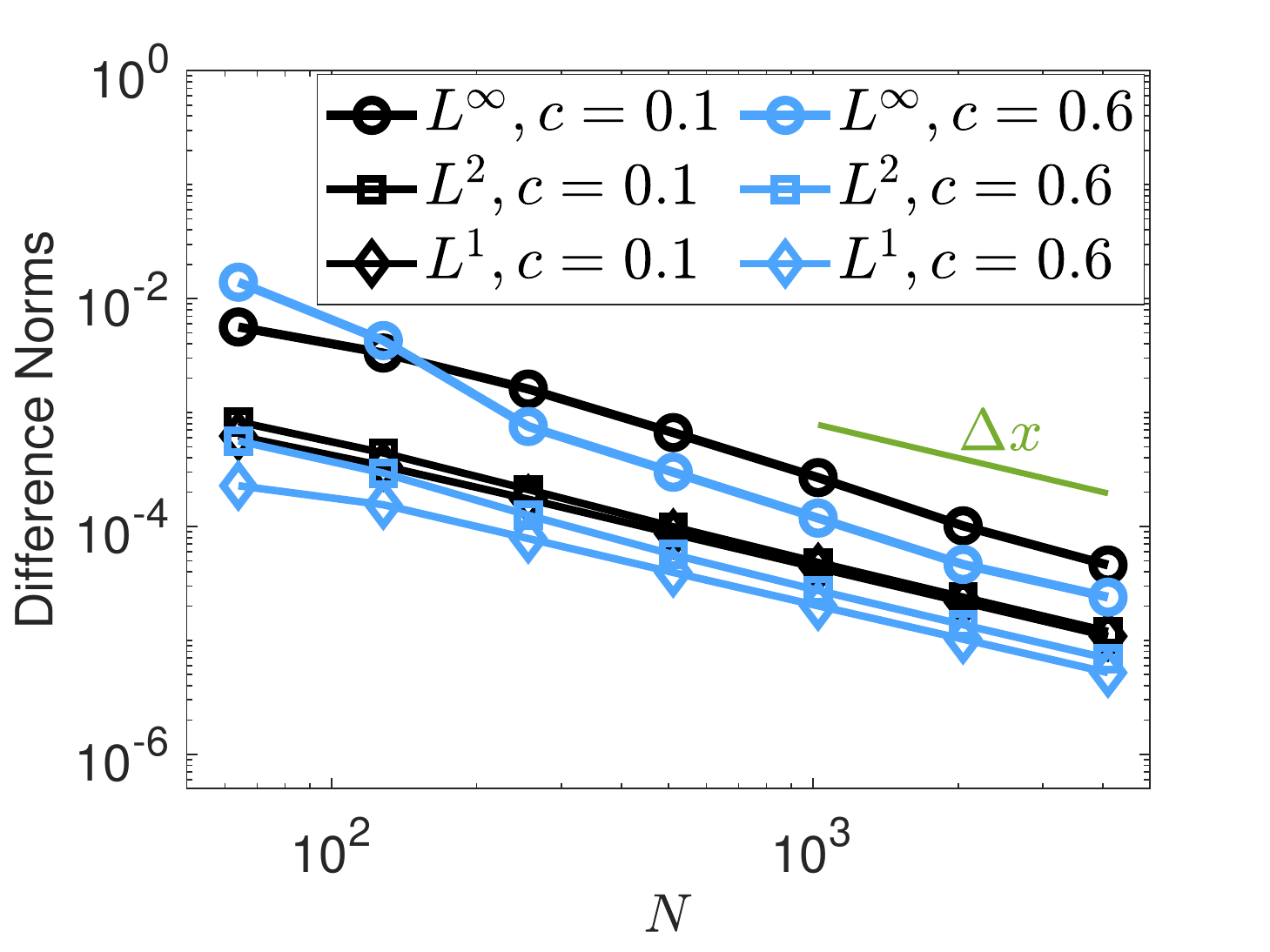}
\caption{\normalsize $u$}
\label{flow past cylinder u ref}
\end{subfigure}
\begin{subfigure}{0.495\textwidth}
\centering
\includegraphics[width=\textwidth]{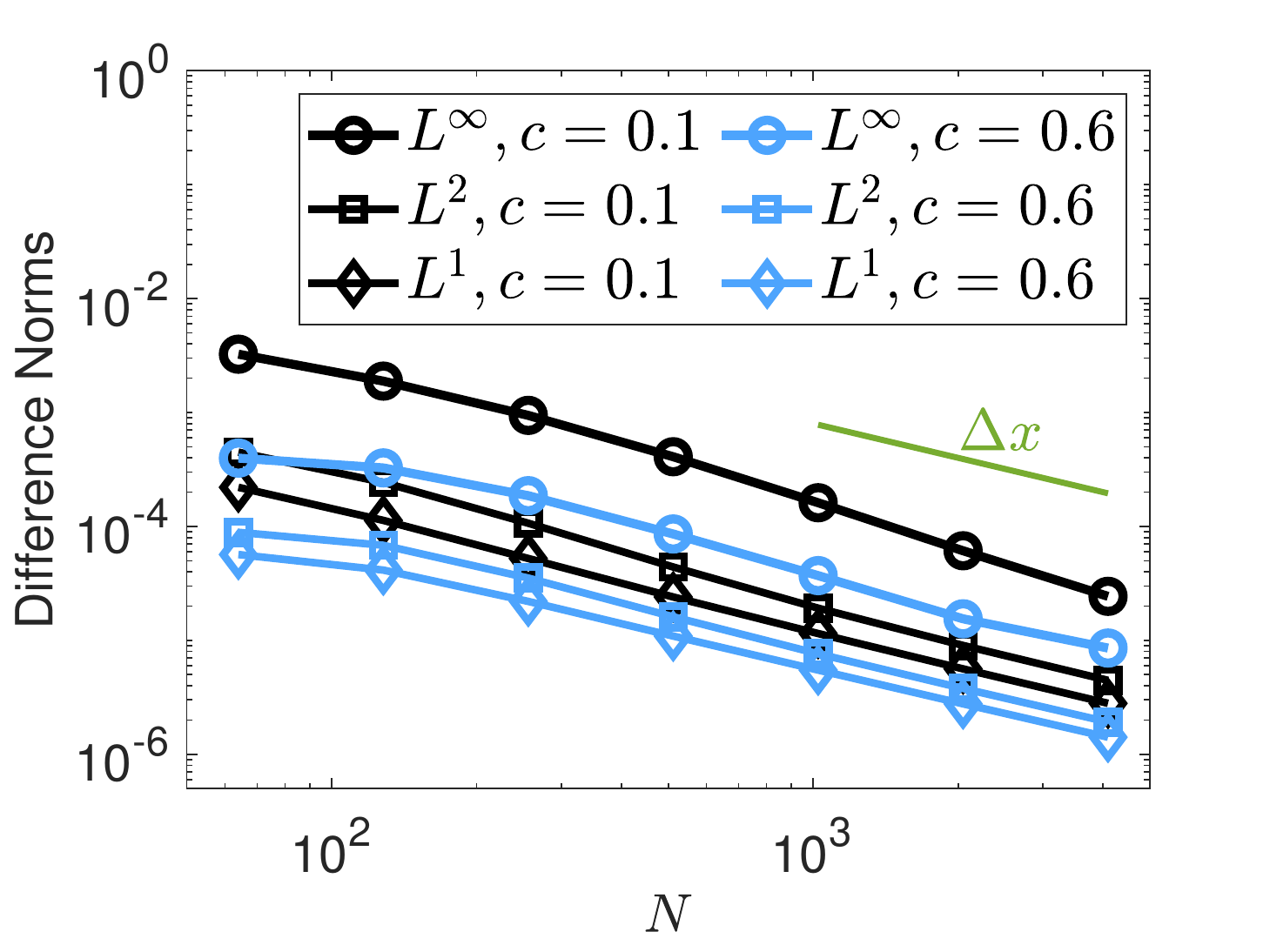}
\caption{\normalsize $v$}
\label{flow past cylinder v ref}
\end{subfigure}
\caption[Solution plot and solution and drag refinement studies for flow past a periodic array of cylinders]{Solution plot and solution and drag refinement studies for flow past a periodic array of cylinders. A Fourier spectral method is used. The computational domain is the periodic box $[-0.5, 0.5]^2$, the boundary point spacing is $ \Delta s  \approx \Delta x$, the interpolation width is $m_1=2(\log_2{N}-4)$, and $\eta=10$. Figure \ref{flow past cylinder solution} shows the solution for grid size $2^{10}$ and $c=0.1$. Figure \ref{flow past cylinder drag ref} shows the refinement study for the dimensionless drag, $D$, using the difference of successive values, relative to the value of the drag. Figures \ref{flow past cylinder u ref}-\ref{flow past cylinder v ref} show refinement studies for the velocity components, using the absolute differences of successive solutions. The refinement studies in Figures \ref{flow past cylinder drag ref}-\ref{flow past cylinder v ref} are for $c=0.1$ and $c=0.6$. }\label{flow past cylinder pics}
\end{figure}

Table \ref{drag table} gives the dimensionless drag forces for the three methods, and Figure \ref{drag plot} plots $D$ versus $c$ for the IBDL method and the asymptotic approximations. For the IBDL and IBSL methods, we use a Fourier spectral method with a boundary point spacing of $\Delta s \approx \Delta x$ and a mesh size of $N=2^{13}$. We see very good agreement with the literature and asymptotic results, especially for the dilute case. For $c\leq0.4$, the relative errors for the IBDL method compared to the Greengard and Kropinski values are less than $0.1\%$. For $c\leq 0.6$, the relative error is less than $1\%$, and for $c=0.7$, less than $2\%$. The IBSL method does outperform our method for a very dense periodic array of cylinders, but we are able to obtain reasonable results for up to $c\approx 0.72$, with a relative error to the asymptotic approximation under $5\%$. The relative gap between cylinders is then approximately $0.04$. Kallemov et al. report a bound on the relative gap for their preconditioned IBSL-like method of approximately $0.02$ \cite{GriffithDonev}. Again, it is worth noting that the method of Kallemov et al. and the IBSL method are able to achieve a much more accurate drag in the case of $c=0.75$, or a relative gap of $0.02$. The poor performance of the IBDL method when the cylinders are this close is worth further investigation, as there are applications in which we may wish to simulate objects with this close interaction. 

Figure \ref{flow past cylinder solution} shows a plot of the IBDL solution for the case of $c=0.1$. The color illustrates the pressure that has mean $0$ on $\C$. This solution was obtained on a grid size of $N=2^{10}$, and the velocity vectors are plotted every $2^6$ grid points. Figures \ref{flow past cylinder drag ref}-\ref{flow past cylinder v ref} then show refinement studies for the dimensionless drag and velocity solutions for two representative flows, with $c=0.1$ and $c=0.6$. The refinement studies are done using the difference of successive solutions, and the drag refinement study uses these differences, relative to the drag values. We see first order convergence in all refinement studies. As usual, the number of iterations of \texttt{gmres} needed for these solutions is essentially independent of grid size. For $0.05 \leq c\leq 0.55$, about $9-15$ iterations were needed. For $c=0.6$, about $16-18$ iterations were needed. For $c=0.7$, about $28$ iterations were needed, and for $c=0.75$, about $44$ iterations were needed.

\subsection{Stokes flow around many objects}\label{6.5 many}

\begin{figure}
\begin{center}
 \begin{tabular}{||c | c | c | c | c | c | c ||} 
 \hline
 \multicolumn{7}{||c||}{Iteration Counts - $9$ Ellipses} \\
 \hline
   &\multicolumn{2}{|c|}{$\Delta s \approx  2\Delta x $}  &\multicolumn{2}{|c|}{$\Delta s \approx  1.5\Delta x $} &\multicolumn{2}{|c||}{$\Delta s \approx  \Delta x $}\\ 
 \hline
 $\Delta x$&\textcolor{blue}{ IBSL} & \textcolor{cyan}{IBDL}&\textcolor{blue}{ IBSL} & \textcolor{cyan}{IBDL} &\textcolor{blue}{ IBSL} & \textcolor{cyan}{IBDL}  \\
 \hline
$2^{-4}$ &  \textcolor{blue}{ 420 } & \textcolor{cyan}{25}&  \textcolor{blue}{ 1058} & \textcolor{cyan}{ 25}
&  \textcolor{blue}{4159 } & \textcolor{cyan}{ 25}  \\
$2^{-5} $&   \textcolor{blue}{509}& \textcolor{cyan}{ 25} &  \textcolor{blue}{1540} & \textcolor{cyan}{25}  &  \textcolor{blue}{5906 } & \textcolor{cyan}{ 25}  \\
$2^{-6 }$&   \textcolor{blue}{599}& \textcolor{cyan}{25} &  \textcolor{blue}{1490 } & \textcolor{cyan}{ 26}&  \textcolor{blue}{9052 } & \textcolor{cyan}{ 26}   \\
$2^{-7}$& \textcolor{blue}{641 }&\textcolor{cyan}{26}   &  \textcolor{blue}{1605} & \textcolor{cyan}{26 } &  \textcolor{blue}{16251 } & \textcolor{cyan}{ 26}   \\
$ 2^{-8}$& \textcolor{blue}{650}&\textcolor{cyan}{26}  &  \textcolor{blue}{1957} & \textcolor{cyan}{26 } &  \textcolor{blue}{ 24436} & \textcolor{cyan}{ 26}   \\
$ 2^{-9}$ & \textcolor{blue}{845}&\textcolor{cyan}{26}   &  \textcolor{blue}{2065} & \textcolor{cyan}{27 } &  \textcolor{blue}{ 26790} & \textcolor{cyan}{ 27}  \\
$ 2^{-10}$  & \textcolor{blue}{1117}&\textcolor{cyan}{27} &  \textcolor{blue}{ 2724} & \textcolor{cyan}{ 27}&  \textcolor{blue}{34571 } & \textcolor{cyan}{ 27} \\
  \hline
  \end{tabular}
\captionof{table}[Number of iterations of \texttt{minres} and \texttt{gmres} for Stokes flow past periodic array of $9$ randomly placed ellipses using the IBSL and IBDL methods]{Number of iterations of \texttt{minres} and \texttt{gmres}, with tolerance $10^{-8}$, needed to solve Equation \eqref{stokes pde flow} using the IBSL and IBDL methods, respectively, for $\Omega$, the region exterior to $9$ randomly placed ellipses, and $\C=[-2, 2]^2$. Both methods use a Fourier spectral method, and the IBDL method uses $\eta=10$.  }
\label{many table}
\end{center}
\end{figure}

\begin{figure}
\begin{subfigure}{0.55\textwidth}
\centering
\includegraphics[width=\textwidth]{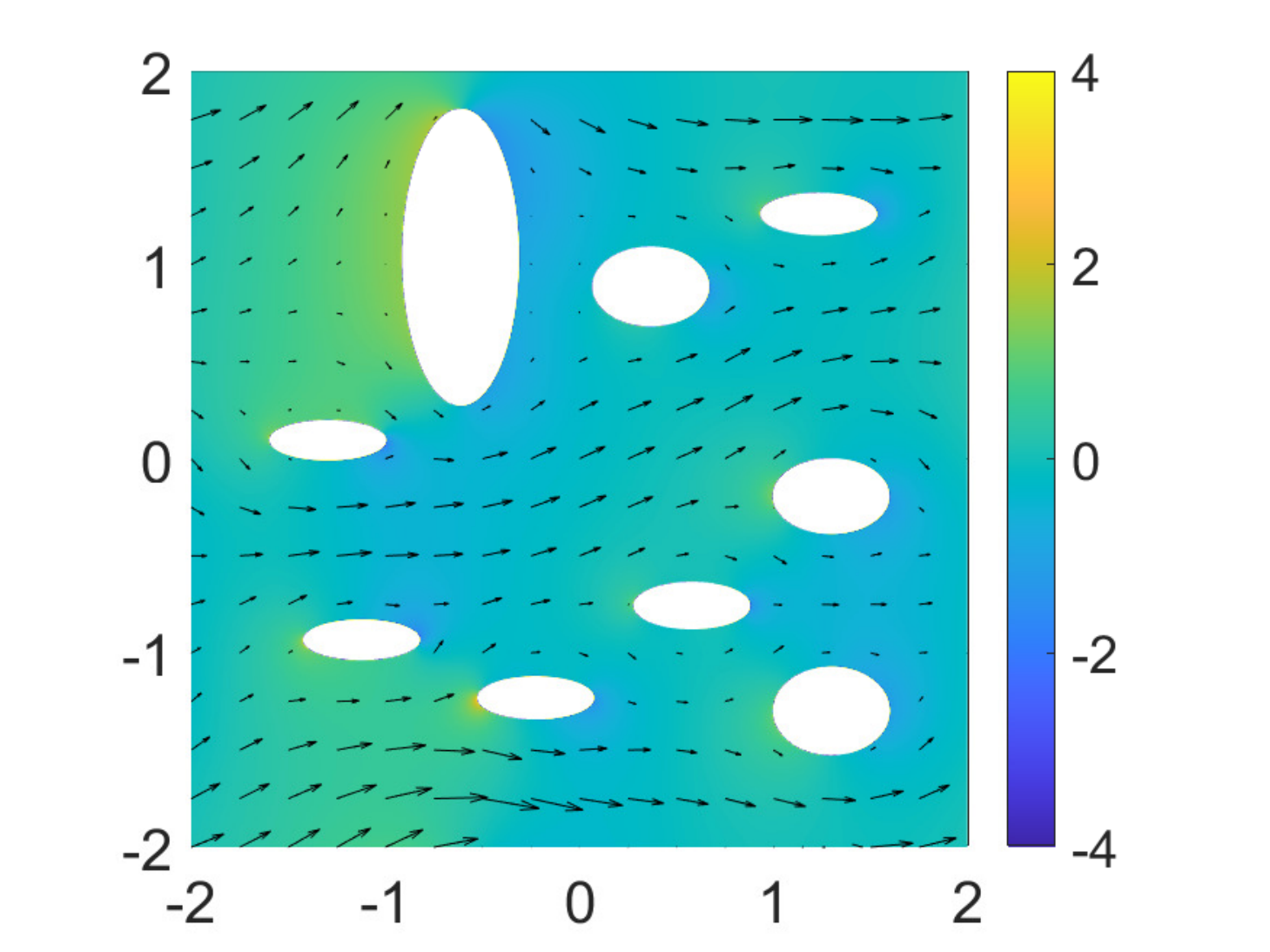}
\caption{\normalsize  Solution plot}
\label{flow past many solution}
\end{subfigure}

\begin{subfigure}{0.495\textwidth}
\centering
\includegraphics[width=\textwidth]{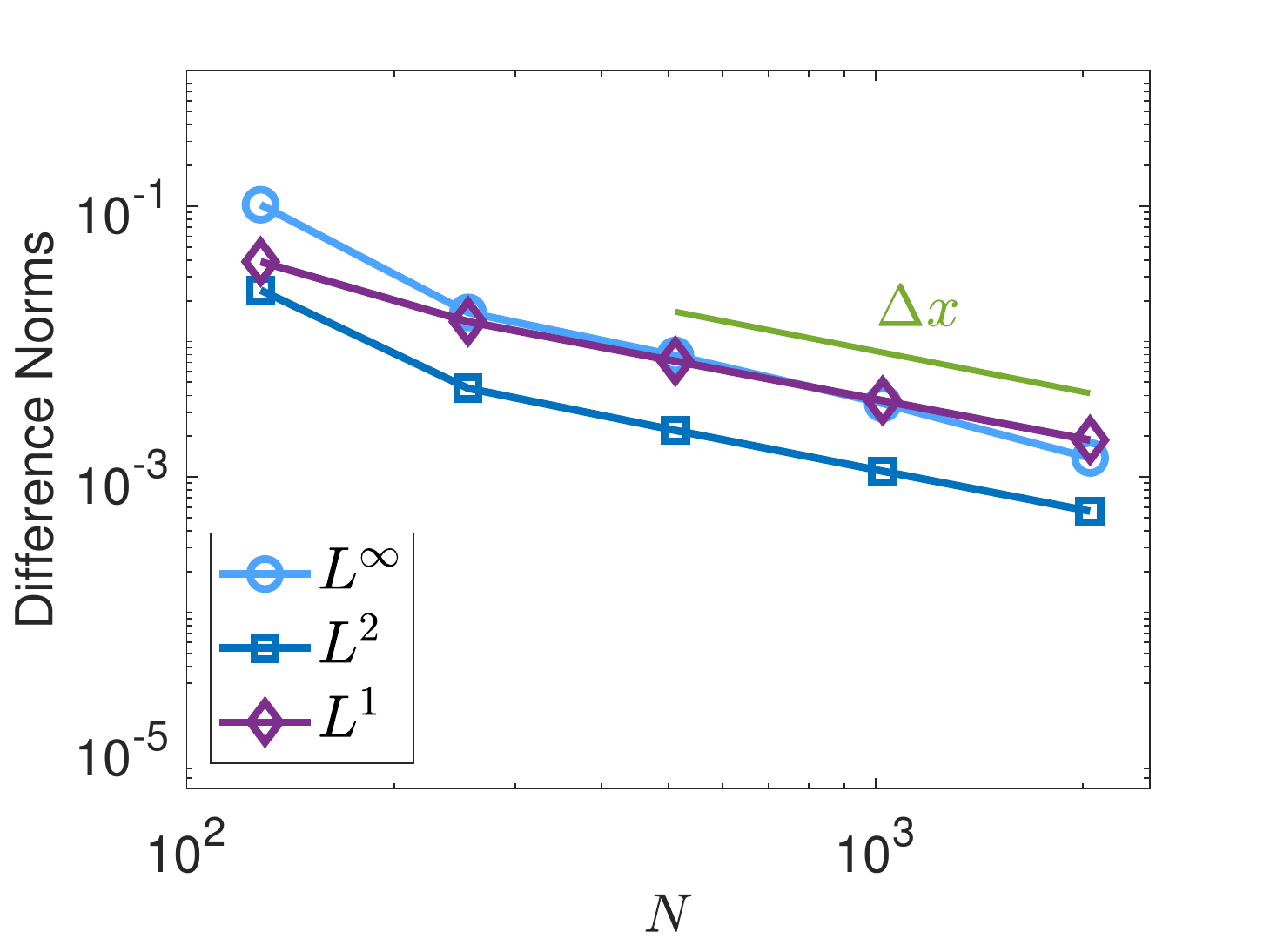}
\caption{\normalsize $u$}
\label{flow past many u ref}
\end{subfigure}
\begin{subfigure}{0.495\textwidth}
\centering
\includegraphics[width=\textwidth]{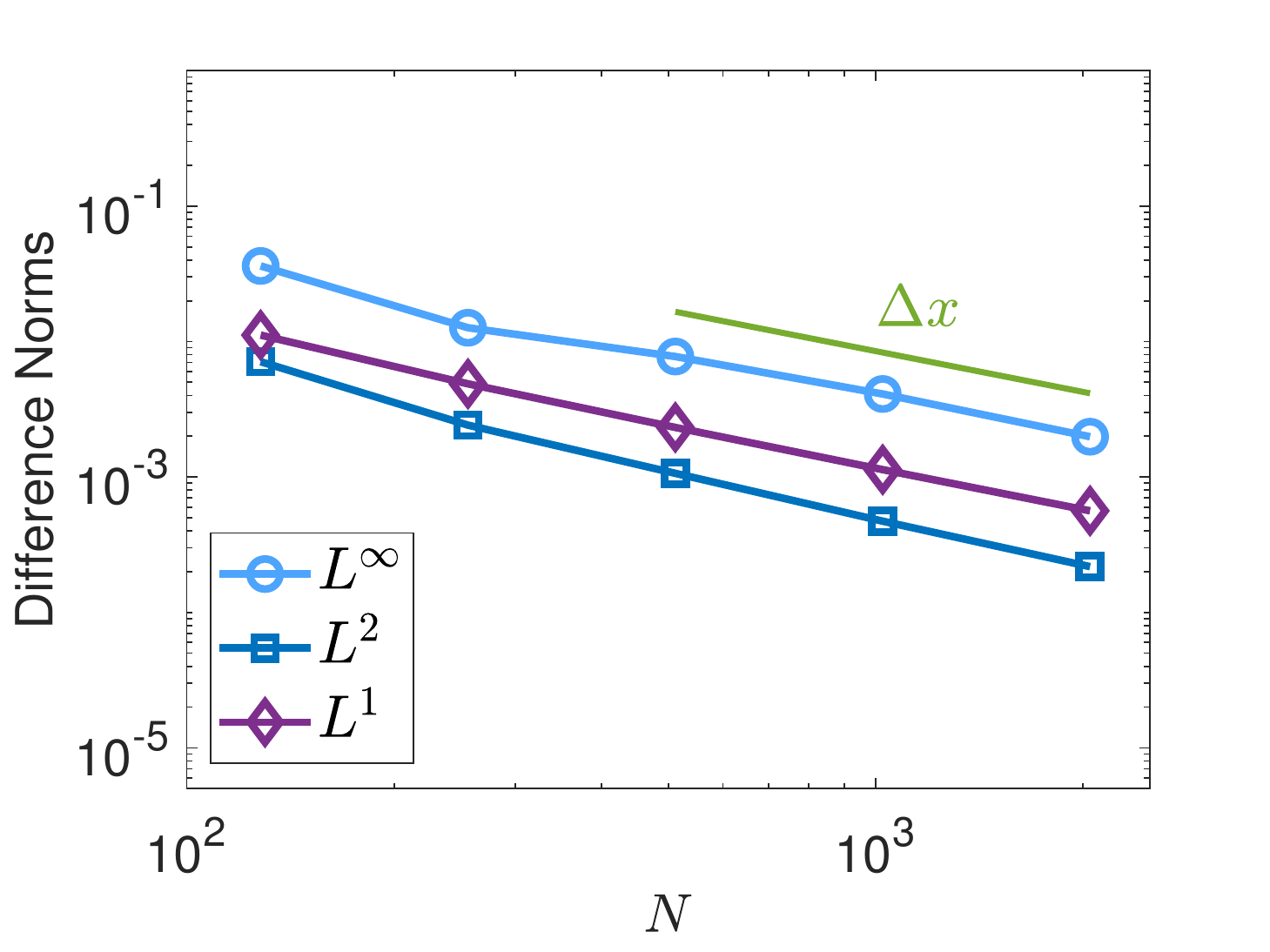}
\caption{\normalsize $v$}
\label{flow past many v ref}
\end{subfigure}
\caption[Solution plot and refinement studies for Stokes flow past periodic array of $9$ randomly placed ellipses]{Solution plot and refinement studies for Stokes flow past a periodic array of $9$ randomly placed ellipses. The computational domain is the periodic box $[-2, 2]^2$, the boundary point spacing is $ \Delta s  \approx  \Delta x$, the interpolation width is $m_1=\log_2{N}-4$, $\eta=10$, and a Fourier spectral method is used. Figure \ref{flow past many solution} shows the solution for grid size $2^{10}$, and Figures \ref{flow past cylinder u ref}-\ref{flow past cylinder v ref} show refinement studies for the velocity components, using the absolute difference of successive solutions.  }\label{many pics}
\end{figure}

We next examine Stokes flow past a periodic array of $9$ arbitrarily placed ellipses that differ in size. We again solve Equation \eqref{stokes pde flow}, with $\Omega$ the region in the computational domain, $\C=[-2,2]^2$, that is exterior to the $9$ ellipses. For the IBDL method, we use an interpolation width of $m_1=\log_2{N}-4$ and $\eta=10$. We solve the equation using both the IBDL and IBSL methods, with a Fourier spectral discretization, using a boundary point spacing of $\Delta s \approx \alpha \Delta x$ for various $\alpha$ values. Table \ref{many table} gives the iteration counts of the Krylov methods needed to solve with these methods. The wide boundary point spacing of $\alpha=2$ is generally used for the IBSL method in order to control the conditioning. However, we see that for this complex domain, the IBSL method requires high iteration counts, even for this wide spacing. When solving this PDE for the region exterior to one circle of radius $0.15$, the IBSL method requires 58 iterations for $\alpha=2$ and $\Delta x= 2^{-10}$. Moving to this more complicated domain, the number of iterations required is multiplied by more than $19$. On the other hand, the iteration count for the IBDL method gets multiplied by less than $3$. Therefore, we can see that the iteration counts resulting from the better conditioning of the IBDL method can become even more favorable when the domain is more complicated, such as when it contains more obstacles. We can see in Table \ref{many table} that the IBSL iteration counts are of course even more extreme for tighter boundary point spacing. In the case of moving objects or Navier-Stokes, where this solve would be required at each step in time, IBSL would be impractical without preconditioning. Figure \ref{many pics} shows the solution plot for $N=2^{10}$ and the IBDL velocity refinement studies. The color on the solution plot illustrates the pressure solution that has mean $0$ on $\C$.

\section{Results: Navier-Stokes equation}\label{ch 6 Navier Stokes}

In this section, we implement the IBDL method for the Navier-Stokes equation, given by 
\begin{subequations} \label{NS}
\begin{alignat}{2}
&\rho\Big( \partial_t \u + \u\dotp \grad \u \Big)=\mu \Delta \u  -\grad p  \qquad && \text{in } \Omega\label{NS 1}  \\
& \Div \u=0 \qquad && \text{in } \Omega  \\
&\u=\U_b \qquad && \text{on } \Gamma.
\end{alignat}
\end{subequations}
We utilize the completed IBDL method described in Section \ref{6.1 Hebeker} and test the performance of the IBDL method by solving for a flow around a cylinder for Reynolds numbers $10$ and $100$. For $Re=10$, we  and compare the drag force against the IBSL method, and for $Re=100$, we compare Strouhal number against the IB method presented by Lai and Peskin \cite{formal}. 

\subsection{Numerical implementation}\label{6.5 implementation}

For time discretization, we use an implicit-explicit scheme analyzed in \cite{timediscret2} and used in an immersed boundary framework in \cite{timediscret1}. In this scheme, the nonlinear terms are treated explicitly in time and the viscous and boundary distributions are treated implicitly in time. The PDE on $\C$ resulting from the completed IBDL formulation is 
\begin{subequations} \label{NS with ibdl}
\begin{alignat}{2}
&\rho\Big( \partial_t \u + \u\dotp \grad \u\Big) =\mu \Delta \u  -\grad p +\eta S \Q + \mu \Div (S\A) \qquad && \text{in } \Omega \\
& \Div \u+S(\Q\dotp \n)=0 \qquad && \text{in } \Omega  \\
&S^*\u+\frac12 \Q=\U_b \qquad && \text{on } \Gamma.
\end{alignat}
\end{subequations}
Using the IMEX scheme, the discretized system is 
\begin{subequations} \label{NS with ibdl discret}
\begin{multline}
\rho\Big(\frac{3\u^{n+1}-4\u^n+\u^{n-1}}{2\Delta t} +2\u^{n} \dotp \grad \u^{n}-\u^{n-1}\dotp\grad\u^{n-1}\Big) \\
 = -\grad p^{n+1}+\mu \Delta \u^{n+1}+\eta S\Q^{n+1} + \mu \Div (S \A^{n+1})\end{multline}
 \begin{alignat}{2}
& \Div \u^{n+1} + S(\Q^{n+1} \dotp \n^{n+1})=0\\
& S^*\u^{n+1}+\frac12\Q^{n+1}=\U_b^{n+1}.
\end{alignat}
\end{subequations}
The first equation can be rewritten as 
\begin{multline}
\Big(\mu\Delta - \frac{3\rho}{2\Delta t}\mathds{I}\Big) \u^{n+1} - \grad p^{n+1}+\eta S\Q^{n+1} + \mu \Div (S \A^{n+1})\\
= \rho \Big(\frac{-4\u^n+\u^{n-1}}{2\Delta t} +2\u^{n} \dotp \grad \u^{n}-\u^{n-1}\dotp\grad\u^{n-1}\Big) , 
\end{multline}
where the right-hand side is known when advancing the system to the next time step. Therefore, by treating the nonlinearity explicitly in time, we recover the Brinkman equation, and we can use the completed IBDL method.

For the application to flow past a cylinder, we use the parameters from \cite{formal}. The computational domain is the periodic box $\C=[0,8]^2$, and the cylinder of radius $0.15$ is located at $(1.85, 4.0)$. To drive the flow, at each step in time, we reset the horizontal velocity to be $u_{\infty}=1$ and the vertical velocity to be $0$ for a strip along the left side of $\C$ that is 4 meshwidths wide. We define the Reynolds number as 
\begin{equation}
Re=\frac{2r \rho u_{\infty} }{\mu},
\end{equation}
where $r$ is the radius of the cylinder. 
We use a time step of $\Delta t= 1.8(10)^{-3}$, and a mesh size of $\Delta x = 2^{-7}$. We use the completed IBDL method with a Fourier spectral discretization and $\eta=10$. For our implementation, we calculate the derivatives for the nonlinear term after completing the interpolation step of the IBDL method. We may instead choose to calculate these derivatives prior to the interpolation and/or by setting the nonphysical domain solutions equal to $0$ to enforce continuity. Initial observations show that these options give comparable results, but more investigation is needed. 

In the case of Navier-Stokes, we examine the dimensionless drag coefficient, defined as
\begin{equation}
C_D=\frac{B_1}{u^2_{\infty} R}
\end{equation}
where $\bm{B}$ is the net force on the cylinder. In the case of Stokes flow, we derived equations for the net force on the cylinder using the boundary distributions of $\F$ and $\Q$ for IBSL and IBDL, respectively. In the case of Navier-Stokes, at steady state, this equation can still be utilized for the IBSL method, but we have found this not to be the case for the IBDL method. More investigation into this matter is needed. Therefore, in the case of Navier-Stokes, we calculate $\bm{B}$ using a more general formula, which we now derive. 

Let us define $R$ as the region inside a box that encircles $\Gamma$. Then let us integrate Equation \eqref{NS 1} over $R'\equiv R\cap \Omega$. Then, using incompressibility to rewrite the second term and replacing the right-hand side with $\Div \sig$, we get
\begin{equation}
0 = \int_{R'} \rho \u_t d\x +\int_{R'}\Div \u\u^{\intercal} d\x- \int_{ R'} \Div \sig d\x.
\end{equation}
Then, using the divergence theorem and the fact that the boundary of $R'$ is given by $\partial R' = \partial R \cup \Gamma$, we get
\begin{equation}
0 = \int_{R'} \rho \u_t d\x+\int_{\partial R} (\u\u^{\intercal}-\sig)\dotp \n_R  dl(\x) +\int_{\Gamma } (\u\u^{\intercal}-\sig)\dotp \n dl(\x),
\end{equation}
where $\n_R$ and $\n$ point out of $R'$. Then since $\u\u^{\intercal}=0$ on $\Gamma$ for flow past a cylinder, we get 
\begin{equation}
\bm{B}=- \int_{\Gamma} \sig \dotp \n dl(\x)= -\int_{R'} \rho \u_t d\x + \int_{\partial R} (\sig - \u\u^{\intercal})\dotp \n_{R} dl(\x).
\end{equation}
In the case of steady flow, the first term vanishes. We use this expression to calculate the force $\bm{B}$ for Navier-Stokes. In this dissertation, we use the box $R=[0.5, 3.203125] \times [2.8125, 5.1875]$.
\subsection{Results: $Re=10$}\label{6.5 cylinder Re10}

\begin{figure}
\begin{subfigure}{0.495\textwidth}
\centering
\includegraphics[width=\textwidth]{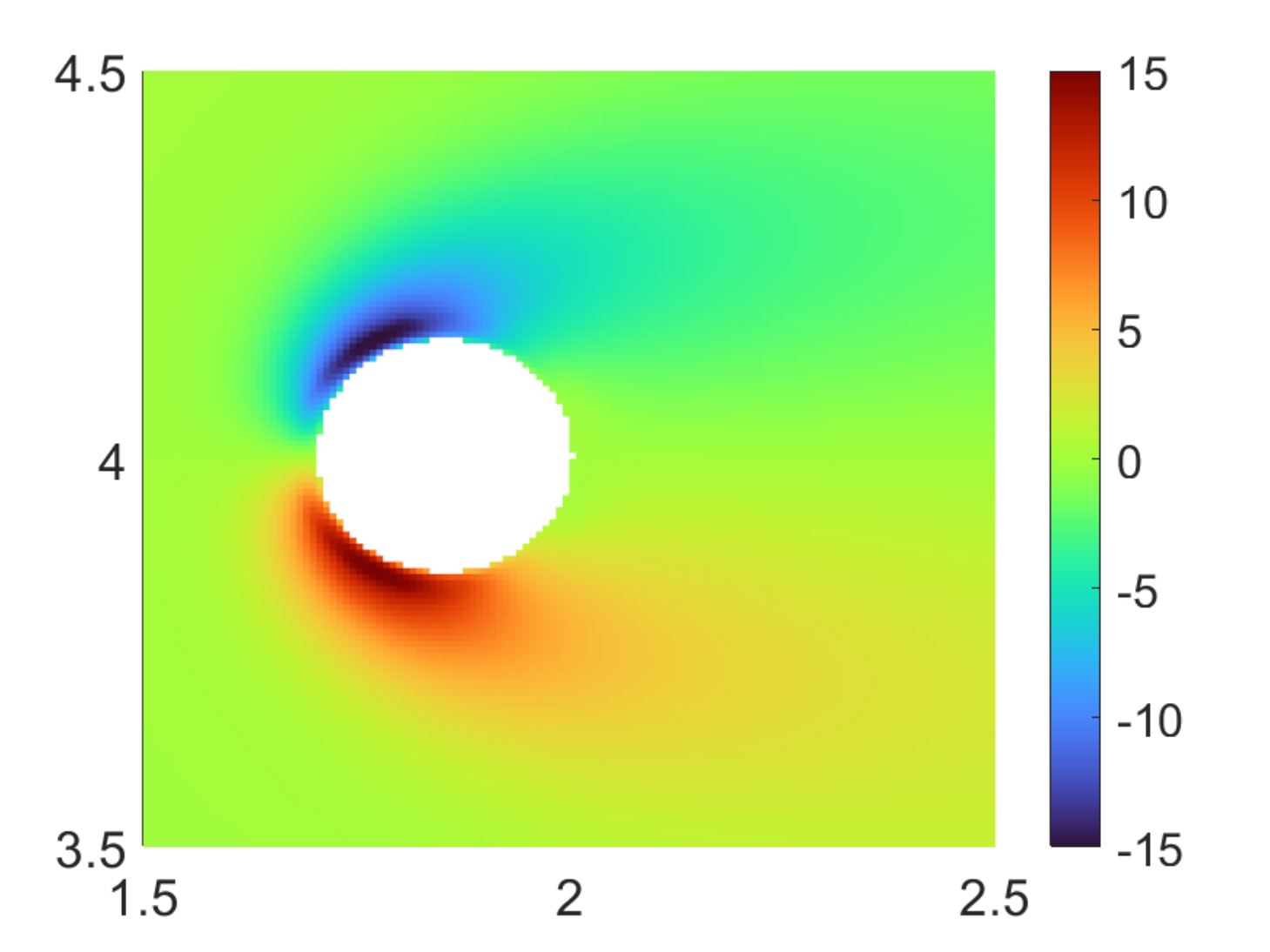}
\caption{\normalsize IBSL}
\label{re10ibsl}
\end{subfigure}
\begin{subfigure}{0.495\textwidth}
\centering
\includegraphics[width=\textwidth]{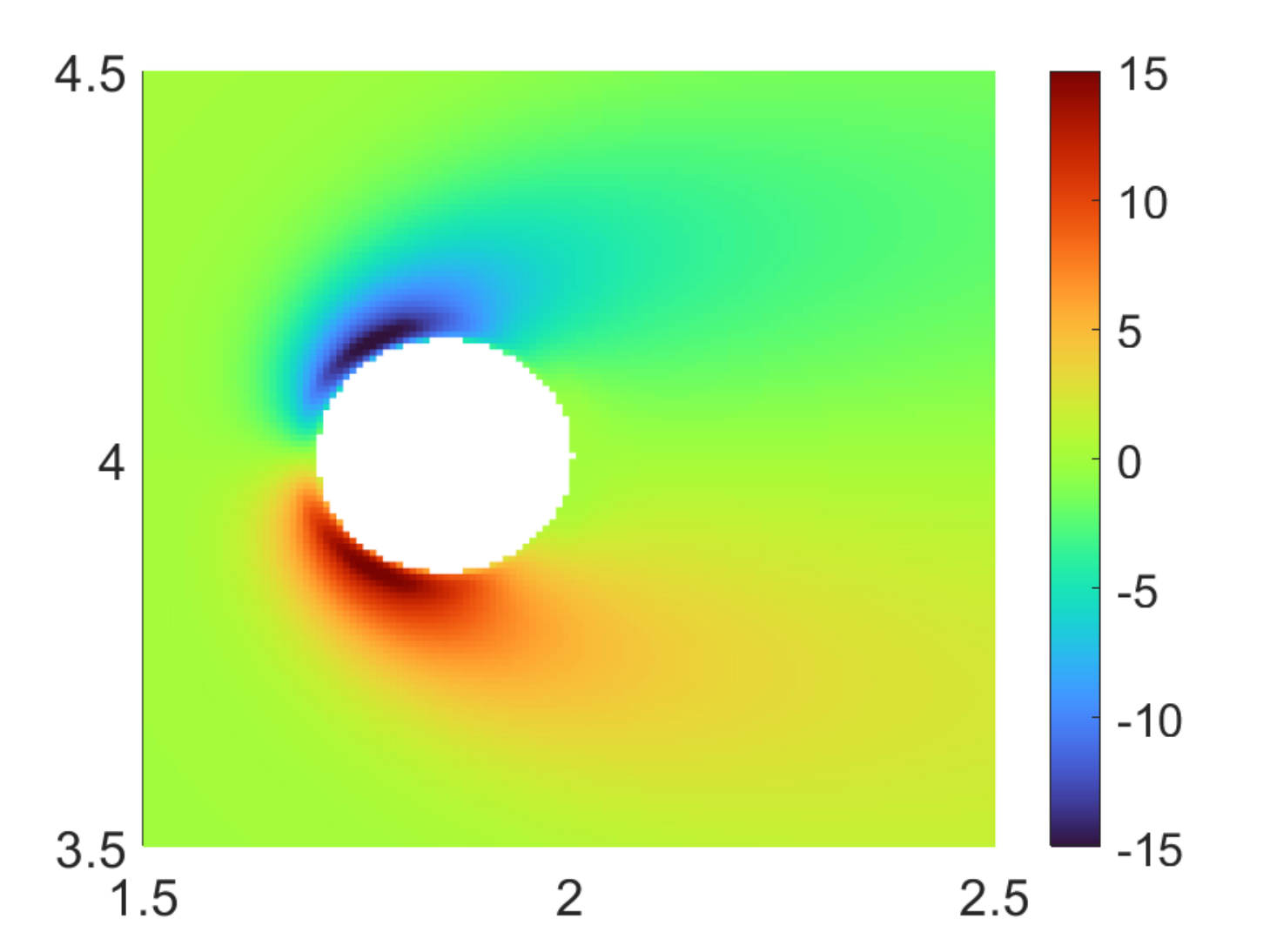}
\caption{\normalsize  IBDL}
\label{re1034}
\end{subfigure}
\caption[Vorticity plots for Navier-Stokes flow past a cylinder for $Re=10$, using IBSL and IBDL methods]{Vorticity plots for Navier-Stokes flow past a cylinder for $Re=10$ at $t=144$, using the IBSL and IBDL methods. The computational domain is the periodic box $[0, 8]^2$, the boundary and grid point spacings are $ \Delta s  \approx  \Delta x=2^{-7}$, the cylinder has radius $0.15$. Both methods use a Fourier spectral method. }\label{vorticityRe10}
\end{figure}

We first solve the PDE for $Re=10$ up to time $t= 144 $ and compare the results from the IBSL method and the IBDL method, with an interpolation width of $m_1=3$ and interior interpolation points $m_2=4$ meshwidths from the boundary. Figure \ref{vorticityRe10} shows the vorticity plots for the IBSL and IBDL methods, and to the eye, we see good vorticity agreement. Table \ref{Re10table} shows the dimensionless drag coefficients for the two methods, and they are within $1\% $ of each other. Further explorations will be made into the affect of the choices of $\eta$, interpolation width, and the method for calculation of the nonlinearity. 

\begin{figure}
\begin{center}
 \begin{tabular}{|| c | c | c ||} 
 \hline
 \multicolumn{3}{||c||}{Drag} \\
 \hline
&IBSL & IBDL\\
\hline
$C_D$ &0.3160    & 0.3137    \\
  \hline
  \end{tabular}
\captionof{table}[Dimensionless drag coefficients for $Re=10$ flow past a cylinder, calcuated with IBSL and IBDL methods]{Dimensionless drag coefficients for $Re=10$ flow past a cylinder, calcuated with IBSL and IBDL methods. The values are averaged over the time period $t=54$ to $t=144$.  }
\label{Re10table}
\end{center}
\end{figure}
\subsection{Results: $Re=100$}\label{6.5 cylinder Re100}

\begin{figure}
\centering
\includegraphics[width=\textwidth]{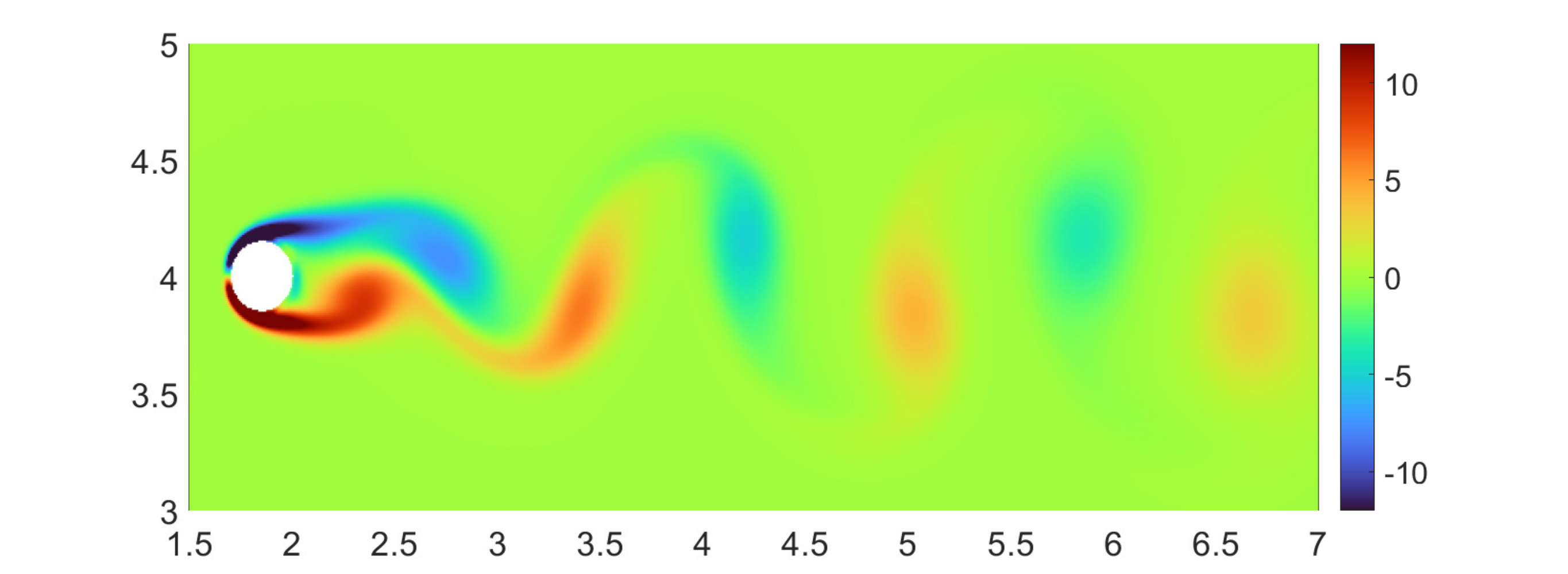}
\caption[Plot of vorticity for Navier-Stokes flow past cylinder for $Re=100$]{Plot of vorticity for Navier-Stokes flow past a cylinder for $Re=100$, using the IBDL method with a Fourier spectral discretization. The computational domain is the periodic box $[0, 8]^2$, the boundary and grid point spacings are $ \Delta s  \approx  \Delta x=2^{-7}$, the cylinder has radius $0.15$, and the method uses an interpolation width of $m_1=6$ and second interpolation points $m_2=8$ meshwidths from the boundary. }\label{Re100plots}
\end{figure}

We next solve the PDE for $Re=100$ up to time $t=234 $ and compare the results from the IBDL method, with $m_1=6$ and $m_2=8$, to that of Lai and Peskin \cite{formal}. Figure \ref{Re100plots} shows the vorticity plot for the IBDL method, demonstrating the vortex shedding that we expect with this Reynolds number. We then report the Strouhal number, which gives the dimensionless frequency with which the vortices are shed. It is defined by 
\begin{equation}
St=\frac{2R f_{vs}}{u_{\infty}},
\end{equation} 
where $f_{vs}$ is the frequency for vortex shedding, which we can measure using $f_{vs}=1/t_{vs}$, where $t_{vs}$ is the average time between peaks in the lift coefficient. We compare to the Strouhal numbers in \cite{formal}, found using two different combinations of time steps, mesh sizes and stiffness parameters. We see that the Strouhal number for the IBDL method matches one of the values from \cite{formal} to three digits. We do not report drag and lift coefficients as more exploration is needed into the sensitivity of these quantities to various IBDL parameters that may be exaggerated in the case of unsteady flow. 

\begin{figure}
\begin{center}
 \begin{tabular}{|| c | c | c | c ||} 
 \hline
 \multicolumn{4}{||c||}{Strouhal Number} \\
 \hline
&IBDL & $L\&P_1$ & $L\&P_2$   \\
\hline
$St$ & 0.155&0.155& 0.165 \\
  \hline
  \end{tabular}
\captionof{table}[Strouhal numbers for $Re=100$ flow past a cylinder]{Strouhal numbers for $Re=100$ flow past a cylinder. The first column is calcluated with the IBDL method using a Fourier spectral discretization and an interpolation width of $m_1=6$, a time step of $\Delta t = 1.8(10^{-3})$, and a grid point spacing of $\Delta x = 2^{-7}$. The second column gives the results from Scheme 2 in Peskin and Lai \cite{formal}, using $\Delta t = 1.8(10^{-3})$, $\Delta x = 2^{-6}$, and stiffness $\kappa=4.8(10^{-4})$. The third column gives the results from Scheme 2, using $\Delta t= 0.9(10^{-3})$ and $\Delta x =2^{-7}$, and  $\kappa=9.6(10^{-4})$.   }
\label{Re100table}
\end{center}
\end{figure}

   \chapter{Discussion}
   \label{chapter 7}

We have developed the Immersed Boundary Double Layer method, a numerical method for linear PDE on complex domains with Dirichlet or Neumann boundary conditions that is more efficient than the constraint formulation of the Immersed Boundary method. We achieved this greater efficiency by reformulating the constraint IB method to correspond to a regularized double layer integral equation, which has better conditioning than the single layer integral equation to which the original constraint formulation corresponds. With this better conditioning, we can solve for the  Lagrangian dipole force distribution, $Q$, in a small number of iterations of a Krylov method, and the iteration count does not increase as we refine the mesh. Furthermore, the computed $Q$ is relatively smooth and converges reasonably to the exact potential strength. Both this convergence and the iteration counts are in stark contrast to the IB constraint method, in which the Lagrangian delta force distribution, $F$, requires many iterations to compute, the number of which increases with finer meshes and tighter boundary point spacing. Additionally, this distribution often fails to converge due to high-frequency noise. 

The IBDL method is related to the method of immersed layers, presented in the recent work of Eldredge \cite{eldredge}. In this work, Eldredge uses the indicator function discussed in Section \ref{5.3 flagging} to develop extended forms of PDEs that govern variables that define different functions on either side of the boundary. Jump quantities corresponding to the strengths of single and double layer potentials naturally emerge in the PDEs. The goal of this method is to be able to enforce different constraints on either side of the immersed boundary in order to obtain solutions on both sides and to accurately predict the surface traction on one side of the boundary. With the inclusion of jumps in both the solution and derivative, this method is able to achieve these goals. We should also note here that our application of the IBDL method to a Neumann problem can be seen as a specific case of the method of immersed layers \cite{eldredge}. However, if we focus on Dirichlet problems, the method of immersed layers again requires inverting an operator corresponding to a first-kind integral equation that is poorly conditioned. 

However, the IBDL method is used for Dirichlet problems for which the solution is only desired on one side of the immersed boundary, and it uses a double layer potential whose strength is the unknown jump in solution values. This allows us to develop a method that corresponds to a second-kind integral equation that can be solved much more efficiently. One can see the IB constraint method of Taira and Colonius \cite{Taira} as corresponding to a single layer integral representation, the IBDL method presented in this dissertation as corresponding to a double layer integral representation, and the method of immersed layers of Eldredge \cite{eldredge} as corresponding to an integral representation that includes both single and double layers. And in the case of Dirichlet boundary conditions, the better conditioning of the IBDL operator gives a more efficient method. 

We also make note of another difference between the IBDL method and the method of immersed layers. While \cite{eldredge} uses a lattice Green's function for unbounded external flows, the IBDL method does not involve an explicit Green's function. Instead the method obtains the convolution of the boundary potential with a regularized Green's function by solving the PDE after spreading the boundary potential with the regularized delta function. Thus the IBDL method can be used on more general domains. In this paper we use a periodic computational domain, but one could also easily use other boundary conditions and domains.

The IBDL method retains much of the flexibility and robustness of the IB method. The communication between the Lagrangian coordinate system and the underlying Cartesian grid is achieved with convolutions with regularized delta functions, and no analytical Green's functions are needed in the computation. Furthermore, minimal geometric information is needed for the immersed boundary. We only need boundary points, unit normals, an indicator function flagging points as inside the PDE domain $\Omega$, and an indicator function flagging points as near-boundary points. Furthermore, we can reduce this input list to just the set of points. The unit normals can be estimated as described in Section \ref{5.3 spread}, and the interior indicator function can be calculated as described in Section \ref{5.3 flagging}. Since such minimal geometric information is needed, this method can be easily used for complex PDE domains. While we focused on two dimensions in this dissertation, nothing in the method formulation requires this. We chose to use an arclength parametrization in our method development, and in three dimensions, this would simply require that we estimate the areas of the discretized boundary elements. Lastly, we would again require unit normal approximations.

There are a few ways in which the method is not as flexible and easily generalized as the IBSL method. Unlike the original IB constraint method, the IBDL method produces a discontinuous solution across the boundary, requiring careful handling of near-boundary points in order to maintain pointwise convergence in this region. We have found that in the case of a scalar PDE, a finite difference discretization often requires less interpolation. However, a Fourier spectral method can still be utilized by interpolating a region for which the size approaches $0$ as we refine the mesh. In the case of Stokes and Brinkman equations, we have found that the form of the IBDL method complicates the implementation of a finite difference method. We have achieved first-order accuracy by approximating the discrete wide Laplacian with the regular Laplacian and utilizing a smoother delta function. However, one area of future research is in developing a more suitable finite difference discretization. 

Additionally, we have found that in the case of the exterior Poisson and Stokes equations, the constraint matrix is singular with a one-dimensional nullspace. It therefore requires careful handling of the inversion process. However, by adding a single layer term to form the completed IBDL method, we get a constraint matrix that is invertible, as in the IBSL case. Another area of future research is in better understanding this nullspace and its role in the IBDL method.  

The increased efficiency of this method makes it practical for use on time-dependent PDEs. For example, with an implicit time-discretization of the diffusion equation, we could use the IBDL method to solve a Helmholtz equation at each time step. We can also use the IBDL method for a PDE with a nonlinearity by using an implicit-explicit time-stepping scheme, as is done for Navier-Stokes equation in Section \ref{ch 6 Navier Stokes}. 

Another benefit of the IBDL method is the convergence of the potential strength. This convergence provides the opportunity to use it to obtain the normal derivative of the solution on the boundary in the case of a scalar PDE or boundary traction in the case of Stokes. We would use analytical expressions for derivatives of double layer potentials in order to achieve this. This opens the possibility of improving the IBDL method by obtaining a continuous extension of the solution across the boundary, possibly eliminating the need for the interpolation step of the method. Another area of future work is in conceivably using such normal derivatives on the boundary to create a method like the Immersed Boundary Smooth Extension (IBSE) method \cite{Stein} to obtain higher accuracy in the IB method. The IBSE method can achieve this higher accuracy, but the poor conditioning of the constraint method is exacerbated, making the method impractical for three dimensions or moving domains. With the greater efficiency of the IBDL method, we may be able to achieve a higher order convergence without this limitation. We therefore have several paths forward for improving the IBDL method, using it to achieve higher accuracy, and applying it to a wider range of problems.


       
   \backmatter
   
   \bibliographystyle{siam}
   
      \let\oldaddcontentsline\addcontentsline
\renewcommand{\addcontentsline}[3]{}
   \bibliography{DissertationBibliography}
   \let\addcontentsline\oldaddcontentsline

   \addcontentsline{toc}{chapter}{\textbf{Bibliography}}
   
\end{document}